# n-LINEAR ALGEBRA OF TYPE II

**W. B. Vasantha Kandasamy**
e-mail: vasanthakandasamy@gmail.com
web: http://mat.iitm.ac.in/~wbv
www.vasantha.in

**Florentin Smarandache**
e-mail: smarand@unm.edu

**2008**

# n-LINEAR ALGEBRA OF TYPE II

**W. B. Vasantha Kandasamy**
**Florentin Smarandache**

# CONTENTS

## PREFACE

This book is a continuation of the book n-linear algebra of type I and its applications. Most of the properties that could not be derived or defined for n-linear algebra of type I is made possible in this new structure: n-linear algebra of type II which is introduced in this book. In case of n-linear algebra of type II, we are in a position to define linear functionals which is one of the marked difference between the n-vector spaces of type I and II. However all the applications mentioned in n-linear algebras of type I can be appropriately extended to n-linear algebras of type II. Another use of n-linear algebra (n-vector spaces) of type II is that when this structure is used in coding theory we can have different types of codes built over different finite fields whereas this is not possible in the case of n-vector spaces of type I. Finally in the case of n-vector spaces of type II we can obtain n-eigen values from distinct fields; hence, the n-characteristic polynomials formed in them are in distinct different fields.

An attractive feature of this book is that the authors have suggested 120 problems for the reader to pursue in order to understand this new notion. This book has three chapters. In the first chapter the notion of n-vector spaces of type II are introduced. This chapter gives over 50 theorems. Chapter two introduces the notion of n-inner product vector spaces of type II, n-bilinear forms and n-linear functionals. The final chapter suggests over a hundred problems. It is important that the reader should be well versed with not only linear algebra but also n-linear algebras of type I.

The authors deeply acknowledge the unflinching support of Dr.K.Kandasamy, Meena and Kama.

W.B.VASANTHA KANDASAMY
FLORENTIN SMARANDACHE

**Chapter One**

# n-VECTOR SPACES OF TYPE II AND THEIR PROPERTIES

In this chapter we for the first time introduce the notion of n-vector space of type II. These n-vector spaces of type II are different from the n-vector spaces of type I because the n-vector spaces of type I are defined over a field F where as the n-vector spaces of type II are defined over n-fields. Some properties enjoyed by n-vector spaces of type II cannot be enjoyed by n-vector spaces of type I. To this; we for the sake of completeness just recall the definition of n-fields in section one and n-vector spaces of type II are defined in section two and some important properties are enumerated.

## 1.1 n-Fields

In this section we define n-field and illustrate it by examples.

**DEFINITION 1.1.1:** *Let $F = F_1 \cup F_2 \cup \ldots \cup F_n$ where each $F_i$ is a field such that $F_i \not\subset F_j$ or $F_j \not\subset F_i$ if $i \neq j$, $1 \leq i, j \leq n$, we call F a n-field.*

We illustrate this by the following example.

***Example 1.1.1:*** Let $F = R \cup Z_3 \cup Z_5 \cup Z_{17}$ be a 4-field.

Now how to define the characteristic of any n-field, $n \geq 2$.

**DEFINITION 1.1.2:** *Let $F = F_1 \cup F_2 \cup \ldots \cup F_n$ be a n-field, we say F is a n-field of n-characteristic zero if each field $F_i$ is of characteristic zero, $1 \leq i \leq n$.*

***Example 1.1.2:*** Let $F = F_1 \cup F_2 \cup F_3 \cup F_4 \cup F_5 \cup F_6$, where $F_1 = Q(\sqrt{2})$, $F_2 = Q(\sqrt{7})$, $F_3 = Q(\sqrt{3}, \sqrt{5})$, $F_4 = Q(\sqrt{11})$, $F_5 = Q(\sqrt{3}, \sqrt{19})$ and $F_6 = Q(\sqrt{5}, \sqrt{17})$; we see all the fields $F_1$, $F_2$, …, $F_6$ are of characteristic zero thus F is a 6-field of characteristic 0.

Now we proceed on to define an n-field of finite characteristic.

**DEFINITION 1.1.3:** *Let $F = F_1 \cup F_2 \cup \ldots \cup F_n$ $(n \geq 2)$, be a n-field. If each of the fields $F_i$ is of finite characteristic and not zero characteristic for i = 1, 2, …, n then we call F to be a n-field of finite characteristic.*

***Example 1.1.3:*** Let $F = F_1 \cup F_2 \cup F_3 \cup F_4 = Z_5 \cup Z_7 \cup Z_{17} \cup Z_{31}$, F is a 4-field of finite characteristic.

*Note:* It may so happen that in a n-field $F = F_1 \cup F_2 \cup \ldots \cup F_n$, $n \geq 2$ some fields $F_i$ are of characteristic zero and some of the fields $F_j$ are of characteristic a prime or a power of a prime. Then how to define such n-fields.

**DEFINITION 1.1.4:** *Let $F = F_1 \cup F_2 \cup \ldots \cup F_n$ be a n-field $(n \geq 2)$, if some of the $F_i$'s are fields of characteristic zero and some of the $F_j$'s are fields of finite characteristic $i \neq j$, $1 \leq i, j \leq n$ then we define the characteristic of F to be a mixed characteristic.*

***Example 1.1.4:*** Let $F = F_1 \cup F_2 \cup F_3 \cup F_4 \cup F_5$ where $F_1 = Z_2$, $F_2 = Z_7$, $F_3 = Q(\sqrt{7})$, $F_4 = Q(\sqrt{3}, \sqrt{5})$ and $F_5 = Q(\sqrt{3}, \sqrt{23}, \sqrt{2})$ then F is a 5-field of mixed characteristic; as $F_1$ is of characteristic two, $F_2$ is a field of characteristic 7, $F_3$, $F_4$ and $F_5$ are fields of characteristic zero.

Now we define the notion of n-subfields.

**DEFINITION 1.1.5:** *Let $F = F_1 \cup F_2 \cup \ldots \cup F_n$ be a n-field ($n \geq 2$). $K = K_1 \cup K_2 \cup \ldots \cup K_n$ is said to be a n-subfield of F if each $K_i$ is a proper subfield of $F_i$, i = 1, 2, …, n and $K_i \not\subset K_j$ or $K_j \not\subset K_i$ if $i \neq j$, $1 \leq i, j \leq n$.*

We now give an example of an n-subfield.

***Example 1.1.5:*** Let $F = F_1 \cup F_2 \cup F_3 \cup F_4$ where

$$F_1 = Q(\sqrt{2}, \sqrt{3}),\ F_2 = Q(\sqrt{7}, \sqrt{5}),$$

$$F_3 = \left(\frac{Z_2[x]}{\langle x^2 + x + 1\rangle}\right) \text{ and } F_4 = \left(\frac{Z_{11}[x]}{\langle x^2 + x + 1\rangle}\right)$$

be a 4-field. Take $K = K_1 \cup K_2 \cup K_3 \cup K_4 = (Q(\sqrt{2}) \cup Q(\sqrt{7}) \cup Z_2 \cup Z_{11} \subseteq F = F_1 \cup F_2 \cup F_3 \cup F_4$. Clearly K is a 4-subfield of F.

It may so happen for some n-field F, we see it has no n-subfield so we call such n-fields to be prime n-fields.

***Example 1.1.6:*** Let $F = F_1 \cup F_2 \cup F_3 \cup F_4 = Z_7 \cup Z_{23} \cup Z_2 \cup Z_{17}$ be a 4-field. We see each of the field $F_i$'s are prime, so F is a n-prime field (n = 4).

**DEFINITION 1.1.6:** *Let $F = F_1 \cup F_2 \cup \ldots \cup F_n$ be a n-field ($n \geq 2$) if each of the $F_i$'s is a prime field then we call F to be prime n-field.*

It may so happen that some of the fields may be prime and others non primes in such cases we call F to be a semiprime n-field. If all the fields $F_i$ in F are non prime i.e., are not prime fields then we call F to be a non prime field. Now in case of n-semiprime field or semiprime n-field if we have an m-subfield m < n then we call it as a quasi m-subfield of F m < n.

We illustrate this situation by the following example.

***Example 1.1.7:*** Let

$$F = Q \cup Z_7 \cup \frac{Z_2[x]}{\langle x^2 + x + 1 \rangle} \cup \frac{Z_7[x]}{\langle x^2 + 1 \rangle} \cup \frac{Z_3[x]}{\langle x^3 + 9 \rangle};$$

be a 5-field of mixed characteristic; clearly F is a 5-semiprime field. For Q and $Z_7$ are prime fields and the other three fields are non-prime; take $K = K_1 \cup K_2 \cup K_3 \cup K_4 \cup K_5 = \phi \cup \phi \cup Z_2 \cup Z_7 \cup Z_3 \subset F = F_1 \cup F_2 \cup F_3 \cup F_4 \cup F_5$. Clearly K is a quasi 3-subfield of F.

## 1.2 n-Vector Spaces of Type II

In this section we proceed on to define n-vector spaces of type II and give some basic properties about them.

**DEFINITION 1.2.1:** *Let $V = V_1 \cup V_2 \cup \ldots \cup V_n$ where each $V_i$ is a distinct vector space defined over a distinct field $F_i$ for each i, i = 1, 2, ..., n; i.e., $V = V_1 \cup V_2 \cup \ldots \cup V_n$ is defined over the n-field $F = F_1 \cup F_2 \cup \ldots \cup F_n$. Then we call V to be a n-vector space of type II.*

We illustrate this by the following example.

***Example 1.2.1:*** Let $V = V_1 \cup V_2 \cup V_3 \cup V_4$ be a 4-vector space defined over the 4-field $Q(\sqrt{2}) \cup Q(\sqrt{3}) \cup Q(\sqrt{5}) \cup Q(\sqrt{7})$. V is a 4-vector space of type II.

Unless mention is made specifically we would by an n-vector space over F mean only n-vector space of type II, in this book.

***Example 1.2.2:*** Let $F = F_1 \cup F_2 \cup F_3$ where $F_1 = Q(\sqrt{3}) \times Q(\sqrt{3}) \times Q(\sqrt{3})$ a vector space of dimension 3 over $Q(\sqrt{3})$. $F_2 = Q(\sqrt{7})$ [x] be the polynomial ring with coefficients from $Q(\sqrt{7})$, $F_2$ is a vector space of infinite dimension over $Q(\sqrt{7})$ and

$$F_3 = \left\{ \begin{bmatrix} a & b \\ c & d \end{bmatrix} \Big/ a, b, c, d \in Q(\sqrt{5}) \right\},$$

$F_3$ is a vector space of dimension 4 over $Q(\sqrt{5})$. Thus $F = F_1 \cup F_2 \cup F_3$ is a 3-vector space over the 3-field $Q(\sqrt{3}) \cup Q(\sqrt{7}) \cup Q(\sqrt{5})$ of type II.

Thus we have seen two examples of n-vector spaces of type II. Now we will proceed on to define the linearly independent elements of the n-vector space of type II. Any element $\alpha \in V = V_1 \cup V_2 \cup \ldots \cup V_n$ is a n-vector of the form $(\alpha_1 \cup \alpha_2 \cup \ldots \cup \alpha_n)$ or $(\alpha^1 \cup \alpha^2 \cup \ldots \cup \alpha^n)$ where $\alpha_i$ or $\alpha^i \in V_i$, i = 1, 2, … , n and each $\alpha_i$ is itself a m row vector if dimension of $V_i$ is m.

We call $S = S_1 \cup S_2 \cup \ldots \cup S_n$ where each $S_i$ is a proper subset of $V_i$ for i = 1, 2, …, n as the n-set or n-subset of the n-vector space V over the n-field F. Any element $\alpha = (\alpha^1 \cup \alpha^2 \cup \ldots \cup \alpha^n) = \{\alpha^1_1, \alpha^1_2, \ldots, \alpha^1_{n_1}\} \cup \{\alpha^2_1, \alpha^2_2, \ldots, \alpha^2_{n_2}\} \cup \ldots \cup \{\alpha^n_1, \alpha^n_2, \ldots, \alpha^n_{n_n}\}$ of $V = V_1 \cup V_2 \cup \ldots \cup V_n$ where $\{\alpha^i_1, \alpha^i_2, \ldots, \alpha^i_{n_i}\} \subseteq V_i$ for i = 1, 2, …, n and $\alpha^i_j \in V_i$, j = 1, 2, …, $n_i$ and $1 \le i \le n$.

**DEFINITION 1.2.2:** *Let $V = V_1 \cup V_2 \cup \ldots \cup V_n$ be a n-vector space of type II defined over the n-field $F = F_1 \cup F_2 \cup \ldots \cup F_n$, $(n \ge 2)$. Let $S = \{\alpha^1_1, \alpha^1_2, \ldots, \alpha^1_{k_i}\} \cup \{\alpha^2_1, \alpha^2_2, \ldots, \alpha^2_{k_2}\} \cup \ldots \cup \{\alpha^n_1, \alpha^n_2, \ldots, \alpha^n_{k_n}\}$ be a proper n-subset of $V = V_1 \cup V_2 \cup \ldots \cup V_n$ we say the n-set S is a n-linearly independent n-subset of V over $F = F_1 \cup F_2 \cup \ldots \cup F_n$ if and only if each subset $\{\alpha^i_1, \alpha^i_2, \ldots, \alpha^i_{k_i}\}$* *is a linearly independent subset of $V_i$ over $F_i$, this must be true for each i, i = 1, 2, ... , n. If even one of the subsets of the n-subset S; say $\{\alpha_1^j, \alpha_2^j, \ldots, \alpha_{k_j}^j\}$ is not a linearly independent subset of $V_j$ over $F_j$, then the subset S of V is not a n-linearly independent subset of V, $1 \le j \le n$. If in the n-subset S of V every subset is not linearly independent subset we say S is a n-linearly dependent subset of V over $F = F_1 \cup F_2 \cup \ldots \cup F_n$.*

*Now if the n-subset $S = S_1 \cup S_2 \cup \ldots \cup S_n \subseteq V$; $S_i \subset V_i$ is a subset of $V_i$, i = 1, 2, ..., n is such that some of the subsets $S_j$ are linearly independent subsets over $F_j$ and some of the subsets $S_i$ of $V_i$ are linearly dependent subsets of $S_i$ over $F_i$, $1 \le i, j \le n$ then we call S to be a semi n-linearly independent n-subset over F or equivalently S is a semi n-linearly dependent n-subset of V over F.*

We illustrate all these situations by the following examples.

***Example 1.2.3:*** Let $V = V_1 \cup V_2 \cup V_3 = \{Q(\sqrt{3}) \times Q(\sqrt{3})\} \cup \{Q(\sqrt{7})[x]_5$ | this contains only polynomials of degree less than or equal to 5$\} \cup$

$$\left\{\begin{bmatrix} a & b \\ c & d \end{bmatrix} \middle| a, b, c, d \in Q(\sqrt{2})\right\}$$

is a 3-vector space over the 3-field $F = Q(\sqrt{3}) \cup Q(\sqrt{7}) \cup Q(\sqrt{2})$.

Take

$$S = \{(1, 0), (5, 7)\} \cup \{x^3, x, 1\} \cup \left\{\begin{bmatrix} 1 & 0 \\ 0 & 2 \end{bmatrix}, \begin{bmatrix} 0 & 1 \\ 0 & 0 \end{bmatrix}\right\}$$

$$= S_1 \cup S_2 \cup S_3 \subseteq V_1 \cup V_2 \cup V_3;$$

a proper 3-subset of V. It is easily verified S is a 3-linearly independent 3-subset of V over F. Take

$$T = \{(1, 3)\ (0, 2)\ (5, 1)\} \cup \{x^3, 1, x^2 + 1, x^3 + x^3 + x^5\} \cup \left\{\begin{bmatrix} 0 & 1 \\ 1 & 0 \end{bmatrix}, \begin{bmatrix} 1 & 0 \\ 0 & 0 \end{bmatrix}\right\}$$

$$= T_1 \cup T_2 \cup T_3 \subseteq V_1 \cup V_2 \cup V_3$$

is a proper 3-subset of V. Clearly T is only a semi 3-dependent 3-subset of V over F or equivalently 3-semi dependent 3-subset of V, over F for $T_1$ and $T_2$ are linearly dependent subsets of $V_1$ and $V_2$ respectively over the fields $Q(\sqrt{3})$ and $Q(\sqrt{7})$ respectively and $T_3$ is a linearly independent subset of $V_3$ over $Q(\sqrt{2})$. Take

$$P = P_1 \cup P_2 \cup P_3$$

$$= \{(1, 2)\ (2, 5), (5, 4), (-1, 0)\} \cup \{1, x^2 + x, x, x^2\} \cup \left\{\begin{bmatrix} 1 & 0 \\ 0 & 1 \end{bmatrix}, \begin{bmatrix} 0 & 1 \\ 1 & 0 \end{bmatrix}, \begin{bmatrix} 1 & 1 \\ 1 & 1 \end{bmatrix}, \begin{bmatrix} 1 & 1 \\ 1 & 0 \end{bmatrix}\right\}$$

$$\subset V_1 \cup V_2 \cup V_3$$

is a 3-subset of V. Clearly P is a 3-dependent 3-subset of V over F.

Now we have seen the notion of n-independent n-subset; n-dependent n-subset and semi n-dependent n-subset of V over F. We would proceed to define the notion of n-basis of V.

**DEFINITION 1.2.3:** *Let $V = V_1 \cup V_2 \cup \ldots \cup V_n$ be a n-vector space over the n-field $F = F_1 \cup F_2 \cup \ldots \cup F_n$. A n-subset $B = B_1 \cup B_2 \cup \ldots \cup B_n$ is a n-basis of V over F if and only if each $B_i$ is a basis of $V_i$ for every i = 1, 2, ..., n. If each basis $B_i$ of $V_i$ is finite for every i = 1, 2, ..., n then we say V is a finite n-dimensional n-vector space over the n-field F. If even one of the basis $B_i$ of $V_i$ is infinite then we call V to be an infinite dimensional n-vector space over the n-field F.*

Now we shall illustrate by an example a finite n-dimensional n-vector space over the n-field $F = F_1 \cup F_2 \cup \ldots \cup F_n$.

***Example 1.2.4:*** Let $V = V_1 \cup V_2 \cup \ldots \cup V_4$ be a 4-vector space over the field $F = Q(\sqrt{2}) \cup Z_7 \cup Z_5 \cup Q(\sqrt{3})$ where $V_1 = Q(\sqrt{2}) \times Q(\sqrt{2})$, a vector space of dimension two over $Q(\sqrt{2})$ $V_2 = Z_7 \times Z_7 \times Z_7$ a vector space of dimension 3 over $Z_7$,

$$V_3 = \left\{ \begin{bmatrix} a & b & c \\ d & e & f \end{bmatrix} \middle| a, b, c, d, e, f \in Z_5 \right\}$$

a vector space of dimension 6 over $Z_5$ and

$$V_4 = Q(\sqrt{3}) \times Q(\sqrt{3}) \times Q(\sqrt{3}) \times Q(\sqrt{3})$$

a vector space of dimension 4 over $Q(\sqrt{3})$. Clearly V is of (2, 3, 6, 4) dimension over F. Since the dimension of each $V_i$ is finite we see V is a finite dimensional 4-vector space over F.

Now we give an example of an infinite dimensional n-vector space over the n-field F.

***Example 1.2.5:*** Let $V = V_1 \cup V_2 \cup V_3$ be a 3-vector space over the 3-field $F = Q(\sqrt{3}) \cup Z_3 \cup Q(\sqrt{2})$ where $V_1 = Q(\sqrt{3})[x]$ is a vector space of infinite dimension over $Q(\sqrt{3})$, $V_2 = Z_3 \times Z_3 \times Z_3$ is a vector space of dimension 3 over $Z_3$ and

$$V_3 = \left\{ \begin{bmatrix} a & b & c \\ d & e & f \end{bmatrix} \middle| a, b, c, d, e, f \in Q\left(\sqrt{2}\right) \right\}$$

a vector space of dimension 6 over $Q(\sqrt{2})$. Clearly V is a 3-vector space of infinite dimension over F.

We now proceed on to define the simple notion of n-linear combination of n-vectors in V.

**DEFINITION 1.2.4:** *Let $V = V_1 \cup V_2 \cup \ldots \cup V_n$ be a n-vector space over the n-field $F = F_1 \cup F_2 \cup \ldots \cup F_n$. Let $\beta = \beta_1 \cup \beta_2 \cup \ldots \cup \beta_n$ be a n-vector in V where $\beta_i \in V_i$ for $i = 1, 2, \ldots, n$. We say $\beta$ is a n-linear combination of the n-vectors*

$$\alpha = \left\{\alpha_1^1, \alpha_2^1, \ldots, \alpha_{n_1}^1\right\} \cup \left\{\alpha_1^2, \alpha_2^2, \ldots, \alpha_{n_2}^2\right\} \cup \ldots \cup \left\{\alpha_1^n, \alpha_2^n, \ldots, \alpha_{n_n}^n\right\}$$

*in V, provided there exists n-scalars*

$$C = \left\{C_1^1, C_2^1, \ldots, C_{n_1}^1\right\} \cup \left\{C_1^2, C_2^2, \ldots, C_{n_2}^2\right\} \cup \ldots \cup \left\{C_1^n, C_2^n, \ldots, C_{n_n}^n\right\}$$

*such that*

$$\beta = (\beta_1 \cup \beta_2 \cup \ldots \cup \beta_n) =$$
$$\left\{C_1^1\alpha_1^1 + C_2^1\alpha_2^1 + \ldots + C_{n_1}^1\alpha_{n_1}^1\right\} \cup \left\{C_1^2\alpha_1^2 + C_2^2\alpha_2^2 + \ldots + C_{n_2}^2\alpha_{n_2}^2\right\} \cup \ldots$$
$$\cup \left\{C_1^n\alpha_1^n + C_2^n\alpha_2^n + \ldots + C_{n_n}^n\alpha_{n_n}^n\right\}.$$

We just recall this for it may be useful in case of representation of elements of a n-vector α in V, V a n-vector space over the n-field F.

**DEFINITION 1.2.5:** *Let S be a n-set of n-vectors in a n-vector space V. The n-subspace spanned by the n-set* $S = S_1 \cup S_2 \cup \ldots \cup S_n$*, where* $S_i \subseteq V_i$*, i = 1, 2, …, n is defined to be the intersection W of all n-subspaces of V which contain the n-set S. When the n-set S is finite n-set of n-vectors,*

$$S = \left\{\alpha_1^1, \alpha_2^1, \ldots, \alpha_{n_1}^1\right\} \cup \left\{\alpha_1^2, \alpha_2^2, \ldots, \alpha_{n_2}^2\right\} \cup \ldots \cup \left\{\alpha_1^n, \alpha_2^n, \ldots, \alpha_{n_n}^n\right\}$$

*we shall simply call W the n-subspace spanned by the n-vectors S.*

In view of this definition it is left as an exercise for the reader to prove the following theorem.

**THEOREM 1.2.1:** *The n-subspace spanned by a non empty n-subset* $S = S_1 \cup S_2 \cup \ldots \cup S_n$ *of a n-vector space* $V = V_1 \cup V_2 \cup \ldots \cup V_n$ *over the n-field* $F = F_1 \cup F_2 \cup \ldots \cup F_n$ *is the set of all linear combinations of n-vectors of S.*

We define now the n-sum of the n-subsets of V.

**DEFINITION 1.2.6:** *Let* $S_1, S_2, \ldots, S_k$ *be n-subsets of a n-vector space* $V = V_1 \cup V_2 \cup \ldots \cup V_n$ *where* $S_i = S_1^i \cup S_2^i \cup \ldots \cup S_n^i$ *for i = 1, 2, …, k. The set of all n-sums*

$$\left\{\alpha_1^1 + \alpha_2^1 + \ldots + \alpha_{k_1}^1\right\} \cup \left\{\alpha_1^2 + \alpha_2^2 + \ldots + \alpha_{k_2}^2\right\} \cup \ldots$$
$$\cup \left\{\alpha_1^n + \alpha_2^n + \ldots + \alpha_{k_n}^n\right\}$$

*where* $\alpha_j^i \in S_i$ *for each i = 1, 2, …, n and* $i \le j \le k_i$ *is denoted by*

$$\left(S_1^1 + S_1^2 + \ldots + S_1^{k_1}\right) \cup \left(S_2^1 + S_2^2 + \ldots + S_2^{k_2}\right) \cup \ldots \cup$$
$$\left(S_n^1 + S_n^2 + \ldots + S_n^{k_n}\right).$$

*If* $W_1, \ldots, W_k$ *are n-subspaces of* $V = V_1 \cup V_2 \cup \ldots \cup V_n$ *where* $W_i = W_1^i \cup W_2^i \cup \ldots \cup W_n^i$ *for i = 1, 2, …, n, then*

$$W = \left(W_1^1 + W_1^2 + \ldots + W_1^{k_1}\right) \cup \left(W_2^1 + W_2^2 + \ldots + W_2^{k_2}\right) \cup \ldots$$
$$\cup \left(W_n^1 + W_n^2 + \ldots + W_n^{k_n}\right)$$

*is a n-subspace of V which contains each of the subspaces* $W = W_1 \cup W_2 \cup \ldots \cup W_n$.

**THEOREM 1.2.2:** *Let* $V = V_1 \cup V_2 \cup \ldots \cup V_n$ *be an n-vector space over the n-field* $F = F_1 \cup F_2 \cup \ldots \cup F_n$. *Let V be spanned by a finite n-set of n-vectors* $\left\{\beta_1^1, \beta_2^1, \ldots, \beta_{m_1}^1\right\} \cup \left\{\beta_1^2, \beta_2^2, \ldots, \beta_{m_2}^2\right\} \cup \ldots \cup \left\{\beta_1^n, \beta_2^n, \ldots, \beta_{m_n}^n\right\}$. *Then any independent set of n-vectors in V is finite and contains no more than* $(m_1, m_2, \ldots, m_n)$ *elements.*

*Proof:* To prove the theorem it is sufficient if we prove it for one of the component spaces $V_i$ of V with associated subset $S_i$ which contains more than $m_i$ elements, this will be true for every i; i = 1, 2, …, n. Let $S = S_1 \cup S_2 \cup \ldots \cup S_n$ be a n-subset of V where each $S_i$ contains more than $m_i$ vectors. Let $S_i$ contain $n_i$ distinct vectors $\alpha_1, \alpha_2, \ldots, \alpha_{n_i}$, $n_i > m_i$, since $\left\{\beta_1^i, \beta_2^i, \ldots, \beta_{m_i}^i\right\}$ spans $V_i$, there exists scalars $A_{jk}^i$ in $F_i$ such that $\alpha_k^i = \sum_{j=1}^{m_i} A_{jk}^i \beta_k^i$.

For any $n_i$ scalars $x_1^i, x_2^i, \ldots, x_{n_i}^i$ we have $x_1^i \alpha_1^i + \ldots + x_{n_i}^i \alpha_{n_i}^i$

$$= \sum_{j=1}^{n_i} x_j^i \alpha_j^i$$

$$= \sum_{j=1}^{n_i} x_j^i \sum_{k=1}^{m_i} A_{kj}^i \beta_k^i$$

$$= \sum_{j=1}^{n_i} \sum_{k=1}^{m_i} (A_{kj}^i x_j^i) \beta_k^i$$

$$= \sum_{k=1}^{m_i} \sum_{j=1}^{n_i} (A_{kj}^i x_j^i) \beta_k^i \;.$$

Since $n_i > m_i$ this imply that there exists scalars $x_1^i, x_2^i, ..., x_{n_i}^i$ not all zero such that

$$\sum_{j=1}^{m_i} A_{kj}^i x_j^i = 0 \;;\; 1 \leq k \leq m_i.$$

Hence $x_1^i \alpha_1^i + x_2^i \alpha_2^i + ... + x_{n_i}^i \alpha_{n_i}^i = 0$. This shows that $S_i$ is a linearly dependent set. This is true for each i. Hence the result holds good for $S = S_1 \cup S_2 \cup \ldots \cup S_n$.

The reader is expected to prove the following theorems.

**THEOREM 1.2.3:** *If V is a finite dimensional n-vector space over the n-field F, then any two n-basis of V have the same number of n-elements.*

**THEOREM 1.2.4:** *Let $V = V_1 \cup V_2 \cup ... \cup V_n$ be a n-vector space over the n-field $F = F_1 \cup F_2 \cup ... \cup F_n$ and if dim V = $(n_1, n_2, ..., n_n)$, then*

1. *any n-subset of V which contains more than n-vectors is n-linearly dependent;*
2. *no n-subset of V which contains less than $(n_1, n_2, ..., n_n)$ vectors can span V.*

**THEOREM 1.2.5:** *Let $S = S_1 \cup S_2 \cup ... \cup S_n$ be a n-linearly independent n-subset of a n-vector space V. Suppose $\beta = (\beta_1 \cup \beta_2 \cup ... \cup \beta_n)$ is a vector in V which is not in the n-subspace* *spanned by S. Then the n-subset obtained by adjoining β to S is n-linearly independent.*

**THEOREM 1.2.6:** *Let $W = W_1 \cup W_2 \cup \ldots \cup W_n$ be a n-subspace of a finite dimensional n-vector space $V = V_1 \cup V_2 \cup \ldots \cup V_n$, every n-linearly independent n-subset of W is finite and is a part of a (finite) n-basis for W.*

The proof of the following corollaries is left as an exercise for the reader.

**COROLLARY 1.2.1:** *If $W = W_1 \cup W_2 \cup \ldots \cup W_n$ is a proper n-subspace of a finite dimensional n-vector space V, then W is finite dimensional and dim W < dim V; i.e. $(m_1, m_2, \ldots, m_n) < (n_1, n_2, \ldots, n_n)$ each $m_i < n_i$ for $i = 1, 2, \ldots, n$.*

**COROLLARY 1.2.2:** *In a finite dimensional n-vector space $V = V_1 \cup V_2 \cup \ldots \cup V_n$ every non-empty n-linearly independent set of n-vectors is part of a n-basis.*

**COROLLARY 1.2.3:** *Let $A = A_1 \cup A_2 \cup \ldots \cup A_n$, be a n-vector space where $A_i$ is a $n_i \times n_i$ matrix over the field $F_i$ and suppose the n-row vectors of A form a n-linearly independent set of n-vectors in $F_1^{n_1} \cup F_2^{n_2} \cup \ldots \cup F_n^{n_n}$. Then A is n-invertible i.e., each $A_i$ is invertible, $i = 1, 2, \ldots, n$.*

**THEOREM 1.2.7:** *Let $W_1 = W_1^1 \cup W_2^1 \cup \ldots \cup W_n^1$ and $W_2 = W_1^2 \cup W_2^2 \cup \ldots \cup W_n^2$ be finite dimensional n-subspaces of a n-vector space V, then*

$$W_1 + W_2 = W_1^1 + W_1^2 \cup W_2^1 + W_2^2 \cup \ldots \cup W_n^1 + W_n^2$$

*is finite dimensional and*

$$\dim W_1 + \dim W_2 = \dim(W_1 \cap W_2) + \dim(W_1 + W_2)$$

*where $\dim W_1 = (m_1^1, m_2^1, \ldots, m_{n_1}^1)$ and $\dim W_2 = (m_1^2, m_2^2, \ldots, m_{n_2}^2)$*

$$\dim W_1 + \dim W_2 = (m_1^1 + m_1^2, m_2^1 + m_2^2, \ldots, m_{n_1}^1 + m_{n_2}^2).$$

*Proof:* By the above results we have

$$W_1 \cap W_2 = W_1^1 \cap W_1^2 \cup W_2^1 \cap W_2^2 \cup \ldots \cup W_n^1 \cap W_n^2$$
has a finite n-basis.

$$\{ \alpha_1^1, \alpha_2^1, \ldots, \alpha_{k_1}^1 \} \cup \{ \alpha_1^2, \alpha_2^2, \ldots, \alpha_{k_2}^2 ) \cup \ldots \cup \{ \alpha_1^n, \alpha_2^n, \ldots, \alpha_{k_n}^n \}$$
which is part of a n-basis.

$$\{ \alpha_1^1, \alpha_2^1, \ldots, \alpha_{k_1}^1, \beta_1^1, \ldots, \beta_{n_1-k_1}^1 \} \cup \{ \alpha_1^2, \alpha_2^2, \ldots, \alpha_{k_2}^2, \beta_1^2, \beta_2^2, \ldots, \beta_{n_2-k_2}^2 \}$$
$$\cup \{ \alpha_1^n, \alpha_2^n, \ldots, \alpha_{k_n}^n, \beta_1^n, \beta_2^n, \ldots, \beta_{n_n-k_n}^n \}$$
is a n-basis of $W_1$ and

$$\{ \alpha_1^1, \alpha_2^1, \ldots, \alpha_{k_1}^1, \gamma_1^1, \gamma_2^1, \ldots, \gamma_{m_1-k_1}^1 \} \cup$$
$$\{ \alpha_1^2, \alpha_2^2, \ldots, \alpha_{k_2}^2, \gamma_1^2, \gamma_2^2, \ldots, \gamma_{m_2-k_2}^2 \} \cup \ldots \cup$$
$$\{ \alpha_1^n, \alpha_2^n, \ldots, \alpha_{k_n}^n, \gamma_1^n, \gamma_2^n, \ldots, \gamma_{m_n-k_n}^n \}$$

is a n-basis of $W_2$. Then the subspace $W_1 + W_2$ is spanned by the n-vectors

$$\{ \alpha_1^1, \ldots, \alpha_{k_1}^1, \beta_1^1, \beta_2^1, \ldots, \beta_{n_1-k_1}^n, \gamma_1^1, \gamma_2^2, \ldots, \gamma_{m_1-k_1}^1 \} \cup$$
$$\{ \alpha_1^2, \alpha_2^2, \ldots, \alpha_{k_2}^2, \beta_1^2, \beta_2^2, \ldots, \beta_{n_2-k_2}^2, \gamma_1^2, \gamma_2^2, \ldots, \gamma_{m_2-k_2}^2 \} \cup \ldots \cup$$
$$\{ \alpha_1^n, \alpha_2^n, \ldots, \alpha_{k_n}^n, \beta_1^n, \beta_2^n, \ldots, \beta_{n_n-k_n}^n, \gamma_1^n, \gamma_2^n, \ldots, \gamma_{m_n-k_n}^n \}$$

and these n-vectors form an n-independent n-set . For suppose

$$\sum x_i^k \alpha_i^k + \sum y_j^k \beta_j^k + \sum z_r^k \gamma_r^k = 0$$
true for k = 1, 2, …, n.
$$-\sum z_r^k \gamma_r^k = \sum x_i^k \alpha_i^k + \sum y_i^k \beta_i^k$$

which shows that $\sum z_r^k \gamma_r^k$ belongs $W_1^k$, true for k = 1, 2, …, n. As already $\sum z_r^k \gamma_r^k$ belongs to $W_2^k$ for k = 1, 2, …, n it follows;

$$\sum z_r^k \gamma_r^k = \sum C_i^k \alpha_i^k$$

for certain scalars $C_1^k, C_2^k, ..., C_{n_k}^k$ true for k = 1, 2, …, n.

But the n-set $\{ \alpha_1^i, \alpha_2^i, \ldots, \alpha_{k_i}^i, \gamma_1^i, \gamma_2^i, \ldots, \gamma_{n_i}^i \}$, i = 1, 2, …, n is n-independent each of the scalars $z_r^k = 0$ true for each k = 1, 2, …, n. Thus $\sum x_i^k \alpha_i^k + \sum y_j^k \beta_j^k = 0$ for k = 1, 2, …, n, and since $\{ \alpha_1^k, \alpha_2^k, \ldots, \alpha_{k_k}^k, \beta_1^k, \beta_2^k, \ldots, \beta_{n_k - k_k}^k \}$ is also an independent set each $x_i^k = 0$ and each $y_j^k = 0$ for k = 1, 2, …, n. Thus

$$\{ \alpha_1^1, \alpha_2^1, \ldots, \alpha_{k_1}^1, \beta_1^1, \beta_2^1, \ldots, \beta_{n_1 - k_1}^1, \gamma_1^1, \gamma_2^1, \ldots, \gamma_{m_1 - k_1}^1 \} \cup$$
$$\{ \alpha_1^2, \alpha_2^2, \ldots, \alpha_{k_2}^2, \beta_1^2, \beta_2^2, \ldots, \beta_{n_2 - k_2}^2, \gamma_1^2, \gamma_2^2, \ldots, \gamma_{m_2 - k_2}^2 \} \cup \ldots \cup$$
$$\{ \alpha_1^n, \alpha_2^n, \ldots, \alpha_{k_n}^n, \beta_1^n, \beta_2^n, \ldots, \beta_{n_n - k_n}^n, \gamma_1^n, \gamma_2^n, \ldots, \gamma_{m_n - k_n}^n \}$$

is a n-basis for

$$W_1 + W_2 = \{ W_1^1 + W_1^2 \} \cup \{ W_2^1 + W_2^2 \} \cup \ldots \cup \{ W_n^1 + W_n^2 \}.$$

Finally

dim $W_1$ + dim $W_2$

= $\{k_1 + m_1 + k_1 + n_1, k_2 + m_2 + k_2 + n_2, \ldots, k_n + m_n + k_n + n_n\}$

= $\{k_1 + (m_1 + k_1 + n_1), k_2 + (m_2 + k_2 + n_2), \ldots, k_n + (m_n + k_n + n_n)\}$

= $\dim(W_1 \cap W_2) + \dim(W_1 + W_2)$.

Suppose $V = V_1 \cup V_2 \cup \ldots \cup V_n$ be a $(n_1, n_2, \ldots, n_n)$ finite dimensional n-vector space over $F_1 \cup F_2 \cup \ldots \cup F_n$. Suppose

$$B = \{ \alpha_1^1, \alpha_2^1, ..., \alpha_{n_1}^1 \} \cup \{ \alpha_1^2, \alpha_2^2, ..., \alpha_{n_2}^2 \} \cup \ldots \cup \{ \alpha_1^n, \alpha_2^n, ..., \alpha_{n_n}^n \}$$

and

$$C = \{ \beta_1^1, \beta_2^1, \ldots, \beta_{n_1}^1 \} \cup \{ \beta_1^2, \beta_2^2, \ldots, \beta_{n_2}^2 \} \cup \ldots \cup \{ \beta_1^n, \beta_2^n, ..., \beta_{n_n}^n \}$$

be two ordered n-basis for V. There are unique n-scalars $P_{ij}^k$ such that $\beta_j^k = \sum_{i=1}^{n_i} P_{ij}^k \alpha_i^k$, k = 1, 2, …, n; $1 \le j \le n_i$.

Let $\{ x_1^1, x_2^1, \ldots, x_{n_1}^1 \} \cup \{ x_1^2, x_2^2, \ldots, x_{n_2}^2 \} \cup \ldots \cup \{ x_1^n, x_2^n, \ldots, x_{n_n}^n \}$ be the n-coordinates of a given n-vector α in the ordered n-basis C.

$$\alpha = \{ x_1^1\beta_1^1 + \ldots + x_{n_1}^1 \beta_{n_1}^1 \} \cup \{ x_1^2\beta_1^2 + \ldots + x_{n_2}^2 \beta_{n_2}^2 \} \cup \ldots \cup \{ x_1^n\beta_1^n + \ldots + x_{n_n}^n \beta_{n_n}^n \}$$

$$= \sum_{j=1}^{n_1} x_j^1\beta_j^1 \cup \sum_{j=1}^{n_2} x_j^2\beta_j^2 \cup \ldots \cup \sum_{j=1}^{n_n} x_j^n\beta_j^n$$

$$= \sum_{j=1}^{n_1} x_j^1 \sum_{i=1}^{n_1} P_{ij}^1\alpha_i^i \cup \sum_{j=1}^{n_2} x_j^2 \sum_{i=1}^{n_2} P_{ij}^2\alpha_i^2 \cup \ldots \cup \sum_{j=1}^{n_n} x_j^n \sum_{i=1}^{n_n} P_{ij}^n\alpha_i^n$$

$$= \sum_{j=1}^{n_1}\sum_{i=1}^{n_1} (P_{ij}^1 x_j^1)\alpha_i^1 \cup \sum_{j=1}^{n_2}\sum_{i=1}^{n_2} (P_{ij}^2 x_j^2)\alpha_i^2 \cup \ldots \cup \sum_{j=1}^{n_n}\sum_{i=1}^{n_n} (P_{ij}^n x_j^n)\alpha_i^n .$$

Thus we obtain the relation

$$\alpha = \sum_{j=1}^{n_1}\sum_{i=1}^{n_1} (P_{ij}^1 x_j^1)\alpha_i^1 \cup \sum_{j=1}^{n_2}\sum_{i=1}^{n_2} (P_{ij}^2 x_j^2)\alpha_i^2 \cup \ldots \cup \sum_{j=1}^{n_n}\sum_{i=1}^{n_n} (P_{ij}^n x_j^n)\alpha_i^n .$$

Since the n-coordinates $(y_1^1, y_2^1, \ldots, y_{n_1}^1) \cup (y_1^2, y_2^2, \ldots, y_{n_2}^2) \cup \ldots \cup (y_1^n, y_2^n, \ldots, y_{n_n}^n)$ of the n-basis B are uniquely determined it follows $y_i^k = \sum_{j=1}^{n_k} (P_{ij}^k x_j^k)$; $1 \le i \le n_k$. Let $P^i$ be the $n_i \times n_i$ matrix whose i, j entry is the scalar $P_{ij}^k$ and let X and Y be the n-coordinate matrices of the n-vector in the ordered n-basis B and C. Thus we get

$$\begin{aligned} Y &= Y^1 \cup Y^2 \cup \ldots \cup Y^n \\ &= P^1X^1 \cup P^2X^2 \cup \ldots \cup P^nX^n \\ &= PX. \end{aligned}$$

Since B and C are n-linearly independent n-sets, Y = 0 if and only if X = 0. Thus we see P is n-invertible i.e.

$P^{-1} = (P^1)^{-1} \cup (P^2)^{-1} \cup \ldots \cup (P^n)^{-1}$ i.e.

$X = P^{-1}Y$
$X^1 \cup X^2 \cup \ldots \cup X^n = (P^1)^{-1}Y^1 \cup (P^2)^{-1}Y^2 \cup \ldots \cup (P^n)^{-1} Y^n$.

This can be put with some new notational convenience as

$$[\alpha]_B = P[\alpha]_C$$
$$[\alpha]_C = P^{-1}[\alpha]_B.$$

In view of the above statements and results, we have proved the following theorem.

**THEOREM 1.2.8:** *Let $V = V_1 \cup V_2 \cup \ldots \cup V_n$ be a finite $(n_1, n_2, \ldots, n_n)$ dimensional n-vector space over the n-field F and let B and C be any two n-basis of V. Then there is a unique necessarily invertible n-matrix $P = P_1 \cup P_2 \cup \ldots \cup P_n$ of order $(n_1 \times n_1) \cup (n_2 \times n_2) \cup \ldots \cup (n_n \times n_n)$ with entries from $F_1 \cup F_2 \cup \ldots \cup F_n$ i.e. entries of $P_i$ are from $F_i$, i = 1, 2, …, n; such that*

$$[\alpha]_B = P[\alpha]_C$$
*and*
$$[\alpha]_C = P^{-1}[\alpha]_B$$

*for every n-vector α in V. The n-columns of P are given by $P_j^i = [\beta_j^i]_B$; j = 1, 2, …, n and for each i; i = 1, 2, …, n.*

Now we prove yet another interesting result.

**THEOREM 1.2.9:** *Suppose $P = P_1 \cup P_2 \cup \ldots \cup P_n$ is a $(n_1 \times n_1) \cup (n_2 \times n_2) \cup \ldots \cup (n_n \times n_n)$ n-invertible matrix over the n-field $F = F_1 \cup F_2 \cup \ldots \cup F_n$, i.e. $P_i$ takes its entries from $F_i$ true for each i = 1, 2, …, n. Let $V = V_1 \cup V_2 \cup \ldots \cup V_n$ be a finite $(n_1, n_2, \ldots, n_n)$ dimensional n-vector space over the n-field F. Let B be an n-basis of V. Then there is a unique n-basis C of V such that*

*(i)* $[\alpha]_B = P[\alpha]_C$
*and*
*(ii)* $[\alpha]_C = P^{-1}[\alpha]_B$

*for every n-vector $\alpha \in V$.*

*Proof:* Let $V = V_1 \cup V_2 \cup \ldots \cup V_n$ be a n-vector space of $(n_1, n_2, \ldots, n_n)$ dimension over the n-field $F = F_1 \cup F_2 \cup \ldots \cup F_n$. Let $B = \{(\alpha_1^1, \alpha_2^1, \ldots, \alpha_{n_1}^1) \cup (\alpha_1^2, \alpha_2^2, \ldots, \alpha_{n_2}^2) \cup \ldots \cup (\alpha_1^n, \alpha_2^n, \ldots, \alpha_{n_n}^n)\}$ be a set of n-vectors in V. If $C = \{(\beta_1^1, \beta_2^1, \ldots, \beta_{n_1}^1) \cup (\beta_1^2, \beta_2^2, \ldots, \beta_{n_2}^2) \cup \ldots \cup (\beta_1^n, \beta_2^n, \ldots, \beta_{n_n}^n)\}$ is an ordered n-basis of V for which (i) is valid, it is clear that for each $V_i$ which has $(\beta_1^i, \beta_2^i, \ldots, \beta_{n_i}^i)$ as its basis we have $\beta_j^i = \sum_{k=1}^{n_i} P_{kj} \alpha_k^i$ – I, true for every i, i = 1, 2, …,n.

Now we need only show that the vectors $\beta_j^i \in V_i$ defined by these equations I form a basis, true for each i, i = 1, 2, …, n. Let $Q = P^{-1}$. Then

$$\sum_j Q_{j_r}^i \beta_j^i = \sum_j Q_{j_r}^i \sum_k P_{kj}^i \alpha_k^i$$

(true for every i = 1, 2, …, n)

$$= \sum_j \sum_k P_{kj}^i Q_{jr}^i \alpha_k^i$$

$$= \sum_k \left( \sum_j P_{kj}^i Q_{jr}^i \right) \alpha_k^i$$

$$= \alpha_k^i$$

true for each i = 1, 2, …, n. Thus the subspace spanned by the set $\{\beta_1^i, \beta_2^i, \ldots, \beta_{n_i}^i\}$ contain $\{\alpha_1^i, \alpha_2^i, \ldots, \alpha_{n_i}^i\}$ and hence equals $V_i$; this is true for each i. Thus

$$B = \{\alpha_1^1, \alpha_2^1, \ldots, \alpha_{n_1}^1\} \cup \{\alpha_1^2, \alpha_2^2, \ldots, \alpha_{n_2}^2\} \cup \ldots \cup \{\alpha_1^n, \alpha_2^n, \ldots, \alpha_{n_n}^n\}$$

is contained in $V = V_1 \cup V_2 \cup \ldots \cup V_n$.

Thus

$$C = \{\beta_1^1, \beta_2^1, \ldots, \beta_{n_1}^1\} \cup \{\beta_1^2, \beta_2^2, \ldots, \beta_{n_2}^2\} \cup \ldots \cup \{\beta_1^n, \beta_2^n, \ldots, \beta_{n_n}^n\}$$

is a n-basis and from its definition and by the above theorem (i) is valid hence also (ii).

Now we proceed onto define the notion of n-linear transformation for n-vector spaces of type II.

**DEFINITION 1.2.7:** *Let $V = V_1 \cup V_2 \cup \ldots \cup V_n$ be a n-vector space of type II over the n-field $F = F_1 \cup F_2 \cup \ldots \cup F_n$ and $W = W_1 \cup W_2 \cup \ldots \cup W_n$ be a n-space over the same n-field $F = F_1 \cup F_2 \cup \ldots \cup F_n$. A n-linear transformation of type II is a n-function $T = T_1 \cup T_2 \cup \ldots \cup T_n$ from V into W such that $T_i(C^i\alpha^i + \beta^i) = C^iT_i\alpha^i + T_i\beta^i$ for all $\alpha^i$, $\beta^i$ in $V_i$ and for all scalars $C^i \in F_i$. This is true for each i; i = 1, 2, …, n. Thus*

$$
\begin{aligned}
T[C\alpha + \beta] &= (T_1 \cup T_2 \cup \ldots \cup T_n) \\
&\quad [C^1\alpha^1 + \beta^1] \cup [C^2\alpha^2 + \beta^2] \cup \ldots \cup [C^n\alpha^n + \beta^n] \\
&= T_1(C^1\alpha^1 + \beta^1) \cup T_2(C^2\alpha^2 + \beta^2) \cup \ldots \cup \\
&\quad T_n(C^n\alpha^n + \beta^n)
\end{aligned}
$$

*for all $\alpha^i \in V_i$ or $\alpha^1 \cup \alpha^2 \cup \ldots \cup \alpha^n \in V$ and $C^i \in F_i$ or $C^1 \cup C^2 \cup \ldots \cup C^n \in F_1 \cup F_2 \cup \ldots \cup F_n$. We say $I = I_1 \cup I_2 \cup \ldots \cup I_n$ is the n-identity transformation of type II from V to V if*

$$
\begin{aligned}
I\alpha &= I(\alpha^1 \cup \alpha^2 \cup \ldots \cup \alpha^n) \\
&= I_1(\alpha^1) \cup I_2(\alpha^2) \cup \ldots \cup I_n(\alpha^n) \\
&= \alpha^1 \cup \alpha^2 \cup \ldots \cup \alpha^n
\end{aligned}
$$

*for all $\alpha^1 \cup \alpha^2 \cup \ldots \cup \alpha^n \in V = V_1 \cup V_2 \cup \ldots \cup V_n$. Similarly the n-zero transformation $0 = 0 \cup 0 \cup \ldots \cup 0$ of type II from V into V is given by $0\alpha = 0\alpha^1 \cup 0\alpha^2 \cup \ldots \cup 0\alpha^n$ $= 0 \cup 0 \cup \ldots \cup 0$, for all $\alpha^1 \cup \alpha^2 \cup \ldots \cup \alpha^n \in V$.*

Now we sketch the proof of the following interesting and important theorem.

**THEOREM 1.2.10:** *Let $V = V_1 \cup V_2 \cup \ldots \cup V_n$ be a $(n_1, n_2, \ldots, n_n)$ dimensional finite n-vector space over the n-field $F = F_1 \cup F_2 \cup \ldots \cup F_n$ of type II. Let $\{\alpha_1^1, \alpha_2^1, \ldots, \alpha_{n_1}^1\} \cup \{\alpha_1^2, \alpha_2^2, \ldots, \alpha_{n_2}^2\} \cup \ldots \cup \{\alpha_1^n, \alpha_2^n, \ldots, \alpha_{n_n}^n\}$ be an n-ordered basis for V. Let W =*

*$W_1 \cup W_2 \cup \ldots \cup W_n$ be a n-vector space over the same n-field F = $F_1 \cup F_2 \cup \ldots \cup F_n$.*

*Let $\{\beta_1^1, \beta_2^1, \ldots, \beta_{n_1}^1\} \cup \{\beta_1^2, \beta_2^2, \ldots, \beta_{n_2}^2\} \cup \ldots \cup \{\beta_1^n, \beta_2^n, \ldots, \beta_{n_n}^n\}$ be any n-vector in W. Then there is precisely a n-linear transformation T = $T_1 \cup T_2 \cup \ldots \cup T_n$ from V into W such that $T_i\alpha_j^i = \beta_j^i$ for i = 1, 2, …, n and j = 1, 2, …, $n_i$.*

*Proof:* To prove there is some n-linear transformation T = $T_1 \cup T_2 \cup \ldots \cup T_n$ with $T_i\alpha_j^i = \beta_j^i$ for each j = 1, 2, …, $n_i$ and for each i = 1, 2, …, n. For every $\alpha = \alpha^1 \cup \alpha^2 \cup \ldots \cup \alpha^n$ in V we have for every $\alpha^i \in V_i$ a unique $x_1^i, x_2^i, \ldots, x_{n_i}^n$ such that α = $x_1^i\alpha_2^i + \ldots + x_{n_i}^i\alpha_{n_i}^i$. This is true for every i; i = 1, 2, …, n. For this vector $\alpha^i$ we define $T_i\alpha^i$ = $x_1^i\beta_1^i + \ldots + x_{n_i}^i\beta_{n_i}^i$ true for i = 1, 2, …, n. Thus $T_i$ is well defined for associating with each vector $\alpha^i$ in $V_i$ a vector $T_i\alpha^i$ in $W_i$ this is well defined rule for T = $T_1 \cup T_2 \cup \ldots \cup T_n$ as it is well defined rule for each $T_i$: $V_i \to W_i$ , i = 1, 2, …, n.

From the definition it is clear that $T_i\alpha_j^i = \beta_j^i$ for each j. To see that T is n-linear. Let $\beta^i$ = $y_1^i\alpha_2^i + \ldots + y_{n_i}^i\alpha_{n_i}^i$ be in V and let $C^i$ be any scalar from $F_i$. Now

$$C^i\alpha^i + \beta^i = (C^ix_1^i + y_1^i)\beta_1^i + \ldots + (C^ix_{n_i}^i + y_{n_i}^i)\beta_{n_i}^i$$

true for every i, i = 1, 2, …, n. On the other hand

$$T_i(C^i\alpha^i + \beta^i) = C^i\sum_{j=1}^{n_i} x_j^i\beta_j^i + \sum_{j=1}^{n_i} y_j^i\beta_j^i$$

true for i = 1, 2, …, n i.e. true for every linear transformation $T_i$ in T. $T_i(C^i\alpha^i + \beta^i)$ = $C^iT_i(\alpha^i) + T_i(\beta^i)$ true for every i. Thus $T(C\alpha + \beta) = T_1(C^1\alpha^1 + \beta^1) \cup T_2(C^2\alpha^2 + \beta^2) \cup \ldots \cup T_n(C^n\alpha^n + \beta^n)$. If U = $U_1 \cup U_2 \cup \ldots \cup U_n$ is a n-linear transformation from V into W with $U_i\alpha_j^i = \beta_j^i$ for j = 1, 2, …, $n_i$ and i = 1, 2, …, n then for the n-vector $\alpha = \alpha^1 \cup \alpha^2 \cup \ldots \cup \alpha^n$ we have for every $\alpha^i$ in α

$$\alpha^i = \sum_{j=1}^{n_i} x_j^i \alpha_j^i$$

we have

$$U_i\alpha^i = U_i \sum_{j=1}^{n_i} x_j^i \alpha_j^i$$

$$= \sum_{j=1}^{n_i} x_j^i (U_i \alpha_j^i)$$

$$= \sum_{j=1}^{n_i} x_j^i \beta_j^i$$

so that U is exactly the rule T which we define. This proves $T\alpha = \beta$, i.e. if $\alpha = \alpha^1 \cup \alpha^2 \cup \ldots \cup \alpha^n$ and $\beta = \beta^1 \cup \beta^2 \cup \ldots \cup \beta^n$ then $T_i \alpha_j^i = \beta_j^i$, $i \leq j \leq n_i$ and i = 1, 2, …, n.

Now we proceed on to give a more explicit description of the n-linear transformation. For this we first recall if $V_i$ is any $n_i$ dimensional vector space over $F_i$ then $V_i \cong (F_i)^{n_i} = F_i^{n_i}$. Further if $W_i$ is also a vector space of dimension $m_i$ over $F_i$, the same field then $W_i \cong (F_i)^{m_i}$.

Now let $V = V_1 \cup V_2 \cup \ldots \cup V_n$ be a n-vector space over the n-field $F = F_1 \cup F_2 \cup \ldots \cup F_n$ of $(n_1, n_2, \ldots, n_n)$ dimension over F, i.e., each $V_i$ is a vector space over $F_i$ of dimension $n_i$ over $F_i$ for i = 1, 2, …, n. Thus $V_i \cong F_i$. Hence $V = V_1 \cup V_2 \cup \ldots \cup V_n \cong F_1^{n_1} \cup F_2^{n_2} \cup \ldots \cup F_n^{n_n}$. Similarly if $W = W_1 \cup W_2 \cup \ldots \cup W_n$ is a n-vector space over the same n-field F and $(m_1, m_2, \ldots, m_n)$ is the dimension of W over F then $W_1 \cup W_2 \cup \ldots \cup W_n \cong F_1^{m_1} \cup F_2^{m_2} \cup \ldots \cup F_n^{m_n}$.

$T = T_1 \cup T_2 \cup \ldots \cup T_n$ is uniquely determined by a sequence of n-linear n-vectors

$$(\beta_1^1, \beta_2^1, \ldots, \beta_{n_1}^1) \cup (\beta_1^2, \beta_2^2, \ldots, \beta_{n_2}^2) \cup \ldots \cup (\beta_1^n, \beta_2^n, \ldots, \beta_{n_n}^n)$$

where $\beta_j^i = T_i E_j$, j = 1, 2, …, $n_i$ and i = 1, 2, …, n. In short T is uniquely determined by the n-images of the standard n-basis vectors. The determinate is

$$\alpha = \{(x_1^1, x_2^1, \ldots, x_{n_1}^1) \cup (x_1^2, x_2^2, \ldots, x_{n_2}^2) \cup \ldots \cup (x_1^n, x_2^n, \ldots, x_{n_n}^n)\}.$$

$$T\alpha = (x_1^1\beta_1^1 + \ldots + x_{n_1}^1\beta_{n_1}^1) \cup (x_1^2\beta_1^2 + \ldots + x_{n_2}^2\beta_{n_2}^2) \cup \ldots \cup (x_1^n\beta_1^n + \ldots + x_{n_n}^n\beta_{n_n}^n).$$

If $B = B^1 \cup B^2 \cup \ldots \cup B^n$ is a $(n_1 \times n_1, n_2 \times n_2, \ldots, n_n \times n_n)$, n-matrix which has n-row vector

$$(\beta_1^1, \beta_2^1, \ldots, \beta_{n_1}^1) \cup (\beta_1^2, \beta_2^2, \ldots, \beta_{n_2}^2) \cup \ldots \cup (\beta_1^n, \beta_2^n, \ldots, \beta_{n_n}^n)$$

then $T\alpha = \alpha B$. In other words if $\beta_k^i = (B_{k1}^i, B_{k2}^i, \ldots, B_{kn_i}^i)$ ; i = 1, 2, …, n, then

$$T\{(x_1^1, x_2^1, \ldots, x_{n_1}^1) \cup \ldots \cup (x_1^n, x_2^n, \ldots, x_{n_n}^n)\}$$
$$= T_1(x_1^1, x_2^1, \ldots, x_{n_1}^1) \cup T_2(x_1^2, x_2^2, \ldots, x_{n_2}^2) \cup \ldots \cup T_n(x_1^n, x_2^n, \ldots, x_{n_n}^n) =$$

$$[x_1^1, x_2^1, \ldots, x_{n_1}^1] \begin{bmatrix} B_{11}^1 & \cdots & B_{1n_1}^1 \\ \vdots & & \vdots \\ B_{m_1 1}^1 & \cdots & B_{m_1 n_1}^1 \end{bmatrix} \cup$$

$$[x_1^2, x_2^2, \ldots, x_{n_2}^2] \begin{bmatrix} B_{11}^2 & \cdots & B_{1n_2}^2 \\ \vdots & & \vdots \\ B_{m_2 1}^2 & \cdots & B_{m_2 n_2}^2 \end{bmatrix} \cup \ldots \cup$$

$$[x_1^n, x_2^n, \ldots, x_{n_n}^n] \begin{bmatrix} B_{11}^n & \cdots & B_{1n_n}^n \\ \vdots & & \vdots \\ B_{m_n 1}^n & \cdots & B_{m_n n_n}^n \end{bmatrix}.$$

This is the very explicit description of the n-linear transformation. If T is a n-linear transformation from V into W then n-range of T is not only a n-subset of W it a n-subspace of

$W = W_1 \cup W_2 \cup \ldots \cup W_n$. Let $R_T = R^1_{T_1} \cup R^2_{T_2} \cup \ldots \cup R^n_{T_n}$ be the n-range of $T = T_1 \cup T_2 \cup \ldots \cup T_n$ that is the set of all n-vectors $\beta = (\beta^1 \cup \beta^2 \cup \ldots \cup \beta^n)$ in $W = W_1 \cup W_2 \cup \ldots \cup W_n$ such that $B = T\alpha$ i.e., $\beta^i_j = T_i\alpha^i_j$ for each i = 1, 2, …, n and for some $\alpha$ in V. Let $\beta^i_1, \beta^i_2 \in R^i_{T_i}$ and $C^i \in F_i$. There are vectors $\alpha^i_1, \alpha^i_2 \in V_i$ such that $T_i\alpha^i_1 = \beta^i_1$ and $T_i\alpha^i_2 = \beta^i_2$. Since $T_i$ is linear for each i; $T_i(C^i\alpha^i_1 + \alpha^i_2)$ = $C^iT_i\alpha^i_1 + T_i\alpha^i_2$ = $C^i\beta^i_1 + \beta^i_2$ which shows that $C^i\beta^i_1 + \beta^i_2$ is also in $R^i_{T_i}$. Since this is true for every i we have

$(C^1\beta^1_1 + \beta^1_2) \cup \ldots \cup (C^n\beta^n_1 + \beta^n_2) \in R_T = R^1_{T_1} \cup R^2_{T_2} \cup \ldots \cup R^n_{T_n}$.

Another interesting n-subspace associated with the n-linear transformation T is the n-set $N = N_1 \cup N_2 \cup \ldots \cup N_n$ consisting of the n-vectors $\alpha \in V$ such that $T\alpha = 0$. It is a n-subspace of V because

(1) $T(0) = 0$ so $N = N_1 \cup N_2 \cup \ldots \cup N_n$ is non empty. If $T\alpha = T\beta = 0$ then $\alpha, \beta \in V = V_1 \cup V_2 \cup \ldots \cup V_n$. i.e., $\alpha = \alpha^1 \cup \alpha^2 \cup \ldots \cup \alpha^n$ and $\beta = \beta^1 \cup \beta^2 \cup \ldots \cup \beta^n$ then

$$\begin{aligned} T(C\alpha + \beta) &= CT\alpha + T\beta \\ &= C.0 + 0 \\ &= 0 \\ &= 0 \cup 0 \cup \ldots \cup 0 \end{aligned}$$

so $C\alpha + \beta \in N = N_1 \cup N_2 \cup \ldots \cup N_n$.

**DEFINITION 1.2.8:** *Let V and W be two n-vector spaces over the same n-field $F = F_1 \cup F_2 \cup \ldots \cup F_n$ of dimension $(n_1, n_2, \ldots, n_n)$ and $(m_1, m_2, \ldots, m_n)$ respectively. Let $T: V \to W$ be a n-linear transformation. The n-null space of T is the set of all n-vectors $\alpha = \alpha^1 \cup \alpha^2 \cup \ldots \cup \alpha^n$ in V such that $T_i\alpha^i = 0$, i = 1, 2, …, n.*

*If V is finite dimensional, the n-rank of T is the dimension of the n-range of T and the n-nullity of T is the dimension of the n-null space of T.*

The following interesting result on the relation between n-rank T and n-nullity T is as follows:

**THEOREM 1.2.11:** *Let V and W be n-vector spaces over the same n-field $F = F_1 \cup F_2 \cup \ldots \cup F_n$ and suppose V is finite say $(n_1, n_2, \ldots, n_n)$ dimensional. T is a n-linear transformation from V into W. Then n-rank T + n-nullity T = n-dim V = $(n_1, n_2, \ldots, n_n)$.*

*Proof:* Given $V = V_1 \cup V_2 \cup \ldots \cup V_n$ to be a $(n_1, n_2, \ldots, n_n)$ dimensional n-vector space over the n-field $F = F_1 \cup F_2 \cup \ldots \cup F_n$. $W = W_1 \cup W_2 \cup \ldots \cup W_n$ is a n-vector space over the same n-field F. Let T be a n-linear transformation from V into W given by $T = T_1 \cup T_2 \cup \ldots \cup T_n$ where $T_i: V_i \to W_i$ is a linear transformation, and $V_i$ is of dimension $n_i$ and both the vector spaces $V_i$ and $W_i$ are defined over $F_i$; i = 1, 2, …, n.
n-rank T = rank $T_1 \cup$ rank $T_2 \cup \ldots \cup$ rank $T_n$; n-nullity T = nullity $T_1 \cup$ nullity $T_2 \cup \ldots \cup$ nullity $T_n$. So n-rank T + n-nullity T = dim V = $(n_1, n_2, \ldots, n_n)$ i.e. rank $T_i$ + nullity $T_i$ = dim $V_i = n_i$.

We shall prove the result for one $T_i$ and it is true for every i, i = 1, 2, …, n. Suppose $\{\alpha_1^i, \alpha_2^i, \ldots, \alpha_{k_i}^i\}$ is a basis of $N_i$; the null space of $T_i$. There are vectors $\alpha_{k_i+1}^i, \alpha_{k_i+2}^i, \ldots, \alpha_{n_i}^i$ in $V_i$ such that $\{\alpha_1^i, \alpha_2^i, \ldots, \alpha_{n_i}^i\}$ is a basis of $V_i$, true for every i, i = 1, 2, …, n.

We shall now prove that $\{T_i\alpha_{k_i+1}^i, \ldots, T_i\alpha_{n_i}^i\}$ is a basis for the range of $T_i$. The vector $T_i\alpha_1^i, T_i\alpha_2^i, \ldots, T_i\alpha_{n_i}^i$ certainly span the range of $T_i$ and since $T_i\alpha_j^i = 0$ for $j \leq k_i$ we see that $\{T_i\alpha_{k_i+1}^i, \ldots, T_i\alpha_{n_i}^i\}$ span the range. To see that these vectors are independent, suppose we have scalar $C_r^i$ such that

$$\sum_{r=k_i+1}^{n_i} C_r^i T_i \alpha_r^i = 0 .$$

This says that

$$T_i\left(\sum_{j=k_i+1}^{n_i} C_r^i \alpha_r^i\right) = 0$$

and accordingly the vector $\alpha^i = \sum_{j=k_i+1}^{n_i} C_r^i \alpha_r^i$ is in the null space of $T_i$. Since $\alpha_1^i, \alpha_2^i, ..., \alpha_{k_i}^i$ form a basis for $N_i$, there must be scalar $\beta_1^i, \beta_2^i, ..., \beta_{k_i}^i$ in $F_i$ such that $\alpha^i = \sum_{r=1}^{k_i} \beta_n^i \alpha_r^i$. Thus

$$\sum_{r=1}^{k_i} \beta_n^i \alpha_r^i - \sum_{j=k_i+1}^{n_i} C_j^i \alpha_j^i = 0$$

and since $\alpha_1^i, \alpha_2^i, ..., \alpha_{n_i}^i$ are linearly independent we must have $\beta_1^i = \beta_2^i = \ldots = \beta_{k_i}^i = C_{k_i+1}^i = \ldots = C_{n_i}^i = 0$. If $r_i$ is the rank of $T_i$, the fact that $T_i\alpha_{k_i+1}, \ldots, T_i\alpha_{n_i}$ form a basis for the range of $T_i$ tells us that $r_i = n_i - k_i$. Since $k_i$ is the nullity of $T_i$ and $n_i$ is the dimension of $V_i$ we have rank $T_i$ + nullity $T_i$ = dim $V_i = n_i$. We see this above equality is true for every i, i = 1, 2, …, n. We have rank $T_i$ + nullity $T_i$ = dim$V_i$ = $n_i$; we see this above equality is true for every i, i = 1, 2, …, n. We have (rank $T_1$ + nullity $T_1$) ∪ (rank $T_2$ + nullity $T_2$) ∪ … ∪ (rank $T_n$ + nullity $T_n$) = dim V = $(n_1, n_2, \ldots, n_n)$.

Now we proceed on to introduce and study the algebra of n-linear transformations.

**THEOREM 1.2.12:** *Let V and W be a n-vector spaces over the same n-field, $F = F_1 \cup F_2 \cup \ldots \cup F_n$. Let T and U be n-linear transformation from V into W. The function (T + U) is defined by (T + U)α = Tα + Uα is a n-linear transformation from V into W. If C is a n-scalar from $F = F_1 \cup F_2 \cup \ldots \cup F_n$, the n-function CT is defined by (CT)α = C(Tα) is a n-linear transformation from V into W. The set of all n-linear transformations form V into W together with addition and scalar multiplication defined above is a n-vector space over the same n-field $F = F_1 \cup F_2 \cup \ldots \cup F_n$.*

*Proof:* Given $V = V_1 \cup V_2 \cup \ldots \cup V_n$ and $W = W_1 \cup W_2 \cup \ldots \cup W_n$ are two n-vector spaces over the same n-field $F = F_1 \cup F_2$ $\cup \ldots \cup F_n$ . Let $U = U_1 \cup U_2 \cup \ldots \cup U_n$ and $T = T_1 \cup T_2 \cup \ldots \cup T_n$ be two n-linear transformations from V into W, to show the function (T + U) defined by $(T + U)(\alpha) = T\alpha + U\alpha$ is a n-linear transformation from V into W.

For $\alpha = \alpha^1 \cup \alpha^2 \cup \ldots \cup \alpha^n \in V = V_1 \cup V_2 \cup \ldots \cup V_n$, $C = C^1 \cup C^2 \cup \ldots \cup C^n$; consider

$$\begin{aligned}
&(T + U)(C\alpha + \beta) \\
&\quad = T(C\alpha + \beta) + U(C\alpha + \beta) \\
&\quad = CT(\alpha) + T\beta + CU(\alpha) + U(\beta) \\
&\quad = CT(\alpha) + CU(\alpha) + T\beta + U\beta \\
&\quad = C(T + U)\alpha + (T + U)\beta
\end{aligned}$$

which shows T + U is a n-linear transformation of V into W given by

$$\begin{aligned}
T + U &= (T_1 \cup T_2 \cup \ldots \cup T_n) + (U_1 \cup U_2 \cup \ldots \cup U_n) \\
&= (T_1 + U_1) \cup (T_2 + U_2) \cup \ldots \cup (T_n + U_n)
\end{aligned}$$

since each $(T_i + U_i)$ is a linear transformation from $V_i$ into $W_i$ true for each i, i = 1, 2, …, m; hence we see T + U is a n-linear transformation from V into W.

Similarly CT is also a n-linear transformation. For

$$\begin{aligned}
CT &= (C^1 \cup C^2 \cup \ldots \cup C^n)\ (T_1 \cup T_2 \cup \ldots \cup T_n) \\
&= C^1T_1 \cup C^2T_2 \cup \ldots \cup C^nT_n.
\end{aligned}$$

Now for $(d\alpha + \beta)$ where $d = d^1 \cup d^2 \cup \ldots \cup d^n$ and $d \in F_1 \cup F_2 \cup \ldots \cup F_n$ i.e., $d^i \in F_i$ for i = 1, 2, …, n and $\alpha = \alpha^1 \cup \alpha^2 \cup \ldots \cup \alpha^n$ and $\beta = \beta^1 \cup \beta^2 \cup \ldots \cup \beta^n$ consider

$$\begin{aligned}
T(d\alpha + \beta) &= (T_1 \cup T_2 \cup \ldots \cup T_n)\ \{[d^1\alpha^1 \cup d^2\alpha^2 \cup \ldots \cup d^n\alpha^n] \\
&\qquad + (\beta^1 \cup \beta^2 \cup \ldots \cup \beta^n)\} \\
&= T_1(d^1\alpha^1 + \beta^1) \cup T_2(d^2\alpha^2 + \beta^2) \cup \ldots \cup T_n(d^n\alpha^n + \\
&\qquad \beta^n)
\end{aligned}$$

now

$$\begin{aligned}
CT &= (C^1 \cup C^2 \cup \ldots \cup C^n)(T_1 \cup T_2 \cup \ldots \cup T_n) \\
&= C^1T_1 \cup C^2T_2 \cup \ldots \cup C^nT_n
\end{aligned}$$

(we know from properties of linear transformation each $C^iT_i$ is a linear transformation for i = 1, 2, …, n) so $CT = C^1T_1 \cup C^2T_2 \cup \ldots \cup C^nT_n$ is a n-linear transformation as

$$CT(d\alpha + \beta) = C^1T_1(d^1\alpha^1 + \beta^1) \cup C^2T_2(d^2\alpha^2 + \beta^2) \cup \ldots \cup C^nT_n(d^n\alpha^n + \beta^n).$$

Hence the claim. CT is a n-linear transformation from V into W.

Now we prove about the properties enjoyed by the collection of all n-linear transformation of V into W. Let $L^n$ (V, W) denote the collection of all linear transformations of V into W, to prove $L^n$ (V, W) is a n-vector space over the n-field; $F = F_1 \cup F_2 \cup \ldots \cup F_n$, where V and W are n-vector spaces defined over the same n-field F. Just now we have proved $L^n(V, W)$ is closed under sum i.e. addition and also $L^n$ (V, W) is closed under the n-scalar from the n-field i.e. we have proved by defining $(T + U)(\alpha) = T\alpha + U\alpha$ for all $\alpha = \alpha^1 \cup \alpha^2 \cup \ldots \cup \alpha^n \in V$, T + U is again a n-linear transformation from V into W. i.e. $L^n(V, W)$ is closed under addition. We have also proved for every n-scalar $C = C^1 \cup C^2 \cup \ldots \cup C^n \in F = F_1 \cup F_2 \cup \ldots \cup F_n$ and $T = T_1 \cup T_2 \cup \ldots \cup T_n$; CT is also a n-linear transformation of V into W; i.e. for every T, U $\in L^n$ (V,W), T + U $\in L^n(V,W)$ and for every $C \in F = F_1 \cup F_2 \cup \ldots \cup F_n$ and for every T $\in L^n(V, W)$, CT $\in L^n(V, W)$. Trivially $0\alpha = 0$ for every $\alpha \in V$ will serve as the n-zero transformation of V into W.

Thus $L^n$ (V, W) is a n-vector space over the same n-field $F = F_1 \cup F_2 \cup \ldots \cup F_n$.

Now we study the dimension of $L^n$ (V, W), when both the n-vector spaces are of finite dimension.

**THEOREM 1.2.13:** *Let $V = V_1 \cup V_2 \cup \ldots \cup V_n$ be a finite $(n_1, n_2, \ldots, n_n)$-dimensional n-vector space over the n-field $F = F_1 \cup F_2 \cup \ldots \cup F_n$. Let $W = W_1 \cup W_2 \cup \ldots \cup W_n$ be a finite $(m_1, m_2, \ldots, m_n)$-dimensional n-vector space over the same n-field F. Then the space $L^n$ (V, W) is of finite dimension and has $(m_1n_1, m_2n_2, \ldots, m_nn_n)$-dimension over the same n-field F.*

*Proof:* Given $V = V_1 \cup V_2 \cup \ldots \cup V_n$ is a n-vector space over the n-field $F = F_1 \cup F_2 \cup \ldots \cup F_n$ of $(n_1, n_2, \ldots, n_n)$-dimension over F. Also $W = W_1 \cup W_2 \cup \ldots \cup W_n$ is a n-vector space over the n-field F of $(m_1, m_2, \ldots, m_n)$ dimension over F. Let $T = T_1 \cup T_2 \cup \ldots \cup T_n$ be any n-linear transformation of type II i.e. $T_i: V_i \to W_i$, this is the only way the n-linear transformation of type II is defined because $V_i$ is defined over $F_i$ and $W_i$ is also defined over $F_i$. Clearly $W_j$ ($i \neq j$) is defined over $F_j$ and $F_i \neq F_j$ so one cannot imagine of defining any n-linear transformation. Let B = $\{\alpha_1^1, \alpha_2^1, \ldots, \alpha_{n_1}^1\} \cup \{\alpha_1^2, \alpha_2^2, \ldots, \alpha_{n_2}^2\} \cup \ldots \cup \{\alpha_1^n, \alpha_2^n, \ldots, \alpha_{n_n}^n\}$ be a n-ordered basis for V over the n-field F. We say n-ordered basis if each basis of $V_i$ is ordered and $B^1 = \{\beta_1^1, \beta_2^1, \ldots, \beta_{m_1}^1\} \cup \{\beta_1^2, \beta_2^2, \ldots, \beta_{m_2}^2\} \cup \ldots \cup \{\beta_1^n, \beta_2^n, \ldots, \beta_{m_n}^n\}$ be a n-ordered basis of W over the n-field F.

For every pair of integers $(p^i, q^i)$; $1 \leq p^i \leq m_i$ and $1 \leq q^i \leq m_i$, $i = 1, 2, \ldots, n$, we define a n-linear transformation $E^{p,q} = E_1^{p^1,q^1} \cup E_2^{p^2,q^2} \cup \ldots \cup E_n^{p^n,q^n}$ from V into W by

$$E_i^{p^i,q^i}(\alpha_j^i) = \begin{cases} 0 \text{ if } j \neq q^i \\ \beta_{p^i}^i \text{ if } j = q^i \end{cases} = \delta_{jq^i}\beta_{p^i}^i,$$

true for $i = 1, 2, \ldots, n_i$, $i = 1, 2, \ldots, m_i$. We have by earlier theorems a unique n-linear transformation from $V_i$ into $W_i$ satisfying these conditions. The claim is that the $m_i n_i$ transformations $E_i^{p^i,q^i}$ form a basis of L $(V_i, W_i)$. This is true for each i. So $L^n (V, W) = L (V_1, W_1) \cup L (V_2, W_2) \cup \ldots \cup L (V_n, W_n)$ is a n-vector space over the n-field of $(m_1n_1, m_2n_2, \ldots, m_nn_n)$ dimension over F. Now suppose $T = T_1 \cup T_2 \cup \ldots \cup T_n$ is a n-linear transformation from V into W.

Now for each j, $1 \leq j \leq n_i$, $i = 1, 2, \ldots, n$ let $A_{kj}^i, \ldots, A_{m_ij}^i$ be the co ordinates of the vector $T_i\alpha_j^i$ in the ordered basis $(\beta_1^i, \beta_2^i, \ldots, \beta_{m_i}^i)$ of $B^i$, $i = 1, 2, \ldots, n$. i.e.

$$T_i\alpha^i_j = \sum_{p^i=1}^{m_i} A^i_{p^ij}\beta^i_{p^i}\ .$$

We wish to show that $T = \sum_{p^i=1}^{m_i}\sum_{q^i=1}^{n_i} A^i_{p^iq^i} E_i^{p^i,q^i}$ (I).

Let U be a n-linear transformation in the right hand member of (I) then for each j

$$\begin{aligned} U_i\alpha^i_j &= \sum_{p^i}\sum_{q^i} A^i_{p^iq^i} E_i^{p^i,q^i}(\alpha^i_j) \\ &= \sum_{p^i}\sum_{q^i} A^i_{p^iq^i}\delta_{jq^i}\beta^i_{p^i} \\ &= \sum_{p^i=1}^{m_i} A^i_{p^i_j}\beta^i_{p^i} \\ &= T_i\,\alpha^i_j \end{aligned}$$

and consequently $U_i = T_i$ as this is true for each i, i = 1, 2, …, n. We see T = U. Now I shows that the $E^{p,q}$ spans $L^n$ (V, W) as each $E_i^{p^i,q^i}$ spans $L(V_i, W_i)$; i = 1, 2, …, n. Now it remains to show that they are n-linearly independent. But this is clear from what we did above for if in the n-transformation U, we have

$$U_i = \sum_{p^i}\sum_{q^i} A^i_{p^iq^i} E_i^{p^i,q^i}$$

to be the zero transformation then $U_i\alpha^i_j = 0$ for each j, $1 \leq j \leq n_i$ so $\sum_{p^i=1}^{m_i} A^i_{p^i_j}\beta^i_{p^i} = 0$ and the independence of the $\beta^i_{p^i}$ implies that $A^i_{p^i_j} = 0$ for every $p^i$ and j. This is true for every i = 1, 2, …, n.

Next we proceed on to prove yet another interesting result about n-linear transformations, on the n-vector spaces of type II.

**THEOREM 1.2.14:** *Let V, W and Z be n-vector spaces over the same n-field $F = F_1 \cup F_2 \cup \ldots \cup F_n$. Let T be a n-linear transformation from V into W and U a n-linear transformation* *from W into Z. Then the n-composed n-function UT defined by UT(α) = U(T(α)) is a n-linear transformation from V into Z.*

*Proof:* Given $V = V_1 \cup V_2 \cup \ldots \cup V_n$ , $W = W_1 \cup W_2 \cup \ldots \cup W_n$ and $Z = Z_1 \cup Z_2 \cup \ldots \cup Z_n$ are n-vector spaces over the n-field $F = F_1 \cup F_2 \cup \ldots \cup F_n$ where $V_i$, $W_i$ and $Z_i$ are vector spaces defined over the same field $F_i$, this is true for i = 1, 2, …, n. Let $T = T_1 \cup T_2 \cup \ldots \cup T_n$ be a n-linear transformation from V into W i.e. each $T_i$ maps $V_i$ into $W_i$ for i = 1, 2, …, n. In no other way there can be a n-linear transformation of type II from V into W. $U = U_1 \cup U_2 \cup \ldots \cup U_n$ is a n-linear transformation of type II from W into Z such that for every $U_i$ is a linear transformation from $W_i$ into $Z_i$; (for $W_i$ and $Z_i$ alone are vector spaces defined over the same field $F_i$); this is true for each i, i = 1, 2, …, n. Suppose the n-composed function UT defined by (UT)(α) = U(T(α)) where $\alpha = \alpha^1 \cup \alpha^2 \cup \ldots \cup \alpha^n \in V$ i.e.

$$
\begin{aligned}
UT &= U_1T_1 \cup U_2T_2 \cup \ldots \cup U_nT_n \text{ and} \\
UT(\alpha) &= (U_1T_1 \cup U_2T_2 \cup \ldots \cup U_nT_n)(\alpha^1 \cup \alpha^2 \cup \ldots \cup \alpha^n) \\
&= (U_1T_1)(\alpha^1) \cup (U_2T_2)(\alpha^2) \cup \ldots \cup (U_nT_n)(\alpha^n).
\end{aligned}
$$

Now from results on linear transformations we know each $U_iT_i$ defined by $(U_iT_i)(\alpha^i) = U_i(T_i(\alpha^i)$ is a linear transformation from $V_i$ into $Z_i$; This is true for each i, i = 1, 2, …, n. Hence UT is a n-linear transformation of type II from V into Z.

We now define a n-linear operator of type II for n-vector spaces V over the n-field F.

**DEFINITION 1.2.9:** *Let $V = V_1 \cup V_2 \cup \ldots \cup V_n$ be a n-vector space defined over the n-field $F = F_1 \cup F_2 \cup \ldots \cup F_n$. A n-linear operator $T = T_1 \cup T_2 \cup \ldots \cup T_n$ on V of type II is a n-linear transformation from V into V i.e. $T_i: V_i \to V_i$ for i = 1, 2, …, n.*

We prove the following lemma.

**LEMMA 1.2.1:** *Let $V = V_1 \cup V_2 \cup \ldots \cup V_n$ be a n-vector space over the n-field $F = F_1 \cup F_2 \cup \ldots \cup F_n$ of type II. Let U, $T^1$ and $T^2$ be the n-linear operators on V and let $C = C^1 \cup C^2 \cup \ldots \cup C^n$ be an element of $F = F_1 \cup F_2 \cup \ldots \cup F_n$. Then*

a. *$IU = UI = U$ where $I = I_1 \cup I_2 \cup \ldots \cup I_n$ is the n-identity linear operator on V i.e. $I_i (v_i) = v_i$ for every $v_i \in V_i$ ; $I_i$; $V_i \to V_i$, i = 1, 2, ..., n.*
b. *$U(T^1 + T^2) = UT^1 + UT^2$*
   *$(T^1 + T^2)U = T^1U + T^2U$*
c. *$C(UT^1) = (CU)T^1 = U(CT^1)$.*

*Proof:* (a) Given $V = V_1 \cup V_2 \cup \ldots \cup V_n$, a n-vector space defined over the n-field of type II for $F = F_1 \cup F_2 \cup \ldots \cup F_n$ and each $V_i$ of V is defined over the field $F_i$ of F. If $I = I_1 \cup I_2 \cup \ldots \cup I_n$**:** $V \to V$ such that $I_i (\upsilon_i) = \upsilon_i$ for every $\upsilon_i \in V_i$; true for i = 1, 2, …, n then $I(\alpha) = \alpha$ for every $\alpha \in V_1 \cup V_2 \cup \ldots \cup V_n$ in V and $UI = (U_1 \cup U_2 \cup \ldots \cup U_n) (I_1 \cup \ldots \cup I_n) = U_1I_1 \cup U_2I_2 \cup \ldots \cup U_nI_n$. Since each $U_i I_i = I_i U_i$ for i = 1, 2, …, n, we have $IU = I_1U_1 \cup I_2U_2 \cup \ldots \cup I_nU_n$. Hence (a) is proved.

$$\begin{aligned}[U(T^1 + T^2)](\alpha) &= (U_1 \cup U_2 \cup \ldots \cup U_n) [(T_1^1 \cup T_2^1 \cup \ldots \cup T_n^1) + (T_1^2 \cup T_2^2 \cup \ldots \cup T_n^2)] [\alpha^1 \cup \alpha^2 \cup \ldots \cup \alpha^n] \\ &= (U_1 \cup U_2 \cup \ldots \cup U_n)[(T_1^1 + T_1^2) \cup (T_2^1 + T_2^2) \cup \ldots \cup (T_n^1 + T_n^2)] (\alpha^1 \cup \alpha^2 \cup \ldots \cup \alpha^n) \\ &= U_1(T_1^1 + T_1^2)\alpha^1 \cup U_2(T_2^1 + T_2^2)(\alpha^2) \cup \ldots \cup U_n(T_n^1 + T_n^2)(\alpha^n).\end{aligned}$$

Since for each i we have

$$U_i(T_i^1 + T_i^2)(\alpha^i) = U_iT_i^1\alpha^i + U_iT_i^2\alpha^i$$

for every $\alpha^i \in V_i$; i = 1, 2, …, n.
Thus $U_i(T_i^1 + T_i^2) = U_iT_i^1 + U_iT_i^2$, for every i, i = 1, 2, …, n hence $U(T^1 + T^2) = UT^1 + UT^2$, similarly $(T^1 + T^2)U = T^1U + T^2U$. Proof of (c) is left as an exercise for the reader. Suppose $L^n(V, V)$ denote the set of all n-linear operators from V to V of type II. Then by the above lemma we see $L^n(V, V)$ is a n-vector space over the n-field of type II.

Clearly $L^n(V, V) = L^n(V_1, V_1) \cup L^n(V_2, V_2) \cup \ldots \cup L^n(V_n, V_n)$ we see each $L^n(V_i, V_i)$ is a vector space over the field $F_i$, true for $i = 1, 2, \ldots, n$. Hence $L^n(V, V)$ is a n-vector space over the n-field $F = F_1 \cup F_2 \cup \ldots \cup F_n$ of type II.

Since we have composition of any two n-linear operators to be contained in $L^n(V, V)$ and $L^n(V, V)$ contains the n-identity we see $L^n(V, V)$ is a n-linear algebra over the n-field F of type II. To this end we just recall the definition of n-linear algebra over the n-field of type II.

**DEFINITION 1.2.10:** *Let $V = V_1 \cup V_2 \cup \ldots \cup V_n$ be a n-vector space over the n-field $F = F_1 \cup F_2 \cup \ldots \cup F_n$ of type II, where $V_i$ is a vector space over the field $F_i$. If each of the $V_i$ is a linear algebra over $F_i$ for every i, i = 1, 2, ..., n then we call V to be a n-linear algebra over the n-field F of type II.*

We illustrate this by the following example.

***Example 1.2.6:*** Given $V = V_1 \cup V_2 \cup V_3 \cup V_4$ where

$$V_1 = \left\{ \begin{bmatrix} a & b \\ c & d \end{bmatrix} \middle| a, b, c, d \in Z_2 \right\},$$

a linear algebra over $Z_2$; $V_2$ = {All polynomials with coefficients from $Z_7$} a linear algebra over $Z_7$. $V_3 = Q(\sqrt{3}) \times Q(\sqrt{3})$ a linear algebra over $Q(\sqrt{3})$ and $V_4 = Z_5 \times Z_5 \times Z_5$ a linear algebra over $Z_5$. Thus $V = V_1 \cup V_2 \cup V_3 \cup V_4$ is a 4-linear algebra over the 4-field $F = Z_2 \cup Z_7 \cup Q(\sqrt{3}) \cup Z_5$ of type II.

All n-vector spaces over a n-field F need not in general be a n-linear algebra over the n-field of type II.
We illustrate this situation by an example.

***Example 1.2.7:*** Let $V = V_1 \cup V_2 \cup V_3$ be a 3-vector space over the 3-field $F = Z_2 \cup Q \cup Z_7$. Here

$$V_1 = \left\{ \begin{bmatrix} a & b & c \\ d & e & f \end{bmatrix} \middle| a, b, c, d, e, f \in Z_2 \right\}$$

a vector space over the field $Z_2$. Clearly $V_1$ is only a vector space over $Z_2$ and never a linear algebra over $Z_2$. $V_2$ = {All polynomials in x with coefficients from Q}; $V_2$ is a linear algebra over Q.

$$V_3 = \left\{ \begin{bmatrix} a & b \\ c & d \\ e & f \end{bmatrix} \middle| a, b, c, d, e, f \in Z_7 \right\};$$

$V_3$ is a vector space over $Z_7$ and not a linear algebra. Thus V is only a 3-vector space over the 3-field F and not a 3-linear algebra over F.

We as in case of finite n-vector spaces V and W over n-field F of type II, one can associate with every n-linear transformation T of type II a n-matrix. We in case of n-vector space V over the n-field F for every n-operator on V associate a n-matrix which will always be a n-square mixed matrix. This is obvious by taking W = V, then instead of getting a n-matrix of n-order $(m_1n_2, m_2n_2, \ldots, m_nn_n)$ we will have $m_i = n_i$ for every i, i = 1, 2, …, n. So the corresponding n-matrix would be a mixed n-square matrix of n-order $(n_1^2, n_2^2, ..., n_n^2)$ provided n-dim(V) = $(n_1, n_2, \ldots, n_n)$. But in case of n-linear transformation of type II we would not be in a position to talk about invertible n-linear transformation. But in case of n-linear operators we can define invertible n-linear operators of a n-vector space over the n-field F.

**DEFINITION 1.2.11:** *Let $V = V_1 \cup V_2 \cup \ldots \cup V_n$ and $W = W_1 \cup W_2 \cup \ldots \cup W_n$ be two n-vector spaces defined over the same n-field $F = F_1 \cup F_2 \cup \ldots \cup F_n$ of type II. A n-linear transformation $T = T_1 \cup T_2 \cup \ldots \cup T_n$ from V into W is n-invertible if and only if*

1. *T is one to one i.e. each $T_i$ is one to one from $V_i$ into $W_i$ such that $T_i\alpha^i = T_i\beta^i$ implies $\alpha^i = \beta^i$ true for each i, i = 1, 2, …, n.*
2. *T is onto that is range of T is (all of) W, i.e. each $T_i$: $V_i \to W_i$ is onto and range of $T_i$ is all of $W_i$, true for every i, i = 1, 2, …, n.*

The following theorem is immediate.

**THEOREM 1.2.15:** *Let $V = V_1 \cup V_2 \cup \ldots \cup V_n$ and $W = W_1 \cup W_2 \cup \ldots \cup W_n$ be two n-vector spaces over the n-field $F = F_1 \cup F_2 \cup \ldots \cup F_n$ of type II. Let $T = T_1 \cup T_2 \cup \ldots \cup T_n$ be a n-linear transformation of V into W of type II. If T is n-invertible then the n-inverse n-function $T^{-1} = T_1^{-1} \cup T_2^{-1} \cup \ldots \cup T_n^{-1}$ is a n-linear transformation from W onto V.*

*Proof:* We know if $T = T_1 \cup T_2 \cup \ldots \cup T_n$ is a n-linear transformation from V into W and when T is one to one and onto then there is a uniquely determined n-inverse function $T^{-1} = T_1^{-1} \cup T_2^{-1} \cup \ldots \cup T_n^{-1}$ which map W onto V such that $T^{-1}T$ is the n-identity n-function on V and $TT^{-1}$ is the n-identity n-function on W. We need to prove only if a n-linear n-function T is n-invertible then the n-inverse $T^{-1}$ is also n-linear. Let $\beta_1 = \beta_1^1 \cup \beta_1^2 \cup \ldots \cup \beta_1^n$, $\beta_2 = \beta_2^1 \cup \beta_2^2 \cup \ldots \cup \beta_2^n$ be two n-vectors in W and let $C = C^1 \cup C^2 \cup \ldots \cup C^n$ be a n-scalar of the n-field $F = F_1 \cup F_2 \cup \ldots \cup F_n$. Consider

$$\begin{aligned} T^{-1}(C\beta_1 + \beta_2) &= (T^{-1} = T_1^{-1} \cup T_2^{-1} \cup \ldots \cup T_n^{-1}) \\ &\quad [C^1\beta_1^1 + \beta_2^1 \cup C^2\beta_1^2 + \beta_2^2 \cup \ldots \cup C^n\beta_1^n + \beta_2^n] \\ &= T_1^{-1}(C_1^1\beta_1^1 + \beta_2^1) \cup T_2^{-1}(C_1^2\beta_1^2 + \beta_2^2) \cup \ldots \cup \\ &\quad T_n^{-1}(C_1^n\beta_1^n + \beta_2^n). \end{aligned}$$

Since from usual linear transformation properties we see each $T_i$ has $T_i^{-1}$ and $T_i^{-1}(C_1^i\beta_1^i + \beta_2^i) = C_1^i(T_i^{-1}\beta_1^i + T_i^{-1}\beta_2^i)$ true for each i, i = 1, 2, …, n. Let $\alpha_j^i = T_i^{-1}\beta_j^i$, j = 1, 2, and i = 1, 2, …, n that is let $\alpha_j^i$ be the unique vector in $V_i$ such that $T_i\alpha_j^i = \beta_j^i$; since $T_i$ is linear

$$\begin{aligned} T_i(C^i\alpha_1^i + \alpha_2^i) &= C^iT_i\alpha_1^i + T_i\alpha_2^i \\ &= C^i\beta_1^i + \beta_2^i. \end{aligned}$$

Thus $C^i\alpha_1^i + \alpha_2^i$ is the unique vector in $V_i$ which is sent by $T_i$ into $T_i$ ( $C^i\beta_1^i + \beta_2^i$ ) and so

$$T_i^{-1}(C^i\beta_1^i + \beta_2^i) = C^i\alpha_1^i + \alpha_2^i = C^iT_i^{-1}\beta_1^i + T_i^{-1}\beta_2^i$$

and $T_i^{-1}$ is linear. This is true for every $T_i$ i.e., for i, i = 1, 2, …, n. So $T^{-1} = T_1^{-1} \cup T_2^{-1} \cup ... \cup T_n^{-1}$ is n-linear. Hence the claim.

We say a n-linear transformation T is n-non-singular if $T\gamma = 0$ implies $\gamma = (0 \cup 0 \cup \ldots \cup 0)$; i.e., each $T_i$ in T is non singular i.e. $T_i\gamma_i = 0$ implies $\gamma_i = 0$ for every i; i.e., if $\gamma = \gamma_1 \cup \gamma_2 \cup \ldots \cup \gamma_n$ then $T\gamma = 0 \cup \ldots \cup 0$ implies $\gamma = 0 \cup 0 \cup \ldots \cup 0$.

Thus T is one to one if and only if T is n-non singular.

The following theorem is proved for the n-linear transformation from V into W of type II.

**THEOREM 1.2.16:** *Let $T = T_1 \cup T_2 \cup \ldots \cup T_n$ be a n-linear transformation from $V = V_1 \cup V_2 \cup \ldots \cup V_n$ into $W = W_1 \cup W_2 \cup \ldots \cup W_n$. Then T is n-non singular if and only if T carries each n-linearly independent n-subset of V into a n-linearly independent n-subset of W.*

*Proof:* First suppose we assume $T = T_1 \cup T_2 \cup \ldots \cup T_n$ is n-non singular i.e. each $T_i: V_i \to W_i$ is non singular for i = 1, 2, …, n. Let $S = S_1 \cup S_2 \cup \ldots \cup S_n$ be a n-linearly, n-independent n-subset of V. If $\{\alpha_1^1, \alpha_2^1, ..., \alpha_{k_1}^1\} \cup \{\alpha_1^2, \alpha_2^2, ..., \alpha_{k_2}^2\} \cup \ldots \cup \{\alpha_1^n, \alpha_2^n, ..., \alpha_{k_n}^n\}$ are n-vectors in $S = S_1 \cup S_2 \cup \ldots \cup S_n$ then the n-vectors $\{T_1\alpha_1^1, T_1\alpha_2^1, ..., T_1\alpha_{k_1}^1\} \cup \{T_2\alpha_1^2, T_2\alpha_2^2, ..., T_2\alpha_{k_2}^2\} \cup \ldots \cup \{T_n\alpha_1^n, T_n\alpha_2^n, ..., T_n\alpha_{k_n}^n\}$ are n-linearly n-independent; for if $C_1^i(T_i\alpha_1^i) + \ldots + C_k^i(T_i\alpha_{k_i}^i) = 0$ for each i, i = 1, 2, …, n. then $T_i(C_1^i\alpha_1^i + ... + C_{k_i}^i\alpha_{k_i}^i) = 0$ and since each $T_i$ is singular we have $C_1^i\alpha_1^i + ... + C_{k_i}^i\alpha_{k_i}^i = 0$ this is true for each i, i = 1, 2, …, n. Thus each $C_t^i = 0$; t = 1, 2, …, $k_i$ and i = 1, 2, …, n; because $S_i$ is an independent subset of the n-set $S = S_1 \cup S_2 \cup \ldots \cup S_n$. This shows image of $S_i$ under $T_i$ is independent. Thus $T = T_1 \cup T_2 \cup \ldots \cup T_n$ is n- independent as each $T_i$ is independent.

Suppose $T = T_1 \cup T_2 \cup \ldots \cup T_n$ is such that it carriers independent n-subsets onto independent n-subsets. Let $\alpha = \alpha^1 \cup \alpha^2 \cup \ldots \cup \alpha^n$ be a non zero n-vector in V i.e. each $\alpha^i$ is a non zero vector of $V_i$, i = 1, 2, …, n. Then the set $S = S_1 \cup S_2 \cup \ldots \cup S_n$ consisting of the one vector $\alpha = \alpha^1 \cup \alpha^2 \cup \ldots \cup \alpha^n$ is n-independent each $\alpha^i$ is independent for i = 1, 2, …, n. The n-image of S is the n-set consisting of the one n-vector $T\alpha = T_1\alpha^1 \cup T_2\alpha^2 \cup \ldots \cup T_n\alpha^n$ and this is n-independent, therefore $T\alpha \neq 0$ because the n-set consisting of the zero n-vector alone is dependent. This shows that the n-nullspace of $T = T_1 \cup T_2 \cup \ldots \cup T_n$ is the zero subspace as each $T_i\alpha^i \neq 0$ implies each of the zero vector alone in $V_i$ is dependent for i = 1, 2, …, n. Thus each $T_i$ is non singular so T is n-non singular.

We prove yet another nice theorem.

**THEOREM 1.2.17:** *Let $V = V_1 \cup V_2 \cup \ldots \cup V_n$ and $W = W_1 \cup W_2 \cup \ldots \cup W_n$ be n-vector spaces over the same n-field $F = F_1 \cup F_2 \cup \ldots \cup F_n$ of type II. If T is a n-linear transformation of type II from V into W, the following are equivalent.*

- *(i) T is n-invertible*
- *(ii) T is n-non-singular*
- *(iii) T is onto that is the n-range of T is W.*

*Proof:* Given $V = V_1 \cup V_2 \cup \ldots \cup V_n$ is a n-vector space over the n-field F of $(n_1, n_2, \ldots, n_n)$ dimension over the n-field F; i.e. dimension of $V_i$ over the field $F_i$ is $n_i$ for i = 1, 2, … , n. Let us further assume n-dim(V) = $(n_1, n_2, \ldots, n_n)$ = n-dim W. We know for a n-linear transformation $T = T_1 \cup T_2 \cup \ldots \cup T_n$ from V into W, n-rank T + n-nullity T = $(n_1, n_2, \ldots, n_n)$ i.e., (rank $T_1$ + nullity $T_1$) $\cup$ (rank $T_2$ + nullity $T_2$) $\cup \ldots \cup$ (rank $T_n$ + nullity $T_n$) = $(n_1, n_2, \ldots, n_n)$. Now $T = T_1 \cup T_2 \cup \ldots \cup T_n$ is n-non singular if and only if n-nullity T = $(0 \cup 0 \cup \ldots \cup 0)$ and since dim W = $(n_1, n_2, \ldots, n_n)$ the n-range of T is W if and only if n-rank T = $(n_1, n_2, \ldots, n_n)$. Since n-rank T + n-nullity T = $(n_1, n_2, \ldots, n_n)$, the n-nullity T is $(0 \cup 0 \cup \ldots \cup 0)$ precisely when the n-rank is $(n_1, n_2, \ldots, n_n)$. Therefore T is n-nonsingular if and only if T(V) = W i.e. $T_1(V_1) \cup T_2(V_2) \cup \ldots \cup T_n(V_n) = W_1 \cup W_2 \cup \ldots \cup W_n$. So if either condition (ii) or (iii) hold good, the other is satisfied as well and T is n-invertible. We further see if the 3 conditions of the theorem are also equivalent to (iv) and (v).

(iv). If $\{\alpha_1^1, \alpha_2^1, \ldots, \alpha_{n_1}^1\} \cup \{\alpha_1^2, \alpha_2^2, \ldots, \alpha_{n_2}^2\} \cup \ldots \cup \{\alpha_1^n, \alpha_2^n, \ldots, \alpha_{n_n}^n\}$ is a n-basis of V then $\{T_1\alpha_1^1, \ldots, T_1\alpha_{n_1}^1\} \cup \{T_2\alpha_1^2, \ldots, T_2\alpha_{n_2}^2\} \cup \ldots \cup \{T_n\alpha_1^n, \ldots, T_n\alpha_{n_n}^n\}$ is a n-basis for W.

(v). There is some n-basis $\{\alpha_1^1, \alpha_2^1, \ldots, \alpha_{n_1}^1\} \cup \{\alpha_1^2, \alpha_2^2, \ldots, \alpha_{n_2}^2\} \cup \ldots \cup \{\alpha_1^n, \alpha_2^n, \ldots, \alpha_{n_n}^n\}$ for V such that $\{T_1\alpha_1^1, T_1\alpha_2^1, \ldots, T_1\alpha_{n_1}^1\} \cup \{T_2\alpha_1^2, T_2\alpha_2^2, \ldots, T_2\alpha_{n_2}^2\} \cup \ldots \cup \{T_n\alpha_1^n, T_n\alpha_2^n, \ldots, T_n\alpha_{n_n}^n\}$ is an n-basis for W.

The proof of equivalence of these conditions (i) to (v) mentioned above is left as an exercise for the reader.

**THEOREM 1.2.18:** *Every $(n_1, n_2, \ldots, n_n)$ dimensional n-vector space V over the n-field $F = F_1 \cup \ldots \cup F_n$ is isomorphic to $F_1^{n_1} \cup \ldots \cup F_n^{n_n}$.*

*Proof:* Let $V = V_1 \cup V_2 \cup \ldots \cup V_n$ be a $(n_1, n_2, \ldots, n_n)$ dimensional space over the n-field $F = F_1 \cup F_2 \cup \ldots \cup F_n$ of type II. Let $B = \{\alpha_1^1, \alpha_2^1, \ldots, \alpha_{n_1}^1\} \cup \{\alpha_1^2, \alpha_2^2, \ldots, \alpha_{n_2}^2\} \cup \ldots \cup \{\alpha_1^n, \alpha_2^n, \ldots, \alpha_{n_n}^n\}$ be an n-ordered n-basis of V. We define an n-function $T = T_1 \cup T_2 \cup \ldots \cup T_n$ from V into $F_1^{n_1} \cup F_2^{n_2} \cup \ldots \cup F_n^{n_n}$ as follows: If $\alpha = \alpha^1 \cup \alpha^2 \cup \ldots \cup \alpha^n$ is in V, let $T\alpha$ be the $(n_1, n_2, \ldots, n_n)$ tuple;

$$(x_1^1, x_2^1, \ldots, x_{n_1}^1) \cup (x_1^2, x_2^2, \ldots, x_{n_2}^2) \cup \ldots \cup (x_1^n, x_2^n, \ldots, x_{n_n}^n)$$

of the n-coordinate of $\alpha = \alpha^1 \cup \alpha^2 \cup \ldots \cup \alpha^n$ relative to the n-ordered n-basis B, i.e. the $(n_1, n_2, \ldots, n_n)$ tuple such that

$$\alpha = (x_1^1\alpha_1^1 + \ldots + x_{n_1}^1\alpha_{n_1}^1) \cup (x_1^2\alpha_1^2 + \ldots + x_{n_2}^2\alpha_{n_2}^2) \cup \ldots \cup (x_1^n\alpha_1^n + \ldots + x_{n_n}^n\alpha_{n_n}^n).$$

Clearly T is a n-linear transformation of type II, one to one and maps V onto $F_1^{n_1} \cup F_2^{n_2} \cup \ldots \cup F_n^{n_n}$ or each $T_i$ is linear and one to one and maps $V_i$ to $F_i^{n_i}$ for every i, i = 1, 2, …, n. Thus we can as in case of vector spaces transformation by matrices give a representation of n-transformations by n-matrices.

Let $V = V_1 \cup V_2 \cup \ldots \cup V_n$ be a n-vector space over the n-field $F = F_1 \cup F_2 \cup \ldots \cup F_n$ of $(n_1, n_2, \ldots, n_n)$ dimension over F. Let $W = W_1 \cup \ldots \cup W_n$ be a n-vector space over the same n-field F of $(m_1, m_2, \ldots, m_n)$ dimension over F. Let B = $\{\alpha_1^1, \alpha_2^1, \ldots, \alpha_{n_1}^1\} \cup \{\alpha_1^2, \alpha_2^2, \ldots, \alpha_{n_2}^2\} \cup \ldots \cup \{\alpha_1^n, \alpha_2^n, \ldots, \alpha_{n_n}^n\}$ be a n-basis for V and C= $(\beta_1^1, \beta_2^1, \ldots, \beta_{m_1}^1) \cup (\beta_1^2, \beta_2^2, \ldots, \beta_{m_2}^2) \cup \ldots \cup (\beta_1^n, \beta_2^n, \ldots, \beta_{m_n}^n)$ be an n-ordered basis for W. If T is any n-linear transformation of type II from V into W then T is determined by its action on the vectors $\alpha^i$. Each of the $(n_1, n_2, \ldots, n_n)$ tuple vector; $T_i\alpha_j^i$ is uniquely expressible as a linear combination; $T_i\alpha_j^i = \sum_{k_i=1}^{m_i} A_{k_i j}^i \beta_{k_i}^i$ . This is true for every i, i = 1, 2, …, n and $1 \leq j \leq n_i$ of $\beta_{k_i}^i$; the scalars $A_{ij}^i, \ldots, A_{mj}^i$ being the coordinates of $T_i\alpha_j^i$ in the ordered basis $\{\beta_1^i, \beta_2^i, \ldots, \beta_{m_i}^i\}$ of C. This is true for each i, i = 1, 2, …, n.

Accordingly the n-transformation $T = T_1 \cup T_2 \cup \ldots \cup T_n$ is determined by the $(m_1n_1, m_2n_2, \ldots, m_nn_n)$ scalars $A_{k_i j}^i$. The $m_i \times n_i$ matrix $A^i$ defined by $A_{(k_i j)}^i = A_{k_i j}^i$ is called the component matrix $T_i$ of T relative to the component basis $\{\alpha_1^i, \alpha_2^i, \ldots, \alpha_{n_i}^i\}$ and $\{\beta_1^i, \beta_2^i, \ldots, \beta_{m_i}^i\}$ of B and C respectively. Since this is true for every i, we have A = $A_{(k_1,j)}^1 \cup A_{(k_2,j)}^2 \cup \ldots \cup A_{(k_n,j)}^n$ = $A^1 \cup A^2 \cup \ldots \cup A^n$ (for simplicity of notation) the n-matrix associated with $T = T_1 \cup T_2 \cup \ldots \cup T_n$. Each $A^i$ determines the linear transformation $T_i$ for i = 1, 2, …, n. If $\alpha^i = x_1^i\alpha_1^i + \ldots + x_{n_i}^i\alpha_{n_i}^i$ is a vector in $V_i$ then

$$T_i\alpha^i = (T_i \sum_{j=1}^{n_i} x_j^i \alpha_j^i)$$

$$= \sum_{j=1}^{n_i} x_j^i (T_i \alpha_j^i)$$

$$= \sum_{j=1}^{n_i} x_j^i \sum_{k=1}^{m_i} A_{k_i j}^i \beta_{k_i}^i$$

$$= \sum_{k=1}^{m_i} \sum_{j=1}^{n_i} (A_{k_i j}^i x_j^i) \beta_{k_i}^i .$$

If $X^i$ is the coordinate matrix of $\alpha^i$ in the component n-basis of B then the above computation shows that $A^i X^i$ is the coordinate matrix of the vector $T_i\alpha^i$; that is the component of the n-basis C because the scalar $\sum_{j=1}^{n_i} A_{k_i j}^i x_{k_i}^i$ is the entry in the $k^{th}$ row of the column matrix $A^i X^i$. This is true for every i, i = 1, 2, … , n. Let us also observe that if $A^i$ is any $m_i \times n_i$ matrix over the field $F_i$ then

$$T_i(\sum_{j=1}^{n_i} x_j^i \alpha_j^i) = \sum_{k=1}^{m_i} (\sum_{j=1}^{n_i} A_{k_i j}^i x_j^i) \beta_{k_i}^i$$

defines a linear transformation $T_i$ from $V_i$ into $W_i$, the matrix of which is $A^i$ relative to $\{\alpha_1^i, \ldots, \alpha_{n_i}^i\}$ and $\{\beta_1^i, \beta_2^i, \ldots, \beta_{m_i}^i\}$, this is true of every i. Hence $T = T_1 \cup T_2 \cup \ldots \cup T_n$ is a linear n-transformation from $V = V_1 \cup V_2 \cup \ldots \cup V_n$ into $W = W_1 \cup W_2 \cup \ldots \cup W_n$, the n-matrix of which is $A = A^1 \cup A^2 \ldots \cup A^n$ relative to the n-basis $B = \{\alpha_1^1, \alpha_2^1, \ldots, \alpha_{n_1}^1\} \cup \{\alpha_1^2, \alpha_2^2, \ldots, \alpha_{n_2}^2\} \cup \ldots \cup \{\alpha_1^n, \alpha_2^n, \ldots, \alpha_{n_n}^n\}$ and $C = \{\beta_1^1, \beta_2^1, \ldots, \beta_{m_1}^1\} \cup \{\beta_1^2, \beta_2^2, \ldots, \beta_{m_2}^2\} \cup \ldots \cup \{\beta_1^n, \beta_2^n, \ldots, \beta_{m_n}^n\}$.

In view of this we have the following interesting theorem.

**THEOREM 1.2.19:** *Let $V = V_1 \cup \ldots \cup V_n$ be a $(n_1, n_2, \ldots, n_n)$-dimension n-vector space over the n-field $F = F_1 \cup \ldots \cup F_n$. Let*

$$B = \{\alpha_1^1, \alpha_2^1, ..., \alpha_{n_1}^1\} \cup \{\alpha_1^2, \alpha_2^2, ..., \alpha_{n_2}^2\} \cup ... \cup \{\alpha_1^n, \alpha_2^n, ..., \alpha_{n_n}^n\}$$

*$= B^1 \cup B^2 \cup ... \cup B^n$ be the n-basis of V over F and*

$$C = \{\beta_1^1, \beta_2^1, ..., \beta_{m_1}^1\} \cup \{\beta_1^2, \beta_2^2, ..., \beta_{m_2}^2\} \cup ... \cup \{\beta_1^n, \beta_2^n, ..., \beta_{m_n}^n\}$$
$$= C^1 \cup C^2 \cup ... \cup C^n,$$

*n-basis for W, where W is a n-vector space over the same n-field F of ($m_1$, $m_2$, ..., $m_n$) dimension over F. For each n-linear transformation $T = T_1 \cup T_2 \cup ... \cup T_n$ of type II from V into W there is a n-matrix $A = A^1 \cup A^2 \cup ... \cup A^n$ where each $A^i$ is a $m_i \times n_i$ matrix with entries in $F_i$ such that $\left[T\alpha\right]_C = A[\alpha]_B$ where for every n-vector $\alpha = \alpha^1 \cup \alpha^2 \cup ... \cup \alpha^n \in V$ we have*

$$\left[T_1\alpha^1\right]_{C^1} \cup ... \cup \left[T_n\alpha^n\right]_{C^n} = A^1\left[\alpha^1\right]_{B^1} \cup ... \cup A^n\left[\alpha^n\right]_{B^n}.$$

*Further more $T \to A$ is a one to one correspondence between the set of all n-linear transformation from V into W of type II and the set of all ($m_1 \times n_1$, $m_2 \times n_2$, ..., $m_n \times n_n$) n-matrices over the n-field F. The n-matrix $A = A^1 \cup A^2 \cup ... \cup A^n$ which is associated with $T = T_1 \cup T_2 \cup ... \cup T_n$ is called the n-matrix of T relative to the n-ordered basis B and C. From the equality*

$$T\alpha = T_1\alpha_j^i \cup T_2\alpha_J^2 \cup ... \cup T_n\alpha_j^n =$$
$$\sum_{k_1=1}^{m_1} A_{k_1 j}^1 \beta_j^1 \cup \sum_{k_2=1}^{m_2} A_{k_2 j}^2 \beta_j^2 \cup ... \cup \sum_{k_n=1}^{m_n} A_{k_n j}^n \beta_j^n$$

*says that $A = A^1 \cup A^2 \cup ... \cup A^n$ is the n-matrix whose n-columns $\{A_1^1, A_2^1, ..., A_{n_1}^1\} \cup \{A_1^2, A_2^2, ..., A_{n_2}^2\} \cup ... \cup \{A_1^n, A_2^n, ..., A_{n_n}^n\}$ are given by $A = [T\alpha]_C$ i.e. $A_j^i = [T_i\alpha_j^i]_{C^i}$, i = 1, 2, ... , n and j = 1, 2, ..., $n_i$. If $U = U_1 \cup U_2 \cup ... \cup U_n$ is another n-linear transformation from $V = V_1 \cup V_2 \cup ... \cup V_n$ into $W = W_1 \cup W_2 \cup ... \cup W_n$ and $P = \{P_1^1, ..., P_{n_1}^1\} \cup \{P_1^2, ..., P_{n_2}^2\} \cup \{P_1^n, ..., P_{n_n}^n\}$ is the n-matrix of U relative to the n-ordered basis B and C then CA + P is the matrix CT + U where $C = C^1 \cup ... \cup C^n \in F_1 \cup F_2 \cup ... \cup F_n$ where $C^i \in F_i$, for i = 1, 2, ...,n. That is clear because* $C^i A^i_j + P^i_j = C^i[T_i\alpha^i_j]_{C^i} + [U_i\alpha^i_j]_{C^i} = C^i[\left(T_i\alpha^i_j + U_i\right)\alpha^i_j]_{C^i}$. *This is true for every i, i = 1, 2, … , n and* $1 \le j \le n_i$.

Several interesting results analogous to usual vector spaces can be derived for n-vector spaces over n-field of type II and its related n-linear transformation of type II.

The reader is expected to prove the following theorem.

**THEOREM 1.2.20:** *Let* $V = V_1 \cup V_2 \cup \ldots \cup V_n$ *be a* $(n_1, n_2, \ldots, n_n)$*, n-vector space over the n-field* $F = F_1 \cup F_2 \cup \ldots \cup F_n$ *and let* $W = W_1 \cup W_2 \cup \ldots \cup W_n$ *be a* $(m_1, m_2, \ldots, m_n)$ *dimensional vector space over the same n-field* $F = F_1 \cup \ldots \cup F_n$*. For each pair of ordered n-basis B and C for V and W respectively the n-function which assigns to a n-linear transformation T of type II, its n-matrix relative to B, C is an n-isomorphism between the n-space* $L^n(V, W)$ *and the space of all n-matrices of n-order* $(m_1 \times n_1, m_2 \times n_2, \ldots, m_n \times n_n)$ *over the same n-field F.*

Since we have the result to be true for every pair of vector spaces $V_i$ and $W_i$ over $F_i$, we can appropriately extend the result for $V = V_1 \cup V_2 \cup \ldots \cup V_n$ and $W = W_1 \cup W_2 \cup \ldots \cup W_n$ over $F = F_1 \cup F_2 \cup \ldots \cup F_n$ as the result is true for every i.

Yet another theorem of interest is left for the reader to prove.

**THEOREM 1.2.21:** *Let* $V = V_1 \cup V_2 \cup \ldots \cup V_n$*,* $W = W_1 \cup W_2 \cup \ldots \cup W_n$ *and* $Z = Z_1 \cup Z_2 \cup \ldots \cup Z_n$ *be three* $(n_1, n_2, \ldots, n_n)$*,* $(m_1, m_2, \ldots, m_n)$*, and* $(p_1, p_2, \ldots, p_n)$ *dimensional n-vector spaces respectively defined over the n-field* $F = F_1 \cup F_2 \cup \ldots \cup F_n$*. Let* $T = T_1 \cup T_2 \cup \ldots \cup T_n$ *be a n-linear transformation of type II from V into W and* $U = U_1 \cup U_2 \cup \ldots \cup U_n$ *a n-linear transformation of type II from W into Z. Let* $B = B^1 \cup B^2 \cup \ldots \cup B^n = \{\alpha^1_1, \alpha^1_2, \ldots, \alpha^1_{n_1}\} \cup \{\alpha^2_1, \alpha^2_2, \ldots, \alpha^2_{n_2}\} \cup \ldots \cup \{\alpha^n_1, \alpha^n_2, \ldots, \alpha^n_{n_n}\}$*,* $C = C^1 \cup C^2 \cup \ldots \cup C^n = \{\beta^1_1, \beta^1_2, \ldots, \beta^1_{m_1}\} \cup \{\beta^2_1, \beta^2_2, \ldots, \beta^2_{m_2}\} \cup \ldots \cup \{\beta^n_1, \beta^n_2, \ldots, \beta^n_{m_n}\}$ *and* $D = D^1 \cup D^2 \cup \ldots \cup D^n = \{\gamma^1_1, \gamma^1_2, \ldots, \gamma^1_{p_1}\} \cup \{\gamma^2_1, \gamma^2_2, \ldots, \gamma^2_{p_2}\} \cup \ldots \cup \{\gamma^n_1, \gamma^n_2, \ldots, \gamma^n_{p_n}\}$ *be n-ordered n-basis* *for the n-spaces V, W and Z respectively, if* $A = A^1 \cup A^2 \cup \ldots \cup A^n$ *is a n-matrix relative to T to the pair B, C and* $R = R^1 \cup R^2 \cup \ldots \cup R^n$ *is a n-matrix of U relative to the pair C and D then the n-matrix of the composition UT relative to the pair B, D is the product n-matrix E = RA i.e., if* $E = E^1 \cup E^2 \cup \ldots \cup E^n = R^1A^1 \cup R^2A^2 \cup \ldots \cup R^nA^n$.

**THEOREM 1.2.22:** *Let* $V = V_1 \cup V_2 \cup \ldots \cup V_n$ *be a* $(n_1, n_2, \ldots, n_n)$ *finite dimensional n-vector space over the n-field* $F = F_1 \cup F_2 \cup \ldots \cup F_n$ *and let* $B = B^1 \cup B^2 \cup \ldots \cup B^n = \{\alpha_1^1, \alpha_2^1, \ldots, \alpha_{n_1}^1\} \cup \{\alpha_1^2, \alpha_2^2, \ldots, \alpha_{n_2}^2\} \cup \ldots \cup \{\alpha_1^n, \alpha_2^n, \ldots, \alpha_{n_n}^n\}$ *and* $C = C^1 \cup C^2 \cup \ldots \cup C^n = \{\beta_1^1, \beta_2^1, \ldots, \beta_{m_1}^1\} \cup \{\beta_1^2, \beta_2^2, \ldots, \beta_{m_2}^2\} \cup \ldots \cup \{\beta_1^n, \beta_2^n, \ldots, \beta_{m_n}^n\}$ *be an n-ordered n-basis for V. Suppose T is a n-linear operator on V.*

*If* $P = P^1 \cup P^2 \cup \ldots \cup P^n = \{P_1^1, \ldots, P_{n_1}^1\} \cup \{P_1^2, \ldots, P_{n_2}^2\} \cup \ldots \cup \{P_1^n, P_2^n, \ldots, P_{n_n}^n\}$ *is a* $(n_1 \times n_1, n_2 \times n_2, \ldots, n_n \times n_n)$ *n-matrix with j* $^{th}$ *component of the n-columns* $P_j^i = [\beta_j^i]_B$; $1 \le j \le n_i$, $i = 1, 2, \ldots, n$ *then* $[T]_C = P^{-1}[T]_B P$ *i.e.* $[T_1]_{C^1} \cup [T_2]_{C^2} \cup \ldots \cup [T_n]_{C^n} = P_1^{-1}[T_1]_{B^1} P_1 \cup P_2^{-1}[T_2]_{B^2} P_2 \cup \ldots \cup P_n^{-1}[T_n]_{B^n} P_n$. *Alternatively if U is the n-invertible operator on V defined by* $U_i\alpha_j^i = \beta_j^i$, $j = 1, 2, \ldots, n_i$, $i = 1, 2, \ldots, n$ *then* $[T]_C = [U]_B^{-1} B[T]_B[U]_B$ *i.e.*

$$[T_1]_{C^1} \cup [T_2]_{C^2} \cup \ldots \cup [T_n]_{C^n} =$$
$$[U_1]_{B^1}^{-1}[T_1]_{B^1}[U_1]_{B^1} \cup [U_2]_{B^2}^{-1}[T_2]_{B^2}[U_2]_{B^2} \cup \ldots \cup [U_n]_{B^n}^{-1}[T_n]_{B^n}[U_n]_{B^n}.$$

The proof of the above theorem is also left for the reader.

We now define similar space n-matrices.

**DEFINITION 1.2.12:** *Let* $A = A_1 \cup A_2 \cup \ldots \cup A_n$ *be a n-mixed square matrix of n-order* $(n_1 \times n_1, n_2 \times n_2, \ldots, n_n \times n_n)$ *over the n-field* $F = F_1 \cup F_2 \cup \ldots \cup F_n$ *i.e., each* $A_i$ *takes its entries from the field* $F_i$, $i = 1,2, \ldots, n$. $B = B^1 \cup B^2 \cup \ldots \cup B^n$ *is a n-mixed* *square matrix of n-order ($n_1 \times n_1$, $n_2 \times n_2$, ..., $n_n \times n_n$) over the same n-field F. We say that B is n-similar to A over the n-field F if there is an invertible n-matrix $P = P_1 \cup P_2 \cup \ldots \cup P_n$ of n-order ($n_1 \times n_1$, $n_2 \times n_2$, ..., $n_n \times n_n$) over the n-field F such that $B = P^{-1} A P$ i.e.*

$$B_1 \cup B_2 \cup \ldots \cup B_n = P_1^{-1} A_1 P_1 \cup P_2^{-1} A_2 P_2 \cup \ldots \cup P_n^{-1} A_n P_n \ .$$

Now we proceed on to define the new notion of n-linear functionals or linear n-functionals. We know we were not in a position to define n-linear functionals in case of n-vector space over the field F i.e. for n-vector spaces of type I.

Only in case of n-vector spaces defined over the n-field of type II we are in a position to define n-linear functionals on the n-vector space V.

Let $V = V_1 \cup V_2 \cup \ldots \cup V_n$ be a n-vector space over the n-field $F = F_1 \cup F_2 \cup \ldots \cup F_n$ of type II. A n-linear transformation $f = f_1 \cup f_2 \cup \ldots \cup f_n$ from V into the n-field of n-scalars is also called a n-linear functional on V i.e., if f is a n-function from V into F such that $f(C\alpha + \beta) = Cf(\alpha) + f(\beta)$ where $C = C^1 \cup C^2 \cup \ldots \cup C^n$; i = 1, 2, …, n, $\alpha = \alpha^1 \cup \alpha^2 \cup \ldots \cup \alpha^n$ and $\beta = \beta^1 \cup \beta^2 \cup \ldots \cup \beta^n$ where $\alpha^i, \beta^i \in V_i$ for each i, i = 1,2, …, n. $f = f_1 \cup f_2 \cup \ldots \cup f_n$ where each $f_i$ is a linear functional on $V_i$; i = 1, 2, …, n. i.e.,

$$\begin{aligned} f(C\alpha + \beta) &= (f_1 \cup f_2 \cup \ldots \cup f_n)\,(C^1\alpha^1 + \beta^1) \cup (C^2\alpha^2 + \beta^2) \cup \ldots \cup (C^n\alpha^n + \beta^n) \\ &= (C^1 f_1(\alpha^1) + f_1(\beta^1)) \cup (C^2 f_2(\alpha^2) + f_2(\beta^2)) \cup \ldots \cup (C^n f_n(\alpha^n) + f_n(\beta^n))] \\ &= f_1(C^1\alpha^1 + \beta^1) \cup f_2(C^2\alpha^2 + \beta^2) \cup \ldots \cup f_n(C^n\alpha^n + \beta^n) \\ &= (C^1 f_1(\alpha^1) + f_1(\beta^1)) \cup (C^2 f_2(\alpha^2) + f_2(\beta^2)) \cup \ldots \cup (C^n f_n(\alpha^n) + f_n(\beta^n)) \end{aligned}$$

for all n-vectors $\alpha, \beta \in V$ and $C \in F$.

We make the following observation.

Let $F = F_1 \cup F_2 \cup \ldots \cup F_n$ be a n-field and let $F_1^{n_1} \cup F_2^{n_2} \cup \ldots \cup F_n^{n_n}$ be a n-vector space over the n-field F of type II. A n-linear function $f = f_1 \cup f_2 \cup \ldots \cup f_n$ from $F_1^{n_1} \cup F_2^{n_2} \cup \ldots \cup F_n^{n_n}$ to $F_1 \cup F_2 \cup \ldots \cup F_n$ given by $f_1(x_1^1, \ldots, x_{n_1}^1) \cup f_2(x_1^2, \ldots, x_{n_2}^2) \cup \ldots \cup f_n(x_1^n, \ldots, x_{n_n}^n) = (\alpha_1^1 x_1^1 + \ldots + \alpha_{n_1}^1 x_{n_1}^1) \cup (\alpha_1^2 x_1^2 + \ldots + \alpha_{n_2}^2 x_{n_2}^2) \cup \ldots \cup (\alpha_1^n x_1^n + \ldots + \alpha_{n_2}^n x_{n_n}^n)$ where $\alpha_j^i$ are in $F_i$, $1 \leq j \leq n_i$; for each $i = 1, 2, \ldots, n$; is a n-linear functional on $F_1^{n_1} \cup F_2^{n_2} \cup \ldots \cup F_n^{n_n}$. It is the n-linear functional which is represented by the n-matrix $[\alpha_1^1, \alpha_2^1, \ldots, \alpha_{n_1}^1] \cup [\alpha_1^2, \alpha_2^2, \ldots, \alpha_{n_2}^2] \cup \ldots \cup [\alpha_1^n, \alpha_2^n, \ldots, \alpha_{n_n}^n]$ relative to the standard ordered n-basis for $F_1^{n_1} \cup F_2^{n_2} \cup \ldots \cup F_n^{n_n}$ and the n-basis $\{1\} \cup \{1\} \cup \ldots \cup \{1\}$ for $F = F_1 \cup F_2 \cup \ldots \cup F_n$. $\alpha_j^i = f_i(E_j^i)$; $j = 1, 2, \ldots, n_i$ for every $i = 1, 2, \ldots, n$. Every n-linear functional on $F_1^{n_1} \cup F_2^{n_2} \cup \ldots \cup F_n^{n_n}$ is of this form for some n-scalars $(\alpha_1^1 \ldots \alpha_{n_1}^1) \cup (\alpha_1^2 \ldots \alpha_{n_2}^2) \cup \ldots \cup (\alpha_1^n \alpha_2^n \ldots \alpha_{n_n}^n)$. That is immediate from the definition of the n-linear functional of type II because we define $\alpha_j^i = f_i(E_j^i)$.

$$\therefore f_1(x_1^1, \ldots, x_{n_1}^1) \cup f_2(x_1^2, \ldots, x_{n_2}^2) \cup \ldots \cup f_n(x_1^n, x_2^n, \ldots, x_{n_n}^n)$$

$$= f_1\left(\sum_{j=1}^{n_1} x_j^1 E_j^1\right) \cup f_2\left(\sum_{j=1}^{n_2} x_j^2 E_j^2\right) \cup \ldots \cup f_n\left(\sum_{j=1}^{n_n} x_j^n E_j^n\right)$$

$$= \sum_j x_j^1 f_1(e_j^1) \cup \sum_j x_j^2 f_2(e_j^2) \cup \ldots \cup \sum_j x_j^n f_n(e_j^n)$$

$$= \sum_j^{n_1} a_j^1 x_j^1 \cup \sum_{j=1}^{n_2} a_j^2 x_j^2 \cup \ldots \cup \sum_{j=1}^{n_n} a_j^n x_j^n .$$

Now we can proceed onto define the new notion of n-dual space or equivalently dual n-space of a n-space $V = V_1 \cup \ldots \cup V_n$ defined over the n-field $F = F_1 \cup F_2 \cup \ldots \cup F_n$ of type II. Now as in case of $L^n(V,W) = L(V_1,W_1) \cup L(V_2,W_2) \cup \ldots \cup L(V_n, W_n)$ we in case of linear functional have $L^n(V,F) = L(V_1,F_1) \cup L(V_2, F_2) \cup \ldots \cup L(V_n,F_n)$. Now $V^* = L^n(V,F) = V_1^* \cup V_2^* \cup \ldots \cup V_n^*$ i.e. each $V_i^*$ is the dual space of $V_i$, $V_i$ defined over the field $F_i$; $i = 1, 2, \ldots, n$. We know in case of vector space $V_i$, dim $V_i^*$ =

dim $V_i$ for every i. Thus dim V = dim $V^*$ = dim $V_1^* \cup$ dim $V_2^* \cup \ldots \cup$ dim $V_n^*$. If B = $\{\alpha_1^1, \alpha_2^1, \ldots, \alpha_{n_1}^1\} \cup \{\alpha_1^2, \alpha_2^2, \ldots, \alpha_{n_2}^2\} \cup \ldots \cup \{\alpha_1^n, \alpha_2^n, \ldots, \alpha_{n_n}^n\}$ is a n-basis for V, then we know for a n-linear functional of type II. $f = f_1 \cup f_2 \cup \ldots \cup f_n$ we have $f_k$ on $V_k$ is such that $f_i^k(\alpha_j^k) = \delta_{ij}^k$ true for k = 1, 2, …, n. In this way we obtain from B a set of n-tuple of $n_i$ sets of distinct n-linear functionals $\{f_1^1, f_2^1, \ldots, f_{n_1}^1\} \cup \{f_1^2, f_2^2, \ldots, f_{n_2}^2\} \cup \ldots \cup \{f_1^n, f_2^n, \ldots, f_{n_k}^n\}$ = $f^1 \cup f^2 \cup \ldots \cup f^n$ on V. These n-functionals are also n-linearly independent over the n-field $F = F_1 \cup F_2 \cup \ldots \cup F_n$ i.e., $\{f_1^i, f_2^i, \ldots, f_{n_i}^i\}$ is linearly independent on $V_i$ over the field $F_i$; this is true for each i, i = 1, 2, …, n.

Thus $f^i = \sum_{j=1}^{n_i} c_j^i f_j^i$, i = 1, 2, …, n i.e.

$$f = \sum_{j=1}^{n_1} c_j^1 f_j^1 \cup \sum_{j=1}^{n_2} c_j^2 f_j^2 \cup \ldots \cup \sum_{j=1}^{n_n} c_j^n f_j^n .$$

$$\begin{aligned} f^i(\alpha_j^i) &= \sum_{k=1}^{n_i} c_k^i f_k^i(\alpha_j^k) \\ &= \sum_{k=1}^{n_i} c_k^i \delta_{kj} \\ &= c_j^i . \end{aligned}$$

This is true for every i = 1, 2, …, n and $1 \leq j \leq n_i$. In particular if each $f_i$ is a zero functional $f^i(\alpha_j^i) = 0$ for each j and hence the scalar $c_j^i$ are all zero. Thus $f_1^i, f_2^i, \ldots, f_{n_i}^i$ are $n_i$ linearly independent functionals of $V_i$ defined on $F_i$, which is true for each i, i = 1, 2, … , n. Since we know $V_i^*$ is of dimension $n_i$ it must be that $B_i^* = \{f_1^i, f_2^i, \ldots, f_{n_i}^i\}$ is a basis of $V_i^*$ which we know is a dual basis of B. Thus $B^* = B_1^* \cup B_2^* \cup \ldots \cup B_n^*$ = $\{f_1^1, f_2^1, \ldots, f_{n_1}^1\} \cup \{f_1^2, f_2^2, \ldots, f_{n_2}^2\} \cup \ldots \cup \{f_1^n, f_2^n, \ldots, f_{n_n}^n\}$ is the n-dual basis of B = $\{\alpha_1^1, \alpha_2^1, \ldots, \alpha_{n_1}^1\} \cup \{\alpha_1^2, \alpha_2^2, \ldots, \alpha_{n_2}^2\} \cup \ldots \cup$ $\{ \alpha_1^n, \alpha_2^n, \ldots, \alpha_{n_n}^n \}$. B* forms the n-basis of $V^* = V_1^* \cup V_2^* \cup \ldots \cup V_n^*$.

We prove the following interesting theorem.

**THEOREM 1.2.23:** *Let $V = V_1 \cup \ldots \cup V_n$ be a finite ($n_1, n_2, \ldots, n_n$) dimensional n-vector space over the n-field $F = F_1 \cup F_2 \cup \ldots \cup F_n$ and let $B = \{ \alpha_1^1, \alpha_2^1, \ldots, \alpha_{n_1}^1 \} \cup \{ \alpha_1^2, \alpha_2^2, \ldots, \alpha_{n_2}^2 \} \cup \ldots \cup \{ \alpha_1^n, \alpha_2^n, \ldots, \alpha_{n_n}^n \}$ be a n-basis for V. Then there is a unique n-dual basis*

$$B^* = \{ f_1^1, f_2^1, \ldots, f_{n_1}^1 \} \cup \{ f_1^2, f_2^2, \ldots, f_{n_2}^2 \} \cup \ldots \cup \{ f_1^n, f_2^n, \ldots, f_{n_n}^n \}$$

*for $V^* = V_1^* \cup V_2^* \cup \ldots \cup V_n^*$ such that $f_i^k(\alpha_j) = \delta_{ij}^k$. For each n-linear functional $f = f^1 \cup f^2 \cup \ldots \cup f^n$ we have*

$$f = \sum_{k=1}^{n_i} f^i(\alpha_k^i) f_k^i \,;$$

*i.e.,*

$$f = \sum_{k=1}^{n_1} f^1(\alpha_k^1) f_k^1 \cup \sum_{k=1}^{n_2} f^2(\alpha_k^2) f_k^2 \cup \ldots \cup \sum_{k=1}^{n_n} f^n(\alpha_k^n) f_k^n$$

*and for each n-vector $\alpha = \alpha^1 \cup \alpha^2 \cup \ldots \cup \alpha^n$ in V we have*

$$\alpha^i = \sum_{k=1}^{n_i} f_k^i(\alpha^i) \alpha_k^i$$

*i.e.,*

$$\alpha = \sum_{k=1}^{n_1} f_k^1(\alpha^1) \alpha_k^1 \cup \sum_{k=1}^{n_2} f_k^2(\alpha^2) \alpha_k^2 \cup \ldots \cup \sum_{k=1}^{n_n} f_k^n(\alpha^n) \alpha_k^n \,.$$

*Proof:* We have shown above that there is a unique n-basis which is dual to B. If f is a n-linear functional on V then f is some n-linear combination of $f_j^i$, $i \leq j \leq k_j$ and i = 1, 2, …, n; and from earlier results observed the scalars $C_j^i$ must be given by $C_j^i = f^i(\alpha_j^i)$, $1 \leq j \leq k_i$, i =1, 2, …, n. Similarly if

$$\alpha = \sum_{i=1}^{n_1} x_i^1 \alpha_i^1 \cup \sum_{i=1}^{n_2} x_i^2 \alpha_i^2 \cup \ldots \cup \sum_{i=1}^{n_n} x_i^n \alpha_i^n$$

is a n-vector in V then

$$f_j(\alpha) = \sum_{i=1}^{n_1} x_i^1 f_j^1(\alpha_i^1) \cup \sum_{i=1}^{n_2} x_i^2 f_j^2(\alpha_i^2) \cup \ldots \cup \sum_{i=1}^{n_n} x_i^n f_j^n(\alpha_i^n)$$

$$= \sum_{i=1}^{n_1} x_i^1 \delta_{ij}^1 \cup \sum_{i=1}^{n_2} x_i^2 \delta_{ij}^2 \cup \ldots \cup \sum_{i=1}^{n_n} x_i^n \delta_{ij}^n$$

$$= x_j^1 \cup x_j^2 \cup \ldots \cup x_j^n;$$

so that the unique expression for α as a n-linear combination of $\alpha_j^k$, k = 1, 2, …, n; $1 \leq j \leq n_j$ i.e.

$$\alpha = \sum_{i=1}^{n_1} f_i^1(\alpha^1)\alpha_i^1 \cup \sum_{i=1}^{n_2} f_i^2(\alpha^2)\alpha_i^2 \cup \ldots \cup \sum_{i=1}^{n_n} f_i^n(\alpha^n)\alpha_i^n .$$

Suppose $N_f = N_{f^1}^1 \cup N_{f^2}^2 \cup \ldots \cup N_{f^n}^n$ denote the n-null space of f the n-dim $N_f$ = dim $N_{f^n}^n \cup \ldots \cup \dim N_{f^n}^n$ but dim $N_{f^i}^i$ = dim $V_i$–1 = $n_i$–1 so n-dim $N_f$ = (dim $V_1$–1) ∪ dim $V_2$–1 ∪ … ∪ dim $V_n$–1 = $n_1$ – 1 ∪ $n_2$ – 1 ∪ … ∪ $n_n$ – 1. In a vector space of dimension n, a subspace of dimension n – 1 is called a hyperspace like wise in a n-vector space V = $V_1 \cup \ldots \cup V_n$ of ($n_1$, $n_2$, … , $n_n$) dimension over the n-field F = $F_1 \cup F_2 \cup \ldots \cup F_n$ then the n-subspace has dimension ($n_1$ – 1, $n_2$ – 1, …, $n_n$ – 1) we call that n-subspace to be a n-hyperspace of V. Clearly $N_f$ is a n-hyper subspace of V.

Now this notion cannot be defined in case of n-vector spaces of type I. This is also one of the marked differences between the n-vector spaces of type I and n-vector spaces of type II.

Now we proceed onto define yet a special feature of a n-vector space of type II, the n-annihilator of V.

**DEFINITION 1.2.13:** *Let V = $V_1$, ∪ … ∪ $V_n$ be a n-vector space over the n-field F = $F_1 \cup F_2 \cup \ldots \cup F_n$ of type II. Let S = $S_1 \cup S_2 \cup \ldots \cup S_n$ be a n-subset of V = $V_1 \cup \ldots \cup V_n$ (i.e. $S_i \subseteq V_i$, i =*

*1, 2, ..., n); the n-annihilator of S is* $S^o = S_1^o \cup \ldots \cup S_n^o$ *of n-linear functionals on V such that f (α) = 0 ∪ 0 ∪ ... ∪ 0 i.e. if f* $= f^1 \cup f^2 \cup \ldots \cup f^n$ *for every* $\alpha \in S$ *i.e.,* $\alpha = \alpha^1 \cup \alpha^2 \cup \ldots \cup \alpha^n \in S_1 \cup S_2 \cup \ldots \cup S_n$*. i.e.,* $f_i(\alpha^i) = 0$ *for every* $\alpha^i \in S_i$*; i = 1, 2, ..., n. It is interesting and important to note that* $S^o = S_1^o \cup \ldots \cup S_n^o$ *is an n-subspace of* $V^* = V_1^* \cup V_2^* \cup \ldots \cup V_n^*$ *whether S is an n-subspace of V or only just a n-subset of V. If S = (0 ∪ 0 ∪ ... ∪ 0) then* $S^o = V^*$*. If S = V then* $S^o$ *is just the zero n-subspace of* $V^*$*.*

Now we prove an interesting result in case of finite dimensional n-vector space of type II.

**THEOREM 1.2.24:** *Let* $V = V_1 \cup V_2 \cup \ldots \cup V_n$ *be a n-vector space of* $(n_1, n_2, \ldots, n_n)$ *dimension over the n-field* $F = F_1 \cup F_2 \cup \ldots \cup F_n$ *of type II. Let* $W = W_1 \cup W_2 \cup \ldots \cup W_n$ *be a n-subspace of V. Then dim W + dim W° = dim V (i.e., if dim W is* $(k_1, k_2, \ldots, k_n)$*) that is* $(k_1, k_2, \ldots, k_n) + (n_1 - k_1, n_2 - k_2, \ldots, n_n - k_n) = (n_1, n_2, \ldots, n_n)$*.*

*Proof:* Let $(k_1, k_2, \ldots, k_n)$ be the n-dimension of $W = W_1 \cup W_2 \cup \ldots \cup W_n$. Let

$P = \{\alpha_1^1, \alpha_2^1, \ldots, \alpha_{k_1}^1\} \cup \{\alpha_1^2, \alpha_2^2, \ldots, \alpha_{k_2}^2\} \cup \ldots \cup \{\alpha_1^n, \alpha_2^n, \ldots, \alpha_{k_n}^n\}$

$= P^1 \cup P^2 \cup \ldots \cup P^n$ be a n-basis of W. Choose n-vectors

$$\{\alpha_{k_1+1}^1, \ldots, \alpha_{n_1}^1\} \cup \{\alpha_{k_2+1}^2, \ldots, \alpha_{n_2}^2\} \cup \ldots \cup \{\alpha_{k_n+1}^n, \ldots, \alpha_{n_n}^n\}$$

in V such that

$B = \{\alpha_1^1, \alpha_2^1, \ldots, \alpha_{n_1}^1\} \cup \{\alpha_1^2, \alpha_2^2, \ldots, \alpha_{n_2}^2\} \cup \ldots \cup \{\alpha_1^n, \alpha_2^n, \ldots, \alpha_{n_n}^n\}$

$= B^1 \cup B^2 \cup \ldots \cup B^n$ is a n-basis for V.

Let

$$f = \{f_1^1, f_2^1, \ldots, f_{n_1}^1\} \cup \{f_1^2, f_2^2, \ldots, f_{n_2}^2\} \cup \ldots \cup \{f_1^n, f_2^n, \ldots, f_{n_n}^n\}$$

$= f^1 \cup f^2 \cup \ldots \cup f^n$ be a n-basis for V* which is the dual n-basis for V. The claim is that $\{f_{k_1+1}^1, \ldots, f_{n_1}^1\} \cup \{f_{k_2+1}^2, \ldots, f_{n_2}^2\} \cup \ldots \cup \{f_{k_n+1}^n, \ldots, f_{n_n}^n\}$ is a n-basis for the n-annihilator of W° $= W_1^o \cup W_2^o \cup \ldots \cup W_n^o$. Certainly $f_i^r$ belongs to $W_r^o$ and $i \geq k_r$

+ 1; r = 1, 2, …, n, $1 \leq i \leq k_r$, because $f_i^r(\alpha_j^r) = \delta_{ij}$ and $\delta_{ij} = 0$ if $i \geq k_r + 1$ and $j \geq k_r + 1$, from this it follows that for $i \leq k_r$. $f_i^r(\alpha^r) = 0$ whenever $\alpha^r$ is a linear combination $a_1^r, a_2^r, ..., a_{k_r}^r$. The functionals $a_{k_r+1}^r, a_{k_r+2}^r, ..., a_{n_r}^r$ are independent for every r = 1, 2, …, n. Thus we must show that they span $W_r^o$, for r = 1, 2, …, n. Suppose $f^r \in V_r^*$. Now $f^r = \sum_{i=1}^{n_r} f^r(\alpha_i^r) f_i^r$ so that if $f^r$ is in $W_r^o$ we have $f^r(\alpha_i^r) = 0$ for $i \leq k_r$ and $f^r = \sum_{i=k_r+1}^{n_r} f^r(\alpha_i^r) f_i^r$.

We have to show if dim $W_r = k_r$ and dim $V_r = n_r$ then dim $W_r^o = n_r - k_r$; this is true for every r. Hence the theorem.

**COROLLARY 1.2.4:** *If $W = W_1 \cup W_2 \cup \ldots \cup W_n$ is a $(k_1, k_2, \ldots, k_n)$ dimensional n-subspace of a $(n_1, n_2, \ldots, n_n)$ dimensional n-vector space $V = V_1 \cup V_2 \cup \ldots \cup V_n$ over the n-field $F = F_1 \cup F_2 \cup \ldots \cup F_n$ then W is the intersection of $(n_1 - k_1) \cup (n_2 - k_2) \cup \ldots \cup (n_n - k_n)$ tuple n-hyper subspaces of V.*

*Proof:* From the notations given in the above theorem in $W = W_1 \cup W_2 \cup \ldots \cup W_n$, each $W_r$ is the set of vector $\alpha^r$ such that $f_i(\alpha^r) = 0$, $i = k_r + 1, \ldots, n_r$. In case $k_r = n_r - 1$ we see $W_r$ is the null space of $f_{n_r}$. This is true for every r, r = 1, 2, …, n. Hence the claim.

**COROLLARY 1.2.5:** *If $W_1 = W_1^1 \cup W_1^2 \cup \ldots \cup W_1^n$ and $W_2 = W_2^1 \cup W_2^2 \cup \ldots \cup W_2^n$ are n-subspace of the n-vector space $V = V_1 \cup V_2 \cup \ldots \cup V_n$ over the n-field $F = F_1 \cup F_2 \cup \ldots \cup F_n$ of dimension $(n_1, n_2, \ldots, n_n)$ then $W_1 = W_2$ if and only if $W_1^o = W_2^o$ i.e., if and only if $\left(W_1^i\right)^o = \left(W_2^i\right)^o$ for every i = 1, 2, …, n.*

*Proof:* If $W_1 = W_2$ i.e., if

$$W_1^1 \cup W_1^2 \cup \ldots \cup W_1^n = W_2^1 \cup W_2^2 \cup \ldots \cup W_2^n$$

then each $W_1^j = W_2^j$ for j = 1, 2, …, n, so that $\left(W_1^j\right)^o = \left(W_2^j\right)^o$ for every j = 1, 2, …,n. Thus $W_1^o = W_2^o$. If on the other hand $W_1 \neq W_2$ i.e. $W_1^1 \cup W_1^2 \cup \ldots \cup W_1^n \neq W_2^1 \cup W_2^2 \cup \ldots \cup W_2^n$ then one of the two n-subspaces contains a n-vector which is not in the other. Suppose there is a n-vector $\alpha_2 = \alpha^1 \cup \alpha^2 \cup \ldots \cup \alpha^n$ which is in $W_2$ and not in $W_1$ i.e. $\alpha_2 \in W_2$ and $\alpha_2 \notin W_1$ by the earlier corollary just proved and the theorem there is a n-linear functional $f = f^1 \cup f^2 \cup \ldots \cup f^n$ such that $f(\beta) = 0$ for all $\beta$ in W, but $f(\alpha_2) \neq 0$. Then f is in $W_1^o$ but not in $W_2^o$ and $W_1^o \neq W_2^o$. Hence the claim.

Next we show a systematic method of finding the n-annihilator n-subspace spanned by a given finite n-set of n-vectors in $F_1^{n_1} \cup F_1^{n_2} \cup \ldots \cup F_n^{n_n}$. Consider a n-system of n-homogeneous linear n-equations, which will be from the point of view of the n-linear functionals. Suppose we have n-system of n-linear equations.

$$\begin{array}{l} A_{11}^1 x_1^1 + \ \ldots \ + A_{1n_1}^1 x_{n_1}^1 = 0 \\ \quad \vdots \qquad\qquad \vdots \\ A_{m_1 1}^1 x_1^1 + \ \ldots \ + A_{m_1 n_1}^1 x_{n_1}^1 = 0, \end{array}$$

$$\begin{array}{l} A_{11}^2 x_1^2 + \ \ldots \ + A_{1n_2}^2 x_{n_2}^2 = 0 \\ \quad \vdots \qquad\qquad \vdots \\ A_{m_2 1}^2 x_1^2 + \ \ldots \ + A_{m_2 n_2}^2 x_{n_2}^2 = 0, \end{array}$$

so on

$$\begin{array}{l} A_{11}^n x_1^n + \ \ldots \ + A_{1n_n}^n x_{n_n}^n = 0 \\ \quad \vdots \qquad\qquad \vdots \\ A_{m_n 1}^n x_1^n + \ \ldots \ + A_{m_n n_n}^n x_{n_n}^n = 0 \end{array}$$

for which we wish to find the solutions. If we let $f_i^k$, i = 1, 2, …n, $m_k$, k = 1, 2, …,n be the linear function on $F_k^{n_k}$, defined by $f_i^k(x_1^k,...,x_{n_k}^k) = A_{i1}^k x_1^k + ... + A_{in_k}^k x_{n_k}^k$, this is true for every k, k = 1, 2, …, n; then we are seeking the n-subspace of $F_1^{n_1} \cup F_2^{n_2} \cup ... \cup F_n^{n_n}$; that is all $\alpha = \alpha^1 \cup \alpha^2 \cup \ldots \cup \alpha^n$ such that $f_i^k(\alpha^k) = 0$, i = 1, 2, …, $m_k$ and k = 1, 2, …, n. In other words we are seeking the n-subspace annihilated by $\{f_1^1, f_2^1, ..., f_{m_1}^1\} \cup \{f_1^2, f_2^2, ..., f_{m_2}^2\} \cup \ldots \cup \{f_1^n, f_2^n, ..., f_{m_n}^n\}$. Row reduction of each of the coefficient matrix of the n-matrix provides us with a systematic method of finding this n-subspace. The ($n_1$, $n_2$, … , $n_n$) tuple $(A_{i1}^1,...,A_{in_1}^1) \cup (A_{i1}^2,...,A_{in_2}^2) \cup \ldots \cup (A_{i1}^n,...,A_{in_n}^n)$ gives the coordinates of the n-linear functional $f_i^k$, k = 1, 2, …, n relative to the n-basis which is n-dual to the standard n-basis of $F_1^{n_1} \cup F_2^{n_2} \cup ... \cup F_n^{n_n}$. The n-row space of the n-coefficient n-matrix may thus be regarded as the n-space of n-linear functionals spanned by $(f_1^1, f_2^1, ..., f_{m_1}^1) \cup (f_1^2, f_2^2, ..., f_{m_2}^2) \cup \ldots \cup (f_1^n, f_2^n, ..., f_{m_n}^n)$. The solution n-space is the n-subspace n-annihilated by this space of n-functionals.

Now one may find the n-system of equations from the n-dual point of view. This is suppose that we are given, $m_i$ n-vectors in $F_1^{n_1} \cup F_1^{n_2} \cup ... \cup F_n^{n_n}$;

$\alpha = \alpha_i^1 \cup ... \cup \alpha_i^n = (A_{i1}^1 A_{i2}^1 ... A_{in_1}^1) \cup \ldots \cup (A_{i1}^n A_{i2}^n ... A_{in_n}^n)$ and we find the n-annihilator of the n-subspace spanned by these vectors. A typical n-linear functional on $F_1^{n_1} \cup F_1^{n_2} \cup ... \cup F_n^{n_n}$ has the form $f^1(x_1^1,...,x_{n_1}^1) \cup f^2(x_1^2,...,x_{n_2}^2) \cup \ldots \cup f^n(x_1^n,...,x_{n_n}^n)$ = $((c_1^1 x_1^1 + ... + c_{n_1}^1 x_{n_1}^1) \cup (c_1^2 x_1^2 + ... + c_{n_2}^2 x_{n_2}^2) \cup \ldots \cup (c_1^n x_1^n + ... + c_{n_n}^n x_{n_n}^n)$ and the condition that $f^2 \cup \ldots \cup f^n$ be in this n-annihilator; that is $\sum_{j=1}^{n_1} A_{i_1 j}^1 c_j^1 \cup \sum_{j=1}^{n_2} A_{i_2 j}^2 c_j^2 \cup ... \cup \sum_{j=1}^{n_n} A_{i_n j}^n c_j^n = 0 \cup 0 \cup \ldots \cup 0$, and $1 \leq i_1 \leq m_1$ ,…, $1 \leq i_n \leq m_n$ that is $(c_1^1, \ldots, c_{n_1}^1) \cup \ldots \cup (c_1^n, \ldots, c_{n_n}^n)$ be the n-solution of the system $A^1X^1 \cup A^2X^2 \cup \ldots \cup A^nX^n = 0 \cup 0 \cup \ldots \cup 0$.

As in case of usual vector space we in case of the n-vector spaces of type II define the double dual. Here also it is important to mention in case of n-vector spaces of type I we cannot define dual or double dual. This is yet another difference between the n-vector spaces of type I and type II.

**DEFINITION 1.2.14:** *Let $V = V_1 \cup V_2 \cup \ldots \cup V_n$ be a n-vector space over the n-field $F = F_1 \cup F_2 \cup \ldots \cup F_n$ (i.e. each $V_i$ is a vector space over $F_i$). Let $V^* = V_1^* \cup V_2^* \cup \ldots \cup V_n^*$ be the n-vector space which is the n-dual of V over the same n-field $F = F_1 \cup F_2 \cup \ldots \cup F_n$. The n-dual of the n-dual space $V^*$ i.e., $V^{**}$ in terms of n-basis and n-dual basis is given in the following.*

*Let $\alpha = \alpha^1 \cup \alpha^2 \cup \ldots \cup \alpha^n$ be a n-vector in V, then x induces a n-linear functional $L_\alpha = L^1_{\alpha^1} \cup L^2_{\alpha^2} \cup \ldots \cup L^n_{\alpha^n}$ defined by $L_\alpha(f) = L^1_{\alpha^1}(f^1) \cup \ldots \cup L^n_{\alpha^n}(f^n)$ (where $f = f^1 \cup f^2 \cup \ldots \cup f^n$) i.e. $L_\alpha(f) = L^1_{\alpha^1}(f^1) \cup \ldots \cup L^n_{\alpha^n}(f^n) = f(\alpha) = f^1(\alpha^1) \cup f^2(\alpha^2) \cup \ldots \cup f^n(\alpha^n)$, $f \in V^* = V_1^* \cup \ldots \cup V_n^*$ $f^i \in V_i^*$ for i = 1, 2, …, n. The fact that each $L^i_{\alpha^i}$ is linear is just a reformulation of the definition of the linear operators in $V_i^*$ for each i = 1, 2, …, n.*

$$\begin{aligned} L_\alpha(cf + g) &= L^1_{\alpha^1}(c^1 f^1 + g^1) \cup \ldots \cup L^n_{\alpha^n}(c^n f^n + g^n) \\ &= (c^1 f^1 + g^1)(\alpha^1) \cup \ldots \cup (c^n f^n + g^n)(\alpha^n) \\ &= (c^1 f^1(\alpha^1) + g^1(\alpha^1)) \cup \ldots \cup (c^n f^n(\alpha^n) + g^n(\alpha^n)) \\ &= (c^1 L^1_{\alpha^1}(f^1) + L^1_{\alpha^1}(g^1)) \cup \ldots \cup (c^n L^n_{\alpha^n}(f^n) + L^n_{\alpha^n}(g^n)) \\ &= cL_\alpha(f) + L_\alpha(g). \end{aligned}$$

*If $V = V_1 \cup V_2 \cup \ldots \cup V_n$ is a finite $(n_1, n_2, \ldots, n_n)$ dimensional and $\alpha \neq 0$ then $L_\alpha \neq 0$, in other words there exits a n-linear functional $f = f^1 \cup f^2 \cup \ldots \cup f^n$ such that $f(x) \neq 0$ i.e. $f(\alpha) = f^1(\alpha^1) \cup \ldots \cup f^n(\alpha^n)$ for each $f^i(\alpha^i) \neq 0$, i = 1, 2, …,n.*

The proof is left for the reader, using the fact if we choose a ordered n-basis $B = \{\alpha^1_1, \alpha^1_2, \ldots, \alpha^1_{n_1}\} \cup \{\alpha^2_1, \alpha^2_2, \ldots, \alpha^2_{n_2}\} \cup \ldots \cup$ $\{\alpha_1^n, \alpha_2^n, \ldots, \alpha_{n_n}^n\}$ for V = $V_1 \cup V_2 \cup \ldots \cup V_n$ such that α = $\alpha_1^1, \alpha_1^2, \ldots, \alpha_1^n$ and let f be a n-linear functional which assigns to each n-vector in V its first coordinate in the n-ordered basis B.

We prove the following interesting theorem.

**THEOREM 1.2.25:** *Let $V = V_1 \cup V_2 \cup \ldots \cup V_n$ be a finite ($n_1, n_2, \ldots, n_n$) dimensional n-vector space over the n-field $F = F_1 \cup F_2 \cup \ldots \cup F_n$. For each n-vector $\alpha = \alpha^1 \cup \alpha^2 \cup \ldots \cup \alpha^n$ in V define $L_\alpha(f) = L^1_{\alpha^1}(f^1) \cup \ldots \cup L^n_{\alpha^n}(f^n) = f(\alpha) = f^1(\alpha^1) \cup \ldots \cup f^n(\alpha^n)$; f in V*. The mapping $\alpha \to L_\alpha$, is then an n-isomorphism of V onto V**.*

*Proof:* We showed that for each $\alpha = \alpha^1 \cup \alpha^2 \cup \ldots \cup \alpha^n$ in V, the n-function $L_\alpha = L^1_{\alpha^1} \cup \ldots \cup L^n_{\alpha^n}$ is n-linear. Suppose $\alpha = \alpha^1 \cup \alpha^2 \cup \ldots \cup \alpha^n$ and $\beta = \beta^1 \cup \beta^2 \cup \ldots \cup \beta^n$ are in $V_1 \cup V_2 \cup \ldots \cup V_n$ and $c = c^1 \cup c^2 \cup \ldots \cup c^n$ is in $F = F_1 \cup F_2 \cup \ldots \cup F_n$ and let $\gamma = c\alpha + \beta$ i.e., $\gamma_1 \cup \gamma_2 \cup \ldots \cup \gamma_n = (c^1\alpha^1 + \beta^1) \cup (c^2\alpha^2 + \beta^2) \cup \ldots \cup (c^n\alpha^n + \beta^n)$. Thus for each f in V*

$$
\begin{aligned}
L_\gamma(f) &= f(\gamma) \\
&= f(c\alpha + \beta) \\
&= f^1(c^1\alpha^1 + \beta^1) \cup \ldots \cup f^n(c^n\alpha^n + \beta^n) \\
&= [c^1f^1(\alpha^1) + f^1(\beta^1)] \cup \ldots \cup [c^nf^n(\alpha^n) + f^n(\beta^n)] \\
&= c^1L^1_{\alpha^1}(f^1) + L^1_{\beta^1}(f^1) \cup \ldots \cup c^nL^n_{\alpha^n}(f^n) + L^n_{\beta^n}(f^n) \\
&= cL_\alpha(f) + L_\beta(f) \\
L_\gamma &= cL_\alpha + L_\beta.
\end{aligned}
$$

This proves the n-map $\alpha \to L_\alpha$ is an n-linear transformation from V into V**. This n-transformation is n-nonsingular i.e. $L_\alpha = 0$ if and only if α = 0. Hence dim V** = dim V* = dim V= ($n_1, n_2, \ldots, n_n$).

The following are the two immediate corollaries to the theorem.

**COROLLARY 1.2.6:** *Let $V = V_1 \cup V_2 \cup \ldots \cup V_n$ be a ($n_1, n_2, \ldots, n_n$) finite dimensional n-vector space over the n-field $F = F_1 \cup$*

*$F_2 \cup \ldots \cup F_n$. If L is a n-linear functional on the dual space V* of V then there is a unique n-vector $\alpha = \alpha^1 \cup \alpha^2 \cup \ldots \cup \alpha^n$ in V such that L(f) = f (α) for every f in V*.*

The proof is left for the reader to prove.

**COROLLARY 1.2.7:** *If $V = V_1 \cup V_2 \cup \ldots \cup V_n$ is a $(n_1, n_2, \ldots, n_n)$ dimensional n-vector space over the n-field $F = F_1 \cup F_2 \cup \ldots \cup F_n$. Each n-basis for V* is the dual of some n-basis for V.*

*Proof:* Given $V = V_1 \cup V_2 \cup \ldots \cup V_n$ is a $(n_1, n_2, \ldots, n_n)$ vector space over the n-field $F = F_1 \cup F_2 \cup \ldots \cup F_n$. $V^* = V_1^* \cup \ldots \cup V_n^*$ is the dual n-space of V over the n-field F.
Let $B^* = \{f_1^1, \ldots, f_{n_1}^1\} \cup \ldots \cup \{f_1^n, \ldots, f_{n_n}^n\}$ be a n-basis for V* by an earlier theorem there is a n-basis $\{L_1^1, L_2^1, \ldots, L_{n_1}^1\} \cup \{L_1^2, L_2^2, \ldots, L_{n_2}^2\} \cup \ldots \cup \{L_1^n, L_2^n, \ldots, L_{n_n}^n\}$ for $V^{**} = V_1^{**} \cup \ldots \cup V_n^{**}$ such that $L_i^k(f_j^k) = \delta_{ij}$ for k = 1, 2, …, n, $1 \le i \le n_k$. Using the above corollary for each i, there is a vector $\alpha_i^k$ in $V_k$ such that $L_i^k(f^k) = f^k(\alpha_i^k)$ for every $f^k$ in $V_k^*$, such that $L_i^k = L_{\alpha_i^k}^k$. It follows immediately that $\{\alpha_1^1, \alpha_2^1, \ldots, \alpha_{n_1}^1\} \cup \{\alpha_1^2, \alpha_2^2, \ldots, \alpha_{n_2}^2\} \cup \ldots \cup \{\alpha_1^n, \alpha_2^n, \ldots, \alpha_{n_n}^n\}$ is a n-basis for V and B* is the n-dual of this n-basis.

In view of this we can say $(W^o)^o = W$.

Now we prove yet another interesting theorem.

**THEOREM 1.2.26:** *If S is any n-set of a finite $(n_1, n_2, \ldots, n_n)$ dimensional n-vector space $V = V_1 \cup V_2 \cup \ldots \cup V_n$ then $(S^o)^o$ is the n-subspace spanned by S.*

*Proof:* Let $W = W_1 \cup W_2 \cup \ldots \cup W_n$ be a n-subspace spanned by the n-set $S = S_1 \cup \ldots \cup S_n$ i.e. each $S_i$ spans $W_i$, i = 1, 2, … , n. Clearly $W^o = S^o$. Therefore what is left over to prove is that $W = W^{oo}$. We can prove this yet in another way. We know dim W + dim $W^o$ = dim V. dim $W^o$ + dim $(W^o)^o$ = dim V* since dim V = dim V* we have dim W + dim $W^o$ = dim $W^o$ + dim $W^{oo}$ which implies dim W = dim $W^{oo}$. Since W is a subspace of $W^{oo}$ we see that $W = W^{oo}$.

Let V be an n-vector space over the n-field of type II. We define an n-hyper subspace or n-hyperspace of V. We assume V is $(n_1, n_2, \ldots, n_n)$ dimension over $F = F_1 \cup F_2 \cup \ldots \cup F_n$. If N is a n-hyperspace of V i.e. N is of $(n_1 - 1, n_2 - 1, \ldots, n_n - 1)$ dimension over F then we can define N to be a n-hyperspace of V if

1. N is a proper n-subspace of V
2. If W is a n-subspace of V which contains N then either W = N or W = V.

Conditions (1) and (2) together say that N is a proper n-subspace and there is no larger proper n-subspace in short N is a maximal proper n-subspace of V.

Now we define n-hyperspace of a n-vector space.

**DEFINITION 1.2.15:** *If $V = V_1 \cup V_2 \cup \ldots \cup V_n$ is a n-vector space over the n-field $F = F_1 \cup F_2 \cup \ldots \cup F_n$ a n-hyperspace in V is a maximal proper n-subspace of V.*

We prove the following theorem on n-hyperspace of V.

**THEOREM 1.2.27:** *If $f = f^1 \cup \ldots \cup f^n$ is a non-zero n-linear functional on the n-vector space $V = V_1 \cup V_2 \cup \ldots \cup V_n$ of type II over the n-field $F = F_1 \cup F_2 \cup \ldots \cup F_n$, then the n-hyperspace in V is the n-null space of a (not unique) non-zero n-linear functional on V.*

*Proof:* Let $f = f^1 \cup f^2 \cup \ldots \cup f^n$ be a non zero n-linear functional on $V = V_1 \cup V_2 \cup \ldots \cup V_n$ and $N_f = N^1_{f^1} \cup N^2_{f^2} \cup \ldots \cup N^n_{f^n}$ its n-null space. Let $\alpha = \alpha^1 \cup \alpha^2 \cup \ldots \cup \alpha^n$ be a n-vector in $V = V^1 \cup V^2 \cup \ldots \cup V^n$ which is not in $N_f$ i.e., a n-vector such that $f(\alpha) \neq 0 \cup \ldots \cup 0$. We shall show that every n-vector in V is in the n-subspace spanned by $N_f$ and α. That n-subspace consist of all n-vector $\gamma + c\alpha$ where $\gamma = \gamma^1 \cup \ldots \cup \gamma^n$ and $c = c^1 \cup c^2 \cup \ldots \cup c^n$, $\gamma$ in $N_f$, c in $F = F^1 \cup F^2 \cup \ldots \cup F^n$. Let $\beta = \beta^1 \cup \beta^2 \cup \ldots \cup \beta^n$ be in V. Define

$$c = \frac{f(\beta)}{f(\alpha)} = \frac{f^1(\beta^1)}{f^1(\alpha^1)} \cup \ldots \cup \frac{f^n(\beta^n)}{f^n(\alpha^n)},$$

i.e., each

$$c^i = \frac{f^i(\beta^i)}{f^i(\alpha^i)}, \; i = 1, 2, \ldots, n;$$

which makes sense because each $f^i(\alpha^i) \neq 0$, i = 1, 2,…, n i.e., $f(\alpha) \neq 0$. Then the n-vector $\gamma = \beta - c\alpha$ is in $N_f$ since $f(\gamma) = f(\beta - c\alpha) = f(\beta) - cf(\alpha) = 0$. So $\beta$ is in the n-subspace spanned by $N_f$ and $\alpha$. Now let N be the n-hyperspace in V. For some n-vector $\alpha = \alpha^1 \cup \alpha^2 \cup \ldots \cup \alpha^n$ which is not in N. Since N is a maximal proper n-subspace, the n-subspace spanned by N and $\alpha$ is the entire n-space V. Therefore each n-vector $\beta$ in V has the form $\beta = \gamma + c\alpha$, $\gamma = \gamma_1 \cup \ldots \cup \gamma_n$ in N. $c = c_1 \cup \ldots \cup c_n$ in $F = F_1 \cup \ldots \cup F_n$ and $N = N_1 \cup N_2 \cup \ldots \cup N_n$ where each $N_i$ is a maximal proper subspace of $V_i$ for i = 1, 2, …, n. The n-vector $\gamma$ and the n-scalars c are uniquely determined by $\beta$. If we have also $\beta = \gamma^1 + c_1\alpha$, $\gamma^1$ in N, $c_1$ in F then $(c_1 - c)\alpha = \gamma - \gamma^1$ if $c_1 - c \neq 0$ then $\alpha$ would be in N, hence $c = c_1$ and $\gamma = \gamma^1$, i.e., if $\beta$ is in V there is a unique n-scalar c such that $\beta - c\alpha$ is in N. Call the n-scalar $g(\beta)$. It is easy to see g is an n-linear functional on V and that N is a n-null space of g.

Now we state a lemma and the proof is left for the reader.

**LEMMA 1.2.2:** *If g and h are n-linear functionals on a n-vector space V then g is a n-scalar space V then g is a n-scalar multiple of f if and only if the n-null space of g contains the n-null space of f that is if and only if $f(x) = 0 \cup \ldots \cup 0$ implies $g(x) = 0 \cup \ldots \cup 0$.*

We prove the following interesting theorem for n-linear functional, on an n-vector space V over the n-field F.

**THEOREM 1.2.28:** *Let $V = V_1 \cup V_2 \cup \ldots \cup V_n$ be an n-vector space over the n-field $F = F^1 \cup F^2 \cup \ldots \cup F^n$. If $\{g^1, f_1^1, \ldots, f_{r_1}^1\} \cup \{g^2, f_1^2, \ldots, f_{r_2}^2\} \cup \ldots \cup \{g^n, f_1^n, \ldots, f_{r_n}^n\}$ be n-linear functionals on the space V with respective n-null spaces $\{N^1, N_1^1, \ldots, N_{r_1}^1\} \cup \{N^2, N_1^2, \ldots, N_{r_2}^2\} \cup \ldots \cup \{N^n, N_1^n, \ldots, N_{r_n}^n\}$. Then $(g^1 \cup g^2 \cup \ldots \cup g^n)$ is a linear combination of $\{f_1^1, \ldots, f_{r_1}^1\} \cup \ldots \cup \{f_1^n, f_2^n, \ldots, f_{r_n}^n\}$ if and only if $N^1 \cup N^2 \cup \ldots \cup N^n$ contains the n-intersection $\{N_1^1 \cap \ldots \cap N_{r_1}^1\} \cup \{N_1^2 \cap \ldots \cap N_{r_2}^2\} \cup \ldots \cup \{N_1^n \cap \ldots \cap N_{r_n}^n\}$.*

*Proof:* We shall prove the result for $V_i$ of $V_1 \cup V_2 \cup \ldots \cup V_n$ and since $V_i$ is arbitrary, the result we prove is true for every i, i = 1, 2, ..., n. Let $g^i = c_1^i f_1^i + \ldots + c_{r_i}^i f_{r_i}^i$ and $f_j^i(\alpha^i) = 0$ for each j the clearly $g^i(\alpha) = 0$. Therefore $N^i$ contains $N_1^i \cap N_2^i \cap \ldots \cap N_{r_i}^i$.

We shall prove the converse by induction on the number $r_i$. The proceeding lemma handles the case $r_i = 1$. Suppose we know the result for $r_i = k_i-1$ and let $f_1^i, f_2^i, \ldots, f_{k_i}^i$ be the linear functionals with null spaces $N_1^i, \ldots, N_{k_i}^i$ such that $N_1^i \cap \ldots \cap N_{k_i}^i$ is contained in $N_i$, the null space of $g^i$. Let $(g^i)'$ $(f_1^i)'$ … $(f_{k_i-1}^i)'$ be the restrictions of $g^i$, $f_1^i$, … , $f_{k_i-1}^i$ to the subspace $N_{k_i}^i$. Then $(g^i)'$, $(f_1^i)'$, … , $(f_{k_i-1}^i)'$ are linear functionals on the vector space $N_{k_i}^i$. Further more if $\alpha^i$ is a vector in $N_{k_i}^i$ and $((f_j^i(\alpha^i))' = 0$, j = 1, …, $k_i - 1$ then $\alpha^i$ is in $N_1^i \cap N_2^i \cap \ldots \cap N_{k_i}^i$ and so $(g^i(\alpha^i))' = 0$. By induction hypothesis (the case $r_i = k_i - 1$) there are scalars $c_j^i$ such that $(g^i)' = c_1^i(f_1^i)' + \ldots + c_{k_i-1}^i(f_{k_i-1}^i)'$. Now let $h^i = g^i - \sum_{j=1}^{k_i-1} c_j^i f_j^i$. Then $h^i$ is a linear functional on $V_i$ and this tells $h^i(\alpha^i) = 0$ for every $\alpha^i$ in $N_{k_i}^i$. By the proceeding lemma $h^i$ is a scalar multiple of $f_{k_i}$ if $h^i = c_{k_i}^i f_{k_i}^i$ then $g^i = \sum_{j=1}^{k_i} c_j^i f_j^i$.

Now the result is true for each i, i = 1, 2, …, n hence the theorem.

Now we proceed onto define the notion of n-transpose of a n-linear transformation $T = T_n \cup \ldots \cup T_n$ of the n-vector spaces V and W of type II.

Suppose we have two n-vector spaces $V = V^1 \cup V^2 \cup \ldots \cup V^n$ and $W = W^1 \cup W^2 \cup \ldots \cup W^n$ over the n-field $F = F^1 \cup F^2 \cup \ldots \cup F^n$. Let $T = T^1 \cup \ldots \cup T^n$ be a n-linear transformation from V into W. Then T induces a n-linear transformation from W* into V* as follows:

Suppose $g = g^1 \cup \ldots \cup g^n$ is a n-linear functional on $W = W^1 \cup \ldots \cup W^n$ and let $f(\alpha) = f_1(\alpha^1) \cup \ldots \cup f_n(\alpha^n)$ (where $\alpha = \alpha^1 \cup \alpha^2 \cup \ldots \cup \alpha^n \in V$) $f(\alpha) = g(T\alpha)$, that is $f(\alpha) = g^1(T_1\alpha^1) \cup \ldots \cup g^n(T_n\alpha^n)$ ; for each $\alpha^i \in V_i$; i = 1,2, … , n. The above equation defines a n-function f from V into F i.e. $V^1 \cup V^2 \cup \ldots \cup V^n$ into $F^1 \cup \ldots \cup F^n$ namely the n-composition of T, a n-function from V into W with g a n-function from W into $F = F^1 \cup F^2 \cup \ldots \cup F^n$. Since both T and g are n-linear f is also n-linear i.e., f is an n-linear functional on V. This T provides us with a rule $T^t = T_1^t \cup \ldots \cup T_n^t$ which associates with each n-linear functional g on $W = W^1 \cup \ldots \cup W^n$ a n-linear functional $f = T^tg$ i.e. $f_1 \cup f_2 \cup \ldots \cup f_n = T_1g^1 \cup \ldots \cup T_ng^n$ on $V = V^1 \cup V^2 \cup \ldots \cup V^n$, i.e., $f_i = T_i^tg^i$ is a linear functional on $V_i$. $T^t = T_1^t \cup \ldots \cup T_n^t$ is actually a n-linear transformation from $W^* = W_1^* \cup \ldots \cup W_n^*$ into $V^* = V_1^* \cup \ldots \cup V_n^*$ for if $g_1, g_2$, are in W i.e., $g_1 = g_1^1 \cup \ldots \cup g_1^n$ and $g_2 = g_2^1 \cup \ldots \cup g_2^n$ in W* and $c = c^1 \cup \ldots \cup c^n$ is a n-scalar

$$T^t(cg_1 + g_2)(\alpha) = c(T^tg_1)(\alpha) + (T^tg_2)(\alpha)$$

$$\text{i.e. } T_1^t(c^1g_1^1 + g_2^1)\alpha^1 \cup \ldots \cup T_n^t(c^ng_1^n + g_2^n) =$$

$$(c^1(T_1^tg_1^1)\alpha^1 + (T_1^tg_2^1)\alpha^1) \cup \ldots \cup (c^n(T_n^tg_1^n)\alpha^n + (T_n^tg_2^n)\alpha^n),$$ so that $T^t(cg_1 + g_2) = T^t(cg_1 + g_2) = cT_{g_1}^t + T_{g_2}^t$ ; this can be summarized into the following theorem.

**THEOREM 1.2.29:** *Let $V = V^1 \cup V^2 \cup \ldots \cup V^n$ and $W = W^1 \cup W^2 \cup \ldots \cup W^n$ be n-vector spaces over the n-field $F = F^1 \cup F^2$*

*$\cup \ldots \cup F^n$. For each n-linear transformation $T = T_1 \cup \ldots \cup T_n$ from V into W there is unique n-linear transformation $T^t = T_1^t \cup \ldots \cup T_n^t$ from $W^* = W_1^* \cup \ldots \cup W_n^*$ into $V^* = V_1^* \cup \ldots \cup V_n^*$ such that $(T_g^t)\alpha = g(T\alpha)$ for every g in W* and α in V.*

*We call $T^t = T_1^t \cup \ldots \cup T_n^t$ as the n-transpose of T. This n-transformation $T^t$ is often called the n-adjoint of T.*

It is interesting to see that the following important theorem.

**THEOREM 1.2.30:** *Let $V = V^1 \cup V^2 \cup \ldots \cup V^n$ and $W = W^1 \cup W^2 \cup \ldots \cup W^n$ be n-vector spaces over the n-field $F = F^1 \cup F^2 \cup \ldots \cup F^n$ and let $T = T^1 \cup \ldots \cup T^n$ be a n-linear transformation from V into W. The n-null space of $T^t = T_1^t \cup \ldots \cup T_n^t$ is the n-annihilator of the n-range of T. If V and W are finite dimensional then*

*(I) n-rank $(T^t)$ = n-rank T*
*(II) The n-range of $T^t$ is the annihilator of the n-null space of T.*

*Proof:* Let $g = g^1 \cup g^2 \cup \ldots \cup g^n$ be in $W^* = W_1^* \cup \ldots \cup W_n^*$ then by definition $(T^t g)\alpha = g(T\alpha)$ for each $\alpha = \alpha^1 \cup \ldots \cup \alpha^n$ in V. The statement that g is in the n-null space of $T^t$ means that $g(T\alpha) = 0$ i.e., $g^1T_1\alpha^1 \cup \ldots \cup g^nT_n\alpha^n = (0 \cup \ldots \cup 0)$ for every $\alpha \in V = V_1 \cup V_2 \cup \ldots \cup V_n$.

Thus the n-null space of $T^t$ is precisely the n-annihilator of n-range of T. Suppose V and W are finite dimensional, say dim V = $(n_1, n_2, \ldots, n_n)$ and dim W = $(m_1, m_2, \ldots, m_n)$. For (I), let $r = (r_1, r_2, \ldots, r_n)$ be the n-rank of T i.e. the dimension of the n-range of T is $(r_1, r_2, \ldots, r_n)$.

By earlier results the n-annihilator of the n-range of T has dimension $(m_1 - r_1, m_2 - r_2, \ldots, m_n - r_n)$. By the first statement of this theorem, the n-nullity of $T^t$ must be $(m_1 - r_1, \ldots, m_n - r_n)$. But since $T^t$ is a n-linear transformation on an $(m_1, \ldots, m_n)$ dimensional n-space, the n-rank of $T^t$ is $(m_1 - (m_1 - r_1), m_2 - (m_2 - r_2), \ldots, m_n - (m_n - r_n))$ and so T and $T^t$ have the same n-rank.

For (II), let $N = N^1 \cup \ldots \cup N^n$ be the n-null space of T. Every n-function in the n-range of $T^t$ is in the n-annihilator of N, for suppose $f = T^t g$ i.e. $f^1 \cup \ldots \cup f^n = T^t_1 g^1 \cup \ldots \cup T^t_n g^n$ for some g in W* then if α is in N; $f(\alpha) = f^1(\alpha^1) \cup \ldots \cup f^n(\alpha^n) = (T^t g)\alpha = (T^t_1 g^1)\alpha^1 \cup \ldots \cup (T^t_n g^n)\alpha^n = g(T\alpha) = g^1(T_1\alpha^1) \cup \ldots \cup g^n(T_n\alpha^n) = g(0) = g^1(0) \cup \ldots \cup g^n(0) = 0 \cup 0 \cup \ldots \cup 0$. Now the n-range of $T^t$ is an n-subspace of the space $N^o$ and dim $N^o = (n_1\text{-dim } N^1) \cup (n_2\text{-dim}N^2) \cup \ldots \cup (n_n\text{-dim } N^n)$ = n-rank T = n-rank $T^t$ so that the n-range of $T^t$ must be exactly $N^o$.

**THEOREM 1.2.31:** *Let $V = V_1 \cup V_2 \cup \ldots \cup V_n$ and $W = W_1 \cup W_2 \cup \ldots \cup W_n$ be two $(n_1, n_2, \ldots, n_n)$ and $(m_1, m_2, \ldots, m_n)$ dimensional n-vector spaces over the n-field $F = F_1 \cup F_2 \cup \ldots \cup F_n$. Let B be an ordered n-basis of V and B* the dual n-basis of V*. Let C be an ordered n-basis of W with dual n-basis of B*. Let $T = T_1 \cup T_2 \cup \ldots \cup T_n$ be a n-linear transformation from V into W, let A be the n-matrix of T relative to B and C and let B be a n-matrix of $T^t$ relative to B*, C*. Then $B^k_{ij} = A^k_{ij}$, for k = 1, 2, ..., n. i.e. $A^1_{ij} \cup A^2_{ij} \cup \ldots \cup A^n_{ij} = B^1_{ij} \cup B^2_{ij} \cup \ldots \cup B^n_{ij}$.*

*Proof:* Let B = $\{\alpha^1_1, \alpha^1_2, \ldots, \alpha^1_{n_1}\} \cup \{\alpha^2_1, \alpha^2_2, \ldots, \alpha^2_{n_2}\} \cup \ldots \cup \{\alpha^n_1, \alpha^n_2, \ldots, \alpha^n_{n_n}\}$ be a n-basis of V. The dual n-basis of B, B* = $(f^1_1, f^1_2, \ldots, f^1_{n_1}) \cup (f^2_1, f^2_2, \ldots, f^2_{n_2}) \cup \ldots \cup (f^n_1, f^n_2, \ldots, f^n_{n_n})$. Let C = $(\beta^1_1, \beta^1_2, \ldots, \beta^1_{m_1}) \cup (\beta^2_1, \beta^2_2, \ldots, \beta^2_{m_2}) \cup \ldots \cup (\beta^n_1, \beta^n_2, \ldots, \beta^n_{m_n})$ be an n-basis of W and C* = $(g^1_1, g^1_2, \ldots, g^1_{m_1}) \cup (g^2_1, g^2_2, \ldots, g^2_{m_2}) \cup \ldots \cup (g^n_1, g^n_2, \ldots, g^n_{m_n})$ be a dual n-basis of C. Now by definition for $\alpha = \alpha^1 \cup \ldots \cup \alpha^n$

$$T_k \alpha^k_j = \sum_{i=1}^{m_k} A^k_{ij} \beta^k_i \; ; \qquad j = 1, 2, \ldots, n_k;\; k = 1, 2, \ldots, n$$

$$T^t_k g^k_j = \sum_{i=1}^{n_k} B^k_{ij} f^k_i \; ; \qquad j = 1, 2, \ldots, m_k;\; k = 1, 2, \ldots, n.$$

Further

$$(T_k^t g_j^k)(\alpha_i^k) = g_j^k(T_k^t \alpha_i^k)$$

$$= g_j^k\ (\sum_{p=1}^{m_k} A_{pi}^k \beta_p^k)$$

$$= \sum_{p=1}^{m_k} A_{pi}^k g_j^k(\beta_p^k) = \sum_{p=1}^{m_k} A_{pi}\delta_{jp} = A_{ji}^k .$$

For any n-linear functional $f = f^1 \cup \ldots \cup f^n$ on V

$$f^k = \sum_{i=1}^{m_k} f^k(\alpha_i^k) f_i^k\ ; \quad k = 1, 2, \ldots, n.$$

If we apply this formula to the functional $f^k = T_k^t g_j^k$ and use the fact $(T_k^t g_j^k)(\alpha_i^k) = A_{ji}^k$, we have $T_k^t g_j^k = \sum_{i=1}^{n_k} A_{ji}^k f_i^k$ from which it follows $B_{ij}^k = A_{ij}^k$; true for k = 1, 2, …, n i.e., $B_{ij}^1 \cup B_{ij}^2 \cup \ldots \cup B_{ij}^n$ = $A_{ij}^1 \cup A_{ij}^2 \cup \ldots \cup A_{ij}^n$. If $A = A^1 \cup A^2 \cup \ldots \cup A^n$ is a $(m_1 \times n_1, m_2 \times n_2, , \ldots, m_n \times n_n)$ n-matrix over the n-field $F = F_1 \cup F_2 \cup \ldots \cup F_n$, the n-transpose of A is the $(n_1 \times m_1, n_2 \times m_2, \ldots, n_n \times m_n)$ matrix $A^t$ defined by

$$(A_{ij}^1)^t \cup (A_{ij}^2)^t \cup \ldots \cup (A_{ij}^n)^t = A_{ij}^1 \cup A_{ij}^2 \cup \ldots \cup A_{ij}^n .$$

We leave it for the reader to prove the n-row rank of A is equal to the n-column rank of A i.e. for each matrix $A^i$ we have the column rank of $A^i$ to be equal to the row rank of $A^i$; i = 1, 2, …, n.

Now we proceed on to define the notion of n-linear algebra over a n-field of type II.

**DEFINITION 1.2.16:** *Let $F = F_1 \cup F_2 \cup \ldots \cup F_n$, be a n-field. The n-vector space, $A = A^1 \cup A^2 \cup \ldots \cup A^n$ over the n-field F of type II is said to be a n-linear algebra over the n-field F if each $A_i$ is a linear algebra over $F_i$ for i = 1, 2, …, n i.e., for $\alpha, \beta \in A_i$ we have a vector $\alpha\beta \in A_i$, called the product of α and β in such a way that*

a. *multiplication is associative* $\alpha(\beta\gamma) = (\alpha\beta)\gamma$ *for* $\alpha, \beta, \gamma \in A_i$ *for* $i = 1, 2, \ldots, n$.
b. *multiplication is distributive with respect to addition,* $\alpha(\beta + \gamma) = \alpha\beta + \alpha\gamma$ *and* $(\alpha + \beta)\gamma = \alpha\gamma + \beta\gamma$ *for* $\alpha, \beta, \gamma \in A_i$ *for* $i = 1, 2, \ldots, n$.
c. *for each scalar* $c_i \in F_i$, $c_i(\alpha\beta) = (c_i\alpha)\beta = \alpha(c_i\beta)$, *true for* $i = 1, 2, \ldots, n$.

If there is an element $1_n = 1 \cup 1 \cup \ldots \cup 1$ in A such that $1_n\alpha = \alpha 1_n = \alpha$ i.e., $1_n\alpha = (1 \cup 1 \cup \ldots \cup 1)(\alpha_1 \cup \ldots \cup \alpha_n) = \alpha_1 \cup \ldots \cup \alpha_n = \alpha$, $\alpha I_n = (\alpha_1 \cup \ldots \cup \alpha_n)\, 1_n = (\alpha_1 \cup \ldots \cup \alpha_n)(1 \cup 1 \cup \ldots \cup 1) = \alpha_1 \cup \ldots \cup \alpha_n = \alpha$. We call A an n-linear algebra with n-identity over the n-field F. If even one of the $A_i$'s do not contain identity then we say A is a n-linear algebra without an n-identity. $1_n = 1 \cup 1 \cup \ldots \cup 1$ is called the n-identity of A.

The n-algebra A is n-commutative if $\alpha\beta = \beta\alpha$ for all $\alpha, \beta \in A$ i.e. if $\alpha = \alpha_1 \cup \ldots \cup \alpha_n$ and $\beta = \beta_1 \cup \ldots \cup \beta_n$, $\alpha\beta = (\alpha_1 \cup \ldots \cup \alpha_n)(\beta_1 \cup \ldots \cup \beta_n) = \alpha_1\beta_1 \cup \ldots \cup \alpha_n\beta_n$. $\beta\alpha = (\beta_1 \cup \ldots \cup \beta_n)(\alpha_1 \cup \ldots \cup \alpha_n) = \beta_1\alpha_1 \cup \ldots \cup \beta_n\alpha_n$. If each $\alpha_i\beta_i = \beta_i\alpha_i$ for $i = 1, 2, \ldots, n$ then we say $\alpha\beta = \beta\alpha$ for every $\alpha, \beta \in A$. We call A in which $\alpha\beta = \beta\alpha$ to be an n-commutative n-linear algebra. Even if one $A^i$ in A is non commutative we don't call A to be an n-commutative n-linear algebra.

1. All n-linear algebras over the n-field F are n-vector spaces over the n-field F.

2. Every n-linear algebra over an n-field F need not be a n-commutative n-linear algebra over F.

3. Every n-linear algebra A over an n-field F need not be a n-linear algebra with n-identity $1_n$ in A.

Now we proceed on to define the notion of n-polynomial over the n-field $F = F_1 \cup F_2 \cup \ldots \cup F_n$.

**DEFINITION 1.2.17:** *Let* $F[x] = F_1[x] \cup F_2[x] \cup \ldots \cup F_n[x]$ *be such that each* $F_i[x]$ *is a polynomial over* $F_i$, $i = 1, 2, \ldots, n$, *where* $F = F_1 \cup F_2 \cup \ldots \cup F_n$ *is the n-field. We call F[x] the n-polynomial over the n-field* $F = F_1 \cup \ldots \cup F_n$. *Any element p(x)* $\in$ *F[x] will be of the form* $p(x) = p_1(x) \cup p_2(x) \cup \ldots \cup p_n(x)$ *where* $p_i(x)$ *is a polynomial in* $F_i[x]$ *i.e.,* $p_i(x)$ *is a polynomial in the variable x with coefficients from* $F_i$*; i = 1, 2, ..., n.*

*The n-degree of p(x) is a n-tuple given by* $(n_1, n_2, \ldots, n_n)$ *where* $n_i$ *is the degree of the polynomial* $p_i(x)$*; i = 1, 2, ..., n.*

We illustrate the n-polynomial over the n-field by an example.

***Example 1.2.8:*** Let $F = Z_2 \cup Z_3 \cup Z_5 \cup Q \cup Z_{11}$ be a 5-field; $F[x] = Z_2[x] \cup Z_3[x] \cup Z_5[x] \cup Q[x] \cup Z_{11}[x]$ is a 5-polynomial vector space over the 5-field F. $p(x) = x^2 + x + 1 \cup x^3 + 2x^2 + 1 \cup 4x^3 + 2x + 1 \cup 81x^7 + 50x^2 - 3x + 1 \cup 10x^3 + 9x^2 + 7x + 1 \in F[x]$.

**DEFINITION 1.2.18:** *Let* $F[x] = F_1[x] \cup F_2[x] \cup \ldots \cup F_n[x]$ *be a n-polynomial over the n-field* $F = F_1 \cup F_2 \cup \ldots \cup F_n$. *F[x] is a n-linear algebra over the n-field F. Infact F[x] is a n-commutative n-algebra over the n-field F. F[x] the n-linear algebra has the n-identity,* $1_n = 1 \cup 1 \cup \ldots \cup 1$. *We call a n-polynomial* $p(x) = p_1(x) \cup p_2(x) \cup \ldots \cup p_n(x)$ *to be a n-monic n-polynomial if each* $p_i(x)$ *is a monic polynomial in x for i = 1, 2, ..., n.*

The reader is expected to prove the following theorem

**THEOREM 1.2.32:** *Let* $F[x] = F_1[x] \cup F_2[x] \cup \ldots \cup F_n[x]$ *be a n-linear algebra of n-polynomials over the n-field* $F = F_1 \cup F_2 \cup \ldots \cup F_n$. *Then*

a. *For f and g two non zero n-polynomials in F[x] where* $f(x) = f_1(x) \cup f_2(x) \cup \ldots \cup f_n(x)$ *and* $g(x) = g_1(x) \cup g_2(x) \cup \ldots \cup g_n(x)$ *the polynomial* $fg = f_1(x)g_1(x) \cup f_2(x)g_2(x) \cup \ldots \cup f_n(x)g_n(x)$ *is a non zero n-polynomial of F[x]*

b. *n-deg (fg) = n-deg f + n-deg g where n-deg f =* $(n_1, n_2, \ldots, n_n)$ *and n-deg g =* $(m_1, m_2, \ldots, m_n)$

c. *fg is n-monic if both f and g are n-monic*

d. *fg is a n-scalar n-polynomial if and only if f and g are scalar polynomials*

e. *If $f + g \neq 0$, n-deg $(f + g) \leq$ max (n-deg f, n-deg g)*

f. *If f, g, h are n-polynomials over the n-field $F = F_1 \cup F_2 \cup \ldots \cup F_n$ such that $f = f_1(x) \cup f_2(x) \cup \ldots \cup f_n(x) \neq 0 \cup 0 \cup \ldots \cup 0$ and $fg = fh$ then $g = h$, where $g(x) = g_1(x) \cup g_2(x) \cup \ldots \cup g_n(x)$ and $h(x) = h_1(x) \cup h_2(x) \cup \ldots \cup h_n(x)$.*

As in case of usual polynomials we see in case of n-polynomials the following. Let $A = A_1 \cup A_2 \cup \ldots \cup A_n$ be a n-linear algebra with identity $1_n = 1 \cup 1 \cup \ldots \cup 1$ over the n-field $F = F_1 \cup F_2 \cup \ldots \cup F_n$, where we make the convention for any $\alpha = \alpha_1 \cup \alpha_2 \cup \ldots \cup \alpha_n \in A$; $\alpha^{o} = \alpha_1^{o} \cup \alpha_2^{o} \cup \ldots \cup \alpha_n^{o} = 1 \cup 1 \cup \ldots \cup 1 = 1_n$ for each $\alpha \in A$.

Now to each n-polynomial $f(x) \in F[x] = F_1[x] \cup F_2[x] \cup \ldots \cup F_n[x]$ over the n-field $F = F_1 \cup F_2 \cup \ldots \cup F_n$ and $\alpha = \alpha_1 \cup \alpha_2 \cup \ldots \cup \alpha_n$ in A we can associate an element

$$F(\alpha) = \sum_{i=0}^{n_1} f_i^1 \alpha_i^1 \cup \sum_{i=0}^{n_2} f_i^2 \alpha_i^2 \cup \ldots \cup \sum_{i=0}^{n_n} f_i^n \alpha_i^n \ ; \ f_i^k \in F_k \ ;$$

for $k = 1, 2, \ldots, n$; $1 \leq i \leq n_k$.

We can in view of this prove the following theorem.

**THEOREM 1.2.33:** *Let $F = F_1 \cup F_2 \cup \ldots \cup F_n$ be a n-field and $A = A_1 \cup A_2 \cup \ldots \cup A_n$ be a n-linear algebra with identity $1_n = 1 \cup 1 \cup \ldots \cup 1$ over the n-field F.*

*Suppose f(x) and g(x) be n-polynomials in $F[x] = F_1[x] \cup F_2[x] \cup \ldots \cup F_n[x]$ over the n-field $F = F_1 \cup F_2 \cup \ldots \cup F_n$ and that $\alpha = \alpha_1 \cup \ldots \cup \alpha_n \in A$ set by the rule*

$$f(\alpha) = \sum_{i=0}^{n_1} f_i^1(\alpha_1^i) \cup \sum_{i=0}^{n_2} f_i^2(\alpha_2^i) \cup \ldots \cup \sum_{i=0}^{n_n} f_i^n(\alpha_n^i)$$

and for $c = c^1 \cup c^2 \cup \ldots \cup c^n \in F = F_1 \cup F_2 \cup \ldots \cup F_n$, we have $(cf + g)\alpha = cf(\alpha) + g(\alpha)$; $(fg)(\alpha) = f(\alpha)g(\alpha)$.

*Proof:*

$$f(x) = \sum_{i=0}^{m_1} f_i^1 x^i \cup \ldots \cup \sum_{i=0}^{m_n} f_i^n x^i$$

and

$$g(x) = \sum_{j=0}^{n_1} g_j^1 x^j \cup \ldots \cup \sum_{j=0}^{n_n} g_j^n x^j$$

$$fg = \sum_{i,j} f_i^1 g_j^1 x^{i+j} \cup \sum_{i,j} f_i^2 g_j^2 x^{i+j} \cup \ldots \cup \sum_{i,j} f_i^n g_j^n x^{i+j}$$

and hence

$$(fg)\alpha = \sum_{i,j} f_i^1 g_j^1 \alpha_1^{i+j} \cup \sum_{i,j} f_i^2 g_j^2 \alpha_2^{i+j} \cup \ldots \cup \sum_{i,j} f_i^n g_j^n \alpha_n^{i+j}$$

($\alpha = \alpha^1 \cup \alpha^2 \cup \ldots \cup \alpha^n \in A$, $\alpha_i \in A_i$ ; for $i = 1, 2, \ldots, n$)

$$= \left(\sum_{i=0}^{m_1} f_i^1 \alpha_1^i\right)\left(\sum_{j=0}^{n_1} g_j^1 \alpha_1^j\right) \cup \left(\sum_{i=0}^{m_2} f_i^2 \alpha_2^i\right)\left(\sum_{j=0}^{n_2} g_j^2 \alpha_2^j\right) \cup \ldots \cup \left(\sum_{i=0}^{m_n} f_i^n \alpha_n^i\right)\left(\sum_{j=0}^{n_n} g_j^n \alpha_n^j\right)$$

$= f(\alpha)g(\alpha)$

$= f_1(\alpha_1)g_1(\alpha_1) \cup f_2(\alpha_2)g_2(\alpha_2) \cup \ldots \cup f_n(\alpha_n)g_n(\alpha_n)$.

Now we define the Lagranges n-interpolation formula.

Let $F = F_1 \cup F_2 \cup \ldots \cup F_n$ be a n-field and let $t_0^1, t_1^1, \ldots, t_{n_1}^1$ be $n_1 + 1$ distinct elements of $F_1$, $t_0^2, t_1^2, \ldots, t_{n_2}^2$ are $n_2 + 1$ distinct elements of $F_2$ , …, $t_0^n, t_1^n, \ldots, t_{n_n}^n$ are the $n_n + 1$ distinct elements of $F_n$.

Let $V = V_1 \cup V_2 \cup \ldots \cup V_n$ be a n-subspace of $F[x] = F_1[x] \cup F_2[x] \cup \ldots \cup F_n[x]$ consisting of all n-polynomial of n-degree less than or equal to $(n_1, \ldots, n_n)$ together with the n-zero polynomial and let $L = L_i^1, L_i^2, \ldots, L_i^n$ be the n-function from $V = V_1 \cup V_2 \cup \ldots \cup V_n$ into $F = F_1 \cup F_2 \cup \ldots \cup F_n$ defined by $L^i(f) = L_1^{i_1}(f^1) \cup \ldots \cup L_n^{i_n}(f^n) = f^1(t_{i_1}^1) \cup \ldots \cup f^n(t_{i_n}^n)$, $0 \le i_1 \le n_1$, $0 \le i_2 \le n_2, \ldots, 0 \le i_n \le n_n$.

By the property $(cf + g)(\alpha) = cf(\alpha) + g(\alpha)$ i.e., if $c = c^1 \cup \ldots \cup c^n$, $f = f^1 \cup \ldots \cup f^n$; $g = g^1 \cup \ldots \cup g^n$ and $\alpha = \alpha_1 \cup \ldots \cup \alpha_n$ then $c_1f^1(\alpha_1) \cup c_2f^2(\alpha_2) \cup \ldots \cup c_nf^n(\alpha_n) + g^1(\alpha_1) \cup g^2(\alpha_2) \cup \ldots \cup g^n(\alpha_n) = (c_1f^1 + g^1)(\alpha_1) \cup \ldots \cup (c_nf^n + g^n)(\alpha_n)$ each $L^i$ is a n-linear functional on V and one of the things we intend to show is that the set consisting of $L^0, L^1, \ldots, L^n$ is a basis for V* the dual space of V. $L^0 = L_1^0(f^1) \cup \ldots \cup L_n^0(f^n) = f^1(t_0^1) \cup \ldots \cup f^n(t_0^n)$ and so on we know from earlier results $(L^0, L^1, \ldots, L^n)$ that is $\{L_1^0, \ldots, L_1^{n_1}\} \cup \{L_2^0, \ldots, L_2^{n_2}\} \cup \ldots \cup \{L_n^0, \ldots, L_n^{n_n}\}$ is the dual basis of $\{P_1^0, \ldots, P_1^{n_1}\} \cup \{P_2^0, \ldots, P_2^{n_2}\} \cup \ldots \cup \{P_n^0, \ldots, P_n^{n_n}\}$ of V. There is at most one such n-basis and if it exists is characterized by $L_{j_0}^0(P_{i_0}^0) \cup L_{j_1}^1(P_{i_1}^1) \cup L_{j_2}^2(P_{i_2}^2) \cup \ldots \cup L_{j_n}^n(P_{i_n}^n) = P_{i_0}^0(t_{j_0}^0) \cup P_{i_1}^1(t_{j_1}^1) \cup \ldots \cup P_{i_n}^n(t_{j_n}^n) = \delta_{i_0j_0} \cup \delta_{i_1j_1} \cup \ldots \cup \delta_{i_nj_n}$. The n-polynomials $P_i = P_{i_1}^1 \cup \ldots \cup P_{i_n}^n$

$$= \frac{(x - t_0^1)\ldots(x - t_{i_1-1}^1)(x - t_{i_1+1}^1)\ldots(x - t_{n_1}^1)}{(t_{i_1}^1 - t_0^1)\ldots(t_{i_1}^1 - t_{i_1-1}^1)(t_{i_1}^1 - t_{i_1+1}^1)\ldots(t_{i_1}^1 - t_{n_1}^1)}$$

$$\cup \ldots \cup \frac{(x - t_0^n)\ldots(x - t_{i_n-1}^n)(x - t_{i_n+1}^n)\ldots(x - t_{n_n}^n)}{(t_{i_n}^n - t_0^n)\ldots(t_{i_n}^n - t_{i_n-1}^n)(t_{i_n}^n - t_{i_n+1}^n)\ldots(t_{i_n}^n - t_{n_n}^n)}$$

$$= \prod_{j_1 \ne i_1}\left(\frac{x - t_{j_1}^1}{t_{i_1}^1 - t_{j_1}^1}\right) \cup \prod_{j_2 \ne i_2}\left(\frac{x - t_{j_2}^2}{t_{i_2}^2 - t_{j_2}^2}\right) \cup \ldots \cup \prod_{j_n \ne i_n}\left(\frac{x - t_{j_n}^n}{t_{i_n}^n - t_{j_n}^n}\right)$$

are of degree $(n_1, n_2, \ldots, n_n)$, hence belongs to $V = V_1 \cup \ldots \cup V_n$. If

$$f = f^1 \cup \ldots \cup f^n$$

$$= \sum_{i_1=1}^{n_1} c^1_{i_1} P^1_{i_1} \cup \sum_{i_2=1}^{n_2} c^2_{i_2} P^2_{i_2} \cup \ldots \cup \sum_{i_n=1}^{n_n} c^n_{i_n} P^n_{i_n}$$

then for each $j_k$ we have k=1, 2, …, n; $1 \le j_k \le n_k$.

$$f^k(t^k_{j_k}) = \sum_{i_k} c^k_{i_k} P^k_{i_k}(t^k_{j_k}) = c^k_{j_k}$$

true for each k, k = 1, 2, …, n.

Since the 0-polynomial has the properly $0(t_i) = 0$ for each $t = t_1 \cup \ldots \cup t_n \in F_1 \cup \ldots \cup F_n$ it follows from the above relation that the n-polynomials $\{P^1_0, P^1_1, \ldots, P^1_{n_1}\} \cup \{P^2_0, P^2_1, \ldots, P^2_{n_2}\} \cup \ldots \cup \{P^n_0, P^n_1, \ldots, P^n_{n_n}\}$ are n-linearly independent. The polynomials $\{1, x, \ldots, x^{n_1}\} \cup \{1, x, \ldots, x^{n_2}\} \cup \ldots \cup \{1, x, \ldots, x^{n_n}\}$ form a n-basis of V and hence the dimension of V is $\{n_1 + 1, n_2 + 1, \ldots, n_n + 1\}$. So the n-independent set $\{P^1_0, \ldots, P^1_{n_1}\} \cup \{P^2_0, \ldots, P^2_{n_2}\} \cup \ldots \cup \{P^n_0, \ldots, P^n_{n_n}\}$ must form an n-basis for V. Thus for each f in V

$$f = \sum_{i_1=0}^{n_1} f^1(t^1_{i_1}) P^1_{i_1} \cup \ldots \cup \sum_{i_n=0}^{n_n} f^n(t^n_{i_n}) P^n_{i_n} \qquad \text{(I)}$$

I is called the Lagranges' n-interpolation formula. Setting $f^k = x^{j_k}$ in I we obtain

$$x^{j_1} \cup \ldots \cup x^{j_k} = \sum_{i_1=0}^{n_1} (t^1_{i_1})^{\alpha_1} P^1_{i_1} \cup \ldots \cup \sum_{i_n=0}^{n_n} (t^n_{i_n})^{\alpha_n} P^n_{i_n}.$$

Thus the n-matrix

$$T = \begin{bmatrix} 1 & t^1_0 & (t^1_0)^2 \ldots (t^1_0)^{n_1} \\ 1 & t^1_1 & (t^1_1)^2 \ldots (t^1_1)^{n_1} \\ \vdots & \vdots & \vdots \quad \vdots \\ 1 & t^1_{n_1} & (t^1_{n_1})^2 \ldots (t^1_{n_1})^{n_1} \end{bmatrix} \cup \ldots \cup \begin{bmatrix} 1 & t^n_0 & (t^n_0)^2 \ldots (t^n_0)^{n_n} \\ 1 & t^n_1 & (t^n_1)^2 \ldots (t^n_1)^{n_n} \\ \vdots & \vdots & \vdots \quad \vdots \\ 1 & t^n_{n_n} & (t^n_{n_n})^2 \ldots (t^n_{n_n})^{n_n} \end{bmatrix} \text{(II)}$$

$= T_1 \cup \ldots \cup T_n$ is n-invertible. The n-matrix in II is called a Vandermonde n-matrix.

It can also be shown directly that the n-matrix is n-invertible when $(t_0^1, t_1^1, \ldots, t_{n_1}^1) \cup (t_0^2, t_1^2, \ldots, t_{n_2}^2) \cup \ldots \cup (t_0^n, t_1^n, \ldots, t_{n_n}^n)$ are $\{n_1 + 1, n_2 + 1, \ldots, n_n + 1\}$ set of n-distinct elements from the n-field $F = F_1 \cup \ldots \cup F_n$ i.e., each $(t_0^i, t_1^i, \ldots, t_{n_i}^i)$ is the $n_i + 1$ distinct elements of $F_i$; true for i = 1, 2, …, n.

Now we proceed on to define the new notion of n-polynomial function of an n-polynomial over the n-field $F = F_1 \cup F_2 \cup \ldots \cup F_n$. If $f = f_1 \cup \ldots \cup f_n$ be any n-polynomial over the n-field F we shall in the discussion denote by $f^{\sim} = f_1^{\sim} \cup f_2^{\sim} \cup \ldots \cup f_n^{\sim}$ the n-polynomial n-function from F into F taking each $t = t_1 \cup \ldots \cup t_n$ in F into f(t) where $t_i \in F_i$ for i = 1, 2, …, n. We see every polynomial function arises in this way, so analogously every n-polynomial function arises in the same way, however if $f^{\sim} = g^{\sim}$ and $f = f_1 \cup \ldots \cup f_n$ and $g = g_1 \cup g_2 \cup \ldots \cup g_n$ then $f_1^{\sim} = g_1^{\sim}, f_2^{\sim} = g_2^{\sim}, \ldots, f_n^{\sim} = g_n^{\sim}$ for any two equal n-polynomial f and g. So we assume two n-polynomials f and g such that $f \neq g$. However this situation occurs only when the n-field $F = F_1 \cup F_2 \cup \ldots \cup F_n$ is such that each field $F_i$ in F has only finite number of elements in it. Suppose f and g are n-polynomials over the n-field F then the product of $f^{\sim}$ and $g^{\sim}$ is the n-function $f^{\sim}g^{\sim}$ from F into F given by $f^{\sim}g^{\sim}(t) = f^{\sim}(t)\, g^{\sim}(t)$ for every $t = t_1 \cup \ldots \cup t_n \in F_1 \cup F_2 \cup \ldots \cup F_n$.

Further (fg)(x) = f(x)g(x) hence $(f^{\sim}g^{\sim})(x) = f^{\sim}(x)\, g^{\sim}(x)$ for each $x = x_1 \cup \ldots \cup x_n \in F_1 \cup F_2 \cup \ldots \cup F_n$.

Thus we see $\left(f^{\sim}g^{\sim}\right) = (fg)^{\sim}$ and it is also an n-polynomial function. We see that the n-polynomial function over the n-field is in fact an n-linear algebra over the n-field F. We shall denote this n-linear algebra over the n-field F by $A^{\sim} = A_1^{\sim} \cup \ldots \cup A_n^{\sim}$.

**DEFINITION 1.2.19:** *Let $F = F_1 \cup \ldots \cup F_n$ be a n-field. $A = A_1 \cup A_2 \cup \ldots \cup A_n$ be a n-linear algebra over the n-field F. Let $A^{\sim}$ be the n-linear algebra of n-polynomial functions over the same field F. The n-linear algebras A and $A^{\sim}$ are said to be n-isomorphic if there is a one to one n-mapping $\alpha \rightarrow \alpha^{\sim}$ such that*

*$(c\alpha + d\beta)^{\sim} = c\alpha^{\sim} + d\beta^{\sim}$ and $(\alpha\beta)^{\sim} = \alpha^{\sim}\beta^{\sim}$ for all $\alpha, \beta \in A$ and n-scalar c, d in the n-field F. Here $c = c_1 \cup \ldots \cup c_n$ and $d = d_1 \cup \ldots \cup d_n$, $\alpha = \alpha_1 \cup \ldots \cup \alpha_n$ and $\beta = \beta_1 \cup \ldots \cup \beta_n$ then $(c\alpha + d\beta)^{\sim} = (c_1\alpha_1 + d_1\beta_1)^{\sim} \cup (c_2\alpha_2 + d_2\beta_2)^{\sim} \cup \ldots \cup (c_n\alpha_n + d_n\beta_n)^{\sim} = (c_1\alpha^{\sim}_1 + d_1\beta^{\sim}_1) \cup (c_2\alpha^{\sim}_2 + d_2\beta^{\sim}_2) \cup \ldots \cup (c_n\alpha^{\sim}_n + d_n\beta^{\sim}_n)$ here $(c_i\alpha_i + d_i\beta_i)^{\sim} = c_i\alpha^{\sim}_i + d_i\beta^{\sim}_i$ where $\alpha_i, \beta_i \in V_i$, the component vector space of the n-vector space and $\alpha_i, d_i \in F_i$; the component field of the n-field F, true for i = 1, 2, …, n. Further $(\alpha\beta)^{\sim} = (\alpha_1\beta_1)^{\sim} \cup \ldots \cup (\alpha_n\beta_n)^{\sim} = \alpha^{\sim}_1\beta^{\sim}_1 \cup \ldots \cup \alpha^{\sim}_n\beta^{\sim}_n$. The n-mapping $\alpha \to \alpha^{\sim}$ is called an n-isomorphism of A onto $A^{\sim}$. An n-isomorphism of A onto $A^{\sim}$ is thus an n-vector space n-isomorphism of A onto $A^{\sim}$ which has the additional property of preserving products.*

We indicate the proof of the important theorem

**THEOREM 1.2.34:** *Let $F = F_1 \cup F_2 \cup \ldots \cup F_n$ be a n-field containing an infinite number of distinct elements, the n-mapping $f \to f^{\sim}$ is an n-isomorphism of the n-algebra of n-polynomials over the n-field F onto the n-algebra of n-polynomial functions over F.*

*Proof:* By definition the n-mapping is onto, and if $f = f_1 \cup \ldots \cup f_n$ and $g = g_1 \cup \ldots \cup g_n$ belong to $F[x] = F_1[x] \cup \ldots \cup F_n[x]$ the n-algebra of n-polynomial functions over the field F; i.e. $(cf+dg)^{\sim} = cf^{\sim} + dg^{\sim}$ where $c = c_1, \ldots, c_n$ and $d = d_1, \ldots, d_n \in F = F_1 \cup \ldots \cup F_n$; i.e.,

$$(cf + dg)^{\sim} = \left(c_1 f_1^{\sim} + d_1 g_1^{\sim}\right) \cup \left(c_2 f_2^{\sim} + d_2 g_2^{\sim}\right) \cup \ldots \cup \left(c_n f_n^{\sim} + d_n g_n^{\sim}\right)$$

for all n-scalars $c, d \in F$. Since we have already shown that $(fg)^{\sim} = f^{\sim}g^{\sim}$ we need only show the n-mapping is one to one. To do this it suffices by linearity of the n-algebras (i.e., each linear algebra is linear) $f^{\sim} = 0$ implies $f = 0$. Suppose then that $f = f_1 \cup f_2 \cup \ldots \cup f_n$ is a n-polynomial of degree $(n_1, n_2, \ldots, n_n)$ or less such that $f^1 = 0$ ; i.e., $f^1 = f_1^1 \cup f_2^1 \cup \ldots \cup f_n^1 = (0 \cup 0 \cup \ldots \cup 0)$. Let $\{t_0^1, t_1^1, \ldots, t_{n_1}^1\} \cup \{t_0^2, t_1^2, \ldots, t_{n_2}^2\} \cup \ldots \cup \{t_0^n, t_1^n, \ldots, t_{n_n}^n\}$ be any $\{n_1 + 1, n_2 + 1, \ldots, n_n + 1\}$; n-distinct elements of the n-field F. Since $f^{\sim} = 0$, $f^k(t_i) = 0$ for k = 1, 2, …, n; i = 0, 1, 2, …, $n_k$ and it is implies f = 0.

Now we proceed on to define the new notion of n-polynomial n-ideal or we can say as n-polynomial ideals. Throughout this section $F[x] = F_1[x] \cup \ldots \cup F_n[x]$ will denote a n-polynomial over the n-field $F = F_1 \cup \ldots \cup F_n$.

We will first prove a simple lemma.

**LEMMA 1.2.3:** *Suppose $f = f_1 \cup \ldots \cup f_n$ and $d = d_1 \cup \ldots \cup d_n$ are any two non zero n-polynomials over the n-field F such that n-deg d ≤ n-deg f, i.e., n-deg d = $(n_1, \ldots, n_n)$ and n-deg f = $(m_1, \ldots, m_n)$, n-deg d ≤ n-deg f if and only if each $n_i \leq m_i$, i = 1, 2, …, n. Then there exists an n-polynomial $g = g_1 \cup g_2 \cup \ldots \cup g_n$ in $F[x] = F_1[x] \cup \ldots \cup F_n[x]$ such that either f – dg = 0 or n-deg (f – dg) < n-deg f.*

*Proof:* Suppose $f = f_1 \cup \ldots \cup f_n$

$$= (a^1_{m_1} x^{m_1} + \sum_{i_1=0}^{m_1-1} a^1_{i_1} x^{i_1}) \cup (a^2_{m_2} x^{m_2} + \sum_{i_2=0}^{m_2-1} a^2_{i_2} x^{i_2})$$

$$\cup \ldots \cup (a^n_{m_n} x^{m_n} + \sum_{i_n=0}^{m_n-1} a^n_{i_n} x^{i_n});$$

$(a^1_{m_1}, a^2_{m_2}, \ldots, a^n_{m_n}) \neq (0, \ldots, 0)$ i.e., each $a^i_{m_i} \neq 0$ for i = 1, 2, …, n.

$$d = d_1 \cup \ldots \cup d_n =$$

$$\left(b^1_{n_1} x^{n_1} + \sum_{i_1=0}^{n_1-1} b^1_{i_1} x^{i_1}\right) \cup \left(b^2_{n_2} x^{n_2} + \sum_{i_2=0}^{n_2-1} b^2_{i_2} x^{i_2}\right) \cup \ldots \cup$$

$$\left(b^n_{n_n} x^{n_n} + \sum_{i_n=0}^{n_n-1} b^n_{i_n} x^{i_n}\right);$$

with $(b^1_{n_1}, \ldots, b^n_{n_n}) \neq (0, 0, \ldots, 0)$ i.e., $b^i_{n_i} \neq 0$, i = 1, 2, …, n. Then $(m_1, m_2, \ldots, m_n) > (n_1, n_2, \ldots, n_n)$ and

$$\left(f_1 - \left(\frac{a^1_{m_1}}{b^1_{n_1}}\right) x^{m_1 - n_1} d_1\right) \cup \left(f_2 - \left(\frac{a^2_{m_2}}{b^2_{n_2}}\right) x^{m_2 - n_2} d_2\right) \cup \ldots \cup$$

$$\left(f_n - \left(\frac{a^n_{m_n}}{b^n_{n_n}}\right)x^{m_n - n_n}d_n\right)$$

$$= 0 \cup \ldots \cup 0,$$

$$\text{or } \deg\left(f_1 - \left(\frac{a^1_{m_1}}{b^1_{n_1}}\right)x^{m_1 - n_1}d_1\right) \cup \deg\left(f_2 - \left(\frac{a^2_{m_2}}{b^2_{n_2}}\right)x^{m_2 - n_2}d_2\right) \cup \ldots \cup$$

$$\deg\left(f_n - \left(\frac{a^n_{m_n}}{b^n_{n_n}}\right)x^{m_n - n_n}d_n\right) < \deg f_1 \cup \deg f_2 \cup \ldots \cup \deg f_n.$$

Thus we make take

$$g = \left(\frac{a^1_{m_1}}{b^1_{n_1}}\right)x^{m_1 - n_1} \cup \left(\frac{a^2_{m_2}}{b^2_{n_2}}\right)x^{m_2 - n_2} \cup \ldots \cup \left(\frac{a^n_{m_n}}{b^n_{n_n}}\right)x^{m_n - n_n}.$$

Using this lemma we illustrate the usual process of long division of n-polynomials over a n-field $F = F_1 \cup F_2 \cup \ldots \cup F_n$.

**THEOREM 1.2.35:** *Let $f = f_1 \cup \ldots \cup f_n$, $d = d_1 \cup \ldots \cup d_n$ be n-polynomials over the n-field $F = F_1 \cup F_2 \cup \ldots \cup F_n$ and $d = d_1 \cup \ldots \cup d_n$ is different from $0 \cup \ldots \cup 0$; then there exists n-polynomials $q = q_1 \cup q_2 \cup \ldots \cup q_n$ and $r = r_1 \cup \ldots \cup r_n$ in F[x] such that*

1. *$f = dq + r$; i.e., $f = f_1 \cup f_2 \cup \ldots \cup f_n = (d_1q_1 + r_1) \cup \ldots \cup (d_nq_n + r_n)$.*
2. *Either $r = r_1 \cup \ldots \cup r_n = (0 \cup \ldots \cup 0)$ or n-deg r < n-deg d.*

*The n-polynomials q and r satisfying conditions (1) and (2) are unique.*

*Proof:* If $f = 0 \cup \ldots \cup 0$ or n-deg f < n-deg d we make take $q = q_1 \cup \ldots \cup q_n = (0 \cup \ldots \cup 0)$ and $r_1 \cup \ldots \cup r_n = f_1 \cup f_2 \cup \ldots \cup f_n$. In case $f = f_1 \cup \ldots \cup f_n \neq 0 \cup \ldots \cup 0$ and n-deg f $\geq$ n-deg d, then the preceding lemma shows we may choose a n-polynomial $g = g_1 \cup g_2 \cup \ldots \cup g_n$ such that $f - dg = 0 \cup \ldots \cup 0$, i.e. $(f_1 - d_1g_1) \cup (f_2 - d_2g_2) \cup \ldots \cup (f_n - d_ng_n) = 0 \cup \ldots \cup 0$, or n-deg$(f - dg) <$ n-deg f. If $f - dg \neq 0$ and n-deg $(f - dg) >$ n-deg d we choose a n-polynomial h such that $(f - dg) - dh = 0$ or $[(f_1 - d_1g_1) - d_1h_1] \cup [(f_2 - d_2g_2) - d_2h_2 ] \cup \ldots \cup [(f_n - d_ng_n) - d_nh_n] = 0 \cup \ldots \cup 0$ or n-deg$[f - d(g + h)] <$ nd$(f - dg)$ i.e., $\deg(f_1 - d_1(g_1 + h_1)) \cup \deg(f_2 - d_2(g_2 + h_2)) \cup \ldots \cup \deg(f_n - d_n(g_n+h_n)) < \deg(f_1 - d_1g_1) \cup \deg(f_2 - d_2g_2) \cup \ldots \cup \deg(f_n - d_ng_n)$.
Continuing this process as long as necessary we ultimately obtain n-polynomials $q = q_1 \cup q_2 \cup \ldots \cup q_n$ and $r = r_1 \cup \ldots \cup r_n$ such that $r = 0 \cup \ldots \cup 0$ or n-deg r < n-deg d i.e. $\deg r_1 \cup \ldots \cup \deg r_n < \deg d_1 \cup \deg d_2 \cup \ldots \cup \deg d_n$, and $f = dq + r$, i.e., $f_1 \cup f_2 \cup \ldots \cup f_n = (d_1q_1 + r_1) \cup (d_2q_2 + r_2) \cup \ldots \cup (d_nq_n + r_n)$.

Now suppose we also have $f = dq^1 + r^1$, i.e., $f_1 \cup \ldots \cup f_n = (d_1q_1^1 + r_1^1) \cup \ldots \cup (d_nq_n^1 + r_n^1)$ where $r^1 = 0 \cup \ldots \cup 0$ or n-deg $r^1 <$ n-deg d; i.e. $\deg r_1^1 \cup \ldots \cup \deg r_n^1 < \deg d_1 \cup \ldots \cup \deg d_n$. Then $dq + r = dq^1 + r^1$ and $d(q - q^1) = r^1 - r$ if $q - q^1 \neq 0 \cup \ldots \cup 0$ then $d(q - q^1) \neq 0 \cup \ldots \cup 0$. n-deg d + n-deg$(q - q^1)$ = n-deg $(r^1 - r)$; i.e. $(\deg d_1 \cup \deg d_2 \cup \ldots \cup \deg d_n) + \deg(q_1 - q^1_1) \cup \ldots \cup \deg(q_n - q^1_n) = \deg (r^1_1 - r_1) \cup \ldots \cup \deg (r^1_n - r_n)$.

But as n-degree of $r^1 - r$ is less than the n-degree of d this is impossible so $q - q^1 = 0 \cup \ldots \cup 0$ hence $r^1 - r = 0 \cup \ldots \cup 0$.

**DEFINITION 1.2.20:** *Let $d = d_1 \cup d_2 \cup \ldots \cup d_n$ be a non zero n-polynomial over the n-field $F = F_1 \cup \ldots \cup F_n$. If $f = f_1 \cup \ldots \cup f_n$ is in $F[x] = F_1[x] \cup \ldots \cup F_n[x]$, the proceeding theorem shows that there is at most one n-polynomial $q = q_1 \cup \ldots \cup q_n$ in F[x] such that $f = dq$ i.e., $f_1 \cup \ldots \cup f_n = d_1q_1 \cup \ldots \cup d_nq_n$. If such a q exists we say that $d = d_1 \cup \ldots \cup d_n$, n-divides $f = f_1 \cup \ldots \cup f_n$ that f is n-divisible by d and f is a n-multiple of d and we call q the n-quotient of f and d and write $q = f/d$ i.e. $q_1 \cup \ldots \cup q_n = f_1/d_1 \cup \ldots \cup f_n/d_n$.*

The following corollary is direct.

**COROLLARY 1.2.8:** *If $f = f_1 \cup \ldots \cup f_n$ is a n-polynomial over the n-field $F = F_1 \cup F_2 \cup \ldots \cup F_n$ and let $c = c_1 \cup c_2 \cup \ldots \cup c_n$* *be an element of F. Then f is n-divisible by $x - c = (x_1 - c_1) \cup \ldots \cup (x_n - c_n)$ if and only if $f(c) = 0 \cup \ldots \cup 0$ i.e., $f_1(c_1) \cup f_2(c_2) \cup \ldots \cup f_n(c_n) = 0 \cup 0 \cup \ldots \cup 0$.*

*Proof:* By theorem $f = (x-c)q + r$ where r is a n-scalar polynomial i.e., $r = r_1 \cup \ldots \cup r_n \in F$. By a theorem proved earlier we have $f(c) = f_1(c_1) \cup \ldots \cup f_n(c_n) = [0q_1(c_1) + r_1(c_1)] \cup [0q_2(c_2) + r_2(c_2)] \cup \ldots \cup [0q_n(c_n) + r_n(c_n)] = r_1(c_1) \cup \ldots \cup r_n(c_n)$. Hence $r = 0 \cup \ldots \cup 0$ if and only if $f(c) = f_1(c_1) \cup \ldots \cup f_n(c_n) = 0 \cup \ldots \cup 0$.

We know if $F = F_1 \cup \ldots \cup F_n$ is a n-field. An element $c = c_1 \cup c_2 \cup \ldots \cup c_n$ in F is said to be a n-root or a n-zero of a given n-polynomial $f = f_1 \cup \ldots \cup f_n$ over F if $f(c) = 0 \cup \ldots \cup 0$, i.e., $f_1(c_1) \cup \ldots \cup f_n(c_n) = 0 \cup \ldots \cup 0$.

**COROLLARY 1.2.9:** *A n-polynomial $f = f_1 \cup f_2 \cup \ldots \cup f_n$ of degree $(n_1, n_2, \ldots, n_n)$ over a n-field $F = F_1 \cup F_2 \cup \ldots \cup F_n$ has at most $(n_1, n_2, \ldots, n_n)$ roots in F.*

*Proof:* The result is true for n-polynomial of n-degree (1, 1, …, 1). We assume it to be true for n-polynomials of n-degree $(n_1 - 1, n_2 - 1, \ldots, n_n - 1)$. If $a = a_1 \cup a_2 \cup \ldots \cup a_n$ is a n-root of $f = f_1 \cup \ldots \cup f_n$, $f = (x - a)q = (x_1 - a_1)q_1 \cup \ldots \cup (x_n - a_n)q_n$ where $q = q_1 \cup \ldots \cup q_n$ has degree $n - 1$. Since $f(b) = 0$, i.e., $f_1(b_1) \cup \ldots \cup f_n(b_n) = 0 \cup \ldots \cup 0$ if and only if $a = b$ or $q(b) = q_1(b_1) \cup \ldots \cup q_n(b_n) = 0 \cup \ldots \cup 0$; it follows by our inductive assumption that f must have $(n_1, n_2, \ldots, n_n)$ roots.

We now define the n-derivative of a n-polynomial.

Let

$$f(x) = f_1(x) \cup \ldots \cup f_n(x)$$

$$= c_0^1 + c_1^1 x + \ldots + c_{n_1}^1 x^{n_1} \cup \ldots \cup c_0^n + c_1^n x + \ldots + c_{n_n}^n x^{n_n}$$

is a n-polynomial

$$f' = f_1' \cup f_2' \cup \ldots \cup f_n'$$

$$= \left(c_1^1 + 2c_2^1 x + \ldots + n_1 c_{n_1}^1 x^{n_1 - 1}\right) \cup \ldots \cup \left(c_1^n + 2c_2^n x + \ldots + cn_n c_{n_n}^n x^{n_n - 1}\right).$$

Now we will use the notation

$$Df = f' = Df_1 \cup \ldots \cup Df_n = f_1' \cup f_2' \cup \ldots \cup f_n'.$$

$$D^2 f = f'' = f_1'' \cup f_2'' \cup \ldots \cup f_n''$$

and so on. Thus $D_f^k = f_1^k \cup \ldots \cup f_n^k$.

Now we will prove the Taylor's formula for n-polynomials over the n-field F.

**THEOREM 1.2.36:** *Let $F = F_1 \cup \ldots \cup F_n$ be a n-field of n-characteristic (0, …, 0); $c = c_1 \cup c_2 \cup \ldots \cup c_n$ be an element in $F = F_1 \cup \ldots \cup F_n$, and $(n_1, n_2, \ldots, n_n)$ a n-tuple of positive integers. If $f = f_1 \cup \ldots \cup f_n$ is a n-polynomial over f with n-deg f $\leq (n_1, n_2, \ldots, n_n)$ then*

$$f = \sum_{k_1=0}^{n_1} \frac{D^{k_1} f_1}{\underline{|k_1}} c_1 (x - c_1)^{k_1} \cup$$

$$\sum_{k_2=0}^{n_2} \frac{D^{k_2} f_2}{\underline{|k_2}} c_2 (x - c_2)^{k_2} \cup \ldots \cup \sum_{k_n=0}^{n_n} \frac{D^{k_n} f_n}{\underline{|k_n}} c_n (x - c_n)^{k_n}.$$

*Proof:* We know Taylor's theorem is a consequence of the binomial theorem and the linearity of the operators $D^1, D^2, \ldots, D^n$. We know the binomial theorem $(a + b)^m = \sum_{k=0}^{m} \binom{m}{k} a^{m-k} b^k$, where $\binom{m}{k} = \frac{m!}{k!(m-k)!} = \frac{m(m-1)\ldots(m-k+1)}{1.2\ldots k}$ is the familiar binomial coefficient giving the number of combinations of m objects taken k at a time.

Now we apply binomial theorem to the n-tuple of polynomials

$$x^{m_1} \cup \ldots \cup x^{m_n} =$$

$$(c_1 + (x - c_1))^{m_1} \cup (c_2 + (x - c_2))^{m_2} \cup \ldots \cup (c_n + (x - c_n))^{m_n}$$

$$= \sum_{k_1=0}^{m_1} \binom{m_1}{k_1} c_1^{m_1 - k_1} (x - c_1)^{k_1} \cup \ldots \cup \sum_{k_n=0}^{m_n} \binom{m_n}{k_n} c_n^{m_n - k_n} (x - c_n)^{k_n}$$

$$= \{c_1^{m_1} + m_1 c_1^{m_1-1}(x-c_1) + \ldots + (x-c_1)^{m_1}\} \cup$$
$$\{c_2^{m_2} + m_2 c_2^{m_2-1}(x-c_2) + \ldots + (x-c_2)^{m_2}\} \cup \ldots \cup$$
$$\{c_n^{m_n} + m_n c_n^{m_n-1}(x-c_n) + \ldots + (x-c_n)^{m_n}\}$$

and this is the statement of Taylor's n-formula for the case

$$f = x^{m_1} \cup \ldots \cup x^{m_n}.$$

If

$$f = \sum_{m_1=0}^{n_1} a_{m_1}^1 x^{m_1} \cup \sum_{m_2=0}^{n_2} a_{m_2}^2 x^{m_2} \cup \ldots \cup \sum_{m_n=0}^{n_n} a_{m_n}^n x^{m_n}$$

$$D_f^k(c) = \sum_{m_1=0}^{n_1} a_{m_1}^1 D^{k_1} x^{m_1}(c_1) \cup \sum_{m_2=0}^{n_2} a_{m_2}^2 D^{k_2} x^{m_2}(c_2) \cup \ldots \cup$$
$$\sum_{m_n=0}^{n_n} a_{m_n}^n D^{k_n} x^{m_n}(c_n)$$

and

$$\sum_{k_1=0}^{m_1} \frac{D^{k_1} f_1(c_1)(x-c_1)^{k_1}}{\lfloor k_1} \cup$$
$$\sum_{k_2=0}^{m_2} \frac{D^{k_2} f_2(c_2)(x-c_2)^{k_2}}{\lfloor k_2} \cup \ldots \cup \sum_{k_n=0}^{m_n} \frac{D^{k_n} f_n(c_n)(x-c_n)^{k_n}}{\lfloor k_n}$$

$$= \sum_{k_1} \sum_{m_1} a_{m_1}^1 \frac{D^{k_1} x^{m_1}(c_1)(x-c_1)^{k_1}}{\lfloor k_1} \cup$$
$$\sum_{k_2} \sum_{m_2} a_{m_2}^2 \frac{D^{k_2} x^{m_2}(c_2)(x-c_2)^{k_2}}{\lfloor k_2} \cup \ldots \cup$$
$$\sum_{k_n} \sum_{m_n} a_{m_n}^n \frac{D^{k_n} x^{m_n}(c_n)(x-c_n)^{k_n}}{\lfloor k_n}$$

$$= \sum_{m_1} a_{m_1}^1 \sum_{k_1} \frac{D^{k_1} x^{m_1}(c_1)(x-c_1)^{k_1}}{\lfloor k_1} \cup$$
$$\sum_{m_2} a_{m_2}^2 \sum_{k_2} \frac{D^{k_2} x^{m_2}(c_2)(x-c_2)^{k_2}}{\lfloor k_2} \cup \ldots \cup$$

$$\sum_{m_n} a^n_{m_n} \sum_{k_n} \frac{D^{k_n} x^{m_n} (c_n)(x-c_n)^{k_n}}{\lfloor k_n}$$

$$= f_1 \cup \ldots \cup f_n.$$

If $c = c_1 \cup \ldots \cup c_n$ is a n-root of the n-polynomial $f = f_1 \cup \ldots \cup f_n$ the n-multiplicity of $c = c_1 \cup \ldots \cup c_n$ as a n-root of $f = f_1 \cup \ldots \cup f_n$ is the largest n-positive integer $(r_1, r_2, \ldots, r_n)$ such that $(x-c_1)^{r_1} \cup \ldots \cup (x-c_n)^{r_n}$ n-divides $f = f_1 \cup f_2 \cup \ldots \cup f_n$.

**THEOREM 1.2.37:** *Let $F = F_1 \cup \ldots \cup F_n$ be a n-field of (0, ..., 0) characteristic (i.e., each $F_i$ is of characteristic 0) for i = 1, 2, ..., n and $f = f_1 \cup \ldots \cup f_n$ be a n-polynomial over the n-field F with n-deg $f \leq (n_1, n_2, \ldots, n_n)$. Then the n-scalar $c = c_1 \cup \ldots \cup c_n$ is a n-root of f of multiplicity $(r_1, r_2, \ldots, r_n)$ if and only if $(D^{k_1} f_1)(c_1) \cup (D^{k_2} f_2)(c_2) \cup \ldots \cup (D^{k_n} f_n)(c_n) = 0 \cup 0 \cup \ldots \cup 0$; $0 \leq k_i \leq r_i - 1$; i = 1, 2, ...., n. $(D^{r_i} f_i)(c_i) \neq 0$, for every, i = 1, 2, ..., n.*

*Proof:* Suppose that $(r_1, r_2, \ldots, r_n)$ is the n-multiplicity of $c = c_1 \cup c_2 \cup \ldots \cup c_n$ as a n-root of $f = f_1 \cup \ldots \cup f_n$. Then there is a n-polynomial $g = g_1 \cup \ldots \cup g_n$ such that

$$f = (x-c_1)^{r_1} g_1 \cup \ldots \cup (x-c_n)^{r_n} g_n$$

and $g(c) = g_1(c_1) \cup \ldots \cup g_n(c_n) \neq 0 \cup \ldots \cup 0$. For otherwise $f = f_1 \cup \ldots \cup f_n$ would be divisible by $(x-c_1)^{r_1+1} \cup \ldots \cup (x-c_n)^{r_n+1}$. By Taylor's n-formula applied to $g = g_1 \cup \ldots \cup g_n$.

$$f = (x-c_1)^{r_1} \left[ \sum_{m_1=0}^{n_1-r_1} \frac{(D^{m_1} g_1)(c_1)(x-c_1)^{m_1}}{\lfloor m_1} \right] \cup$$

$$(x-c_2)^{r_2} \left[ \sum_{m_2=0}^{n_2-r_2} \frac{(D^{m_2} g_2)(c_2)(x-c_2)^{m_2}}{\lfloor m_2} \right] \cup \ldots \cup$$

$$(x-c_n)^{r_n} \left[ \sum_{m_n=0}^{n_n-r_n} \frac{(D^{m_n} g_n)(c_n)(x-c_n)^{m_n}}{\lfloor m_n} \right]$$

$$= \sum_{m_1=0}^{n_1-r_1} \frac{D^{m_1} g_1 (x-c_1)^{r_1+m_1}}{\lfloor m_1} \cup \ldots \cup \sum_{m_n=0}^{n_n-r_n} \frac{D^{m_n} g_n (x-c_n)^{r_n+m_n}}{\lfloor m_n}.$$

Since there is only one way to write $f = f_1 \cup \ldots \cup f_n$. (i.e., only one way to write each component $f_i$ of f) as a n-linear combination of the n-powers $(x-c_1)^{k_1} \cup \ldots \cup (x-c_n)^{k_n}$; $0 \leq k_i \leq n_i$; i = 1, 2, …, n; it follows that

$$\frac{(D^{k_i} f_i)(c_i)}{\lfloor k_i} = \begin{cases} 0 \text{ if} & 0 \leq k_i \leq r_i - 1 \\ \dfrac{D^{k_i - r_i} g_i(c_i)}{(k_i - r_i)!} & r_i \leq k_i \leq n. \end{cases}$$

This is true for every i, i = 1, 2, …, n. Therefore $D^{k_i} f_i(c_i) = 0$ for $0 \leq k_i \leq r_i - 1$; i = 1, 2, …, n and $D^{r_i} f_i(c_i) \neq g_i(c_i) \neq 0$ ; for every i = 1, 2, …, n. Conversely if these conditions are satisfied, it follows at once from Taylor's n-formula that there is a n-polynomial $g = g_1 \cup \ldots \cup g_n$ such that $f = f_1 \cup \ldots \cup f_n = (x-c_1)^{r_1} g_1 \cup \ldots \cup (x-c_n)^{r_n} g_n$ and $g(c) = g_1(c_1) \cup \ldots \cup g_n(c_n) \neq 0 \cup 0 \cup \ldots \cup 0$.

Now suppose that $(r_1, r_2, \ldots, r_n)$ is not the largest positive n-integer tuple such that $(x-c_1)^{r_1} \cup (x-c_2)^{r_2} \cup \ldots \cup (x-c_n)^{r_n}$ divides $f_1 \cup \ldots \cup f_n$; i.e., each $(x-c_i)^{r_i}$ divides $f_i$ for i = 1, 2, …, n; then there is a n-polynomial $h = h_1 \cup \ldots \cup h_n$ such that $f = (x-c_1)^{r_1+1} h_1 \cup \ldots \cup (x-c_n)^{r_n+1} h_n$. But this implies $g = g_1 \cup g_2 \cup \ldots \cup g_n = (x - c_1)h_1 \cup \ldots \cup (x - c_n)h_n$; hence $g(c) = g_1(c_1) \cup \ldots \cup g_n(c_n) = 0 \cup 0 \cup \ldots \cup 0$; a contradiction, hence the claim.

**DEFINITION 1.2.21:** *Let $F = F_1 \cup \ldots \cup F_n$ be a n-field. An n-ideal in $F[x] = F_1[x] \cup F_2[x] \cup \ldots \cup F_n[x]$ is a n-subspace; $m = m_1 \cup m_2 \cup \ldots \cup m_n$ of $F[x] = F_1[x] \cup \ldots \cup F_n[x]$ such that when $f = f_1 \cup \ldots \cup f_n$ and $g = g_1 \cup \ldots \cup g_n$ then $fg = f_1g_1 \cup f_2g_2$*

*$\cup \ldots \cup f_n g_n$ belongs to $m = m_1 \cup \ldots \cup m_n$; i.e. each $f_i g_i \in m_i$ whenever f is in F[x] and $g \in m$.*

*If in particular the n-ideal $m = dF[x]$ for some polynomial $d = d_1 \cup \ldots \cup d_n \in F[x]$ i.e. $m = m_1 \cup \ldots \cup m_n = d_1F[x] \cup \ldots \cup d_nF[x]$; i.e. the n-set of all n-multiples $d_1f_1 \cup \ldots \cup d_nf_n$ of $d = d_1 \cup \ldots \cup d_n$ by arbitrary $f = f_1 \cup \ldots \cup f_n$ in F[x] is a n-ideal. For m is non empty, m in fact contains d. If f, $g \in F[x]$ and c is a scalar then $c(df) - dg = (c_1d_1f_1 - d_1g_1) \cup \ldots \cup (c_nd_nf_n - d_ng_n) = d_1(c_1f_1 - g_1) \cup \ldots \cup d_n(c_nf_n - g_n)$ belongs to $m = m_1 \cup \ldots \cup m_n$; i.e. $d_i(c_if_i - g_i) \in m_i$ ; i = 1, 2, …., n so that m is a n-subspace. Finally m contains $(df)g = d(fg) = (d_1f_1)g_1 \cup \ldots \cup (d_nf_n)g_n = d_1(f_1g_1) \cup \ldots \cup d_n(f_ng_n)$ as well; m is called the principal n-ideal generated by $d = d_1 \cup \ldots \cup d_n$.*

Now we proceed on to prove an interesting theorem about the n-principal ideal of F[x].

**THEOREM 1.2.38:** *Let $F = F_1 \cup \ldots \cup F_n$ be a n-field and $m = m_1 \cup \ldots \cup m_n$, a non zero n-ideal in $F[x] = F_1[x] \cup \ldots \cup F_n[x]$. Then there is a unique monic n-polynomial $d = d_1 \cup \ldots \cup d_n$ in F[x] where each $d_i$ is a monic polynomial in $F_i[x]$ ; i = 1, 2, …, n such that m is the principal n-ideal generated by d.*

*Proof:* Given $F = F_1 \cup F_2 \cup \ldots \cup F_n$ is a n-field $F[x] = F_1[x] \cup \ldots \cup F_n[x]$ be the n-polynomial over the n-field F. Let $m = m_1 \cup \ldots \cup m_n$ be a non zero n-ideal of F[x]. We call a n-polynomial p(x) to be n-monic i.e. if in $p(x) = p_1(x) \cup \ldots \cup p_n(x)$ every $p_i(x)$ is a monic polynomial for i = 1, 2, …, n. Similarly we call a n-polynomial to be n-minimal if in $p(x) = p_1(x) \cup \ldots \cup p_n(x)$ each polynomial $p_i(x)$ is of minimal degree. Now $m = m_1 \cup \ldots \cup m_n$ contains a non zero n-polynomial $p(x) = p_1(x) \cup \ldots \cup p_n(x)$ where each $p_i(x) \neq 0$ for i = 1, 2, …, n. Among all the non zero n-polynomial in m there is a n-polynomial $d = d_1 \cup \ldots \cup d_n$ of minimal n-degree. Without loss in generality we may assume that minimal n-polynomial is monic i.e., d is monic. Suppose $f = f_1 \cup \ldots \cup f_n$ any n-polynomial in m then we know $f = dq + r$ where r = 0 or n-deg r < n-deg d i.e., $f = f_1 \cup \ldots \cup f_n = (d_1q_1 + r_1) \cup \ldots \cup (d_nq_n + r_n)$. Since d is in m, $dq = d_1q_1 \cup \ldots \cup d_nq_n \in$ m and $f \in m$ so $f - dq = r = r_1 \cup r_2 \cup \ldots \cup r_n \in m$. But since d is an n-polynomial in m of minimal n-degree we cannot have n-deg r < n-deg d so $r = 0 \cup \ldots \cup 0$. Thus $m = dF[x] = d_1F_1[x] \cup \ldots \cup d_nF_n[x]$. If g is any other n-monic polynomial such that $gF[x] = m = g_1F_1[x] \cup \ldots \cup g_nF_n[x]$ then their exists non zero n-polynomial $p = p_1 \cup p_2 \cup \ldots \cup p_n$ and $q = q_1 \cup q_2 \cup \ldots \cup q_n$ such that d = gp and g = dq i.e., $d = d_1 \cup \ldots \cup d_n = g_1p_1 \cup \ldots \cup g_np_n$ and $g_1 \cup \ldots \cup g_n = d_1q_1 \cup \ldots \cup d_nq_n$. Thus $d = dpq = d_1p_1q_1 \cup \ldots \cup d_np_nq_n = d_1 \cup \ldots \cup d_n$ and n-deg d = n-deg d + n-deg p + n-deg q.

Hence n-deg p = n-deg q = (0, 0, …, 0) and as d and g are n-monic p = q = 1. Thus d = g. Hence the claim.

In the n-ideal m we have f = pq + r where p, $f \in m$ i.e. $p = p_1 \cup p_2 \cup \ldots \cup p_n \in m$ and $f = f_1 \cup \ldots \cup f_n \in m$; $f = f_1 \cup \ldots \cup f_n = (p_1q_1 + r_1) \cup \ldots \cup (p_nq_n + r_n)$ where the n-remainder $r = r_1 \cup \ldots \cup r_n \in m$ is different from $0 \cup \ldots \cup 0$ and has smaller n-degree than p.

**COROLLARY 1.2.10:** *If $p^1, p^2, \ldots, p^n$ are n-polynomials over a n-field $F = F_1 \cup \ldots \cup F_n$ not all of which are zero i.e. $0 \cup \ldots \cup 0$, there is a unique n-monic polynomial d in $F[x] = F_1[x] \cup \ldots \cup F_n[x]$ and $d = d_1 \cup \ldots \cup d_n$ such that*

a. *$d = d_1 \cup \ldots \cup d_n$ is in the n-ideal generated by $p^1, \ldots, p^n$, where $p^i = p^i_1 \cup p^i_2 \cup \ldots \cup p^i_n$; i = 1, 2, …, n.*

b. *$d = d_1 \cup \ldots \cup d_n$, n-divides each of the n-polynomials $p^i = p^i_1 \cup \ldots \cup p^i_n$ i.e. $d_j / p^i_j$ for j = 1, 2, …, n true for i = 1, 2, …, n.*

*Any n-polynomial satisfying (a) and (b) necessarily satisfies*

c. *d is n-divisible by every n-polynomial which divides each of the n-polynomials $p^1, p^2, \ldots, p^n$.*

*Proof:* Let $F = F_1 \cup \ldots \cup F_n$ be a n-field, $F[x] = F_1[x] \cup \ldots \cup F_n[x]$ be a n-polynomial ring over the n-field F. Let d be a n-monic generator of the n-ideal

$$\{ p_1^1 F_1[x] \cup p_1^2 F_1[x] \cup \ldots \cup p_1^n F_1[x]\} \cup \{ p_2^1 F_2[x] \cup \ldots \cup p_2^n F_2[x]\} \cup \ldots \cup \{ p_n^1 F_n[x] \cup \ldots \cup p_n^n F_n[x]\}.$$

Every member of this n-ideal is divisible by $d = d_1 \cup \ldots \cup d_n$. Thus each of the n-polynomials $p^i = p_1^i \cup \ldots \cup p_n^i$ is n-divisible by d. Now suppose $f = f_1 \cup \ldots \cup f_n$ is a n-polynomial which n-divides each of the n-polynomial

$$p^1 = p_1^1 \cup \ldots \cup p_n^1$$
$$p^2 = p_1^2 \cup \ldots \cup p_n^2$$

and

$$p^n = p_1^n \cup \ldots \cup p_n^n.$$

Then there exists n-polynomials $g^1, \ldots, g^n$ such that $p^i = fg^i$ i.e., $p_1^i \cup \ldots \cup p_n^i = f_1 g_1^i \cup \ldots \cup f_n g_n^i$ i.e. $p_t^i = f_t g_t^i$ for $t = 1, 2, \ldots, n$; $1 \leq i \leq n$. Also since $d = d_1 \cup \ldots \cup d_n$ is in the n-ideal $\{ p_1^1 F_1[x] \cup p_1^2 F_1[x] \cup \ldots \cup p_1^n F_1[x]\} \cup \{ p_2^1 F_2[x] \cup \ldots \cup p_2^n F_2[x]\} \cup \ldots \cup \{ p_n^1 F_n[x] \cup \ldots \cup p_n^n F_n[x]\}$ there exists n-polynomials $q^1, \ldots, q^n$ in $F[x] = F_1[x] \cup \ldots \cup F_n[x]$ with $q^i = q_1^i \cup \ldots \cup q_n^i$; $i = 1, 2, \ldots, n$, such that $d = \{ p_1^1 q_1^1 \cup p_2^1 q_2^1 \cup \ldots \cup p_n^1 q_n^1 \} \cup \ldots \cup \{ p_1^n q_1^n \cup p_2^n q_2^n \cup \ldots \cup p_n^n q_n^n \} = d_1 \cup \ldots \cup d_n$. Thus $d = f_1( p_1^1 q_1^1 \cup \ldots \cup p_n^1 q_n^1 ) \cup f_2( p_1^2 q_1^2 \cup \ldots \cup p_n^2 q_n^2 ) \cup \ldots \cup f_n( p_1^n q_1^n \cup \ldots \cup p_n^n q_n^n )$.

We have shown that $d = d_1 \cup \ldots \cup d_n$ is a n-monic polynomial satisfying (a), (b) and (c). If $d^1 = d_1^1 \cup \ldots \cup d_n^1$ is any n-polynomial satisfying (a) and (b) it follows from (a) and the definition of d that $d^1$ is a scalar multiple of d and satisfies (c) as well. If $d^1$ is also n-monic then $d = d^1$.

**DEFINITION 1.2.22:** *If $p^1, p^2, \ldots, p^n$ where $p^i = p_1^i \cup p_2^i \cup \ldots \cup p_n^i$ are n-polynomials over the n-field $F = F_1 \cup F_2 \cup \ldots \cup F_n$ for $i = 1, \ldots, n$, such that not all the n-polynomials are $0 \cup \ldots \cup 0$ the monic generator $d = d_1 \cup \ldots \cup d_n$ of the n-ideal $\{ p_1^1 F_1[x] \cup p_1^2 F_1[x] \cup \ldots \cup p_1^n F_1[x]\} \cup \{ p_2^1 F_2[x] \cup p_2^2 F_2[x] \cup \ldots \cup p_2^n F_2[x]\}$*

*$\cup \ldots \cup \{ p_n^1 F_n[x] \cup \ldots \cup p_n^n F_n[x]\}$ is called the greatest common n-divisor or n-greatest common divisor of $p^1, \ldots, p^n$. This terminology is justified by the proceeding corollary. We say the n-polynomials $p^1 = p_1^1 \cup \ldots \cup p_n^1$, $p^2 = p_1^2 \cup \ldots \cup p_n^2, \ldots, p^n = p_1^n \cup \ldots \cup p_n^n$ are n-relatively prime if their n-greatest common divisor is (1, 1, ... , 1) or equivalently if the n-ideal they generate is all of $F[x] = F_1[x] \cup \ldots \cup F_n[x]$.*

Now we talk about the n-factorization, n-irreducible, n-prime polynomial over the n-field F.

**DEFINITION 1.2.23:** *Let F be an n-field i.e. $F = F_1 \cup F_2 \cup \ldots \cup F_n$. A n-polynomial $f = f_1 \cup \ldots \cup f_n$ in $F[x] = F_1[x] \cup \ldots \cup F_n[x]$ is said to be n-reducible over the n-field $F = F_1 \cup F_2 \cup \ldots \cup F_n$, if there exists n-polynomials g, h $\in F[x]$, $g = g_1 \cup \ldots \cup g_n$ and $h = h_1 \cup \ldots \cup h_n$ in F[x] of n-degree $\geq$ (1, 1, ..., 1) such that f = gh, i.e., $f_1 \cup \ldots \cup f_n = g_1h_1 \cup \ldots \cup g_nh_n$ and if such g and h does not exists, $f = f_1 \cup \ldots \cup f_n$ is said to be n-irreducible over the n-field $F = F_1 \cup \ldots \cup F_n$.*

*A non n-scalar, n-irreducible n-polynomial over $F = F_1 \cup \ldots \cup F_n$ is called the n-prime polynomial over the n-field F and we some times say it is n-prime in $F[x] = F_1[x] \cup \ldots \cup F_n[x]$.*

**THEOREM 1.2.39:** *Let $p = p^1 \cup p^2 \cup \ldots \cup p^n$, $f = f^1 \cup f^2 \cup \ldots \cup f^n$ and $g = g^1 \cup g^2 \cup \ldots \cup g^n$ be n-polynomial over the n-field $F = F_1 \cup \ldots \cup F_n$. Suppose that p is a n-prime n-polynomial and that p n-divides the product fg, then either p n-divides f or p n-divides g.*

*Proof:* Without loss of generality let us assume $p = p^1 \cup p^2 \cup \ldots \cup p^n$ is a n-monic n-prime n-polynomial i.e. monic prime n-polynomial. The fact that $p = p^1 \cup p^2 \cup \ldots \cup p^n$ is prime then simply says that only monic n-divisor of p are $1_n$ and p. Let d be the n-gcd or greatest common n-divisor of f and p. $f = f^1 \cup \ldots \cup f^n$ and $p = p^1 \cup \ldots \cup p^n$ since d is a monic n-polynomial which n-divides p. If d = p then p n-divides f and we are done. So suppose $d = 1_n = (1 \cup 1 \cup \ldots \cup 1)$ i.e., suppose f and p are n-relatively prime if f = $f^1 \cup \ldots \cup f^n$ and p = $p^1 \cup \ldots \cup p^n$, then $f^1$ and $p^1$ are relatively prime $f^2$and $p^2$ are relatively prime i.e. $(f^i, p^i) = 1$ for i = 1, 2, …, n. We shall prove that p n-divides g. Since $(f, p) = (f^1, p^1) \cup \ldots \cup (f^n, p^n) = 1 \cup \ldots \cup 1$ there are n-polynomials $f_0 = f_0^1 \cup \ldots \cup f_0^n$ ; $p_0 = p_0^1 \cup \ldots \cup p_0^n$ such that $1 \cup \ldots \cup 1 = f_0f + p_0p$ i.e. $1 \cup \ldots \cup 1 = (f_0^1f^1 + p_0^1p^1) \cup (f_0^2f^2 + p_0^2p^2) \cup \ldots \cup (f_0^nf^n + p_0^np^n)$ n-multiplying by $g = g^1 \cup \ldots \cup g^n$ we get $g = f_0fg + p_0pg$ i.e. $(g^1 \cup \ldots \cup g^n) = (f_0^1f^1g^1 + p_0^1p^1g^1) \cup \ldots \cup (f_0^nf^ng^n + p_0^np^ng^n)$. Since p n-divides fg it n-divides $(fg)f_0$ and certainly p n-divides $pp_0g$. Thus p n-divides g. Hence the claim.

**COROLLARY 1.2.11:** *If $p = p^1 \cup \ldots \cup p^n$ is a n-prime that n-divides a n-product $f_1, \ldots, f_n$ i.e. $f_1^1 \ldots f_n^1 \cup f_1^2 \ldots f_n^2 \cup \ldots \cup f_1^n \ldots f_n^n$. Then p n-divides one of the n-polynomials $f_1^1 \ldots f_n^1 \cup \ldots \cup f_1^n \ldots f_n^n$.*

The proof is left for the reader.

**THEOREM 1.2.40:** *If $F = F_1 \cup F_2 \cup \ldots \cup F_n$ is a n-field, a non n-scalar monic n-polynomial in $F[x] = F_1[x] \cup \ldots \cup F_n[x]$ can be n-factored as a n-product of n-monic primes in F[x] in one and only one way except for the order.*

*Proof:* Given $F = F_1 \cup \ldots \cup F_n$ is a n-field. $F[x] = F_1[x] \cup \ldots \cup F_n[x]$ is a n-polynomial over F[x]. Suppose $f = f^1 \cup \ldots \cup f^n$ is a non scalar monic n-polynomial over the n-field F. As n-polynomial of n-degree (1, 1, …, 1) are irreducible there is nothing to prove if n-deg f = (1, 1, …, 1). Suppose f has n degree $(n_1, n_2, \ldots, n_n) > (1, 1, \ldots, 1)$, by induction we may assume the theorem is true for all non scalar monic n-polynomials of n-degree less than $(n_1, n_2, \ldots, n_n)$. If f is n-irreducible it is already n-factored as a n-product of monic n-primes and otherwise $f = gh = f^1 \cup f^2 \cup \ldots \cup f^n = g^1h^1 \cup \ldots \cup g^nh^n$ where g and h are non scalar monic n-polynomials (i.e. $g = g^1 \cup \ldots \cup g^n$ and $h = h^1 \cup \ldots \cup h^n$) of n-degree less than $(n_1, n_2,$

…, $n_n$). Thus f and g can be n-factored as n-products of monic n-primes in $F[x] = F_1[x] \cup \ldots \cup F_n[x]$ and hence $f = f^1 \cup \ldots \cup f^n = (p_1^1, \ldots, p_{m_1}^1) \cup \ldots \cup (p_1^n, \ldots, p_{m_n}^n) = (q_1^1, \ldots, q_{n_1}^1) \cup \ldots \cup (q_1^n, q_2^n, \ldots, q_{n_n}^n)$ where $(p_1^1, \ldots, p_{m_1}^1) \cup \ldots \cup (p_1^n, p_2^n, \ldots, p_{m_n}^n)$ and $(q_1^1, \ldots, q_{n_1}^1) \cup \ldots \cup (q_1^n, q_2^n, \ldots, q_{n_n}^n)$ are monic n-primes in $F[x] = F_1[x] \cup \ldots \cup F_n[x]$. Then $(p_{m_1}^1 \cup \ldots \cup p_{m_n}^n)$ must n-divide some $(q_{i_1}^i \cup \ldots \cup q_{i_n}^n)$. Since both $(p_{m_1}^1 \cup \ldots \cup p_{m_n}^n)$ and $(q_{i_1}^i \cup \ldots \cup q_{i_n}^n)$ are monic n-primes this means that $q_{i_t}^t = p_{m_t}^t$ for every t = 1, 2, …, n. Thus we see $m_i = n_i = 1$ for each i = 1, 2, …, n, if either $m_i = 1$ or $n_i = 1$ for i = 1, 2, …, n. For

$$\text{n-deg } f = \left( \sum_{i_1=1}^{m_1} \deg p_{i_1}^1 = \sum_{j_1=1}^{n_1} \deg q_{j_1}^1, \ldots, \sum_{i_n=1}^{m_n} \deg p_{i_n}^n = \sum_{j_n=1}^{n_n} \deg q_{j_n}^n \right).$$

In this case we have nothing more to prove. So we may assume $m_i > 1$, i = 1, 2, …, n and $n_j > 1$, j = 1, 2, …, n.

By rearranging $q^i$'s we can assume $p_{m_i}^i = q_{n_i}^i$ and that $\left(p_1^1, \ldots, p_{m_1-1}^1 p_{m_1}^1, \ldots, p_{m_n-1}^n\right) = (q_1^1 q_2^1, \ldots, q_{n_1-1}^1 p_{m_1}^1, \ldots, q_1^n, \ldots, q_{m_n}^n p_{m_n}^n)$. Thus $\left(p_1^1, \ldots, p_{m_1-1}^1, p_{m_1}^1, \ldots, p_{m_n-1}^n\right) = (q_1^1, \ldots, q_{n_1-1}^1, \ldots, q_1^n, \ldots, q_{n_n-1}^n)$. As the n-polynomial has n-degree less than ($n_1$, $n_2$, …, $n_n$) our inductive assumption applies and shows the n-sequence $(q_1^1, \ldots, q_{n_1-1}^1, \ldots, q_1^n, \ldots, q_{n_n-1}^n)$ is at most rearrangement of n-sequence $(p_1^1, \ldots, p_{m_1-1}^1, \ldots, p_1^n, \ldots, p_{m_n-1}^1)$.

This shows that the n-factorization of $f = f_1 \cup \ldots \cup f_n$ as a product of monic n-primes is unique up to order of factors. Several interesting results in this direction can be derived. The reader is expected to define n-primary decomposition of the n-polynomial $f = f^1 \cup \ldots \cup f^n$ in $F[x] = F_1[x] \cup \ldots \cup F_n[x]$. The reader is requested to prove the following theorems

**THEOREM 1.2.41:** *Let $f = f^1 \cup \ldots \cup f^n$ be a non scalar monic n-polynomial over the n-field $F = F_1 \cup \ldots \cup F_n$ and let*

$$f = \left( p_1^{n_1^1} \ldots p_{k_1}^{n_{k_1}^1} \cup \ldots \cup p_1^{n_1^n} \ldots p_{k_n}^{n_{k_n}^n} \right)$$

*be the prime n-factorization of f. For each $j_t$, $1 \le j_t \le k_t$, let*

$$f_j^t = f^t \Big/ p_j^{n_j^t} = \prod_{i \ne j} p_i^{n_{k_i}^t} \,,$$

*$t = 1, 2, \ldots, n$. Then $f_1^t, \ldots, f_{k_t}^t$ are relatively prime for $t = 1, 2, \ldots, n$.*

**THEOREM 1.2.42:** *If $f = f_1 \cup \ldots \cup f_n$ is a n-polynomial over the n-field $F = F_1 \cup F_2 \cup \ldots \cup F_n$ with derivative $f' = f_1' \cup f_2' \cup \ldots \cup f_n'$ Then f is a n-product of distinct irreducible n-polynomials over the n-field F if and only if f and f ' are relatively prime i.e., each $f_i$ and $f_i'$ are relatively prime for $i = 1, 2, \ldots, n$.*

The proof is left as an exercise for the reader.

Next we proceed on to define the notion of n-characteristic value of type II i.e., n-characteristic values for n-linear operator on the n-vector space V.

Throughout this section we assume the n-vector space is defined, over n-field; i.e., n-vector space of type II.

**DEFINITION 1.2.24:** *Let $V = V_1 \cup \ldots \cup V_n$ be a n-vector space over the n-field $F = F_1 \cup F_2 \cup \ldots \cup F_n$ and let T be a n-linear operator on V, i.e., $T = T_1 \cup \ldots \cup T_n$ and $T_i: V_i \to V_i$, $i = 1, 2, \ldots, n$. This is the only way n-linear operator can be defined on V. A n-characteristic value of T is a n-scalar $C = C_1 \cup \ldots \cup C_n$ ($C_i \in F_i$, $i = 1, 2, \ldots, n$) in F such that there is a non zero n-vector $\alpha = \alpha_1 \cup \ldots \cup \alpha_n$ in $V = V_1 \cup \ldots \cup V_n$ with $T\alpha = C\alpha$ i.e., $T = T_1\alpha_1 \cup \ldots \cup T_n\alpha_n = C_1\alpha_1 \cup \ldots \cup C_n\alpha_n$ i.e., $T_i\alpha_i = C_i\,\alpha_i$; $i = 1, 2, \ldots, n$. If C is a n-characteristic value of T then*

- *a. any $\alpha = \alpha_1 \cup \ldots \cup \alpha_n$ such that $T\alpha = C\alpha$ is called the n-characteristic n-vector of T associated with the n-characteristic value $C = C_1 \cup \ldots \cup C_n$.*
- *b. The collection of all $\alpha = \alpha_1 \cup \ldots \cup \alpha_n$ such that $T\alpha = C\alpha$ is called the n-characteristic space associated with C.*

*If $T = T_1 \cup \ldots \cup T_n$ is any n-linear operator on the n-vector space $V = V_1 \cup \ldots \cup V_n$. We call the n-characteristic values associated with T to be n-characteristic roots, n-latent roots, n-eigen values, n-proper values or n-spectral values.*

*If T is any n-linear operator and $C = C_1 \cup \ldots \cup C_n$ in any n-scalar, the set of n-vectors $\alpha = \alpha_1 \cup \ldots \cup \alpha_n$ such that $T\alpha = C\alpha$ is a n-subspace of V. It is in fact the n-null space of the n-linear transformation $(T - CI) = (T_1 - C_1I_1) \cup \ldots \cup (T_n - C_nI_n)$ where $I_j$ denotes a unit matrix for $j = 1, 2, \ldots, n$. We call $C_1 \cup \ldots \cup C_n$ the n-characteristic value of T if this n-subspace is different from the n-zero space $0 = 0 \cup 0 \cup \ldots \cup 0$, i.e.,$(T - CI)$ fails to be a one to one n-linear transformation and if the n-vector space is finite dimensional we see that $(T - CI)$ fails to be one to one and the n-det$(T - CI) = 0 \cup 0 \cup \ldots \cup 0$.*

We have the following theorem in view of these properties.

**THEOREM 1.2.43:** *Let $T = T_1 \cup \ldots \cup T_n$ be a n-linear operator on a finite dimensional n-vector space $V = V_1 \cup \ldots \cup V_n$ and let $C = C_1 \cup \ldots \cup C_n$ be a scalar. The following are equivalent.*

a. *$C = C_1 \cup \ldots \cup C_n$ is a n-characteristic value of $T = T_1 \cup \ldots \cup T_n$.*
b. *The n-operator $(T_1 - C_1I_1) \cup \ldots \cup (T_n - C_nI_n) = T - CI$ is n-singular or (not n-invertible).*
c. *$\det(T - CI) = 0 \cup 0 \cup \ldots \cup 0$ i.e., $\det(T_1 - C_1I) \cup \ldots \cup \det(T_n - C_nI_n) = 0 \cup \ldots \cup 0$.*

Now we define the n-characteristic value of a n-matrix $A = A_1 \cup \ldots \cup A_n$ where each $A_i$ is a $n_i \times n_i$ matrix with entries from the field $F_i$ so that A is a n-matrix defined over the n-field $F = F_1 \cup F_2 \cup \ldots \cup F_n$. A n-characteristic value of A in the n-field $F = F_1 \cup F_2 \cup \ldots \cup F_n$ is a n-scalar $C = C_1 \cup \ldots \cup C_n$ in F such that the n-matrix $A - CI = (A_1 - C_1I_1) \cup (A_2 - C_2I_2) \cup \ldots \cup (A_n - C_nI_n)$ is n-singular or not n-invertible.

Since $C = C_1 \cup \ldots \cup C_n$ is a n-characteristic value of $A = A_1 \cup \ldots \cup A_n$, A a $(n_1 \times n_1, \ldots, n_n \times n_n)$ n-matrix over the n-field, $F = F_1 \cup F_2 \cup \ldots \cup F_n$, if and only if n-det$(A - CI) = 0 \cup 0 \cup \ldots \cup 0$ i.e., $\det(A_1 - C_1I_1) \cup \ldots \cup \det(A_n - C_nI_n) = (0 \cup 0 \cup \ldots \cup 0)$, we form the n-matrix $(xI - A) = (xI_1 - A_1) \cup \ldots \cup (xI_n - A_n)$.

Clearly the n-characteristic values of A in $F = F_1 \cup F_2 \cup \ldots \cup F_n$ are just n-scalars $C = C_1 \cup \ldots \cup C_n$ in $F = F_1 \cup F_2 \cup \ldots \cup F_n$ such that the n-scalars $C = C_1 \cup \ldots \cup C_n$ in $F = F_1 \cup F_2 \cup \ldots \cup F_n$ such that $f(C) = f_1(C_1) \cup \ldots \cup f_n(C_n) = 0 \cup 0 \cup \ldots \cup 0$. For this reason $f = f_1 \cup f_2 \cup \ldots \cup f_n$ is called the n-characteristic polynomial of A. Clearly f is a n-polynomial of different degree in x over different fields. It is important to note that $f = f_1 \cup f_2 \cup \ldots \cup f_n$ is a n-monic n-polynomial which has n-degree exactly $(n_1, n_2, \ldots, n_n)$. The n-monic polynomial is also a n-polynomial over $F = F_1 \cup F_2 \cup \ldots \cup F_n$.

First we illustrate this situation by the following example.

***Example 1.2.9:*** Let

$$A = \begin{bmatrix} 1 & 0 & 1 \\ 0 & 1 & 0 \\ 1 & 0 & 0 \end{bmatrix} \cup \begin{bmatrix} 2 & 1 & 0 & 1 \\ 1 & 1 & 0 & 0 \\ 0 & 2 & 2 & 1 \\ 0 & 0 & 0 & 1 \end{bmatrix} \cup \begin{bmatrix} 0 & 4 \\ 1 & 0 \end{bmatrix} \cup \begin{bmatrix} 1 & 0 & 0 & 0 & 0 \\ 0 & 1 & 0 & 0 & 0 \\ 0 & 4 & 6 & 0 & 1 \\ 0 & 0 & 0 & 5 & 0 \\ 0 & 0 & 0 & 0 & 3 \end{bmatrix}$$

$= A_1 \cup A_2 \cup A_3 \cup A_4$; be a 4-matrix of order $(3 \times 3, 4 \times 4, 2 \times 2, 5 \times 5)$ over the 4-field $F = Z_2 \cup Z_3 \cup Z_5 \cup Z_7$. The 4-characteristic 4-polynomial associated with A is given by

$$(xI - A) = (xI_{3\times3} - A_1) \cup (xI_{4\times4} - A_2) \cup (xI_{2\times2} - A_3) \cup (xI_{5\times5} - A_4)$$

$$= \begin{bmatrix} x+1 & 0 & 1 \\ 0 & x+1 & 0 \\ 1 & 0 & x \end{bmatrix} \cup \begin{bmatrix} x+1 & 2 & 0 & 2 \\ 2 & x+2 & 0 & 0 \\ 0 & 1 & x+1 & 2 \\ 0 & 0 & 0 & x+2 \end{bmatrix} \cup$$

$$\begin{bmatrix} x & 1 \\ 4 & x \end{bmatrix} \cup \begin{bmatrix} x+6 & 0 & 0 & 0 & 0 \\ 0 & x+6 & 0 & 0 & 0 \\ 0 & 3 & x+1 & 0 & 6 \\ 0 & 0 & 0 & x+2 & 0 \\ 0 & 0 & 0 & 0 & x+4 \end{bmatrix}$$

is 4-matrix with polynomial entries.

$$\begin{aligned} f &= f_1 \cup f_2 \cup \ldots \cup f_4 \\ &= \det(xI - A) \\ &= \det(xI_1 - A_1) \cup \det(xI_2 - A_2) \cup \ldots \cup \det(xI - A) \\ &= \{(x+1)^2 x + (x+1)\} \cup \{(x+1)(x+2) \times (x+1)(x+2) \\ &\quad - 2 \times 2(x+1)(x+2)\} \cup (x^2 - 4) \cup \{(x+6)(x+6)(x+ \\ &\quad 1)(x+2)(x+4)\}. \end{aligned}$$

We see det(xI – A) is a 4-polynomial which is a monic 4-polynomial and degree of f is (3, 4, 2, 5) over the 4-field $F = Z_2 \cup Z_3 \cup Z_5 \cup Z_7$.

Now we first define the notion of similar n-matrices when the entries of the n-matrices are from the n-field.

**DEFINITION 1.2.25:** *Let $A = A_1 \cup \ldots \cup A_n$ be a $(n_1 \times n_1, \ldots, n_n \times n_n)$ matrix over the n-field $F = F_1 \cup F_2 \cup \ldots \cup F_n$ i.e., each $A_i$ takes its entries from the field $F_i$, i = 1, 2, ..., n. We say two n-matrices A and B of same order are similar if there exits a n-non invertible n-matrix $P = P_1 \cup \ldots \cup P_n$ of $(n_1 \times n_1, \ldots, n_n \times n_n)$ order such that*

$$B = P^{-1}AP \text{ where } P^{-1} = P_1^{-1} \cup P_2^{-1} \cup \ldots \cup P_n^{-1}.$$

$$B = B_1 \cup B_2 \cup \ldots \cup B_n = P_1^{-1} A_1 P_1 \cup P_2^{-1} A_2 P_2 \cup \ldots \cup P_n^{-1} A_n P_n$$

*then*

$$\begin{aligned} \det(xI - B) &= \det((xI - P^{-1}AP)) \\ &= \det P^{-1}(xI - A) P \\ &= \det P^{-1} \det(xI - A) \det P \\ &= \det(xI - A) \\ &= \det(xI_1 - A_1) \cup \ldots \cup \det(xI_n - A_n). \end{aligned}$$

***Example 1.2.10:*** Let A = $A_1 \cup A_2 \cup A_3$

$$= \begin{bmatrix} 0 & -1 \\ 1 & 0 \end{bmatrix} \cup \begin{bmatrix} 3 & 1 & -1 \\ 2 & 2 & -1 \\ 2 & 2 & 0 \end{bmatrix} \cup \begin{bmatrix} 0 & 1 & 0 & 1 \\ 1 & 0 & 1 & 0 \\ 0 & 1 & 0 & 1 \\ 1 & 0 & 1 & 0 \end{bmatrix}$$

be a 3-matrix over the 3-field F = $Z_2 \cup Z_5 \cup Z_3$. det(xI – A) = det($xI_1 - A_1$) $\cup$ det($xI_2 - A_2$) $\cup$ det($xI_3 - A_3$) is a 3-polynomial of 3-degree(2, 3, 4), can be easily obtained.

We see the 3-polynomial is a monic 3-polynomial.

**DEFINITION 1.2.26:** *Let $T = T_1 \cup T_2 \cup \ldots \cup T_n$ be a non-linear operator over the n-space $V = V_1 \cup V_2 \cup \ldots \cup V_n$. We say T is n-diagonalizable if there is a n-basis for V for each n-vector of which is a n-characteristic vector of T.*

*Suppose $T = T_1 \cup T_2 \cup \ldots \cup T_n$ is a n-diagonalizable n-linear operator. Let $\{C_1^1, \ldots, C_{k_1}^1\} \cup \{C_1^2, \ldots, C_{k_2}^2\} \cup \ldots \cup \{C_1^n, \ldots, C_{k_n}^n\}$ be the n-distinct n-characteristic values of T. Then there is an ordered n-basis $B = B_1 \cup \ldots \cup B_n$ in which T is represented by a n-diagonal matrix which has for its n-diagonal entries the scalars $C_i^t$ each repeated a certain number of times t = 1, 2, ..., n. If $C_i^t$ is repeated $d_i^t$ times then (we may arrange that) the n-matrix has the n-block form*

$$[T]_B = [T_1]_{B_1} \cup \ldots \cup [T_n]_{B_n}$$

$$= \begin{bmatrix} C_1^1 I_1^1 & 0 & \ldots & 0 \\ 0 & C_2^1 I_2^1 & \ldots & 0 \\ \vdots & \vdots & & \vdots \\ 0 & 0 & \ldots & C_{k_1}^1 I_{k_1}^1 \end{bmatrix} \cup \begin{bmatrix} C_1^2 I_1^2 & 0 & \ldots & 0 \\ 0 & C_2^2 I_2^2 & \ldots & 0 \\ \vdots & \vdots & & \vdots \\ 0 & 0 & \ldots & C_{k_2}^2 I_{k_2}^2 \end{bmatrix} \cup$$

$$\cup \ldots \cup \begin{bmatrix} C_1^n I_1^n & 0 & \ldots & 0 \\ 0 & C_2^n I_2^n & \ldots & 0 \\ \vdots & \vdots & & \vdots \\ 0 & 0 & \ldots & C_{k_n}^n I_{k_n}^n \end{bmatrix}$$

*where* $I_j^t$ *is the* $d_j^t \times d_j^t$ *identity matrix.*

*From this n-matrix we make the following observations. First the n-characteristic n-polynomial for* $T = T_1 \cup T_2 \cup \ldots \cup T_n$ *is the n-product of n-linear factors*

$$\begin{aligned} f &= f_1 \cup f_2 \cup \ldots \cup f_n \\ &= (x - C_1^1)^{d_1^1} \ldots (x - C_{k_1}^1)^{d_{k_1}^n} \cup (x - C_1^2)^{d_1^2} \ldots (x - C_{k_2}^2)^{d_{k_2}^n} \\ &\quad \cup \ldots \cup (x - C_1^n)^{d_1^n} \ldots (x - C_{k_n}^n)^{d_{k_n}^n}. \end{aligned}$$

*If the n-scalar field* $F = F_1 \cup F_2 \cup \ldots \cup F_n$ *is n-algebraically closed i.e., if each* $F_i$ *is algebraically closed for i = 1, 2, …, n; then every n-polynomial over* $F = F_1 \cup F_2 \cup \ldots \cup F_n$ *can be n-factored; however if* $F = F_1 \cup F_2 \cup \ldots \cup F_n$ *is not n-algebraically closed we are citing a special property of* $T = T_1 \cup \ldots \cup T_n$ *when we say that its n-characteristic polynomial has such a factorization. The second thing to be noted is that* $d_i^t$ *is the number of times* $C_i^t$ *is repeated as a root of* $f_t$ *which is equal to the dimension of the space in* $V_t$ *of characteristic vectors associated with the characteristic value* $C_i^t$*; i = 1, 2, …,* $k_t$*, t = 1, 2, …, n. This is because the n-nullity of a n-diagonal n-matrix is equal to the number of n-zeros which has on its main n-diagonal and the n-matrix,* $[T - CI]_B = [T_1 - C_{i_1}^1 I]_{B_1} \cup \ldots \cup [T_n - C_{i_n}^n I]_{B_n}$ *has* $(d_{i_1}^1 \ldots d_{i_n}^n)$*, n-zeros on its main n-diagonal.*

This relation between the n-dimension of the n-characteristic space and the n-multiplicity of the n-characteristic value as a n-root of f does not seem exciting at first, however it will provide us with a simpler way of determining whether a given n-operator is n-diagonalizable.

**LEMMA 1.2.4:** *Suppose that $T\alpha = C\alpha$, where $T = T_1 \cup \ldots \cup T_n$ is the n-linear operator $C = C_1 \cup \ldots \cup C_n$ where the n-scalar $C_i \in F_i$; $i = 1, 2, \ldots, n$ and $\alpha = \alpha_1 \cup \ldots \cup \alpha_n$ is a n-vector. If $f = f_1 \cup f_2 \cup \ldots \cup f_n$ is any n-polynomial then $f(T)\alpha = f(C)\alpha$ i.e., $f_1(T_1)\alpha_1 \cup \ldots \cup f_n(T_n)\alpha_n = f_1(C_1)\alpha_1 \cup \ldots \cup f_n(C_n)\alpha_n$.*

It can be easily proved by any reader.

**LEMMA 1.2.5:** *Let $T = T_1 \cup T_2 \cup \ldots \cup T_n$ be a n-linear operator on the finite $(n_1, n_2, \ldots, n_n)$ dimensional n-vector space $V = V_1 \cup \ldots \cup V_n$ over the n-field $F = F_1 \cup \ldots \cup F_n$. If*

$$\left\{C_1^1 \ldots C_{k_1}^1\right\} \cup \left\{C_1^2 \ldots C_{k_2}^2\right\} \cup \ldots \cup \left\{C_1^n C_2^n \ldots C_{k_n}^n\right\}$$

*be the distinct n-characteristic values of T and let $W_i = W_{i_1}^1 \cup W_{i_2}^2 \cup \ldots \cup W_{i_n}^n$ be the n-space of n-characteristic n-vectors associated with the n-characteristic values $C_i = C_{i_1}^1 \cup C_{i_2}^2 \cup \ldots \cup C_{i_n}^n$. If $W = \{W_1^1 + \ldots + W_{k_1}^1\} \cup \{W_1^2 + \ldots + W_{k_2}^2\} \cup \ldots \cup \{W_1^n + \ldots + W_{k_n}^n\}$ the n-dim $W = ((\dim W_1^1 + \ldots + \dim W_{k_1}^1), (\dim W_1^2 + \ldots + \dim W_{k_2}^2), \ldots, (\dim W_1^n + \ldots + W_{k_n}^n)) = \dim W^1 \cup \ldots \cup \dim W^n$. In fact if $B_{i_t}^t$ is the ordered basis of $W_{i_t}^t$, $t = 1, 2, \ldots, k_t$ and $t = 1, 2, \ldots, n$ then $B = \{B_1^1, \ldots, B_{k_1}^1\} \cup \{B_1^2, \ldots, B_{k_2}^2\} \cup \ldots \cup \{B_1^n, \ldots, B_{k_n}^n\}$ is an n-ordered n-basis of W.*

*Proof:* We will prove the result for a $W^t$ for that will hold for every t, t = 1, 2, …, n. Let $W^t = W_1^t \cup \ldots \cup W_{k_t}^t$, $1 \leq t \leq n$, be the subspace spanned by all the characteristic vector of $T_t$. Usually when we form the sum $W^t$ of subspaces $W_i^t$ we see that $\dim W^t < \dim W_1^t + \ldots + \dim W_{k_t}^t$ because of linear relations which may exist between vectors in the various spaces. This result state the characteristic spaces associated with different characteristic values are independent of one another.

Suppose that (for each i) we have a vector $\beta_i^t$ in $W_i^t$ and assume that $\beta_1^t + \ldots + \beta_{k_t}^t = 0$. We shall show that $\beta_i^t = 0$ for each i. Let f be any polynomial. Since $T_t\beta_i^t = C_i^t\beta_i^t$ the proceeding lemma tells us that 0 = $f_t(T_t)$ thus 0 = $f_t(T_t)\beta_1^t + \ldots + f_t(T_t)\beta_{k_t}^t = f_t(C_1^t)\beta_1^t + \ldots + f_t(C_{k_t}^t)\beta_{k_t}^t$. Choose polynomials $f_1^t + \ldots + f_{k_t}^t$ such that

$$f_i^t\left(C_j^t\right) = \delta_{ij}^t = \begin{cases} 1 & i = j \\ 0 & i \neq j \end{cases}.$$

Then $0 = f_i^t\left(T_t\right), 0 = \sum_j \delta_{ij}^t\beta_j^t = \beta_i^t$.

Now let $B_i^t$ be an ordered basis for $W_i^t$ and let $B^t$ be the sequence $B^t = \left(B_1^t, \ldots, B_{k_t}^t\right)$. Then $B^t$ spans the subspace $W^t = W_1^t, \ldots, W_{k_t}^t$. Also $B^t$ is a linearly independent sequence of vectors for the following reason.

Any linear relation between the vectors in $B^t$ will have the form $\beta_1^t + \beta_2^t + \ldots + \beta_{k_t}^t = 0$, where $\beta_i^t$ is some linear combination of the vectors in $B_i^t$. From what we just proved we see that $\beta_i^t = 0$ for each i = 1, 2, …, $k_t$. Since each $B_i^t$ is linearly independent we see that we have only the trivial linear relation between the vectors in $B^t$.

This is true for each t = 1, 2, …, n. Thus we have

$$\begin{aligned} \dim W &= ((\dim W_1^1 + \ldots + \dim W_{k_1}^1), (\dim W_1^2 + \ldots + \dim W_{k_2}^2), \\ &\quad \ldots, (\dim W_1^n + \ldots + \dim W_{k_n}^n)) \\ &= (\dim W^1, \dim W^2, \ldots, \dim W^n). \end{aligned}$$

We leave the proof of the following theorem to the reader.

**THEOREM 1.2.44:** *Let $T = T_1 \cup \ldots \cup T_n$ be a n-linear operator of a finite $(n_1, n_2, \ldots, n_n)$ dimensional n-vector space $V = V_1 \cup V_2 \cup \ldots \cup V_n$ over the n-field $F = F_1 \cup \ldots \cup F_n$. Let $\left\{C_1^1, C_2^1, \ldots, C_{k_1}^1\right\} \cup \ldots \cup \left\{C_1^n, C_2^n, \ldots, C_{k_n}^n\right\}$ be the distinct n-* *characteristic values of T and let* $W_i = W^1_{i_1} \cup W^2_{i_2} \cup \ldots \cup W^n_{i_n}$ *be the n-null space of* $T - C_iI = (T_1 - C^i_1 I) \cup \ldots \cup (T_n - C^i_{k_i} I)$. *The following are equivalent:*

*(i)* *T is n-diagonalizable.*

*(ii)* *The n-characteristic n-polynomial for* $T = T_1 \cup \ldots \cup T_n$ *is* $f = f^1 \cup \ldots \cup f^n$

$$= \left(x - C^1_1\right)^{d^1_1} \ldots \left(x - C^1_{k_1}\right)^{d^1_{k_1}} \cup \left(x - C^2_1\right)^{d^2_1} \ldots \left(x - C^2_{k_2}\right)^{d^2_{k_2}}$$

$$\cup \ldots \cup \left(x - C^n_1\right)^{d^n_1} \ldots \left(x - C^n_{k_n}\right)^{d^n_{k_n}}$$

*and* $\dim W_i = \left(d^1_{i_1}, d^2_{i_2}, \ldots, d^n_{i_n}\right)$, $i = 1, 2, \ldots, k$.

*(iii)* $((\dim W^1_1 + \ldots + \dim W^1_{k_1}), (\dim W^2_1 + \ldots + \dim W^2_{k_2}), \ldots, (\dim W^n_1 + \ldots + \dim W^n_{k_n})) =$ *n-dim* $V = (n_1, n_2, \ldots, n_n)$.

*The n-matrix analogous of the above theorem may be formulated as follows. Let* $A = A_1 \cup \ldots \cup A_n$ *be a* $(n_1\times n_1, n_2\times n_2, \ldots, n_n\times n_n)$ *n-matrix with entries from the n-field,* $F = F_1 \cup \ldots \cup F_n$ *and let* $\left\{C^1_1, \ldots, C^1_{k_1}\right\} \cup \left\{C^1_2, \ldots, C^2_{k_2}\right\} \cup \ldots \cup \left\{C^n_1, \ldots, C^n_{k_n}\right\}$ *be the n-distinct n-characteristic values of A in F. For each i let* $W_i = W^1_{i_1} + \ldots + W^n_{i_n}$ *be the n-subspace of all n-column matrices* $X = X_1 \cup \ldots \cup X_n$ *(with entries from the n-field* $F = F_1 \cup \ldots \cup F_n$*) such that* $(A - C_iI)X = (A_1 - C^1_{i_1} I_{i_1})X_1 \cup \ldots \cup (A_n - C^n_{i_n} I_{i_n})X_n = (0 \cup \ldots \cup 0) = (B^1_{i_1} \cup \ldots \cup B^n_{i_n})$ *and let* $B_i$ *be an n-ordered n-basis for* $W_i = W^1_{i_1} + \ldots + W^n_{i_n}$. *The n-basis* $B = B_1, \ldots, B_n$ *collectively string together to form the n-sequence of n-columns of a n-matrix* $P = [P^1_1, \ldots, P^n_n] = \{B_1, \ldots, B_n\} = \left\{B^1_1, \ldots, B^1_{k_1}\right\} \cup \left\{B^1_2, \ldots, B^1_{k_2}\right\} \cup \ldots \cup \left\{B^n_1, B^n_2, \ldots, B^n_{k_n}\right\}$. *The n-matrix A over the n-field* $F = F_1 \cup \ldots \cup F_n$ *is similar to a n-diagonal n-matrix if and only if P is a* $(n_1\times n_1, \ldots, n_n\times n_n)$ *n-square matrix. When P is square, P is n-invertible and* $P^{-1}AP$ *is n-diagonal i.e., each* $P_i^{-1}A_iP_i$ *is diagonal in* $P^{-1}AP = P_1^{-1}A_1P_1 \cup \ldots \cup P_n^{-1}A_nP_n$.

Now we proceed on to define the new notion of n-annihilating n-polynomials of the n-linear operator T, the n-minimal n-polynomial for T and the analogue of the Cayley-Hamilton theorem for a n-vector space V over n-field of type II and its n-linear operator on V.

In order to know more about the n-linear operator $T = T_1 \cup \ldots \cup T_n$ on $V_1 \cup \ldots \cup V_n$ over the n-field $F = F_1 \cup \ldots \cup F_n$ one of the most useful things to know is the class of n-polynomials which n-annihilate $T = T_1 \cup \ldots \cup T_n$.

To be more precise suppose $T = T_1 \cup \ldots \cup T_n$ is a n-linear operator on V, a n-vector space over the n-field $F = F_1 \cup \ldots \cup F_n$. If $p = p_1 \cup \ldots \cup p_n$ is a n-polynomial over the n-field $F = F_1 \cup \ldots \cup F_n$ then $p(T) = p_1(T_1) \cup \ldots \cup p_n(T_n)$ is again a n-linear operator on $V = V_1 \cup \ldots \cup V_n$.

If $q = q_1 \cup \ldots \cup q_n$ is another n-polynomial over the same n-field $F = F_1 \cup \ldots \cup F_n$ then

$$(p + q)T = (p)T + (q)T$$

i.e.,

$$(p_1 + q_1)T_1 \cup \ldots \cup (p_n + q_n)T_n$$
$$= [p_1(T_1) \cup \ldots \cup p_n(T_n)] + [q_1(T_1) \cup \ldots \cup q_n(T_n)]$$

and

$$(pq)T = (p)T(q)T$$

i.e.,

$$(p_1q_1)T_1 \cup \ldots \cup (p_nq_n)T_n = [p_1T_1q_1T_1 \cup \ldots \cup p_nT_nq_nT_n].$$

Therefore the collection of n-polynomials $P = P_1 \cup \ldots \cup P_n$ which n-annihilate $T = T_1 \cup \ldots \cup T_n$ in the sense that $p(T) = p_1(T_1) \cup \ldots \cup p_n(T_n) = 0 \cup \ldots \cup 0$ is an n-ideal of the n-polynomial n-algebra $F(x) = F_1(x) \cup \ldots \cup F_n(x)$. It may be the zero n-ideal i.e., it may be, that T is not n-annihilated by any non-zero n-polynomial. But that cannot happen if the n-space $V = V_1 \cup \ldots \cup V_n$ is finite dimensional i.e., V is of $(n_1, n_2, \ldots, n_n)$ dimension over the n-field $F = F_1 \cup \ldots \cup F_n$.

Suppose $T = T_1 \cup T_2 \cup \ldots \cup T_n$ is a n-linear operator on the $(n_1, n_2, \ldots, n_n)$ dimension n-space $V = V_1 \cup V_2 \cup \ldots \cup V_n$. The first n, $(n_i^2 + 1)$ operators, (i = 1, 2, …, n) must be n-linearly dependent i.e., the first $(n_1^2 + 1,\ n_2^2 + 1,\ \ldots,\ n_n^2 + 1)$ n-linear operators are n-linearly dependent i.e., we have $C_0^1 I_1 + C_1^1 T_1 + \ldots + C_{n_1^2}^1 T_1^{n_1^2} = 0$, $C_0^2 I_2 + C_1^2 T_2 + \ldots + C_{n_2^2}^2 T_2^{n_2^2} = 0$ and so on, $C_0^n I_n + C_1^n T_n + \ldots + C_{n_n^2}^n T_n^{n_n^2} = 0$;

That is $\{C_0^1 I_1 + C_1^1 T_1 + \ldots + C_{n_1^2}^1 T_1^{n_1^2}\} \cup \{C_0^2 I_2 + C_1^2 T_2 + \ldots + C_{n_2^2}^2 T_2^{n_2^2}\} \cup \ldots \cup \{C_0^n I_n + C_1^n T_n + \ldots + C_{n_n^2}^n T_n^{n_n^2}\} = 0 \cup 0 \cup \ldots \cup 0$ for some n-scalars; $C_{i_1}^1, C_{i_2}^2, \ldots, C_{i_n}^n$ not all zero. $1 \leq i_p \leq n_p$ and p = 1, 2, …, n. Thus the n-ideal of n-polynomials which n-annihilate T contains a non zero n-polynomial of n-degree $(n_1^2, n_2^2, \ldots, n_n^2)$ or less.

We know that every n-polynomial n-ideal consists of all n-multiples of some fixed n-monic n-polynomials which is the n-generator of the n-ideal. Thus there corresponds to the n-operator $T = T_1 \cup \ldots \cup T_n$ a n-monic n-polynomial $p = p_1 \cup \ldots \cup p_n$.

If f is any other n-polynomial over the n-field $F = F_1 \cup \ldots \cup F_n$ then $f(T) = 0 \cup 0 \cup \ldots \cup 0$ i.e., $f_1(T_1) \cup \ldots \cup f_n(T_n) = 0 \cup 0 \cup \ldots \cup 0$ if and only if f = pg where $g = g_1 \cup \ldots \cup g_n$ is some polynomial over the n-field $F = F_1 \cup \ldots \cup F_n$ i.e., $f = f_1 \cup \ldots \cup f_n = p_1 g_1 \cup \ldots \cup p_n g_n$.

Now we define the new notion of n-polynomial for the n-operator $T: V \to V$.

**DEFINITION 1.2.27:** *Let $T = T_1 \cup \ldots \cup T_n$ be a n-linear operator on a finite $(n_1, \ldots, n_n)$ dimensional n-vector space $V = V_1 \cup \ldots \cup V_n$ over the field $F_1 \cup \ldots \cup F_n$. The n-minimal n-polynomial for T is the (unique) monic n-generator of the n-ideal of n-polynomials over the n-field, $F = F_1 \cup \ldots \cup F_n$ which n-annihilate $T = T_1 \cup \ldots \cup T_n$.*

*The n-minimal n-polynomial starts from the fact that the n-generator of a n-polynomial n-ideal is characterised by being the n-monic n-polynomial of n-minimum n-degree in the n-ideal that implies that the n-minimal n-polynomial $p = p_1 \cup \ldots \cup p_n$*

*for the n-linear operator $T = T_1 \cup T_2 \cup \ldots \cup T_n$ is uniquely determined by the following properties.*

1. *p is a n-monic n-polynomial over the n-scalar n-field $F = F_1 \cup \ldots \cup F_n$.*
2. *$p(T) = p_1(T_1) \cup \ldots \cup p_n(T_n) = 0 \cup \ldots \cup 0$.*
3. *No n-polynomial over the n-field $F = F_1 \cup \ldots \cup F_n$ which n-annihilates $T = T_1 \cup T_2 \cup \ldots \cup T_n$ has smaller n-degree than $p = p_1 \cup \ldots \cup p_n$ has.*

*$(n_1 \times n_1, n_2 \times n_2, \ldots, n_n \times n_n)$ is the order of n-matrix $A = A_1 \cup \ldots \cup A_n$ over the n-field $F = F_1 \cup \ldots \cup F_n$ where each $A_i$ has $n_i \times n_i$ matrix with entries from the field $F_i$, i = 1, 2, …, n.*

*The n-minimal n-polynomial for $A = A_1 \cup \ldots \cup A_n$ is defined in an analogous way as the unique n-monic generator of the n-ideal of all n-polynomial over the n-field, $F = F_1 \cup \ldots \cup F_n$ which n-annihilate A.*

*If the n-operator $T = T_1 \cup T_2 \cup \ldots \cup T_n$ is represented by some ordered n-basis by the n-matrix $A = A_1 \cup \ldots \cup A_n$ then T and A have same n-minimal polynomial because $f(T) = f_1(T_1) \cup \ldots \cup f_n(T_n)$ is represented in the n-basis by the n-matrix $f(A) = f_1(A_1) \cup \ldots \cup f_n(A_n)$ so $f(T) = 0 \cup \ldots \cup 0$ if and only if $f(A) = 0 \cup \ldots \cup 0$ i.e., $f_1(A_1) \cup \ldots \cup f_n(A_n) = 0 \cup \ldots \cup 0$ if and only if $f_1(T_1) \cup \ldots \cup f_n(T_n) = 0 \cup \ldots \cup 0$. So $f(P^{-1}AP) = f_1(P^{-1}A_1P_1) \cup \ldots \cup f_n(P_n^{-1}A_nP_n) = P_1^{-1}f_1(A_1)P_1 \cup \ldots \cup P_n^{-1}f_n(A_n)P_n^{-1} = P^{-1}f(A)P$ for every n-polynomial $f = f_1 \cup f_2 \cup \ldots \cup f_n$.*

*Another important feature about the n-minimal polynomials of n-matrices is that suppose $A = A_1 \cup \ldots \cup A_n$ is a $(n_1 \times n_1, \ldots, n_n \times n_n)$ n-matrix with entries from the n-field $F = F_1 \cup \ldots \cup F_n$. Suppose $K = K_1 \cup \ldots \cup K_n$ is n-field which contains the n-field $F = F_1 \cup \ldots \cup F_n$ i.e., $K \supseteq F$ and $K_i \supseteq F_i$ for every i, i = 1, 2, …, n. $A = A_1 \cup \ldots \cup A_n$ is a $(n_1 \times n_1, \ldots, n_n \times n_n)$ n-matrix over F or over K but we do not obtain two n-minimal n-polynomial only one minimal n-polynomial.*

This is left as an exercise for the reader to verify. Now we proceed on to prove an interesting theorem about the n-minimal polynomials for T(or A).

**THEOREM 1.2.45:** *Let $T = T_1 \cup \ldots \cup T_n$ be a n-linear operator on a $(n_1, n_2, \ldots, n_n)$ dimensional n-vector space $V = V_1 \cup \ldots \cup V_n$ [or let A be a $(n_1 \times n_1, \ldots, n_n \times n_n)$n-matrix i.e., $A = A_1 \cup \ldots \cup A_n$ where each $A_i$ is a $n_i \times n_i$ matrix with its entries from the field $F_i$ of $F = F_1 \cup \ldots \cup F_n$, true for i = 1, 2, …, n].*

*The n-characteristic and n-minimal n-polynomial for T [for A] have the same n-roots except for n-multiplicities.*

*Proof:* Let $p = p_1 \cup \ldots \cup p_n$ be a n-minimal n-polynomial for $T = T_1 \cup \ldots \cup T_n$. Let $C = C_1 \cup \ldots \cup C_n$ be a n-scalar of the n-field $F = F_1 \cup \ldots \cup F_n$. To prove $p(C) = p_1(C_1) \cup \ldots \cup p_n(C_n) = 0 \cup \ldots \cup 0$ if and only if $C = C_1 \cup \ldots \cup C_n$ is the n-characteristic value of T. Suppose $p(C) = p_1(C_1) \cup \ldots \cup p_n(C_n) = 0 \cup \ldots \cup 0$; then $p = (x - C_1)q_1 \cup (x - C_2)q_2 \cup \ldots \cup (x - C_n)q_n$ where $q = q_1 \cup \ldots \cup q_n$ is a n-polynomial, since n-deg q < n-deg p, the n-minimal n-polynomial $p = p_1 \cup \ldots \cup p_n$ tells us $q(T) = q_1(T_1) \cup \ldots \cup q_n(T_n) \neq 0 \cup \ldots \cup 0$. Choose a n-vector $\beta = \beta_1 \cup \ldots \cup \beta_n$ such that $q(T)\beta = q_1(T_1)\beta_1 \cup \ldots \cup q_n(T_n)\beta_n \neq 0 \cup \ldots \cup 0$. Let $\alpha = q(T)\beta$ i.e., $\alpha = \alpha_1 \cup \ldots \cup \alpha_n = q_1(T_1)\beta_1 \cup \ldots \cup q(T_n)\beta_n$. Then

$$
\begin{aligned}
0 \cup \ldots \cup 0 &= p(T)\beta = p_1(T_1)\beta_1 \cup \ldots \cup p_n(T_n)\beta_n. \\
&= (T - CI)q(T)\beta \\
&= (T_1 - C_1I_1)q_1(T_1)\beta_1 \cup \ldots \cup (T_n - C_nI_n)q_n(T_n)\beta_n \\
&= (T_1 - C_1I_1)\alpha_1 \cup \ldots \cup (T_n - C_nI_n)\alpha_n
\end{aligned}
$$

and thus $C = C_1 \cup \ldots \cup C_n$ is a n-characteristic value of $T = T_1 \cup \ldots \cup T_n$.

Suppose $C = C_1 \cup \ldots \cup C_n$ is a n-characteristic value of $T = T_1 \cup \ldots \cup T_n$ say $T\alpha = C\alpha$ i.e., $T_1\alpha_1 \cup \ldots \cup T_n\alpha_n = C_1\alpha \cup \ldots \cup C_n\alpha$ with $\alpha \neq 0 \cup \ldots \cup 0$. From the earlier results we have $p(T)\alpha = p(C)\alpha$ i.e., $p_1(T_1)\alpha_1 \cup \ldots \cup p_n(T_n)\alpha_n. = p_1(C_1)\alpha_1 \cup \ldots \cup p_n(C_n)\alpha_n$; since $p(T) = p_1(T_1) \cup \ldots \cup p_n(T_n) = 0 \cup \ldots \cup 0$ and $\alpha = \alpha_1 \cup \ldots \cup \alpha_n \neq 0$ we have $p(C) = p_1(C_1) \cup \ldots \cup p_n(C_n) \neq 0 \cup \ldots \cup 0$.

Let $T = T_1 \cup \ldots \cup T_n$ be a n-diagonalizable n-linear operator and let $\{C_1^1 \ldots C_{k_1}^1\} \cup \{C_1^2 \ldots C_{k_2}^2\} \cup \ldots \cup \{C_1^n \ldots C_{k_n}^n\}$ be the n-distinct n-characteristic values of T. Then the n-minimal n-polynomial for T is the n-polynomial $p = p_1 \cup \ldots \cup p_n = (x - C^1_1) \ldots (x - C^1_{k_1}) \cup (x - C^2_1) \ldots (x - C^2_{k_2}) \cup \ldots \cup (x - C^n_1) \ldots (x - C^n_{k_n})$.

If $\alpha = \alpha_1 \cup \ldots \cup \alpha_n$ is a n-characteristic n-vector then one of the n-operators $\{(T_1 - C^1_1 I_1), \ldots, (T_1 - C^1_{k_1} I_1)\}$, $\{(T_2 - C^2_1 I_2), \ldots, (T_2 - C^2_{k_2} I_2)\}$, …, $\{(T_n - C^n_1 I_n), \ldots, (T_n - C^n_{k_n} I_n)\}$ send $\alpha = \alpha_1 \cup \ldots \cup \alpha_n$ into $0 \cup \ldots \cup 0$, thus resulting in $\{(T_1 - C^1_1 I_1), \ldots, (T_1 - C^1_{k_1} I_1)\}$, $\{(T_2 - C^2_1 I_2), \ldots, (T_2 - C^2_{k_2} I_2)\}$, $\{(T_n - C^n_1 I_n), \ldots, (T_n - C^n_{k_n} I_n)\} = 0 \cup \ldots \cup 0$ for every n-characteristic n-vector $\alpha = \alpha_1 \cup \ldots \cup \alpha_n$.

Hence there exists an n-basis for the underlying n-space which consist of n-characteristic vectors of $T = T_1 \cup \ldots \cup T_n$. Hence $p(T) = p_1(T_1) \cup \ldots \cup p_n(T_n) = \{(T_1 - C^1_1 I_1), \ldots, (T_1 - C^1_{k_1} I_1)\} \cup \ldots \cup \{(T_n - C^n_1 I_n), \ldots, (T_n - C^n_{k_n} I_n)\} = 0 \cup \ldots \cup 0$.

Thus we can conclude if T is n-diagonalizable, n-linear operator then the n-minimal n-polynomial for T is a product of n-distinct n-linear factors.

**THEOREM 1.2.46: (CAYLEY-HAMILTON):** *Let $T = T_1 \cup \ldots \cup T_n$ be a n-linear operator on a finite $(n_1, n_2, \ldots, n_n)$ dimensional vector space $V = V_1 \cup \ldots \cup V_n$ over the n-field $F = F_1 \cup \ldots \cup F_n$. If $f = f_1 \cup \ldots \cup f_n$ is the n-characteristic, n-polynomial for T then $f(T) = f_1(T_1) \cup \ldots \cup f_n(T_n) = 0 \cup \ldots \cup 0$; in otherwords the n-minimal polynomial divides the n-characteristic polynomial for T.*

*Proof:* Let $K = K_1 \cup \ldots \cup K_n$ be a n-commutative ring with n-identity $1_n = (1, \ldots, 1)$ consisting of all n-polynomials in T; K is actually a n-commutative algebra with n-identity over the scalar n-field $F = F_1 \cup \ldots \cup F_n$.

Let $\{\alpha^1_1 \ldots \alpha^1_{n_1}\} \cup \ldots \cup \{\alpha^n_1 \ldots \alpha^n_{n_n}\}$ be an ordered n-basis for V and let $A = A^1 \cup \ldots \cup A^n$ be the n-matrix which represents $T = T_1 \cup \ldots \cup T_n$ in the given n-basis. Then

$$T\alpha_i = T_1\alpha^1_{i_1} \cup \ldots \cup T_n\alpha^n_{i_n}$$

$$= \sum_{j_1=1}^{n_1} A^1_{j_1 i_1}\alpha^1_{j_1} \cup \sum_{j_2=1}^{n_2} A^2_{j_2 i_2}\alpha^2_{j_2} \cup \ldots \cup \sum_{j_n=1}^{n_n} A^n_{j_n i_n}\alpha^n_{j_n};$$

$1 \le j_i \le n_{j_i}$, i = 1, …, n. These n-equations may be equivalently written in the form

$$\sum_{j_1=1}^{n_1}\left(\delta_{j_1 i_1}T_1 - A^1_{j_1 i_1}I_{i_1}\right)\alpha^1_{j_1} \cup \sum_{j_2=1}^{n_2}\left(\delta_{j_2 i_2}T_2 - A^2_{j_2 i_2}I_{i_2}\right)\alpha^2_{j_2.}$$

$$\cup \ldots \cup \sum_{j_n=1}^{n_n}\left(\delta_{j_n i_n}T_n - A^n_{j_n i_n}I_{i_n}\right)\alpha^n_{j_n}$$

$$= 0 \cup \ldots \cup 0,\ 1 \le i_n \le n.$$

Let $B = B^1 \cup \ldots \cup B^n$ denote the element of $K_1^{n_1\times n_1} \cup \ldots \cup K_n^{n_n \times n_n}$ i.e., $B^i$ is an element of $K_i^{n_i \times n_i}$ with entries $B^t_{i_t j_t} = \delta_{i_t j_t}T_t - A_{i_t j_t}I_t$, t = 1, 2, …, n. When $n_t = 2$; $1 \le i_t, j_t \le n_t$.

$$B^t = \begin{bmatrix} T_t - A^t_{11}I_t & -A^t_{21}I_t \\ -A^t_{12}I_t & T_t - A_{22}I_t \end{bmatrix}$$

and det $B^t = (T_t - A^t_{11} I_t)(T_t - A^t_{22} I_t) - (A^t_{12} A^t_{21})I_t = f_t(T_t)$ where $f_t$ is the characteristic polynomial associated with $T_t$, t = 1, 2, …, n. $f_t = x^2$-trace $A^t x$ + det $A^t$. For case $n_t > 2$ it is clear that det $B^t = f_t(T_t)$ since $f_t$ is the determinant of the matrix $xI_t - A_t$ whose entries are polynomial $\left(xI_t - A^t\right)_{i_t j_t} = \delta_{i_t j_t}x - A^t_{i_t j_t}$.

We will show $f_t(T_t) = 0$. In order that $f_t(T_t)$ is a zero operator, it is necessary and sufficient that $(\det B^t)_{\alpha^t_k} = 0$ for $k_t = 0, 1, \ldots, n_t$.

By the definition of $B^t$; the vectors $\alpha^t_1 \cup \ldots \cup \alpha^t_{n_t}$ satisfy the equations;

$$\sum_{j_t=0}^{n_t} B^t_{i_t j_t} \alpha^t_{jt} = 0,$$

$1 \le i_t \le n_t$. When $n_t = 2$ we can write the above equation in the form

$$\begin{bmatrix} T_t - A^t_{11} I_t & -A^t_{21} I_t \\ -A^t_{21} & T_t - A^t_{22} I_t \end{bmatrix} \begin{bmatrix} \alpha^t_1 \\ \alpha^t_2 \end{bmatrix} = \begin{bmatrix} 0 \\ 0 \end{bmatrix}.$$

In this case the usual adjoint $B^t$ is the matrix

$$\tilde{B}^t = \begin{bmatrix} T_t - A^t_{22} I_t & A^t_{21} I \\ A^t_{12} I & T_t - A^t_{11} I \end{bmatrix}$$

and

$$\tilde{B}^t B^t = \begin{bmatrix} \det B^t & 0 \\ 0 & \det B^t \end{bmatrix}.$$

Hence

$$\det B^t \begin{bmatrix} \alpha^t_1 \\ \alpha^t_2 \end{bmatrix} = \left(\tilde{B}^t B^t\right) \begin{bmatrix} \alpha^t_1 \\ \alpha^t_2 \end{bmatrix}$$

$$= \tilde{B}^t B^t \begin{bmatrix} \alpha^t_1 \\ \alpha^t_2 \end{bmatrix} = \tilde{B}^t B^t \begin{bmatrix} \alpha^t_1 \\ \alpha^t_2 \end{bmatrix} = \begin{bmatrix} 0 \\ 0 \end{bmatrix}.$$

In the general case $\tilde{B}^t = \text{adj } B^t$. Then

$$\sum_{j_t=1}^{n_t} \tilde{B}^t_{k_t i_t} B^t_{i_t j_t} \alpha^t_j = 0,$$

for each pair $k_t$, $i_t$ and summing on $i_t$ we have

$$0 = \sum_{i_t=1}^{n_t} \sum_{j_t=1}^{n_t} \tilde{B}^t_{k_t i_t} B^t_{i_t j_t} \alpha^t_{j_t}$$

$$= \sum_{i_t=1}^{n_t} \left( \sum_{j_t=1}^{n_t} \tilde{B}^t_{k_t i_t} B^t_{i_t j_t} \alpha^t_{j_t} \right).$$

Now $\tilde{B}^t B^t = (\det B^t) I_t$ so that

$$\sum_{i_t=1}^{n} \tilde{B}^t_{k_t i_t} B^t_{i_t j_t} = \delta_{k_t j_t} \det B^t .$$

Therefore

$$0 = \sum_{j_t=1}^{n_t} \delta_{k_t j_t} \left(\det B^t\right) \alpha^t_{j_t} = \left(\det B^t\right) \alpha^t_{k_t},\ 1 \le k_t \le n_t.$$

Since this is true for each t, t = 1, 2, …, n we have $0 \cup \ldots \cup 0 = (\det B^t)\,\alpha^1_{k_1} \cup \ldots \cup (\det B^n)\,\alpha^n_{k_n}$, $1 \le k_i \le n_i$, i = 1, 2, …, n.

The Cayley-Hamilton theorem is very important for it is useful in narrowing down the search for the n-minimal n-polynomials of various n-operators.

If we know the n-matrix $A = A^1 \cup \ldots \cup A^n$ which represents $T = T_1 \cup \ldots \cup T_n$ in some ordered n-basis then we can compute the n-characteristic polynomial $f = f_1 \cup \ldots \cup f_n$. We know the n-minimal polynomial $p = p_1 \cup \ldots \cup p_n$ n-divides f i.e., each $p_i/f_i$ for i = 1, 2, …, n (which we call as n-divides f) and that the two n-polynomials have the same n-roots.

However we do not have a method of computing the roots even in case of polynomials so more difficult is finding the n-roots of the n-polynomials. However if $f = f_1 \cup \ldots \cup f_n$ factors as

$$f = (x - C^1_1)^{d^1_1} \ldots (x - C^1_{k_1})^{d^1_{k_1}} \cup (x - C^2_1)^{d^2_1} \ldots (x - C^2_{k_2})^{d^2_{k_2}} \cup \ldots \cup (x - C^n_1)^{d^n_1} \ldots (x - C^n_{k_n})^{d^n_{k_n}}\ \{C^1_1, \ldots, C^1_{k_1}\} \cup \{C^2_1, \ldots, C^2_{k_2}\} \cup \ldots \cup \{C^n_1, \ldots, C^n_{k_n}\}$$

distinct n-sets, $d^t_{i_t} \ge 1$, t = 1, 2, …, $k_t$ then

$$p = p_1 \cup \ldots \cup p_n = (x - C^1_1)^{r^1_1} \ldots (x - C^1_{k_1})^{r^1_{k_1}} \cup \ldots \cup (x - C^n_1)^{r^n_1} \ldots (x - C^n_{k_n})^{r^n_{k_n}};\ 1 \le r^t_j \le d^t_j .$$

Now we illustrate this by a simple example.

***Example 1.2.11:*** Let

$$A = \begin{bmatrix} 1 & 1 & 0 & 0 \\ -1 & -1 & 0 & 0 \\ -2 & -2 & 2 & 1 \\ 1 & 1 & -1 & 0 \end{bmatrix} \cup \begin{bmatrix} 3 & 1 & -1 \\ 2 & 2 & -1 \\ 2 & 2 & 0 \end{bmatrix} \cup \begin{bmatrix} 0 & -1 \\ 1 & 0 \end{bmatrix}$$

be a 3-matrix over the 3-field $F = Z_3 \cup Z_5 \cup Q$. Clearly the 3-characteristic 3-polynomial associated with A is given by $f = f_1 \cup f_2 \cup f_3 = x^2 (x - 1)^2 \cup (x - 1)(x - 2)^2 \cup x^2 + 1$. It is easily verified that $p = p_1 \cup \ldots \cup p_3 = x^2(x - 1)^2 \cup (x - 1)(x - 2)^2 \cup (x^2 + 1)$ is the 3-minimal 3-polynomial of A.

Now we proceed on to define the new notion of n-invariant subspaces or equivalently we may call it as invariant n-subspaces.

**DEFINITION 1.2.28:** *Let $V = V_1 \cup \ldots \cup V_n$ be a n-vector space over the n-field $F = F_1 \cup F_2 \cup \ldots \cup F_n$ of type II. Let $T = T_1 \cup \ldots \cup T_n$ be a n-linear operator on V. If $W = W_1 \cup \ldots \cup W_n$ is a n-subspace of V we say W is n-invariant under T if each of the n-vectors in W, i.e., for the n-vector $\alpha = \alpha_1 \cup \ldots \cup \alpha_n$ in W the n-vector $T\alpha = T_1\alpha_1 \cup \ldots \cup T_n\alpha_n$ is in W i.e., each $T_i\alpha_i \in W_i$ for every $\alpha_i \in W_i$ under the operator $T_i$ for i = 1, 2, ..., n i.e., if T(W) is contained in W i.e., if $T_i(W_i)$ is contained in $W_i$ for i = 1, 2, ..., n i.e., are thus represented as $T_1(W_1) \cup \ldots \cup T_n(W_n) \subseteq W_1 \cup \ldots \cup W_n$.*

This simple example is that we can say V the n-vector space is invariant under a n-linear operator T in V. Similarly the zero n-subspace is invariant under T.

Now we give the n-block matrix associated with T. Let $W = W_1 \cup \ldots \cup W_n$ be a n-subspace of the n-vector space $V = V_1 \cup \ldots \cup V_n$. Let $T = T^1 \cup \ldots \cup T^n$ be a n-operator on V such that W is n-invariant under the n-operator T then T induces a n-linear operator $T_w = T^1_{w_1} \cup \ldots \cup T^n_{w_n}$ on the n-space W. This n-linear operator $T_w$ defined by $T_w(\alpha) = T(\alpha)$ for $\alpha \in W$ i.e., if $\alpha = \alpha_1 \cup$

… $\cup$ $\alpha_n$ then $T_w(\alpha_1 \cup \ldots \cup \alpha_n) = T^1_{w_1}(\alpha_1) \cup \ldots \cup T^n_{w_n}(\alpha_n)$. Clearly $T_w$ is different from T as domain is W and not V. When V is finite dimensional say $(n_1, \ldots, n_n)$ dimensional the n-invariance of W under T has a simple n-matrix interpretation. Let $B = B_1 \cup \ldots \cup B_n = \{\alpha^1_1 \ldots \alpha^1_{r_1}\} \cup \ldots \cup \{\alpha^n_1 \ldots \alpha^n_{r_n}\}$ be a chosen n-basis for V such that $B' = B'_1 \cup \ldots \cup B'_n = \{\alpha^1_1 \ldots \alpha^1_{r_1}\} \cup \ldots \cup \{\alpha^n_1 \ldots \alpha^n_{r_n}\}$ ordered n-basis for W. $r = (r_1, r_2, \ldots, r_n) =$ n dim W.

Let $A = [T]_B$ i.e., if $A = A_1 \cup \ldots \cup A_n$ then $A = A_1 \cup \ldots \cup A_n = [T^1]_{B_1} \cup \ldots \cup [T^n]_{B_n}$ so that

$$T^t_{\alpha^t_{j_1}} = \sum_{i_t=1}^{n_t} A^t_{i_t j_t} \alpha^t_i$$

for t = 1, 2, …, n i.e.,

$$T\alpha_j = T^1 \alpha^1_{j_1} \cup \ldots \cup T^n \alpha^n_{j_n}$$

$$= \sum_{i_1=1}^{n_1} A^1_{i_1 j_1} \alpha^1_{i_1} \cup \sum_{i_2=1}^{n_2} A^2_{i_2 j_2} \alpha^2_{i_2} \cup \ldots \cup \sum_{i_n=1}^{n_n} A^n_{i_n j_n} \alpha^n_{i_n}.$$

Since W is n-invariant under T, the n-vector $T\alpha_j$ belongs to W for $j < r$ i.e., $(j_1 < r_1, \ldots, j_n < r_n)$.

$$T\alpha_j = \sum_{i_1=1}^{r_1} A^1_{i_1 j_1} \alpha^1_{i_1} \cup \ldots \cup \sum_{i_n=1}^{r_n} A^n_{i_n j_n} \alpha^n_{i_n}$$

i.e., $A^k_{i_k j_k} = 0$ if $j_k < r_k$ and $i_k > r_k$ for every k = 1, 2, …, n. Schematically A has the n-block

$$A = \begin{bmatrix} B & C \\ O & D \end{bmatrix} = \begin{bmatrix} B^1 & C^1 \\ O & D^1 \end{bmatrix} \cup \ldots \cup \begin{bmatrix} B^n & C^n \\ O & D^n \end{bmatrix}$$

where $B_t$ is a $r_t \times r_t$ matrix. $C_t$ is a $r_t \times (n_t - r_t)$ matrix and D is an $(n_t - r_t) \times (n_t - r_t)$ matrix for t = 1, 2, …, n. It is $B = B^1 \cup \ldots \cup B^n$ is the n-matrix induced by the n-operator $T_w$ in the n-basis $B' = B'_1 \cup \ldots \cup B'_n$.

**LEMMA 1.2.6:** *Let $W = W_1 \cup \ldots \cup W_n$ be an n-invariant subspace of the n-linear operator $T = T_1 \cup \ldots \cup T_n$ on $V = V_1 \cup \ldots \cup V_n$. The n-characteristic, n-polynomial for the n-restriction operator $T_W = T_{1W_1} \cup \ldots \cup T_{nW_n}$ divides the n-characteristic polynomial for T. The n-minimal polynomial for $T_W$ divides the n-minimal polynomial for T.*

*Proof:* We have

$$A = \begin{bmatrix} B & C \\ O & D \end{bmatrix} = \begin{bmatrix} B^1 & C^1 \\ O & D^1 \end{bmatrix} \cup \ldots \cup \begin{bmatrix} B^n & C^n \\ O & D^n \end{bmatrix}$$

where $A = [T]_B = T_{1B_1} \cup \ldots \cup T_{nB_n}$ and $B = [T_W]_{B'} = B^1 \cup \ldots \cup B^n = [T_{1w_1}]_{B'_1} \cup \ldots \cup [T_{nw_n}]_{B'_n}$. Because of the n-block form of the n-matrix

$\det(xI - A) = \det(xI_1 - A_1) \cup \det(xI_2 - A_2) \cup \ldots \cup \det(xI_n - A_n)$ where

$$\begin{aligned} A &= A_1 \cup \ldots \cup A_n \\ &= \det(xI - B) \det(xI - D) \\ &= \{\det (xI_1 - B^1) \det(xI_1 - D^1) \cup \ldots \cup \det(xI_n - B^n) \det(xI_n - D^n)\}. \end{aligned}$$

That proves the statement about n-characteristic polynomials. Notice that we used $I = I_1 \cup \ldots \cup I_n$ to represent n-identity matrix of these n-tuple of different sizes.

The $k^{th}$ power of the n-matrix A has the n-block form

$$A^k = (A^1)^k \cup \ldots \cup (A^n)^k.$$

$$A^k = \begin{bmatrix} (B^1)^k & (C^1)^k \\ O & (D^1)^k \end{bmatrix} \cup \ldots \cup \begin{bmatrix} (B^n)^k & (C^n)^k \\ O & (D^n)^k \end{bmatrix}$$

where $C^k = (C^1)^k \cup \ldots \cup (C^n)^k$ is $\{r_1 \times (n_1 - r_1), \ldots, r_n \times (n_n - r_n)\}$ n-matrix. Therefore any n-polynomial which n-annihilates A also n-annihilates B (and D too). So the n-minimal polynomial for B n-divides the n-minimal polynomial for A.

Let $T = T_1 \cup \ldots \cup T_n$ be any n-linear operator on a $(n_1, n_2, \ldots, n_n)$finite dimensional n-space $V = V_1 \cup \ldots \cup V_n$. Let $W = W_1 \cup \ldots \cup W_n$ be n-subspace spanned by all of the n-characteristic vectors of $T = T_1 \cup \ldots \cup T_n$. Let $\{C_1^1, \ldots, C_{k_1}^1\} \cup \{C_1^2, \ldots, C_{k_2}^2\} \cup \ldots \cup \{C_1^n, \ldots, C_{k_n}^n\}$ be the n-distinct characteristic values of T. For each i let $W_i = W_{i_1}^1 \cup \ldots \cup W_{i_n}^n$ be the n-space of n-characteristic vectors associated with the n-characteristic value $C_i = C_{i_1}^1 \cup \ldots \cup C_{i_n}^n$ and let $B_i = \{B_i^1 \cup \ldots \cup B_i^n\}$ be the ordered basis of $W_i$ i.e., $B_i^t$ is a basis of $W_i^t$. $B' = \{B_1^1, \ldots, B_{k_1}^1\} \cup \ldots \cup \{B_1^n, \ldots, B_{k_n}^n\}$ is a n-ordered n-basis for $W = (W_1^1 + \ldots + W_{k_1}^1) \cup \ldots \cup (W_1^n + \ldots + W_{k_n}^n) = W = W_1 \cup \ldots \cup W_n$.

In particular n-dimW = $\{(\dim W_1^1 + \ldots + \dim W_{k_1}^1), (\dim W_1^2 + \ldots + \dim W_{k_2}^2), \ldots, (\dim W_1^n + \ldots + \dim W_{k_n}^n)\}$. We prove the result for one particular $W_i = \{W_1^i + \ldots + W_{k_1}^i\}$ and since $W_i$ is arbitrarily chosen the result is true for every i, i = 1, 2, …, n. Let $B_i' = \{\alpha_1^i, \ldots, \alpha_{r_i}^i\}$ so that the first few $\alpha^i$'s form the basis of $B_i'$, the next few $B_2'$ and so on.

Then $T_i\alpha_j^t = t_j^i\alpha_j^t$, j = 1, 2, …, $r_i$ where $(t_1^i, \ldots, t_{r_i}^i) = \{C_1^i, \ldots, C_1^i, \ldots, C_{k_i}^i, C_{k_i}^i, \ldots, C_{k_i}^i)$ where $C_j^i$ is repeated $\dim W_j^i$ times, j = 1, …, $r_i$. Now $W_i$ is invariant under $T_i$ since for each $\alpha^i$ in $W_i$, we have

$$\alpha^i = x_1^i\alpha_1^i + \ldots + x_{r_i}^i\alpha_{r_i}^i$$

$$T_i\alpha^i = t_1^i x_1^i\alpha_1^i + \ldots + t_{r_i}^i x_{r_i}^i\alpha_{r_i}^i.$$

Choose any other vector $\alpha_{r_i+1}^i, \ldots, \alpha_{n_i}^i$ in $V_i$ such that $B_i = \{\alpha_1^i, \ldots, \alpha_{n_i}^i\}$ is a basis for $V_i$. The matrix of $T_i$ relative to $B_i$ has the block form mentioned earlier and the matrix of the restriction operator relative to the basis $B_i'$ is

$$B^i = \begin{bmatrix} t_1^i & 0 & \dots & 0 \\ 0 & t_2^i & \dots & 0 \\ \vdots & \vdots & & \vdots \\ 0 & 0 & \dots & t_{r_i}^i \end{bmatrix}.$$

The characteristic polynomial of $B^i$ i.e., of $T_{iw_i}$ is $g_i = g_i(x - C_1^i)^{e_1^i} \dots (x - C_{k_i}^i)^{e_{k_1}^i}$ where $e_j^i = \dim W_j^i$; j = 1, 2, …, $k_i$. Further more $g_i$ divides $f_i$, the characteristic polynomial for $T_i$. Therefore the multiplicity of $C_j^i$ as a root of $f_i$ is at least $\dim W_j^i$. Thus $T_i$ is diagonalizable if and only if $r_i = n_i$ i.e., if and only if $e_1^i + \dots + e_{k_i}^i = n_i$. Since what we proved for $T_i$ is true for $T = T_1 \cup \dots \cup T_n$. Hence true for every $B^1 \cup \dots \cup B^n$.

We now proceed on to define T-n conductor of $\alpha$ into $W = W_1 \cup \dots \cup W_n \subseteq V_1 \cup \dots \cup V_n$.

**DEFINITION 1.2.29:** *Let $W = W_1 \cup \dots \cup W_n$ be a n-invariant n-subspace for $T = T_1 \cup \dots \cup T_n$ and let $\alpha = \alpha_1 \cup \dots \cup \alpha_n$ be a n-vector in $V = V_1 \cup \dots \cup V_n$. The T-n conductor of $\alpha$ into W is the set $S_T(\alpha; W) = S_{T_1}(\alpha_1; W_1) \cup \dots \cup S_{T_n}(\alpha_n; W_n)$ which consists of all n-polynomials $g = g_1 \cup \dots \cup g_n$ over the n-field $F = F_1 \cup F_2 \cup \dots \cup F_n$ such that $g(T)\alpha$ is in W; that is $g_1(T_1)\alpha_1 \cup g_2(T_2)\alpha_2 \cup \dots \cup g_n(T_n)\alpha_n \in W_1 \cup \dots \cup W_n$.*

*Since the n-operator T will be fixed through out the discussions we shall usually drop the subscript T and write $S(\alpha; W) = S(\alpha_1; W_1) \cup \dots \cup S(\alpha_n; W_n)$. The authors usually call the collection of n-polynomials the n-stuffer. We as in case of vector spaces prefer to call as n-conductor i.e., the n-operator $g(T) = g_1(T_1) \cup \dots \cup g_n(T_n)$; slowly leads to the n-vector $\alpha_1 \cup \dots \cup \alpha_n$ into $W = W_1 \cup \dots \cup W_n$. In the special case when $W = \{0\} \cup \dots \cup \{0\}$, the n-conductor is called the T-annihilator of $\alpha_1 \cup \dots \cup \alpha_n$.*

We prove the following simple lemma.

**LEMMA 1.2.7:** *If $W = W_1 \cup \ldots \cup W_n$ is an n-invariant subspace for $T = T_1 \cup \ldots \cup T_n$, then W is n-invariant under every n-polynomial in $T = T_1 \cup \ldots \cup T_n$. Thus for each $\alpha = \alpha_1 \cup \ldots \cup \alpha_n$ in $V = V_1 \cup \ldots \cup V_n$ the n-conductor $S(\alpha; W) = S(\alpha_1;W_1) \cup \ldots \cup S(\alpha_n; W_n)$ is an n-ideal in the n-polynomial algebra $F[x] = F_1[x] \cup \ldots \cup F_n[x]$.*

*Proof:* Given $W = W_1 \cup \ldots \cup W_n \subseteq V = V_1 \cup \ldots \cup V_n$ a n-vector space over the n-field $F = F_1 \cup \ldots \cup F_n$. If $\beta = \beta_1 \cup \ldots \cup \beta_n$ is in $W = W_1 \cup \ldots \cup W_n$, then $T\beta = T_1\beta_1 \cup \ldots \cup T_n\beta_n$ is in $W = W_1 \cup \ldots \cup W_n$. Thus $T(T\beta) = T^2\beta = T_1^2 \beta_1 \cup \ldots \cup T_n^2 \beta_n$ is in W. By induction $T^k\beta = T_1^{k_1} \beta_1 \cup \ldots \cup T_n^{k_n} \beta_n$ is in W for each k. Take linear combinations to see that $f(T)\beta = f_1(T_1)\beta_1 \cup \ldots \cup f_n(T_n)\beta_n$ is in W for every polynomial $f = f_1 \cup \ldots \cup f_n$.

The definition of $S(\alpha; W) = S(\alpha_1;W_1) \cup \ldots \cup S(\alpha_n; W_n)$ is meaningful if $W = W_1 \cup \ldots \cup W_n$ is any n-subset of V. If W is a n-subspace then $S(\alpha; W)$ is a n-subspace of $F[x] = F_1[x] \cup \ldots \cup F_n[x]$ because $(cf + g)T = cf(T) + g(T)$ i.e., $(c_1f_1 + g_1)T_1 \cup \ldots \cup (c_nf_n + g_n)T_1 = c_1f_1(T_1) + g_1(T_1) \cup \ldots \cup c_nf_n(T_n) + g_n(T_n)$. If $W = W_1 \cup \ldots \cup W_n$ is also n-invariant under $T = T_1 \cup \ldots \cup T_n$ and let $g = g_1 \cup \ldots \cup g_n$ be a n-polynomial in $S(\alpha; W) = S(\alpha_1;W_1) \cup \ldots \cup S(\alpha_n; W_n)$ i.e., let $g(T)\alpha = g_1(T_1)\alpha_1 \cup g_2(T_2)\alpha_2 \cup \ldots \cup g_n(T_n)\alpha_n$ be in $W = W_1 \cup \ldots \cup W_n$ is any n-polynomial then $f(T)g(T)\alpha$ is in $W = W_1 \cup \ldots \cup W_n$ i.e.,

$$f(T)[g(T)\alpha] = f_1(T_1) [g_1(T_1)\alpha_1] \cup \ldots \cup f_n(T_n)[g_n(T_n)\alpha_n]$$

will be in $W = W_1 \cup \ldots \cup W_n$.

Since $(fg)T = f(T)g(T)$ we have

$$(f_1g_1)T_1 \cup \ldots \cup (f_ng_n)T_n = f_1(T_1)\, g_1(T_1) \cup \ldots \cup f_n(T_n)g_n(T_n)\alpha_n$$

be $(fg) \in S(\alpha; W)$ i.e., $(f_ig_i) \in S(\alpha_i; W_i)$; i = 1, 2, …, n. Hence the claim.

The unique n-monic generator of the n-ideal $S(\alpha; W)$ is also called the T-n-conductor of $\alpha$ in W(the T-n annihilator in case $W = \{0\} \cup \{0\} \cup \ldots \cup \{0\}$). The T-n-conductor of $\alpha$ into W is the n-monic polynomial g of least degree such that $g(T)\alpha = g_1(T_1)\alpha_1 \cup \ldots \cup g_n(T_n)\alpha_n$ is in $W = W_1 \cup \ldots \cup W_n$.

A n-polynomial $f = f_1 \cup \ldots \cup f_n$ is in $S(\alpha; W) = S(\alpha_1;W_1) \cup \ldots \cup S(\alpha_n; W_n)$ if and only if g n-divides f. Note the n-conductor S(α; W) always contains the n-minimal polynomial for T, hence every T-n-conductor n-divides the n-minimal polynomial for T.

**LEMMA 1.2.8:** *Let $V = V_1 \cup \ldots \cup V_n$ be a $(n_1, n_2, \ldots, n_n)$ dimensional n-vector space over the n-field $F = F_1 \cup \ldots \cup F_n$. Let $T = T_1 \cup \ldots \cup T_n$ be a n-linear operator on V such that the n-minimal polynomial for T is a product of n-linear factors $p = p_1 \cup \ldots \cup p_n = (x - c_1^1)^{r_1^1} (x - c_2^1)^{r_2^1} \ldots (x - c_{k_1}^1)^{r_{k_1}^1} \cup (x - c_1^2)^{r_1^2} (x - c_2^2)^{r_2^2} \ldots (x - c_{k_2}^2)^{r_{k_2}^2} \cup \ldots \cup (x - c_1^n)^{r_1^n} (x - c_2^n)^{r_2^n} \ldots (x - c_{k_n}^n)^{r_{k_n}^n}$; $c_{t_i}^i \in F_i$, $k_1 \le t_i \le k_n$, $i = 1, 2, \ldots, n$.*

*Let $W_1 \cup \ldots \cup W_n$ be a proper n-subspace of $V (V \neq W)$ which is n-invariant under T. There exist a n-vector $\alpha_1 \cup \ldots \cup \alpha_n$ in V such that*

*(1) α is not in W*

*(2) $(T - cI)\alpha = (T_1 - c_1 I_1)\alpha_1 \cup \ldots \cup (T_n - c_n I_n)\alpha_n$ is in W for some n-characteristic value of the n-operator T.*

*Proof:* (1) and (2) express that T-n conductor of $\alpha = \alpha_1 \cup \ldots \cup \alpha_n$ into $W_1 \cup \ldots \cup W_n$ is a n-linear polynomial. Suppose $\beta = \beta_1 \cup \ldots \cup \beta_n$ is any n-vector in V which is not in W. Let $g = g_1 \cup \ldots \cup g_n$ be the T-n conductor of β in W. Then g n-divides $p = p_1 \cup \ldots \cup p_n$ the n-minimal polynomial for T. Since β is not in W, the n-polynomial g is not constant. Therefore $g = g_1 \cup \ldots \cup g_n = (x - c_1^1)^{e_1^1} \ldots (x - c_{k_1}^1)^{e_{k_1}^1} \cup (x - c_1^2)^{e_1^2} \ldots (x - c_{k_2}^2)^{e_{k_2}^2} \cup \ldots \cup (x - c_1^n)^{e_1^n} \ldots (x - c_{k_n}^n)^{e_{k_n}^n}$ where at least one of the n-tuple of integers $e_i^1 \cup e_i^2 \cup \ldots \cup e_i^n$ is positive. Choose $j_t$ so that $e_{j_t}^t > 0$, then $(x - c_j) = (x - c_{j_1}^1) \cup (x - c_{j_2}^2) \cup \ldots \cup (x - c_{j_n}^n)$ n-divides g. $g = (x - c_j)h$ i.e., $g = g_1 \cup \ldots \cup g_n = (x - c_{j_1}^1)h_1 \cup \ldots \cup (x - c_{j_n}^n)h_n$. But by the definition of g the n-vector $\alpha = \alpha_1 \cup \ldots \cup \alpha_n = h_1(T_1)\beta_1 \cup \ldots \cup h_n(T_n)\beta_n = h(T)\beta$ cannot be in W. But $(T - c_j I)\alpha = (T - c_j I)h(T)\beta = g(T)\beta$ is in W i.e.,

$$(T_1 - c^1_{j_1}\, {}_jI_1)\alpha_1 \cup \ldots \cup (T_n - c^n_{j_n}\, {}_jI_n)\alpha_n$$
$$= (T_1 - c^1_{j_1} I_1)h_1(T_1)\beta_1 \cup \ldots \cup (T_n - c^n_{j_n} I_n)h_n(T_n)\beta_n$$
$$= g_1(T_1)\beta_1 \cup g_2(T_2)\beta_2 \cup \ldots \cup g_n(T_n)\beta_n$$
with $g_i(T_i)\beta_i \in W_i$ for i = 1, 2, …, n.

Now we obtain the condition for T to be n-triangulable.

**THEOREM 1.2.47:** *Let $V = V_1 \cup \ldots \cup V_n$ be a finite $n_1, n_2, \ldots, n_n$ dimensional n-vector space over the n-field $F = F_1 \cup \ldots \cup F_n$ and let $T = T_1 \cup \ldots \cup T_n$ be a n-linear operator on V. Then T is n-triangulable if and only if the n-minimal polynomial for T is a n-product of n-linear polynomials over $F = F_1 \cup \ldots \cup F_n$.*

*Proof:* Suppose the n-minimal polynomial $p = p_1 \cup \ldots \cup p_n$, n-factors as $p = (x - c^1_1)^{r^1_1} \ldots (x - c^1_{k_1})^{r^1_{k_1}} \cup (x - c^2_1)^{r^1_2} \cup (x - c^2_1)^{r^2_1} \ldots (x - c^2_{k_2})^{r^2_{k_2}} \cup \ldots \cup (x - c^n_1)^{r^n_1} \ldots (x - c^n_{k_n})^{r^n_{k_n}}$. By the repeated application of the lemma just proved we arrive at a n-ordered n-basis. $B = \{\alpha^1_1 \ldots \alpha^1_{n_1}\} \cup \{\alpha^2_1 \ldots \alpha^2_{n_2}\} \cup \ldots \cup \{\alpha^n_1 \ldots \alpha^n_{n_n}\} = B_1 \cup B_2 \cup \ldots \cup B_n$ in which the n-matrix representing $T = T_1 \cup \ldots \cup T_n$ is n-upper triangular.

$$[T_1]_{B_1} \cup \ldots \cup [T_n]_{B_n}$$

$$= \begin{bmatrix} a^1_{11} & a^1_{12} & \ldots & a^1_{1n_1} \\ 0 & a^1_{22} & & a^1_{2n_1} \\ \vdots & \vdots & & \vdots \\ 0 & 0 & \ldots & a^1_{n_1n_1} \end{bmatrix} \cup \begin{bmatrix} a^2_{11} & a^2_{12} & \ldots & a^2_{1n_2} \\ 0 & a^2_{22} & & a^2_{2n_2} \\ \vdots & \vdots & & \vdots \\ 0 & 0 & \ldots & a^2_{n_2n_2} \end{bmatrix} \cup \ldots \cup$$

$$\begin{bmatrix} a^n_{11} & a^n_{12} & \ldots & a^n_{1n_n} \\ 0 & a^n_{22} & & a^n_{2n_n} \\ \vdots & \vdots & & \vdots \\ 0 & 0 & \ldots & a^n_{n_nn_n} \end{bmatrix}.$$

Merely $[T]_B$ = the n-triangular matrix of $(n_1 \times n_1, \ldots, n_n \times n_n)$ order shows that

$$T\alpha_j = T_1 \alpha^1_{j_1} \cup \ldots \cup T_n \alpha^n_{j_n}$$
$$= a^1_{1j_1}\alpha^1_1 + \ldots + a^1_{j_1 j_1}\alpha^1_{j_1} \cup a^2_{1j_2}\alpha^2_1 + \ldots + a^2_{j_2 j_2}\alpha^2_{j_2} \cup \ldots \cup a^n_{1j_n}\alpha^n_1 + \ldots + a^n_{j_n j_n}\alpha^n_{j_n} \quad \text{(a)}$$

$1 \leq j_i \leq n_i$; i = 1, 2, …, n, that is $T\alpha_j$ is in the n-subspace spanned by $\{\alpha^1_1 \ldots \alpha^1_{j_1}\} \cup \ldots \cup \{\alpha^n_1 \ldots \alpha^n_{j_n}\}$. To find $\{\alpha^1_1 \ldots \alpha^1_{j_1}\} \cup \ldots \cup \{\alpha^n_1 \ldots \alpha^n_{j_n}\}$; we start by applying the lemma to the n-subspace $W = W_1 \cup \ldots \cup W_n = \{0\} \cup \ldots \cup \{0\}$ to obtain the n-vector $\alpha^1_1 \cup \ldots \cup \alpha^n_1$. Then apply the lemma to $W^1_1 \cup W^2_1 \cup \ldots \cup W^n_1$, the n-space spanned by $\alpha^1 = \alpha^1_1 \cup \ldots \cup \alpha^n_1$, and we obtain $\alpha^2 = \alpha^1_2 \cup \ldots \cup \alpha^n_2$. Next apply lemma to $W_2 = W^1_2 \cup \ldots \cup W^n_2$, the n-space spanned by $\alpha^1_1 \cup \ldots \cup \alpha^n_1$ and $\alpha^1_2 \cup \ldots \cup \alpha^n_2$. Continue in that way. After $\alpha^1, \alpha^2, \ldots, \alpha^i$ we have found it is the triangular type relation given by equation (a) for $j_i = 1, 2, \ldots, n_i$, i = 1, 2, …, n which proves that the n-subspace spanned by $\alpha^1, \alpha^2, \ldots, \alpha^i$ is n-invariant under T.

If T is n-triangulable it is evident that the n-characteristic polynomial for T has the form $f = f_1 \cup \ldots \cup f_n = (x - c^1_1)^{d^1_1} \ldots (x - c^1_{k_1})^{d^1_{k_1}} \cup \ldots \cup (x - c^n_1)^{d^n_1} \ldots (x - c^n_{k_n})^{d^n_{k_n}}$. The n-diagonal entries $(a^1_{11} \ldots a^1_{1n_1}) \cup (a^2_{11} \ldots a^2_{1n_2}) \cup \ldots \cup (a^n_{11} \ldots a^n_{1n_n})$ are the n-characteristic values with $c^t_j$ repeated $d^t_{jt}$ times. But if f can be so n-factored, so can the n-minimal polynomial p because p n-divides f.

We leave the following corollary to be proved by the reader.

**COROLLARY 1.2.12:** *If $F = F_1 \cup \ldots \cup F_n$ is an n-algebraically closed n-field. Every $(n_1 \times n_1, \ldots, n_n \times n_n)$ n-matrix over F is similar over the n-field F to be a n-triangular matrix.*

Now we accept with some deviations in the definition of the n-algebraically closed field when ever $F_i$ is the complex field C.

**THEOREM 1.2.48:** *Let $V = V_1 \cup \ldots \cup V_n$ be a $(n_1, \ldots, n_n)$ dimensional n-vector space over the n-field $F = F_1 \cup \ldots \cup F_n$ and let $T = T_1 \cup \ldots \cup T_n$ be a n-linear operator on $V = V_1 \cup \ldots \cup V_n$. Then T is n-diagonalizable if and only if the n-minimal polynomial for T has the form*

$$\begin{aligned} p &= p_1 \cup \ldots \cup p_n \\ &= (x - c_1^1) \ldots (x - c_{k_1}^1) \cup \ldots \cup (x - c_1^n) \ldots (x - c_{k_n}^n) \end{aligned}$$

*where $\{c_1^1, \ldots, c_{k_1}^1\} \cup \ldots \cup \{c_1^n, \ldots, c_{k_n}^n\}$ are n-distinct elements of $F = F_1 \cup \ldots \cup F_n$.*

*Proof:* We know if $T = T_1 \cup \ldots \cup T_n$ is n-diagonalizable its n-minimal polynomial is a n-product of n-distinct linear factors. Hence one way of the proof is clear.

To prove the converse let $W = W_1 \cup \ldots \cup W_n$ be a subspace spanned by all the n-characteristic n-vectors of T and suppose $W \neq V$. Then we know by the properties of n-linear operator that their exists a n-vector $\alpha = \alpha_1 \cup \ldots \cup \alpha_n$ in V and not in W and the a n-characteristic value $c_j = c_{j_1}^1 \cup \ldots \cup c_{j_n}^n$ of T such that the n-vector

$$\begin{aligned} \beta &= (T - c_j I)\alpha \\ &= (T_1 - c_{j_1}^1 I_1)\alpha_1 \cup \ldots \cup (T_n - c_{j_n}^n I_n)\alpha_n \\ &= \beta_1 \cup \ldots \cup \beta_n \end{aligned}$$

lies in $W = W_1 \cup \ldots \cup W_n$ where each $\beta_i \in W_i$; i = 1, 2, …, n.

Since $\beta = \beta_1 \cup \ldots \cup \beta_n$ is in W; $\beta_i = \beta_i^1 + \ldots + \beta_i^{k_i}$; for each i = 1, 2, …, n with $T_i \beta_i^t = c_i^t \beta_i^t$; t = 1, 2, …, $k_i$, this is true for every i = 1, 2, …, n and hence the n-vector

$$\begin{aligned} h(T)\beta = \{&h^1(c_1^1)\beta_1^1 + \ldots + h^1(c_{k_1}^1)\beta_{k_1}^1\} \\ &\cup \ldots \cup \{h^n(c_1^n)\beta_1^n + \ldots + h^n(c_{k_n}^n)\beta_{k_n}^n\} \end{aligned}$$

for every n-polynomial h. Now

$$\begin{aligned} p &= (x - c_j)q \\ &= p_1 \cup \ldots \cup p_n \\ &= (x - c^1_{j_1})q_1 \cup \ldots \cup (x - c^n_{j_n})q_n \end{aligned}$$

for some n-polynomial $q = q_1 \cup \ldots \cup q_n$. Also

$$q - q(c_j) = (x - c_j)h,$$

i.e., $q_1 - q_1(c^1_{j_1}) = (x - c^1_{j_1})h^1_{j_1}$, …, $q_n - q_n(c^n_j) = (x - c^n_{j_n})h^n_{j_n}$.

We have

$$q(T)\alpha - q(c_j)\alpha = h(T)(T - c_jI)\alpha = h(T)\beta$$

But $h(T)\beta$ is in $W = W_1 \cup \ldots \cup W_n$ and since

$$\begin{aligned} 0 = p(T)\alpha &= (T - c_jI)\, q(T)\alpha \\ &= p_1(T_1)\alpha_1 \cup \ldots \cup p_n(T_n)\alpha_n \\ &= (T_1 - c^1_{j_1}I_1)q_1(T_1)\alpha_1 \cup \ldots \cup (T_n - c^n_{j_n})q_n(T_n)\,\alpha_n \end{aligned}$$

and the n-vector $q(T)\alpha$ is in W i.e., $q_1(T_1)\alpha_1 \cup \ldots \cup q_n(T_n)\alpha_n$ is in W. Therefore

$$q(c_j)\,\alpha = q_1(c^1_{j_1})\alpha_1 \cup \ldots \cup q_n(c^n_{j_n})\alpha_n$$

is in W. Since $\alpha = \alpha_1 \cup \ldots \cup \alpha_n$ is not in W, we have $q(c_j) = q_1(c^1_{j_1}) \cup \ldots \cup q_n(c^n_{j_n}) = 0 \cup \ldots \cup 0$.

This contradicts the fact that p has distinct roots. Hence the claim.

We can now describe this more in terms of how the values are determined and its relation to Cayley Hamilton Theorem for n-vector spaces of type II. Suppose $T = T_1 \cup \ldots \cup T_n$ is a n-linear operator on a n-vector space of type II which is represented by the n-matrix $A = A_1 \cup \ldots \cup A_n$ in some ordered n-basis for which we wish to find out whether T is n-diagonalizable. We compute the n-characteristic polynomial $f = f_1 \cup \ldots \cup f_n$. If we can n-factor $f = f_1 \cup \ldots \cup f_n$ as

$$(x - c^1_1)^{d^1_1} \ldots (x - c^1_{k_1})^{d^1_{k_1}} \cup \ldots \cup (x - c^n_1)^{d^n_1} \ldots (x - c^n_{k_n})^{d^n_{k_n}}$$

we have two different methods for finding whether or not T is n-diagonalizable. One method is to see whether for each i(t) (i(t) means i is independent on t) we can find a $d^t_i$ (t = 1, 2, …, n); $1 \le i \le k_t$ independent characteristic vectors associated with the characteristic value $c^t_j$. The other method is to check whether or not

$$(T - c_1 I) \cup \ldots \cup (T - c_k I) =$$
$$(T_1 - c_1^1 I_1) \ldots (T_1 - c_{k_1}^1 I_1) \cup \ldots \cup (T_n - c_1^n I_n) \ldots (T_n - c_{k_n}^n I_n)$$

is the n-zero operator.

Here we have some problems about speaking of n-algebraically closed n-field in a universal sense; we can only speak of the n-algebraically closed n-field relative to a n-polynomial over the same n-field. Now we proceed on to define or introduce the notion of simultaneous n-diagonalization and simultaneous n-triangulation. Throughout this section $V = V_1 \cup \ldots \cup V_n$ will denote a $(n_1, n_2, \ldots, n_n)$ finite n-vector space over the n-field $F = F_1 \cup \ldots \cup F_n$ and let $\Im = \Im_1 \cup \ldots \cup \Im_n$ be the n-family of n-linear operator on V. We now discuss when one can simultaneously n-triangularize or n-diagonalize the n-operators in F. i.e., to find one n-basis $B = B_1 \cup \ldots \cup B_n$ such that all n-matrices $[T]_B = [T_1]_{B_1} \cup \ldots \cup [T_n]_{B_n}$, T in $\Im$ are n-triangular (or n-diagonal). In case of n-diagonalization it is necessary that f be in the commuting family of n-operators UT = TU, i.e., $U_1T_1 \cup \ldots \cup U_nT_n = T_1U_1 \cup \ldots \cup T_nU_n$ for all T, U in $\Im$. That follows from the simple fact that all n-diagonal n-matrices commute. Of course it is also necessary that each n-operator in $\Im$ be an n-diagonalizable operator. In order to simultaneously n-triangulate each n-operator in $\Im$ we see each n-operator must be n-triangulable. It is not necessary that $\Im$ be a n-commuting family, however that condition is sufficient for simultaneous n-triangulation (if each $T = T_1 \cup \ldots \cup T_n$ can be individually n-triangulated).

We recall a subspace $W = W_1 \cup \ldots \cup W_n$ is n-invariant under $\Im$ if W is n-invariant under each operator in $\Im$ i.e., each $W_i$ in $W = W_1 \cup \ldots \cup W_n$ is invariant under the operator $T_i$ in $T = T_1 \cup \ldots \cup T_i \cup \ldots \cup T_n$; true for every i = 1, 2, …, n.

**LEMMA 1.2.9:** *Let $\Im = \Im_1 \cup \ldots \cup \Im_n$ be a n-commuting family of n-triangulable n-linear operators on the n-vector space $V = V_1 \cup \ldots \cup V_n$. Let $W = W_1 \cup \ldots \cup W_n$ be a proper n-subspace of $V = V_1 \cup \ldots \cup V_n$ which is n-invariant under $\Im$.*
*There exists a n-vector $\alpha = \alpha_1 \cup \ldots \cup \alpha_n$ in V such that*

*a. $\alpha = \alpha_1 \cup \ldots \cup \alpha_n$ is not in $W = W_1 \cup \ldots \cup W_n$*

b. *for each $T = T_1 \cup \ldots \cup T_n$ in $\mathfrak{I} = \mathfrak{I}_1 \cup \ldots \cup \mathfrak{I}_n$ the n-vector $T\alpha = T_1\alpha_1 \cup T_2\alpha_2 \cup \ldots \cup T_n\alpha_n$ is in the n-subspace spanned by $\alpha = \alpha_1 \cup \ldots \cup \alpha_n$ and $W = W_1 \cup \ldots \cup W_n$.*

*Proof:* Without loss in generality let us assume that $\mathfrak{I}$ contains only a finite number of n-operators because of this observations. Let $\{T^1, \ldots, T^r\}$ be a maximal n-linearly independent n-subset of $\mathfrak{I}$; i.e., a n-basis for the n-subspace spanned by $\mathfrak{I}$. If $\alpha = \alpha_1 \cup \ldots \cup \alpha_n$ is a n-vector such that (b) holds for each $T^i = T_1^i \cup \ldots \cup T_n^i, 1 \leq i \leq r$ then (b) will hold for every n-operator which is a n-linear combination of $T^1, \ldots, T^r$.

By earlier results we can find a n-vector $\beta^1 = \beta_1^1 \cup \ldots \cup \beta_n^1$ in V (not in W) and a n-scalar $c^1 = c_1^1 \cup \ldots \cup c_n^1$ such that $(T^1 - c^1 I)\beta^1$ is in W; i.e., $(T_1^1 - c_1^1 I_1)\beta_1^1 \cup \ldots \cup (T_n^1 - c_n^1)\, I_n \beta_n^1 \in W_1 \cup \ldots \cup W_n$ i.e., each $(T_r^1 - c_r^1)\beta_r^1 \in W_r$; r = 1, 2, …, n. Let $V^1 = V_1^1 \cup \ldots \cup V_n^1$ be the collection of all vector $\beta$ in V such that $(T^1 - c^1 I)\beta$ is in W. Then $V^1$ is a n-subspace of V which is properly bigger than W. Further more $V^1$ is n-invariant under $\mathfrak{I}$ for this reason. If $T = T_1 \cup \ldots \cup T_n$ commutes with $T^1 = T_1^1 \cup \ldots \cup T_n^1$ then $(T^1 - c^1 I)T\beta = T(T^1 - c^1 I)\beta$. If $\beta$ is in $V^1$ then $(T_1 - c^1 I)\beta$ is in W i.e., $T\beta$ is in $V^1$ for all $\beta$ in $V^1$ and for all T in $\mathfrak{I}$. Now W is a proper n-subspace of $V^1$. Let $U^2 = U_1^2 \cup \ldots \cup U_n^2$ be the n-linear operator on $V^1$ obtained by restricting $T^2 = T_1^2 \cup \ldots \cup T_n^2$ to the n-subspace $V^1$. The n-minimal polynomial for $U^2$ divides the n-minimal polynomial for $T^2$. Therefore we may apply the earlier results to the n-operator and the n-invariant n-subspace W. We obtain a n-vector $\beta_2$ in $V^1$ (not in W) and a n-scalar $c^2$ such that $(T^2 - c^2 I)\beta_2$ is in W. Note that

a. $\beta_2$ is not in W
b. $(T^1 - c^1 I)\beta_2$ is in W
c. $(T^2 - c^2 I)\beta_2$ is in W.

Let $V^2$ be the set of all n-vectors $\beta$ in $V^1$ such that $(T^2 - c^2 I)\beta$ is in W. Then $V^2$ is n-invariant under $\mathfrak{I}$. Applying earlier results to $U^3 = U_1^3 \cup \ldots \cup U_n^3$ the restriction of $T^3$ to $V^2$. If we continue in this way we shall reach a n-vector $\alpha = \beta^r$ (not in W) $\beta^r = \beta_1^r \cup \ldots \cup \beta_n^r$ such that $(T^j - c^j I)\alpha$ is in W; j = 1, 2, …, r.

The following theorem is left as an exercise for the reader.

**THEOREM 1.2.49:** *Let $V = V_1 \cup \ldots \cup V_n$ be a finite $(n_1, \ldots, n_n)$ dimension n-vector space over the n-field F. Let $\mathfrak{I}$ be a commuting family of n-triangulable, n-linear operators on $V = V_1 \cup \ldots \cup V_n$. There exists an n-ordered basis for $V = V_1 \cup \ldots \cup V_n$ such that every n-operator in $\mathfrak{I}$ is represented by a triangular n-matrix in that n-basis.*

In view of this theorem the following corollary is obvious.

**COROLLARY 1.2.13:** *Let $\mathfrak{I}$ be a commuting family of $(n_1 \times n_1, n_2 \times n_2, \ldots, n_n \times n_n)$ n-matrices (square), i.e., (n-mixed square matrices) over a special algebraically closed n-field $F = F_1 \cup \ldots \cup F_n$. There exists a non-singular $(n_1 \times n_1, n_2 \times n_2, \ldots, n_n \times n_n)$ n-matrix $P = P_1 \cup P_2 \cup \ldots \cup P_n$ with entries in $F = F_1 \cup \ldots \cup F_n$ such that $P^{-1}A\,P = P_1^{-1}A_1P_1 \cup \ldots \cup P_n^{-1}A_nP_n$ in n-upper triangular for every n-matrix $A = A_1 \cup \ldots \cup A_n$ in F.*

Next we prove the following theorem.

**THEOREM 1.2.50:** *Let $\mathfrak{I} = \mathfrak{I}_1 \cup \ldots \cup \mathfrak{I}_n$ be a commuting family of n-diagonalizable n-linear operators on a finite dimensional n-vector space V. There exists an n-ordered basis for V such that every n-operator in $\mathfrak{I}$ is represented in that n-basis by a n-diagonal matrix.*

*Proof:* We give the proof by induction on the $(n_1, n_2, \ldots, n_n)$ dimension of the n-vector space $V = V_1 \cup \ldots \cup V_n$. If n-dim V = (1, 1, …, 1) we have nothing to prove. Assume the theorem for n-dim V < $(n_1, \ldots, n_n)$ where V is given as $(n_1, \ldots, n_n)$ dimensional space. Choose any $T = T_1 \cup \ldots \cup T_n$ in $\mathfrak{I}$ which is not a scalar multiple of n-identity.

Let $\{c_1^1, \ldots, c_{k_1}^1\} \cup \{c_1^2, \ldots, c_{k_2}^2\} \cup \ldots \cup \{c_1^n, \ldots, c_{k_n}^n\}$ be distinct n-character values of T and for each (i(t)) let $W_i^t$ be the null space of the n-null space $W^i = W_i^1 \cup \ldots \cup W_i^n$; t = 1, 2, …, n of $(T - c_iI_i) = (T_1 - C_1I_1) \cup \ldots \cup (T_n - c_nI_n)$. Fix i, then $W_i^t$ is invariant under every operator that commutes with $T_t$. Let $\{\mathfrak{I}_i^t\}$ be the family of linear operators on $W_i^t$ obtained by restricting the operators in $\mathfrak{I}^t$ where $\mathfrak{I} = \mathfrak{I}^1 \cup \ldots \cup \mathfrak{I}^n$; $1 \leq t \leq n$ to the subspace $W_i^t$. Each operator in $\mathfrak{I}_i^t$ is diagonalizable because its minimal polynomial divides the minimal polynomial for the corresponding operator in $\mathfrak{I}^t$ since $\dim W_i^t < \dim V_i$, the operators in $\mathfrak{I}_i^t$ can be simultaneously diagonalized.

In other words $W_i^t$ has a basis $B_i^t$ which consists of vectors which are simultaneously characteristic vectors for every operator $\mathfrak{I}_i^t$. Since this is true for every t, t = 1, 2, …, n we see $T = T_1 \cup \ldots \cup T_n$ is n-diagonalizable and

$$B = \{B_1^1, \ldots, B_{k_1}^1\} \cup \{B_1^2, \ldots, B_{k_2}^2\} \cup \ldots \cup \{B_1^n, \ldots, B_{k_n}^n\}$$

is a n-basis for the n-vector space $V = V_1 \cup \ldots \cup V_n$.

Let $V = V_1 \cup \ldots \cup V_n$ be a n-vector space over the n-field, $F = F_1 \cup \ldots \cup F_n$. Let $W = \{W_1^1, \ldots W_{k_1}^1\} \cup \ldots \cup \{W_1^n, \ldots, W_{k_n}^n\}$ be n-subspaces of the n-vector space V. We say that $W_1^t, \ldots, W_{k_t}^t$ are independent if $\alpha_1^1 + \ldots + \alpha_{k_t}^t = 0, \alpha_i^t$ in $W_i^t$ implies that each $\alpha_i^t = 0$ i.e., if $(\alpha_1^1 + \ldots + \alpha_{k_1}^1) \cup \ldots \cup (\alpha_1^n + \ldots + \alpha_{k_n}^n) = 0 \cup \ldots \cup 0$ implies each $\alpha_i^t = 0, 1 \leq i \leq k_t$; t = 1, 2, …, n; then

$$\{W_1^1, \ldots, W_{k_1}^1\} \cup \ldots \cup \{W_1^n, \ldots, W_{k_n}^n\}$$

is said to be n-independent.

The following lemma would be useful for developing more properties about the n-independent n-subspaces.

**LEMMA 1.2.10:** *Let $V = V_1 \cup \ldots \cup V_n$ be a finite dimensional ($n_1$, $n_2$, …, $n_n$) n-vector space of type II. Let $\{W_1^1, \ldots, W_{k_1}^1\} \cup \ldots \cup \{W_1^n, \ldots, W_{k_n}^n\}$ be n-subspaces of V and let $W = W_1 \cup \ldots \cup W_n = W_1^1 + \ldots + W_{k_1}^1 \cup W_1^2 + \ldots + W_{k_2}^2 \cup \ldots \cup W_1^n + \ldots + W_{k_n}^n$. Then the following are equivalent*

a. *$\{W_1^1, \ldots, W_{k_1}^1\} \cup \ldots \cup \{W_1^n, \ldots, W_{k_n}^n\}$ are n-independent i.e., $\{W_1^t, \ldots, W_{k_t}^t\}$ are independent for t = 1, 2, …, n.*
b. *For each $j_t$, $2 \leq j_t \leq k_t$, t = 1, 2, …, n we have $W_{j_t}^t \cap (W_1^t + \ldots + W_{j_{t-1}}^t) = \{0\}$ for t = 1, 2, …, n.*
c. *If $B_i^t$ is an ordered basis for $W_i^t$, $1 \leq i \leq k_t$; t = 1, 2, …, n then the n-sequence $\{B_1^1, \ldots, B_{k_1}^1\} \cup \ldots \cup \{B_1^n, \ldots, B_{k_n}^n\}$ is an n-ordered basis for n-subspace $W = W_1 \cup \ldots \cup W_n = W_1^1 + \ldots + W_{k_2}^1 \cup \ldots \cup W_1^n + \ldots + W_{k_n}^n$.*

*Proof:* Assume (a) Let $\alpha^t \in W_{j_t}^t \cap (W_1^t + \ldots + W_{j_{t-1}}^t)$, then there are vectors $\alpha_1^t, \ldots, \alpha_{j_{t-1}}^t$ with $\alpha_i^t \in W_i^t$ such that $\alpha_1^t + \ldots + \alpha_{j_{t-1}}^t$ + $\alpha^t$ = 0 + … + 0 = 0 and since $W_1^t, \ldots, W_{k_t}^t$ are independent it must be that $\alpha_1^t = \alpha_2^t = \ldots = \alpha_{j_{t-1}}^t = \alpha^t = 0$. This is true for each t; t = 1, 2, …, n. Now let us observe that (b) implies (a). Suppose 0 = $\alpha_1^t + \ldots + \alpha_{k_t}^t$; $\alpha_i^t \in W_i^t$; i = 1, 2, …, $k_t$. (We denote both the zero vector and zero scalar by 0). Let $j_t$ be the largest integer $i_t$ such that $\alpha_i^t \neq 0$. Then $0 = \alpha_i^t + \ldots + \alpha_{j_t}^t$; $\alpha_{j_t}^t \neq 0$ thus

$$\alpha_{j_t}^t = -\alpha_1^t - \ldots - \alpha_{j_{t-1}}^t$$

is a non zero vector in $W_{j_t}^t \cap (W_1^t + \ldots + W_{j_{t-1}}^t)$. Now that we know (a) and (b) are the same let us see why (a) is equivalent to (c). Assume (a). Let $B_i^t$ be a basis of $W_i^t; 1 \leq i \leq k_t$, and let $B^t = \{B_1^t, \ldots, B_{k_t}^t\}$ true for each t, t = 1, 2, …, n.

Any linear relation between the vector in $B^t$ will have the form $\beta_1^t + \ldots + \beta_{k_t}^t = 0$ where $\beta_i^t$ is some linear combination of vectors in $B_i^t$. Since $W_1^t, \ldots, W_{k_t}^t$ are independent each of $\beta_i^t$ is 0. Since each $B_i^t$ is an independent relation. The relation between the vectors in $B^t$ is trivial. This is true for every t, t = 1, 2, …, n; so in

$$B = B^1 \cup \ldots \cup B^n$$
$$= \{B_1^1, \ldots, B_{k_1}^1\} \cup \ldots \cup \{B_1^n, \ldots, B_{k_n}^n\}$$

every n-relation in n-vectors in B is the trivial n-relation.

It is left for the reader to prove (c) implies (a).

If any of the conditions of the above lemma hold, we say the n-sum W = $(W_1^1 + \ldots + W_{k_1}^1) \cup \ldots \cup (W_1^n + \ldots + W_{k_n}^n)$; n-direct or that W is the n-direct sum of $\{W_1^1, \ldots, W_{k_1}^1\}$ ∪ … ∪ $\{W_1^n, W_2^n, \ldots, W_{k_n}^n\}$ and we write

$$W = (W_1^1 \oplus \ldots \oplus W_{k_1}^1) \cup \ldots \cup (W_1^n \oplus \ldots \oplus W_{k_n}^n).$$

This n-direct sum will also be known as the n-independent sum or the n-interior direct sum of $\{W_1^1, \ldots, W_{k_1}^1\}$ ∪ … ∪ $\{W_1^n, \ldots, W_{k_n}^n\}$.

Let $V = V_1 \cup \ldots \cup V_n$ be a n-vector space over the n-field $F = F_1 \cup \ldots \cup F_n$. A n-projection of V is a n-linear operator $E = E_1 \cup \ldots \cup E_n$ on V such that $E^2 = E_1^2 \cup \ldots \cup E_n^2 = E_1 \cup \ldots \cup E_n = E$.

Since E is a n-projection. Let $R = R_1 \cup \ldots \cup R_n$ be the n-range of E and let $N = N_1 \cup \ldots \cup N_n$ be the null space of E.

1. The n-vector $\beta = \beta_1 \cup \ldots \cup \beta_n$ is the n-range R if and only if $E\beta = \beta$. If $\beta = E\alpha$ then $E\beta = E^2\alpha = E\alpha = \beta$. Conversely if $\beta = E\beta$ then of course $\beta$, is in the n-range of E.

2. $V = R \oplus N$ i.e., $V_1 \cup \ldots \cup V_n = R_1 \oplus N_1 \cup \ldots \cup R_n \oplus N_n$, i.e., each $V_i = R_i \oplus N_i$; i = 1, 2, …, n.

3. The unique expression for $\alpha$ as a sum of n-vector in R and N is $\alpha = E\alpha + (\alpha - E\alpha)$, i.e., $\alpha_1 \cup \ldots \cup \alpha_n = E_1\alpha_1 + (\alpha_1 - E_1\alpha_1) \cup \ldots \cup E_n\alpha_n + (\alpha_n - E_n\alpha_n)$.

From (1), (2) and (3) it is easy to verify. If $R = R_1 \cup \ldots \cup R_n$ and $N = N_1 \cup \ldots \cup N_n$ are n-subspaces of V such that $V = R \oplus N = R_1 \oplus N_1 \cup \ldots \cup R_n \oplus N_n$ there is one and only one n-projection operator $E = E_1 \cup \ldots \cup E_n$ which has n-range R and n-null space N. That operator is called the n-projection on R along N. Any n-projection $E = E_1 \cup\ldots\cup E_n$ is trivially n-diagonalizable. If $\{\alpha_1^1, \ldots, \alpha_{r_1}^1\} \cup \ldots \cup \{\alpha_1^n, \ldots, \alpha_{r_n}^n\}$ is a n-basis of R and $\left\{\alpha_{r_1+1}^1,...,\alpha_{n_1}^1\right\} \cup ... \cup \left\{\alpha_{r_n+1}^n,...,\alpha_{n_n}^n\right\}$ a n-basis for N then the n-basis $B = \{\alpha_1^1, \ldots, \alpha_{n_1}^1\} \cup \ldots \cup \{\alpha_1^n, \ldots, \alpha_{n_n}^n\} = B_1 \cup \ldots \cup B_n$, n-diagonalizes $E = E_1 \cup \ldots \cup E_n$.

$$\begin{aligned}[E]_B &= [E_1]_{B_1} \cup \ldots \cup [E_n]_{B_n} \\ &= \begin{bmatrix} I_1 & 0 \\ 0 & 0 \end{bmatrix} \cup \ldots \cup \begin{bmatrix} I_n & 0 \\ 0 & 0 \end{bmatrix},\end{aligned}$$

where $I_1$ is a $r_1 \times r_1$ identity matrix so on and $I_n$ is a $r_n \times r_n$ identity matrix.

Projections can be used to describe the n-direct sum decomposition of the n-space $V = V_1 \cup \ldots \cup V_n$. For suppose

$$V = (W_1^1 \oplus \ldots \oplus W_{k_1}^1) \cup \ldots \cup \{W_1^n \oplus \ldots \oplus W_{k_n}^n\}$$

for each j(t) we define $E_j^t$ on $V_t$. Let $\alpha$ be in $V = V_1 \cup \ldots \cup V_n$ say

$$\alpha = (\alpha_1^1 + \ldots + \alpha_{k_1}^1) \cup \ldots \cup (\alpha_1^n + \ldots + \alpha_{k_n}^n)$$

with $\alpha_i^t$ in $W_i^t$; $1 \leq i \leq k_t$ and t = 1, 2,…, n. Define $E_j^t\alpha^t = \alpha_j^t$; then $E_j^t$ is a well defined rule. It is easy to see that $E_j^t$ is linear and that range of $E_j^t$ is $W_j^t$ and $(E_j^t)^2 = E_j^t$. The null space of $E_j^t$ is the subspace $W_1^t + \ldots + W_{j-1}^t + W_{j+1}^t + \ldots + W_k^t$ for the statement $E_j^t\alpha^t = 0$ simply means $\alpha_j^t = 0$ i.e., $\alpha$ is actually a sum of vectors from the spaces $W_i^t$ with i $\neq$ j. In terms of the projections $E_j^t$ we have $\alpha^t = E_1^t\alpha^t + \ldots + E_k^t\alpha^t$ for each $\alpha$ in V. The above equation implies; $I_t = E_1^t + \ldots + E_k^t$. Note also that if i $\neq$ j then $E_i^t E_j^t = 0$ because the range of $E_j^t$ is the subspace $W_j^t$ which is contained in the null space of $E_i^t$. This is true for each t, t = 1, 2, …, n. Hence true on the n-vector space

$$V = (W_i^t \oplus \ldots \oplus W_{k_1}^1) \cup \ldots \cup (W_1^n \oplus \ldots \oplus W_{k_n}^n).$$

Now as in case of n-vector spaces of type I we in case of n-vector spaces of type II obtain the proof of the following theorem.

**THEOREM 1.2.51:** *If V = $V_1 \cup \ldots \cup V_n$ is a n-vector space of type II and suppose*

$$V = (W_1^1 \oplus \ldots \oplus W_{k_1}^1) \cup \ldots \cup (W_1^n \oplus \ldots \oplus W_{k_n}^n)$$

*then there exists ($k_1$, …, $k_n$); n-linear operators $\{E_1^1, \ldots, E_{k_1}^1\}$ $\cup$ … $\cup$ $\{E_1^n, \ldots, E_{k_n}^n\}$ on V such that*

a. *each $E_i^t$ is a projection, i.e., $(E_i^t)^2 = E_i^t$ for t = 1, 2, …, n; 1 $\leq$ i $\leq$ $k_t$.*
b. *$E_i^t E_j^t = 0$ if i $\neq$ j; 1 $\leq$ i, j $\leq$ $k_t$; t = 1, 2, …, n.*
c. *$I = I_1 \cup \ldots \cup I_n$*
   *$= (E_1^1 + \ldots + E_{k_1}^1) \cup \ldots \cup (E_1^n + \ldots + E_{k_n}^n)$*
d. *The range of $E_i^t$ is $W_i^t$ i = 1, 2,…, $k_t$ and t = 1, 2, …, n.*

*Proof:* We are primarily interested in n-direct sum n-decompositions $V = (W_1^1 \oplus \ldots \oplus W_{k_1}^1) \cup \ldots \cup (W_1^n \oplus \ldots \oplus W_{k_n}^n)$ $= W_1 \cup \ldots \cup W_n$; where each of the n-subspaces $W_t$ is invariant under some given n-linear operator $T = T_1 \cup \ldots \cup T_n$.

Given such a n-decomposition of V, T induces n-linear operators $\{T_i^1 \cup \ldots \cup T_i^n\}$ on each $W_i^1 \cup \ldots \cup W_i^n$ by restriction. The action of T is, if $\alpha$ is a n-vector in V we have unique n-vectors $\{\alpha_1^1, \ldots, \alpha_{k_1}^1\} \cup \ldots \cup \{\alpha_1^n, \ldots, \alpha_{k_n}^n\}$ with $\alpha_i^t$ in $W_i^t$ such that $\alpha = \alpha_1^1 + \ldots + \alpha_{k_1}^1 \cup \ldots \cup \alpha_1^n + \ldots + \alpha_{k_n}^n$ and then

$$T\alpha = T_1^1\alpha_1^1 + \ldots + T_{k_1}^1\alpha_{k_1}^1 \cup \ldots \cup T_1^n\alpha_1^n + \ldots + T_{k_n}^n\alpha_{k_n}^n .$$

We shall describe this situation by saying that T is the n-direct sum of the operators $\{T_1^1, \ldots, T_{k_1}^1\} \cup \ldots \cup \{T_1^n, \ldots, T_{k_n}^n\}$. It must be remembered in using this terminology that the $T_i^t$ are not n-linear operators on the space V = $V_1 \cup \ldots \cup V_n$ but on the various n-subspaces

$$\begin{aligned} V &= W_1 \cup \ldots \cup W_n. \\ &= (W_1^1 \oplus \ldots \oplus W_{k_1}^1) \cup \ldots \cup (W_1^n \oplus \ldots \oplus W_{k_n}^n) \end{aligned}$$

which enables us to associate with each $\alpha = \alpha_1 \cup \ldots \cup \alpha_n$ in V a unique n, k-tuple $(\alpha_1^1, \ldots, \alpha_{k_1}^1) \cup \ldots \cup (\alpha_1^n, \ldots, \alpha_{k_n}^n)$ of vectors $\alpha_i^t \in W_i^t$; i = 1, 2, …, $k_t$; t = 1, 2, …, n.

$$\alpha = (\alpha_i^t + \ldots + \alpha_{k_1}^1) \cup \ldots \cup (\alpha_1^n + \ldots + \alpha_{k_n}^n)$$

in such a way that we can carry out the n-linear operators on V by working in the individual n-subspaces $W_i = W_i^1 + \ldots + W_i^n$. The fact that each $W_i$ is n-invariant under T enables us to view the action of T as the independent action of the operators $T_i^t$ on the n-subspaces $W_i^t$; i = 1, 2, …, $k_t$ and t = 1, 2, …, n. Our purpose is to study T by finding n-invariant n-direct sum decompositions in which the $T_i^t$ are operators of an elementary nature.

As in case of n-vector spaces of type I the following theorem can be derived verbatim, which is left for the reader.

**THEOREM 1.2.52:** *Let $T = T_1 \cup \ldots \cup T_n$ be a n-linear operator on the n-space $V = V_1 \cup \ldots \cup V_n$ of type II. Let*

$$\{W_1^1, \ldots, W_{k_1}^1\} \cup \ldots \cup \{W_1^n, \ldots, W_{k_n}^n\}$$

*and*

$$\{E_1^1, \ldots, E_{k_1}^1\} \cup \ldots \cup \{E_1^n, \ldots, E_{k_n}^n\}$$

*be as before. Then a necessary and sufficient condition that each n-subspace $W_i^t$ ; to be a n-invariant under $T_t$ for $1 \leq i \leq k_t$ is that $E_i^t T_t = E_i^t T_t$ or TE = ET for every $1 \leq i \leq k_t$ and t = 1, 2, …, n.*

Now we proceed onto define the notion of n-primary decomposition for n-vector spaces of type II over the n-field $F = F_1 \cup \ldots \cup F_n$.

Let T be a n-linear operator of a n-vector space $V = V_1 \cup \ldots \cup V_n$ of $(n_1, n_2, \ldots, n_n)$ dimension over the n-field F. Let $p = p_1 \cup \ldots \cup p_n$ be a n-minimal polynomial for T; i.e.,

$$p = (x - c_1^1)^{r_1^1} \ldots (x - c_{k_1}^1)x^{r_{k_1}^1} \cup \ldots \cup (x - c_1^n)^{r_1^n} \ldots (x - c_{k_n}^n)x^{r_{k_n}^n}$$

where $\{c_1^1, \ldots, c_{k_1}^1\} \cup \ldots \cup \{c_1^n, \ldots, c_{k_n}^n\}$ are distinct elements of F $= F_1 \cup \ldots \cup F_n$, i.e., $\{c_1^t, \ldots, c_{k_t}^t\}$ are distinct elements of $F_t$, t = 1, 2, …, n, then we shall show that the n-space $V = V_1 \cup \ldots \cup V_n$ is the n-direct sum of null spaces $(T_s - c_i^s I_s)^{r_i^s}$, i = 1, 2, …, $k_s$ and s = 1, 2, …, n.

The hypothesis about p is equivalent to the fact that T is n-triangulable.

Now we proceed on to give the primary n-decomposition theorem for a n-linear operator T on the finite dimensional n-vector space $V = V_1 \cup \ldots \cup V_n$ over the n-field $F = F_1 \cup \ldots \cup F_n$ of type II.

**THEOREM 1.2.53:** (PRIMARY N-DECOMPOSITION THEOREM) *Let $T = T_1 \cup \ldots \cup T_n$ be a n-linear operator on a finite $(n_1, n_2, \ldots, n_n)$ dimensional n-vector space $V = V_1 \cup \ldots \cup V_n$ over the n-field $F = F_1 \cup \ldots \cup F_n$. Let $p = p_1 \cup \ldots \cup p_n$ be the n-minimal polynomial for T. $p = p_{11}^{r_1^1} \ldots p_{1k_1}^{r_{k_1}^1} \cup \ldots \cup p_{n1}^{r_1^n} \ldots p_{nk_n}^{r_{k_n}^n}$ where $p_{t_i}^t$ are distinct irreducible monic polynomials over $F_t$, i = 1, 2, …, $k_t$ and t = 1, 2, …, n and the $r_i^t$ are positive integers. Let $W_i = W_i^1 \cup \ldots \cup W_i^n$ be the null space of*

$$p(T) = p_{1i}(T_i^1)^{r_i^1} \cup \ldots \cup p_{ni}(T_i^n)^{r_i^n};$$

*i = 1, 2, …, n. Then*

i. $V = W_1 \cup \ldots \cup W_n$
$= (W_1^1 \oplus \ldots \oplus W_{k_1}^1) \cup \ldots \cup (W_1^n \oplus \ldots \oplus W_{k_n}^n)$

ii. *each* $W_i = W_i^1 + \ldots + W_i^n$ *is n-invariant under T; i = 1, 2, …, n.*

iii. *if* $T_i^r$ *is the operator induced on* $W_i^r$ *by* $T_i$ *then the minimal polynomial for* $T_i^r$ *is r = 1, 2, …,* $k_i$*; i = 1, 2, …, n.*

The proof is similar to that of n-vector spaces of type I and hence left as an exercise for the reader. In view of this theorem we have the following corollary the proof of which is also direct.

**COROLLARY 1.2.14:** *If* $\{E_1^1, \ldots, E_{k_1}^1\} \cup \ldots \cup \{E_1^n, \ldots, E_{k_n}^n\}$ *are the n-projections associated with the n-primary decomposition of T =* $T_1 \cup \ldots \cup T_n$*, then each* $E_i^t$ *is a polynomial in T; 1 ≤ i ≤* $k_t$ *and t = 1, 2, …, n and accordingly if a linear operator U commutes with T then U commutes with each of the* $E_i$*, i.e., each subspace* $W_i$ *is invariant under U.*

*Proof:* We can as in case of n-vector spaces of type I define in case of n-linear operator T of type II, the notion of n-diagonal part of T and n-nilpotent part of T.

Consider the n-minimal polynomial for T, which is the product of first degree polynomials i.e., the case in which each $p_i$ is of the form $p_i^t = x - c_i^t$. Now the range of $E_i^t$ is the null space of $W_i^t$ of $(T_t - c_i^t I_t)^{r_i^t}$; we know by earlier results D is a diagonalizable part of T. Let us look at the n-operator

$$N = T - D$$
$$= (T_1 - D_1) \cup \ldots \cup (T_n - D_n).$$
$$T = (T_1 E_1^1 + \ldots + T_1 E_{k_1}^1) \cup (T_2 E_1^2 + \ldots + T_2 E_{k_2}^2) \cup \ldots \cup (T_n E_1^n + \ldots + T_n E_{k_n}^n)$$

and

$$D = (c_1^1 E_1^1 + \ldots + c_{k_1}^1 E_{k_1}^1) \cup \ldots \cup (c_1^n E_1^n + \ldots + c_{k_n}^n E_{k_n}^n)$$

so

$$N = \{(T_1 - c_1^1 I_1) E_1^1 + \ldots + (T_1 - c_{k_1}^1 I_1) E_{k_1}^1\} \cup \ldots \cup \{(T_n - c_1^n I_n) E_1^n + \ldots + (T_n - c_{k_1}^n I_n) E_{k_n}^n\}.$$

The reader should be familiar enough with n-projections by now so that

$$N^2 = \{(T_1 - c_1^1 I_1)^2 E_1^1 + \ldots + (T_1 - c_{k_1}^1 I_1)^2 E_{k_1}^1\} \cup \ldots \cup \{(T_n - c_1^n I_n)^2 E_1^n + \ldots + (T_n - c_{k_n}^n I_n)^2 E_{k_n}^n\};$$

and in general,

$$N^r = \{(T_1 - c_1^1 I_1)^{r_1} E_1^1 + \ldots + (T_1 - c_{k_1}^1 I_1)^{r_1} E_{k_1}^1\} \cup \ldots \cup \{(T_n - c_1^n I_n)^{r_n} E_1^n + \ldots + (T_n - c_{k_n}^n I_n)^{r_n} E_{k_n}^n\}.$$

When $r \geq r_i$ for each i we have $N^r = 0$ because each of the n-operators $(T_t - c_i^t I_t)^{r_t} = 0$ on the range of $E_i^t$; $1 \leq t \leq k_i$ and i = 1, 2, …, n. Thus $(T - c\, I)^r = 0$ for a suitable r.

Let N be a n-linear operator on the n-vector space $V = V_1 \cup \ldots \cup V_n$. We say that N is n-nilpotent if there is some n-positive integer $(r_1, \ldots, r_n)$ such that $N_i^{r_i} = 0$. We can choose $r > r_i$ for i = 1, 2, …, n then $N^r = 0$, where $N = N_1 \cup \ldots \cup N_n$.

In view of this we have the following theorem for n-vector space type II, which is similar to the proof of n-vector space of type I.

**THEOREM 1.2.54:** *Let $T = T_1 \cup \ldots \cup T_n$ be a n-linear operator on the $(n_1, n_2, \ldots, n_n)$ finite dimensional vector space $V = V_1 \cup \ldots \cup V_n$ over the field F. Suppose that the n-minimal polynomial for T decomposes over $F = F_1 \cup \ldots \cup F_n$ into a n-product of n-linear polynomials. Then there is a n-diagonalizable operator $D = D_1 \cup \ldots \cup D_n$ on $V = V_1 \cup \ldots \cup V_n$ and a nilpotent operator $N = N_1 \cup \ldots \cup N_n$ on $V = V_1 \cup \ldots \cup V_n$ such that*

*(i) $T = D + N$.*

*(ii) $D N = N D$.*

*The n-diagonalizable operator D and the n-nilpotent operator N are uniquely determined by (i) and (ii) and each of them is n-polynomial in $T_1, \ldots, T_n$.*

Consequent of the above theorem the following corollary is direct.

**COROLLARY 1.2.15:** *Let V be a finite dimensional n-vector space over the special algebraically closed field $F = F_1 \cup \ldots \cup F_n$. Then every n-linear operator $T = T_1 \cup \ldots \cup T_n$ on $V = V_1 \cup \ldots \cup V_n$ can be written as the sum of a n-diagonalizable operator $D = D_i \cup \ldots \cup D_n$ and a n-nilpotent operator $N = N_1 \cup \ldots \cup N_n$ which commute. These n-operators D and N are unique and each is a n-polynomial in $(T_1, \ldots, T_n)$.*

Let $V = V_1 \cup \ldots \cup V_n$ be a finite dimensional n-vector space over the n-field $F = F_1 \cup \ldots \cup F_n$ and $T = T_1 \cup \ldots \cup T_n$ be an arbitrary and fixed n-linear operator on V. If $\alpha = \alpha_1 \cup \ldots \cup \alpha_n$ is n-vector in V, there is a smallest n-subspace of V which is n-invariant under T and contains $\alpha$. This n-subspace can be defined as the n-intersection of all T-invariant n-subspaces which contain $\alpha$; however it is more profitable at the moment for us to look at things this way.

If $W = W_1 \cup \ldots \cup W_n$ be any n-subspace of the n-vector space $V = V_1 \cup \ldots \cup V_n$ which is n-invariant under T and contains $\alpha = \alpha_1 \cup \ldots \cup \alpha_n$ i.e., each $T_i$ in T is such that the subspace $W_i$ in $V_i$ is invariant under $T_i$ and contains $\alpha_i$; true for i = 1, 2, …, n.

Then W must also contain $T\alpha$ i.e., $T_i\alpha_i$ is in $W_i$ for each i = 1, 2, …, n, hence $T(T\alpha)$ is in W i.e., $T_i (T_i\alpha_i) = T_i^2\alpha_i$ is in W so on; i.e., $T_i^{m_i}(\alpha_i)$ is in $W_i$ for each i so $T^m(\alpha) \in W$ i.e., W must contain $g(T)\alpha$ for every n-polynomial $g = g_1 \cup \ldots \cup g_n$ over the n-field $F = F_1 \cup \ldots \cup F_n$. The set of all n-vectors of the form $g(T)\alpha = g_1(T_1)\alpha_1 \cup \ldots \cup g_n(T_n)\alpha_n$ with $g \in F[x] = F_1[x] \cup \ldots \cup F_n[x]$, is clearly n-invariant and is thus the smallest n-T-invariant (or T-n-invariant) n-subspace which contains $\alpha$. In view of this we give the following definition for a n-linear operator on V.

**DEFINITION 1.2.30:** *If $\alpha = \alpha_1 \cup \ldots \cup \alpha_n$ is any n-vector in $V = V_1 \cup \ldots \cup V_n$, a n-vector space over the n-field $F = F_1 \cup \ldots \cup F_n$. The T-n-cyclic n-subspace generated by $\alpha$ is a n-subspace $Z(\alpha; T) = Z(\alpha_1; T_1) \cup \ldots \cup Z(\alpha_n; T_n)$ of all n-vectors $g(T)\alpha = g_1(T_1)\alpha_1 \cup \ldots \cup g_n(T_n)\alpha_n$; $g = g_1 \cup \ldots \cup g_n$ in $F[x] = F_1[x] \cup \ldots \cup F_n[x]$. If $Z(\alpha; T) = V$ then $\alpha$ is called a n-cyclic vector for T.*

*Another way of describing this n-subspace $Z(\alpha; T)$ is that $Z(\alpha; T)$ is the n-subspace spanned by the n-vectors $T_\alpha^k$ ; $k \geq 0$ and thus $\alpha$ is a n-cyclic n-vector for T if and only if these n-vectors span V; i.e., each $T_{i\alpha_i}^{k_i}$ span $V_i$, $k_i \geq 0$ and thus $\alpha_i$ is a cyclic vector for $T_i$ if and only if these vectors span $V_i$, true for i = 1, 2, …, n.*

We caution the reader that the general n-operator T has no n-cyclic n-vectors.

For any T, the T n-cyclic n-subspace generated by the n-zero vector is the n-zero n-subspace of V. The n-space Z(α; T) is (1, …, 1) dimensional if and only if α is a n-characteristic vector for T. For the n-identity operator, every non zero n-vector generates a (1, 1, …, 1) dimensional n-cyclic n-subspace; thus if n-dim V > (1, 1, …, 1) the n-identity operator has non-cyclic vector. For any T and α we shall be interested in the n-linear relations $c_0\alpha + c_1T\alpha + \ldots + c_kT\alpha^k = 0$; where $\alpha = \alpha_1 \cup \ldots \cup \alpha_n$ so

$$c_0^1\alpha_1 + c_1^1T_1\alpha_1 + \ldots + c_{k_1}^1 T_1^{k_1}\alpha_1 = 0,$$

$$c_0^2\alpha_2 + c_1^2T_2\alpha_2 + \ldots + c_{k_2}^2 T_2^{k_2}\alpha_2 = 0$$

so on;

$$c_0^n\alpha_n + c_1^nT_n\alpha_n + \ldots + c_{k_n}^n T_n^{k_n}\alpha_n = 0;$$

between the n-vectors $T\alpha^j$, we shall be interested in the n-polynomials $g = g_1 \cup \ldots \cup g_n$ where $g_i = c_0^i + c_1^i x + \ldots + c_{k_i}^i x^{k_i}$ true for i = 1, 2, …, n. which has the property that $g(T)\alpha = 0$.

The set of all g in F[x] = $F_1[x] \cup \ldots \cup F_n[x]$ such that g(T)α =0 is clearly an n-ideal in F[x]. It is also a non zero n-ideal in F[x] because it contains the n-minimal polynomial p = $p_1 \cup \ldots \cup p_n$ of the n-operator T. p(T)α = $0 \cup \ldots \cup 0$, i.e., $p_1(T_1)\alpha_1 \cup \ldots \cup p_n(T_n)\alpha_n = 0 \cup \ldots \cup 0$ for every $\alpha = \alpha_1 \cup \ldots \cup \alpha_n$ in $V = V_1 \cup \ldots \cup V_n$).

**DEFINITION 1.2.31:** *If $\alpha = \alpha_1 \cup \ldots \cup \alpha_n$ is any n-vector in $V = V_1 \cup \ldots \cup V_n$ the T-annihilator ($T = T_1 \cup \ldots \cup T_n$) of $\alpha$ is the n-ideal $M(\alpha; T)$ in $F[x] = F_1[x] \cup \ldots \cup F_n[x]$ consisting of all n-polynomials $g = g_1 \cup \ldots \cup g_n$ over $F = F_1 \cup \ldots \cup F_n$ such that $g(T)\alpha = 0 \cup \ldots \cup 0$ i.e., $g_1(T_1)\alpha_1 \cup \ldots \cup g_n(T_n)\alpha_n = 0 \cup \ldots \cup 0$. The unique monic n-polynomial $p_\alpha = p_1\alpha_1 \cup \ldots \cup p_n\alpha_n$ which n-generates this n-ideal will also be called the n-T annihilator of $\alpha$ or T-n-annihilator of $\alpha$. The n-T-annihilator $p_\alpha$ n-divides the n-minimal n-polynomial of the n-operator T. Clearly n-deg ($p_\alpha$) > (0, 0, …, 0) unless $\alpha$ is the zero n-vector.*

We now prove the following interesting theorem.

**THEOREM 1.2.55:** *Let $\alpha = \alpha_1 \cup \ldots \cup \alpha_n$ be any non zero n-vector in $V = V_1 \cup \ldots \cup V_n$ and let $p_\alpha = p_{1\alpha_1} \cup \ldots \cup p_{n\alpha_n}$ be the n-T-annihilator of $\alpha$.*

*(i) The n-degree of $p_\alpha$ is equal to the n-dimension of the n-cyclic subspace $Z(\alpha; T) = Z(\alpha_1; T_1) \cup \ldots \cup Z(\alpha_n; T_n)$.*

*(ii) If the n-degree of $p_\alpha = p_{1\alpha_1} \cup \ldots \cup p_{n\alpha_n}$ is $(k_1, k_2, \ldots, k_n)$ then the n-vectors $\alpha, T\alpha, T\alpha^2, \ldots, T\alpha^{k-1}$ form a n-basis for $Z(\alpha; T)$ i.e.,*

$$\{\alpha_1, T_1\alpha_1, T_1^2\alpha_1, \ldots, T_1^{k_1-1}\alpha_1\} \cup$$
$$\{\alpha_2, T_2\alpha_2, T_2^2\alpha_2, \ldots, T_2^{k_2-1}\alpha_2\} \cup \ldots \cup$$
$$\{\alpha_n, T_n\alpha_n, T_n^2\alpha_n, \ldots, T_n^{k_n-1}\alpha_n\}$$

*form a n-basis for $Z(\alpha; T) = Z(\alpha_1; T_1) \cup \ldots \cup Z(\alpha_n; T_n)$ i.e., $Z(\alpha_i, T_i)$ has $\{\alpha_i, T_i\alpha_i, \ldots, T_i^{k_i-1}\alpha_i\}$ as its basis; true for every i, i = 1, 2, …, n.*

(iii) *If $U = U_1 \cup \ldots \cup U_n$ is a n-linear operator on $Z(\alpha, T)$ induced by T, then the n-minimal polynomial for U is $p_\alpha$.*

*Proof:* Let $g = g_1 \cup \ldots \cup g_n$ be a n-polynomial over the n-field $F = F_1 \cup \ldots \cup F_n$. Write $g = p_\alpha q + r$, i.e., $g_1 \cup \ldots \cup g_n = p_{1\alpha_1} q_1 + r_1 \cup \ldots \cup p_{n\alpha_n} q_n + r_n$ where $p_\alpha = p_{1\alpha_1} \cup \ldots \cup p_{n\alpha_n}$ for $\alpha = \alpha_1 \cup \ldots \cup \alpha_n$, $q = q_1 \cup \ldots \cup q_n$ and $r = r_1 \cup \ldots \cup r_n$ so $g_i = p_{i\alpha_i} q_i + r_i$ true for i = 1, 2, …, n. Here either $r = 0 \cup \ldots \cup 0$ or n-deg r < n-deg $p_\alpha = (k_1, \ldots, k_n)$.

The n-polynomial $p_\alpha q = p_{1\alpha_1} q_1 \cup \ldots \cup p_{n\alpha_n} q_n$ is in the T-n-annihilator of $\alpha = \alpha_1 \cup \ldots \cup \alpha_n$ and so $g(T)\alpha = r(T)\alpha$. i.e., $g_1(T_1)\alpha_1 \cup \ldots \cup g_n(T_n)\alpha_n = r_1(T_1)\alpha_1 \cup \ldots \cup r_n(T_n)\alpha_n$. Since $r = r_1 \cup \ldots \cup r_n = 0 \cup 0 \cup \ldots \cup 0$ or n-deg $r < (k_1, k_2, \ldots, k_n)$ the n-vector $r(T)\alpha = r_1T_1(\alpha_1) \cup \ldots \cup r_nT_n(\alpha_n)$ is a n-linear combination of the n-vectors $\alpha, T\alpha, \ldots, T^{k-1}\alpha$; i.e., a n-linear combination on n-vectors $\alpha = \alpha_1 \cup \ldots \cup \alpha_n$.

$$T\alpha = T_1\alpha_1 \cup \ldots \cup T_n\alpha_n,$$
$$T^2\alpha = T_1^2\alpha_1 \cup \ldots \cup T_n^2\alpha_n, \ldots,$$
$$T^{k-1}\alpha = T_1^{k_1-1}\alpha_1 \cup \ldots \cup T_n^{k_n-1}\alpha_n$$

and since $g(T)\alpha = g_1(T_1)\alpha_1 \cup \ldots \cup g_n(T_n)\alpha_n$ is a typical n-vector in $Z(\alpha; T)$ i.e., each $g_i(T_i)\alpha_i$ is a typical vector in $Z(\alpha_i; T_i)$ for i = 1, 2, …, n. This shows that these $(k_1, \ldots, k_n)$, n-vectors span $Z(\alpha; T)$.

These n-vectors are certainly n-linearly independent, because any non-trivial n-linear relation between them would give us a non zero n-polynomial g such that $g(T)(\alpha) = 0$ and n-deg g < n-deg $p_\alpha$ which is absurd. This proves (i) and (ii).

Let $U = U_1 \cup \ldots \cup U_n$ be a n-linear operator on $(Z_\alpha; T)$ obtained by restricting T to that n-subspace. If $g = g_1 \cup \ldots \cup g_n$ is any n-polynomial over $F = F_1 \cup \ldots \cup F_n$ then

$$p_\alpha(U)g(T)\alpha = p_\alpha(T)g(T)\alpha$$

i.e., $p_{1\alpha_1}(U_1)\, g_1(T_1)\, \alpha_1 \cup \ldots \cup p_{n\alpha_n}(U_n)\, g_n(T_n)\, \alpha_n$

$= p_{1\alpha_1}(T_1)\, g_1(T_1)\, \alpha_1 \cup \ldots \cup p_{n\alpha_n}(T_n)\, g_n(T_n)\, \alpha_n$

$= \quad g(T)\, p_\alpha(T)\alpha$

$= \quad g_1(T_1)\, p_{1\alpha_1}(T_1)\, \alpha_1 \cup \ldots \cup g_n(T_n)\, p_{n\alpha_n}(T_n)\, \alpha_n$

$= \quad g_1(T_1)0 \cup \ldots \cup g_n(T_n)0$

$= \quad 0 \cup \ldots \cup 0.$

Thus the n-operator $p_\alpha U = p_1\alpha_1(U_1) \cup \ldots \cup p_n\alpha_n(U_n)$ sends every n-vector in $Z(\alpha; T) = Z(\alpha_1; T_1) \cup \ldots \cup Z(\alpha_n; T_n)$ into $0 \cup 0 \cup \ldots \cup 0$ and is the n-zero operator on $Z(\alpha; T)$. Furthermore if $h = h_1 \cup \ldots \cup h_n$ is a n-polynomial of n-degree less than $(k_1, \ldots, k_n)$ we cannot have $h(U) = h_1(U_1) \cup \ldots \cup h_n(U_n) = 0 \cup \ldots \cup 0$ for then $h(U)\alpha = h_1(U_1)\alpha_1 \cup \ldots \cup h_n(U_n)\alpha_n) = h_1(T_1)\alpha_1 \cup \ldots \cup h_n(T_n)\alpha_n = 0 \cup \ldots \cup 0$; contradicting the definition of $p_\alpha$. This show that $p_\alpha$ is the n-minimal polynomial for U.

A particular consequence of this nice theorem is that if $\alpha = \alpha_1 \cup \ldots \cup \alpha_n$ happens to be a n-cyclic vector for $T = T_1 \cup \ldots \cup T_n$ then the n-minimal polynomial for T must have n-degree equal to the n-dimension of the n-space $V = V_1 \cup \ldots \cup V_n$ hence by the Cayley Hamilton theorem for n-vector spaces we have that the n-minimal polynomial for T is the n-characteristic polynomial for T. We shall prove later that for any T there is a n-vector $\alpha = \alpha_1 \cup \ldots \cup \alpha_n$ in $V = V_1 \cup \ldots \cup V_n$ which has the n-minimal polynomial for $T = T_1 \cup \ldots \cup T_n$ for its n-annihilator. It will then follow that T has a n-cyclic vector if and only if the n-minimal and n-characteristic polynomials for T are identical. We now study the general n-operator $T = T_1 \cup \ldots \cup T_n$ by using n-operator vector. Let us consider a n-linear operator $U = U_1 \cup \ldots \cup U_n$ on the n-space $W = W_1 \cup \ldots \cup W_n$ of n-dimension $(k_1, \ldots, k_n)$ which has a cyclic n-vector $\alpha = \alpha_1 \cup \ldots \cup \alpha_n$.

By the above theorem just proved the n-vectors $\alpha, U\alpha, \ldots, U^{k-1}\alpha$ i.e., $\{\alpha_1, U_1\alpha_1, \ldots, U^{k-1}\alpha_1\}$, $\{\alpha_2, U_2\alpha_2, \ldots, U_2^{k_2-1}\alpha_2\}$, …, $\{\alpha_n, U_n\alpha_n, \ldots, U_n^{k_n-1}\alpha_n\}$ form a n-basis for the n-space $W = W_1 \cup \ldots \cup W_n$ and the annihilator $p_\alpha = p_{1\alpha_1} \cup \ldots \cup p_{n\alpha_n}$ of $\alpha = \alpha_1 \cup \ldots \cup \alpha_n$ is the n-minimal polynomial for $U = U_1 \cup \ldots \cup U_n$ (hence also the n-characteristic polynomial for U).

If we let $\alpha^i = U^{i-1}\alpha$; i.e., $\alpha^i = \alpha_1^i \cup \ldots \cup \alpha_n^i$ and $\alpha^i = U^{i-1}\alpha$ implies $\alpha_1^i = U_1^{i_1-1}\alpha_1, \ldots, \alpha_n^i = U_n^{i_n-1}\alpha_n$; $1 \leq i \leq k_i - 1$ then the action of U on the ordered n-basis $\{\alpha_1^1, \ldots, \alpha_{k_1}^1\} \cup \ldots \cup \{\alpha_1^n, \ldots, \alpha_{k_n}^n\}$ is $U\alpha^i = \alpha^{i+1}$ for i = 1, 2, …, k – 1 i.e., $U_t\alpha_t^i = \alpha_t^{i_t+1}$ for i = 1, 2, …, $k_t - 1$ and t = 1, 2, …, n. $U\alpha^k = -c_0\alpha^1 - \ldots - c_{k-1}\alpha^k$ i.e.,

$$U_t\alpha_t^k = -c_0^t\alpha_t^1 - \ldots - c_{k_t-1}^t\alpha_t^k$$

for t = 1, 2, …, n, where

$$\begin{aligned} p_\alpha = \; & \{c_0^1 + c_1^1 x + \ldots + c_{k_1-1}^1 x^{k_1-1} + x^{k_1}\} \cup \ldots \cup \\ & \{c_0^n + c_1^n x + \ldots + c_{k_n-1}^n x^{k_n-1} + x^{k_n}\}. \end{aligned}$$

The n-expression for $U\alpha_k$ follows from the fact $p_\alpha(U)\alpha = 0 \cup \ldots \cup 0$ i.e.,

$$p_{1\alpha_1}(U_1)\alpha_1 \cup \ldots p_{1\alpha_n}(U_n)\alpha_n = 0 \cup \ldots \cup 0.$$

i.e.,

$$U^k\alpha + c_{k-1}U^{k-1}\alpha + \ldots + c_1U\alpha + c_0^\alpha = 0 \cup \ldots \cup 0,$$

i.e.,

$$\begin{aligned} & U_1^{k_1}\alpha_1 + c_{k_1-1}^1U_1^{k_1-1}\alpha_1 + \ldots + c_1^1U_1\alpha_1 + c_0^1\alpha_1 \cup \\ & U_2^{k_2}\alpha_2 + c_{k_2-1}^2U_2^{k_2-1}\alpha_2 + \ldots + c_1^2U_2\alpha_2 + c_0^2\alpha_2 \cup \ldots \cup \\ & U_n^{k_n}\alpha_n + c_{k_n-1}^nU_n^{k_n-1}\alpha_n + \ldots + c_1^nU_n\alpha_n + c_0^n\alpha_n \\ & = 0 \cup \ldots \cup 0. \end{aligned}$$

This is given by the n-matrix of $U = U_1 \cup \ldots \cup U_n$ in the n-ordered basis

$$\begin{aligned} B &= B_1 \cup \ldots \cup B_n \\ &= \{\alpha_1^1 \ldots \alpha_{k_1}^1\} \cup \ldots \cup \{\alpha_1^n \ldots \alpha_{k_n}^n\} \end{aligned}$$

$$= \begin{bmatrix} 0 & 0 & 0 & \ldots & 0 & -c_0^1 \\ 1 & 0 & 0 & \ldots & 0 & -c_1^1 \\ 0 & 1 & 0 & \ldots & 0 & -c_2^1 \\ \vdots & \vdots & \vdots & & \vdots & \vdots \\ 0 & 0 & 0 & \ldots & 1 & -c_{k_1-1}^1 \end{bmatrix} \cup \ldots \cup$$

$$\begin{bmatrix} 0 & 0 & 0 & \dots & 0 & -c_0^n \\ 1 & 0 & 0 & \dots & 0 & -c_1^n \\ 0 & 1 & 0 & \dots & 0 & -c_2^n \\ \vdots & \vdots & \vdots & & \vdots & \vdots \\ 0 & 0 & 0 & \dots & 1 & -c_{k_n-1}^n \end{bmatrix}.$$

The n-matrix is called the n-companion n-matrix of the monic n-polynomial $p_\alpha = p_{1\alpha_1} \cup \ldots \cup p_{n\alpha_n}$ or can also be represented with some flaw in notation as $p_{\alpha_1}^1 \cup \ldots \cup p_{\alpha_n}^n$ where p = $p^1 \cup \ldots \cup p^n$.

Now we prove yet another interesting theorem.

**THEOREM 1.2.56:** *If U = $U_1 \cup \ldots \cup U_n$ is a n-linear operator on a finite ($n_1, n_2, \ldots, n_n$) dimensional n-space W = $W_1 \cup \ldots \cup W_n$, then U has a n-cyclic n-vector if and only if there is some n-ordered n-basis for W in which U is represented by the n-companion n-matrix of the n-minimal polynomial for U.*

*Proof:* We have just noted that if U has a n-cyclic n-vector then there is such an n-ordered n-basis for W = $W_1 \cup \ldots \cup W_n$. Conversely if there is some n-ordered n-basis $\{\alpha_1^1, \ldots, \alpha_{k_1}^1\} \cup \ldots \cup \{\alpha_1^n, \ldots, \alpha_{k_n}^n\}$ for W in which U is represented by the n-companion n-matrix of its n-minimal polynomial it is obvious that $\alpha_1^1 \cup \ldots \cup \alpha_1^n$ is a n-cyclic vector for U.

Now we give another interesting corollary to the above theorem.

**COROLLARY 1.2.16:** *If A = $A_1 \cup \ldots \cup A_n$ is a n-companion n-matrix of a n-monic n-polynomial p = $p_1 \cup \ldots \cup p_n$ (each $p_i$ is monic) then p is both the n-minimal polynomial and the n-characteristic polynomial of A.*

*Proof:* One way to see this is to let U = $U_1 \cup \ldots \cup U_n$ a n-linear operator on $F_1^{k_1} \cup \ldots \cup F_n^{k_n}$ which is represented by A = $A_1 \cup \ldots \cup A_n$ in the n-ordered n-basis and to apply the earlier theorem and the Cayley Hamilton theorem for n-vector spaces.

Yet another method is by a direct calculation.

Now we proceed on to define the notion of n-cyclic decomposition or we can call it as cyclic n-decomposition and its n-rational form or equivalently rational n-form.

Our main aim here is to prove that any n-linear operator T on a finite $(n_1, \ldots, n_n)$ dimensional n-vector space V = $V_1 \cup \ldots \cup V_n$, there exits n-set of n-vectors $\{\alpha_1^1, \ldots, \alpha_{r_1}^1\}$, …, $\{\alpha_1^n, \ldots, \alpha_{r_n}^n\}$ in V such that V = $V_1 \cup \ldots \cup V_n$

$$= Z(\alpha_1^1; T_1) \oplus \ldots \oplus Z(\alpha_{r_1}^1; T_1) \cup Z(\alpha_1^2; T_2) \oplus \ldots \oplus Z(\alpha_{r_2}^2; T_2) \cup \ldots \cup Z(\alpha_1^n; T_n) \oplus \ldots \oplus Z(\alpha_{r_n}^n; T_n).$$

In other words we want to prove that V is a n-direct sum of n-T-cyclic n-subspaces. This will show that T is a n-direct sum of a n-finite number of n-linear operators each of which has a n-cyclic n-vector. The effect of this will be to reduce many problems about the general n-linear operator to similar problems about an n-linear operator which has a n-cyclic n-vector.

This n-cyclic n-decomposition theorem is closely related to the problem in which T n-invariant n-subspaces W have the property that there exists a T- n-invariant n-subspaces W' such that V = W ⊕ W', i.e., V = $V_1 \cup \ldots \cup V_n = W_1 \oplus W'_1 \cup \ldots \cup W_n \oplus W'_n$.

If W = $W_1 \cup \ldots \cup W_n$ is any n-subspace of finite $(n_1, \ldots, n_n)$ dimensional n-space V = $V_1 \cup \ldots \cup V_n$ then there exists a n-subspace $W^1 = W_1^1 + \ldots + W_n^1$ such that V = W ⊕ W'; i.e., V = $V_1 \cup \ldots \cup V_n = (W_1 \oplus W'_1) \cup \ldots \cup (W_n \oplus W'_n)$ i.e., each $V_i$ is a direct sum of $W_i$ and $W'_i$ for i = 1, 2, …, n, i.e., $V_i = W_i \oplus W'_n$.

Usually there are many such n-spaces W' and each of these is called n-complementary to W. We study the problem when a T- n-invariant n-subspace has a complementary n-subspace which is also n-invariant under the same T.

Let us suppose $V = W \oplus W'$; i.e., $V = V_1 \cup \ldots \cup V_n = W_1 \oplus W'_1 \cup \ldots \cup W_n \oplus W'_n$ where both W and W' are n-invariant under T, then we study what special property is enjoyed by the n-subspace W. Each n-vector $\beta = \beta_1 \cup \ldots \cup \beta_n$ in $V = V_1 \cup \ldots \cup V_n$ is of the form $\beta = \gamma + \gamma'$ where $\gamma$ is in W and $\gamma'$ is in W' where $\gamma = \gamma_1 \cup \ldots \cup \gamma_n$ and $\gamma' = \gamma'_1 \cup \ldots \cup \gamma'_n$.

If $f = f_1 \cup \ldots \cup f_n$ any n-polynomial over the scalar n-field $F = F_1 \cup \ldots \cup F_n$ then $f(T)\beta = f(T)\gamma + f(T)\gamma'$; i.e., $f_1(T_1)\beta_1 \cup \ldots \cup f_n(T_n)\beta_n$
$= \{f_1(T_1)\gamma_1 \cup \ldots \cup f_n(T_n)\gamma_n\} + \{f_1(T_1)\gamma'_1 \cup \ldots \cup f_n(T_n)\gamma'_n\}$
$= \{f_1(T_1)\gamma_1 \cup \ldots \cup f_n(T_n)\gamma'_1\} + \{f_1(T_1)\gamma_n \cup \ldots \cup f_n(T_n)\gamma'_n\}$.

Since W and W' are n-invariant under $T = T_1 \cup \ldots \cup T_n$ the n-vector $f(T)\gamma = f_1(T_1)\gamma_1 \cup \ldots \cup f_n(T_n)\gamma_n$ is in $W = W_1 \cup \ldots \cup W_n$ and $f(T)\gamma' = f_1(T_1)\gamma'_1 \cup \ldots \cup f_n(T_n)\gamma'_n$ is in $W' = W'_1 \cup \ldots \cup W'_n$.

Therefore $f(T)\beta = f_1(T_1)\beta_1 \cup \ldots \cup f_n(T_n)\beta_n$ is in W if and only if $f(T)\gamma' = 0 \cup \ldots \cup 0$; i.e., $f_1(T_1)\gamma'_1 \cup \ldots \cup f_n(T_n)\gamma'_n = 0 \cup \ldots \cup 0$. So if $f(T)\beta$ is in W then $f(T)\beta = f(T)\gamma$.

Now we define yet another new notion on the n-linear operator $T = T_1 \cup \ldots \cup T_n$.

**DEFINITION 1.2.32:** *Let $T = T_1 \cup \ldots \cup T_n$ be a n-linear operator on the n-vector space $V = V_1 \cup \ldots \cup V_n$ and let $W = W_1 \cup \ldots \cup W_n$ be a n-subspace of V. We say that W is n- T-admissible if*

1. *W is n-invariant under T*
2. *$f(T)\beta$ is in W, for each $\beta \in V$,*

*there exists a n-vector $\gamma$ in W such that $f(T)\beta = f(T)\gamma$; i.e., if W is n-invariant and has a n-complementary n-invariant n-subspace then W is n-admissible.*

*Thus n-admissibility characterizes those n-invariant n-subspaces which have n-complementary n-invariant n-subspaces.*

We see that the n-admissibility property is involved in the n-decomposition of the n-vector space V = $Z(\alpha_1; T) \oplus \ldots \oplus Z(\alpha_r; T)$; i.e., $V = V_1 \cup \ldots \cup V_n = \{Z(\alpha_1^1; T_1) \oplus \ldots \oplus Z(\alpha_{r_1}^1; T_1)\} \cup \ldots \cup \{Z(\alpha_1^n; T_n) \oplus \ldots \oplus Z(\alpha_{r_n}^n; T_n)\}$.

We arrive at such a n-decomposition by inductively selecting n-vectors $\{\alpha_1^1, \ldots, \alpha_{r_1}^1\} \cup \ldots \cup \{\alpha_1^n, \ldots, \alpha_{r_n}^1\}$.

Suppose that by some method or another we have selected the n-vectors $\{\alpha_1^1, \ldots, \alpha_{r_1}^1\} \cup \ldots \cup \{\alpha_1^n, \ldots, \alpha_{r_n}^n\}$ and n-subspaces which is proper say

$$\begin{aligned} W_j &= W_j^1 \cup \ldots \cup W_j^n \\ &= \{Z(\alpha_1^1; T_1) + \ldots + Z(\alpha_{j_1}^1; T_1)\} \cup \ldots \cup \\ &\quad \{Z(\alpha_1^n; T_n) + \ldots + Z(\alpha_{j_n}^n; T_n)\}. \end{aligned}$$

We find the non zero n-vector $(\alpha_{j_1+1}^1 \cup \ldots \cup \alpha_{j_n+1}^n)$ such that

$$W_j \cap Z(\alpha_{j+1}, T) = 0 \cup \ldots \cup 0$$

i.e.,

$$(W_j^1 \cap Z(\alpha_{j_1+1}^1; T_1)) \cup \ldots \cup (W_j^n \cap Z(\alpha_{j_n+1}^n; T_n))$$
$$= 0 \cup \ldots \cup 0$$

because the n-subspace.

$$W_{j+1} = W_j \oplus Z(\alpha_{j+1}, T)$$

i.e.,

$$W_{j+1} = W_{j_1+1}^1 \cup \ldots \cup W_{j_n+1}^n$$

$$W_j^1 \oplus Z(\alpha_{j_1+1}; T_1) \cup \ldots \cup W_j^n \oplus Z(\alpha_{j_n+1}; T_n)$$

would be at least one n-dimension nearer to exhausting V.

But are we guaranteed of the existence of such

$$\alpha_{j+1} = \alpha_{j_1+1}^1 \cup \ldots \cup \alpha_{j_n+1}^n .$$

If $\{\alpha^1_1, \ldots, \alpha^1_{j_1}\} \cup \ldots \cup \{\alpha^n_1, \ldots, \alpha^n_{j_n}\}$ have been chosen so that $W_j$ is T-n-admissible n-subspace then it is rather easy to find a suitable $\alpha^1_{j+1} \cup \ldots \cup \alpha^n_{j_n+1}$.

Let $W = W_1 \cup \ldots \cup W_n$ be a proper T-n-invariant n-subspace. Let us find a non zero n-vector $\alpha = \alpha_1 \cup \ldots \cup \alpha_n$ such that $W \cap Z(\alpha; T) = \{0\} \cup \ldots \cup \{0\}$, i.e., $W_1 \cap Z(\alpha_1; T_1) \cup \ldots \cup W_n \cap Z(\alpha_n; T_n) = \{0\} \cup \ldots \cup \{0\}$. We can choose some n-vector $\beta = \beta_1 \cup \ldots \cup \beta_n$ which is not in $W = W_1 \cup \ldots \cup W_n$ i.e., each $\beta_i$ is not in $W_i$; i = 1, 2 ,…, n. Consider the T-n-conductor $S(\beta; W) = S(\beta_1; W_1) \cup \ldots \cup S(\beta_n; W_n)$ which consists of all n-polynomials $g = g_1 \cup \ldots \cup g_n$ such that $g(T)\beta = g_1(T_1)\beta_1 \cup \ldots \cup g_n(T_n)\beta_n$ is in $W = W_1 \cup \ldots \cup W_n$. Recall that the n-monic polynomial $f = S(\beta; W)$; i.e., $f_1 \cup \ldots \cup f_n = S(\beta_1; W_1) \cup \ldots \cup S(\beta_n; W_n)$ which n-generate the n-ideals $S(\beta; W) = S(\beta_1; W_1) \cup \ldots \cup S(\beta_n; W_n)$; i.e., each $f_i = S(\beta_I; W_i)$ generate the ideal $S(\beta_I; W_i)$ for i = 1, 2, …, n, i.e., $S(\beta; W)$ is also the T-n-conductor of $\beta$ into W. The n-vector $f(T)\beta = f_1(T_1)\beta_1 \cup \ldots \cup f_n(T_n)\beta_n$ is in $W = W_1 \cup \ldots \cup W_n$. Now if W is T-n-admissible there is a $\gamma = \gamma_1 \cup \ldots \cup \gamma_n$ in W with $f(T)\beta = f(T)\gamma$. Let $\alpha = \beta - \gamma$ and let g be any n-polynomial. Since $\beta - \alpha$ is in W, $g(T)\beta$ will be in W if and only if $g(T)\alpha$ is in W; in other words $S(\alpha; W) = S(\beta; W)$. Thus the n-polynomial f is also the T-n conductor of $\alpha$ into W. But $f(T)\alpha = 0 \cup \ldots \cup 0$. That tells us $f_1(T_1)\alpha_1 \cup \ldots \cup f_n(T_n)\alpha_n = 0 \cup \ldots \cup 0$; i.e., $g(T)\alpha$ is in W if and only if $g(T)\alpha = 0 \cup \ldots \cup 0$ i.e., $g_1(T_1)\alpha_1 \cup \ldots \cup g_n(T_n)\alpha_n = 0 \cup \ldots \cup 0$. The n-subspaces $Z(\alpha; T) = Z(\alpha_1; T_1) \cup \ldots \cup Z(\alpha_n; T_n)$ and $W = W_1 \cup \ldots \cup W_n$ are n-independent and f is the T-n-annihilator of $\alpha$.

Now we prove the cyclic decomposition theorem for n-linear operators.

**THEOREM 1.2.57:** (N-CYCLIC DECOMPOSITION THEOREM): *Let $T = T_1 \cup \ldots \cup T_n$ be a n-linear operator on a finite dimensional $(n_1, \ldots, n_n)$ n-vector space $V = V_1 \cup \ldots \cup V_n$ and let $W_o = W_o^1 \cup \ldots \cup W_o^n$ be a proper T-n-admissible n-subspace of V. There* *exists non zero n-vectors* $\left\{\alpha_1^1 \ldots \alpha_{r_1}^1\right\} \cup \ldots \cup \left\{\alpha_1^n \ldots \alpha_{r_n}^n\right\}$ *in V with respective T-n-annihilators* $\left\{p_1^1 \ldots p_{r_1}^1\right\} \cup \ldots \cup \left\{p_1^n \ldots p_{r_n}^n\right\}$ *such that*

(i) $V = W_o \oplus Z(\alpha_1; T) \oplus \ldots \oplus Z(\alpha_r; T) =$

$$\left\{W_o^1 \oplus Z\left(\alpha_1^1;T_1\right) \oplus \cdots \oplus Z\left(\alpha_{r_1}^1;T_1\right)\right\} \cup$$
$$\left\{W_o^2 \oplus Z\left(\alpha_1^2;T_2\right) \oplus \cdots \oplus Z\left(\alpha_{r_2}^2;T_2\right)\right\} \cup \ldots \cup$$
$$\left\{W_o^n \oplus Z\left(\alpha_1^n;T_n\right) \oplus \cdots \oplus Z\left(\alpha_{r_n}^n;T_n\right)\right\}.$$

(ii) $p_{k_t}^t$ *divides* $p_{k_t-1}^t$, *k = 1, 2,…, r; t = 1, 2, … , n.*

*Further more the integer r and the n-annihilators* $\left\{p_1^1, ..., p_{r_1}^1\right\} \cup \ldots \cup \left\{p_1^n, ..., p_{r_n}^n\right\}$ *are uniquely determined by (i) and (ii) and infact that no* $\alpha_{k_t}^t$ *is zero, t = 1, 2, …, n.*

*Proof:* The proof is rather very long, hence they are given under four steps. For the first reading it may seem easier to take $W_0 = \{0\} \cup \ldots \cup \{0\} = W_0^1 \cup \ldots \cup W_0^n$; i.e., each $W_0^i = \{0\}$ for i = 1, 2, …, n, although it does not produce any substantial simplification. Throughout the proof we shall abbreviate f(T)β to fβ. i.e., $f_1(T_1)\beta_1 \cup \ldots \cup f_n(T_n)\beta_n$ to $f_1\beta_1 \cup \ldots \cup f_n\beta_n$.

STEP I:

There exists non-zero n-vectors $\left\{\beta_1^1 \ldots \beta_{r_1}^1\right\} \cup \ldots \cup \left\{\beta_1^n \ldots \beta_{r_n}^n\right\}$ in $V = V_1 \cup \ldots \cup V_n$ such that

1. $V = W_0 + Z(\beta_1; T) + \ldots + Z(\beta_r; T)$

$$= \left\{W_0^1 + Z\left(\beta_1^1;T_1\right) + \ldots + Z\left(\beta_{r_1}^1;T_1\right)\right\} \cup$$
$$\left\{W_0^2 + Z\left(\beta_1^2;T_2\right) + \ldots + Z\left(\beta_{r_2}^2;T_2\right)\right\} \cup \ldots \cup$$

$$\left\{W_0^n + Z\left(\beta_1^n; T_n\right) + \ldots + Z\left(\beta_{r_n}^n; T_n\right)\right\}$$

2. if $1 \le k_i \le r_i$; i = 1, 2, …, n and

$$\begin{aligned} W_k &= W_{k_1}^1 + \cdots + W_{k_n}^n \\ &= \left\{W_o^1 + Z\left(\beta_1^1; T_1\right) + \cdots + Z\left(\beta_{k_1}^1; T_1\right)\right\} \cup \\ &\quad \left\{W_o^2 + Z\left(\beta_1^2; T_2\right) + \cdots + Z\left(\beta_{k_2}^2; T_2\right)\right\} \cup \ldots \cup \\ &\quad \left\{W_o^n + Z\left(\beta_1^n; T_n\right) + \cdots + Z\left(\beta_{k_n}^n; T_n\right)\right\} \end{aligned}$$

then the n-conductor

$$\begin{aligned} p_k &= p_{k_1}^1 \cup \ldots \cup p_{k_n}^n \\ &= S\left(\beta_{k_1}^1; W_{k_1-1}^1\right) \cup \ldots \cup S\left(\beta_{k_n}^n; W_{k_n-1}^n\right) \end{aligned}$$

has maximum n-degree among all T-n-conductors into the n-subspace $W_{k-1} = \left(W_{k_1-1}^1 \cup \ldots \cup W_{k_n-1}^n\right)$ that is for every $(k_1, \ldots, k_n)$;

$$\begin{aligned} \text{n-deg } p_k &= \max_{\alpha^1 \text{in } V_1} \deg\left(S\left(\alpha^1; W_{k_1-1}^1\right)\right) \cup \ldots \cup \\ &\quad \max_{\alpha^n \text{in } V_n} \deg\left(S\left(\alpha^n; W_{k_n-1}^n\right)\right). \end{aligned}$$

This step depends only upon the fact that $W_0 = \left(W_0^1 \cup \ldots \cup W_0^n\right)$ is an n-invariant n-subspace. If $W = W_1 \cup \ldots \cup W_n$ is a proper T-n-invariant n-subspace then

$$0 < \max_{\alpha} n \deg\left(S\left(\alpha; W\right)\right) \le n \dim V$$

i.e.,

$$\begin{aligned} 0 \cup \ldots \cup 0 &< \max_{\alpha_1} \deg\left(S\left(\alpha_1; W_1\right)\right) \cup \max_{\alpha_2} \deg\left(S\left(\alpha_2; W_2\right)\right) \\ &\quad \cup \ldots \cup \max_{\alpha_n} \deg\left(S\left(\alpha_n; W_n\right)\right) \\ &\le (n_1, \ldots, n_n) \end{aligned}$$

and we can choose a n-vector $\beta = \beta_1 \cup \ldots \cup \beta_n$ so that n-deg(S(β; W)) = deg(S($\beta_1$, $W_1$)) ∪ … ∪ deg(S($\beta_n$; $W_n$)) attains that maximum. The n-subspace W + Z(β; T) = ($W_1$ + Z($\beta_1$; $T_1$)) ∪ … ∪ ($W_n$ + Z($\beta_n$; $W_n$)) is then T-n-invariant and has n-dimension larger than n-dim W. Apply this process to W = $W_0$ to obtain $\beta_1 = \beta_1^1 \cup \cdots \cup \beta_n^1$. If $W_1 = W_0 + Z(\beta^1; T)$, i.e.,

$$W_1^1 \cup \ldots \cup W_n^1 = \left(W_0^1 + Z\left(\beta_1^1; T_1\right)\right) \cup \ldots \cup \left(W_0^n + Z\left(\beta_n^n; T_n\right)\right)$$

is still proper then apply the process to $W_1$ to obtain $\beta_2 = \beta_2^1 \cup \ldots \cup \beta_2^n$. Continue in that manner.

Since n-dim $W_k$ > n-dim $W_{k-1}$ i.e., $\dim W_{k_1}^1 \cup \ldots \cup \dim W_{k_n}^n >$ $\dim W_{k_1-1}^1 \cup \ldots \cup \dim W_{k_n-1}^n$ we must reach $W_r = V$; i.e., $W_{r_1}^1 \cup \ldots \cup W_{r_n}^n = V_1 \cup \ldots \cup V_n$ in not more than n-dim V steps.

STEP 2:

Let $\left\{\beta_1^1, \ldots, \beta_{r_1}^1\right\} \cup \ldots \cup \left\{\beta_1^n, \ldots, \beta_{r_n}^n\right\}$ be a n-set of nonzero n-vectors which satisfy conditions (1) and (2) of Step 1. Fix ($k_1$, …, $k_n$); $1 \le k_i \le r_i$; i = 1, 2, …, n. Let $\beta = \beta_1 \cup \ldots \cup \beta_n$ be any n-vector in $V = V_1 \cup \ldots \cup V_n$ and let $f = S(\beta; W_{k-1})$ i.e., $f_1 \cup \ldots \cup f_n = S(\beta_1; W_{k_1-1}^1) \cup \ldots \cup S(\beta_n; W_{k_n-1}^n)$.

If $f\beta = \beta_0 + \sum_{1\le i\le k} g_i\beta_i$ i.e., $f_1\beta_1 \cup \ldots \cup f_n\beta_n$

$$= \left(\beta_0^1 + \sum_{1\le i_1\le k_1} g_{i_1}^1\beta_{i_1}^1\right) \cup \left(\beta_0^2 + \sum_{1\le i_2\le k_2} g_{i_2}^2\beta_{i_2}^2\right) \cup \ldots \cup \left(\beta_0^n + \sum_{1\le i_n\le k_n} g_{i_n}^n\beta_{i_n}^n\right)$$

$\beta_{i_t} \in W_{i_t}^t$; t = 1, 2 ,…, n, then $f = f_1 \cup \ldots \cup f_n$ n-divides each n-polynomial $g_i = g_i^1 \cup \ldots \cup g_{i_n}^n$ and $\beta_0 = f\gamma_o$ i.e., $\beta_0^1 \cup \ldots \cup \beta_0^n = f_1\gamma_0^1 \cup \ldots \cup f_n\gamma_0^n$ where $\gamma_0 = \gamma_0^1 \cup \ldots \cup \gamma_0^n \in W_0 = W_0^1 \cup \ldots \cup W_0^n$.
If each $k_i$ = 1, i = 1, 2, …, n, this is just the statement that $W_o$ is T-n-admissible. In order to prove this assertion for ($k_1$, …, $k_n$) >

(1, 1, …, 1) apply the n-division algorithms $g_i = fh_i + r_i$ , $r_i = 0$ if n-deg $r_i$ < n-deg f i.e.,

$$g^1_{i_1} \cup \cdots \cup g^n_{i_n} = \left(f_i h^1_{i_1} + r^1_{i_1}\right) \cup \ldots \cup \left(f_n h^n_{i_n} + r^n_{i_n}\right);$$

$r_i = (0 \cup \ldots \cup 0)$ if n-deg $r_i$ < n-deg f.

We wish to show that $r_i = (0 \cup \ldots \cup 0)$ for each $i = (i_1 ,\ldots, i_n)$.

Let

$$\gamma = \beta - \sum_{1}^{k-1} h_i \beta_i ;$$

i.e.,

$$\gamma_1 \cup \ldots \cup \gamma_n = \left(\beta_1 - \sum_{1}^{k_1-1} h^1_{i_1} \beta^1_{i_1}\right) \cup \ldots \cup \left(\beta_n - \sum_{1}^{k_n-1} h^n_{i_n} \beta^n_{i_n}\right).$$

Since $\gamma - \beta$ is in $W_{k-1}$ i.e., $(\gamma_1 - \beta_1) \cup \ldots \cup (\gamma_n - \beta_n)$ is in $W^1_{k_1-1} \cup \ldots \cup W^n_{k_n-1}$. Since

$$S\left(\gamma_i ; W^i_{k_i-1}\right) = S\left(\beta_i ; W^i_{k_i-1}\right)$$

i.e.,

$$S\left(\gamma_1; W^1_{k_1-1}\right) \cup \ldots \cup S\left(\gamma_n ; W^n_{k_n-1}\right)$$

$$= S\left(\beta_1; W^1_{k_1-1}\right) \cup \cdots \cup S\left(\beta_n ; W^n_{k_n-1}\right)$$

$$= f_1 \cup \ldots \cup f_n.$$

$$S(\gamma; W_{k-1}) = S(\beta; W_{k-1}) = f.$$

Further more $f\gamma = \beta_0 + \sum_{1}^{k-1} \gamma_i \beta_i$ i.e.,

$$f_1\gamma_1 \cup \ldots \cup f_n\gamma_n = \left(\beta^1_0 + \sum_{i_1}^{k_1-1} \gamma^1_{i_1} \beta^1_{i_1}\right) \cup \ldots \cup \left(\beta^n_0 + \sum_{i_n}^{k_n-1} \gamma^n_{i_n} \beta^n_{i_n}\right)$$

$r_j = \left(r^1_{j_1}, \ldots, r^n_{j_n}\right) \neq (0, \ldots, 0)$ and n-deg $r_j$ < n-deg f. Let $p = S(\gamma; W_{j-1})$; $j = (j_1, \ldots, j_n)$. i.e.,

$$p_1 \cup \ldots \cup p_n = S\left(\gamma_1, W^1_{j_1-1}\right) \cup \ldots \cup S\left(\gamma_n, W^n_{j_n-1}\right).$$

Since

$$W_{k-1} = W^1_{k_1-1} \cup \ldots \cup W^n_{k_n-1}$$

contains

$$W_{j-1} = W^1_{j_1-1} \cup \ldots \cup W^n_{j_n-1} \;;$$

the n-conductor $f = S(\gamma; W_{k-1})$ i.e., $f_1 \cup \ldots \cup f_n = S(\gamma_1;\ W^1_{k_1-1}) \cup \ldots \cup S(\gamma_n;\ W^n_{k_n-1})$ must n-divide p. $p = fg$ i.e., $p_1 \cup \ldots \cup p_n = f_1g_1 \cup \ldots \cup f_ng_n$. Apply $g(T) = g_1(T_1) \cup \ldots \cup g_n(T_n)$ to both sides i.e., $p\gamma = gf\gamma = gr_j\,\beta_j + g\beta_o + \sum_{1\le i<j} gr_i\beta_i$ i.e.,

$p_1\gamma_1 \cup \ldots \cup p_n\gamma_n = g_1f_1\gamma_1 \cup \ldots \cup g_nf_n\gamma_n$

$$= \left(g_1\gamma^1_{j_1}\beta^1_{j_1} + g_1\beta^1_0 + \sum_{1\le i_1<j_1} g_1r^1_{i_1}\beta^1_{i_1}\right) \cup$$

$$\left(g_2\gamma^2_{j_2}\beta^2_{j_2} + g_2\beta^2_0 + \sum_{1\le i_2<j_2} g_2r^2_{i_2}\beta^2_{i_2}\right) \cup \ldots \cup$$

$$\left(g_n\gamma^n_{j_n}\beta^n_{j_n} + g_n\beta^n_0 + \sum_{1\le i_n<j_n} g_nr^n_{j_n}\beta^n_{j_n}\right).$$

By definition, $p\gamma$ is in $W_{j-1}$ and the last two terms on the right side of the above equation are in $W_{j-1} = W^1_{j_1-1} \cup \ldots \cup W^n_{j_n-1}$. Therefore

$$gr_j\beta_j = g_1r_{j_1}\beta_{j_1} \cup \ldots \cup g_nr_{j_n}\beta_{j_n}$$

is in

$$W_{j-1} = W^1_{j_1-1} \cup \ldots \cup W^n_{j_n-1}.$$

Now using condition (2) of step 1

$$\text{n-deg}(gr_j) \quad \ge \quad \text{n-deg}(S(\beta_j; W_{j-1}))\;;$$

i.e.,

$\deg(g_1 r_{j_1}) \cup \ldots \cup \deg(g_n r_{j_n})$

$$\begin{aligned}
&\ge \quad \deg(S(\beta_{j_1}; W^1_{j_1-1})) \cup \ldots \cup \deg(S(\beta_{j_n};\ W^n_{j_n-1}))\\
&= \quad \text{n-deg } p_j\\
&= \quad \deg p^1_{j_1} \cup \ldots \cup \deg p^n_{j_n}\\
&\ge \quad \text{n-deg } (S(\gamma; W_{j-1})
\end{aligned}$$

$$= \deg (S(\gamma_1; W^1_{j_1-1})) \cup \ldots \cup \deg(S(\gamma_n; W^n_{j_n-1}))$$
$$= \text{n-deg } p$$
$$= \deg p_1 \cup \ldots \cup \deg p_n$$
$$= \text{n-deg } fg$$
$$= \deg f_1g_1 \cup \ldots \cup \deg f_ng_n.$$

Thus n-deg $r_j$ > n-deg f; i.e., $\deg r_{j_1} \cup \ldots \cup \deg r_{j_n} > \deg f_1 \cup \ldots \cup \deg f_n$ and that contradicts the choice of $j = (j_1, \ldots, j_n)$. We now know that $f = f_1 \cup \ldots \cup f_n$, n-divides each $g_i = g^1_{i_1} \cup \ldots \cup g^n_{i_n}$ i.e., $f_{i_t}$ divides $g^t_{i_t}$ for t = 1, 2, …, n and hence that $\beta_0 = f\gamma$ i.e., $\beta^1_0 \cup \ldots \cup \beta^n_0 = f_1\gamma_1 \cup \ldots \cup f_n\gamma_n$. Since $W_0 = W^1_0 \cup \ldots \cup W^n_0$ is T-n-admissible (i.e., each $W^k_0$ is $T_k$-admissible for k = 1, 2, …, n); we have $\beta_0 = f\gamma_0$ where $\gamma_0 = \gamma^1_0 \cup \ldots \cup \gamma^n_0 \in W_0 = W^1_0 \cup \ldots \cup W^n_0$; i.e., $\beta^1_0 \cup \ldots \cup \beta^n_0 = f_1\gamma^1_0 \cup \ldots \cup f_n\gamma^n_0$ where $\gamma_0 \in W_0$. We make a mention that step 2 is a stronger form of the assertion that each of the n-subspaces $W_1 = W^1_1 \cup \ldots \cup W^n_1$, $W_2 = W^1_2 \cup \ldots \cup W^n_2$, …, $W_r = W^1_r \cup \ldots \cup W^n_r$ is T-n admissible.

STEP 3:

There exists non zero n-vectors $(\alpha^1_1, \ldots, \alpha^1_{r_1}) \cup \ldots \cup (\alpha^n_1, \ldots, \alpha^n_{r_n})$ in $V = V_1 \cup \ldots \cup V_n$ which satisfy condition (i) and (ii) of the theorem. Start with n-vectors $\{\beta^1_1, \ldots, \beta^1_{r_1}\} \cup \ldots \cup \{\beta^n_n, \ldots, \beta^n_{r_n}\}$ as in step 1. Fix $k = (k_1, \ldots, k_n)$ as $1 \le k_i \le r_i$. We apply step 2 to the n-vector $\beta = \beta_1 \cup \ldots \cup \beta_n = \beta^1_{k_1} \cup \ldots \cup \beta^n_{k_n} = \beta_k$ and T-n-conductor $f = f_1 \cup \ldots \cup f_n = p^1_{k_1} \cup \ldots \cup p^n_{k_n} = p_k$. We obtain $p_k\beta_k = p_k\gamma_0 + \sum_{1 \le i < k} p_k h_i \beta_i$;

i.e.,

$$p^1_{k_1}\beta^1_{k_1} \cup \ldots \cup p^n_{k_n}\beta^n_{k_n} =$$
$$\left(p^1_{k_1}\gamma^1_0 + \sum_{1 \le i_1 < k_1} p^1_{k_1} h^1_{i_1} \beta^1_{i_1}\right) \cup \ldots \cup \left(p^n_{k_n}\gamma^n_0 + \sum_{1 \le i_n < k_n} p^n_{k_n} h^n_{i_n} \beta^n_{i_n}\right)$$

where $\gamma_0 = \gamma_0^1 \cup \ldots \cup \gamma_0^n$ is in $W_0 = W_0^1 \cup \ldots \cup W_0^n$ and $\left\{h_1^1, \ldots, h_{k_1-1}^1\right\} \cup \ldots \cup \left\{h_1^n, \ldots, h_{k_n-1}^n\right\}$ are n-polynomials. Let

$$\alpha_k = \beta_k - \gamma_0 - \sum_{1 \le i < k} h_i \beta_i$$

i.e.,

$$\left\{\alpha_{k_1}^1 \cup \ldots \cup \alpha_{k_n}^n\right\} =$$

$$\left(\beta_{k_1}^1 - \gamma_0^1 - \sum_{1 \le i_1 < k_1} h_{i_1}^1 \beta_{i_1}^1\right) \cup \ldots \cup \left(\beta_{k_n}^n - \gamma_0^n - \sum_{1 \le i_n < k_n} h_{i_n}^n \beta_{i_n}^n\right).$$

Since $\beta_k - \alpha_k = \left(\beta_{k_1}^1 - \alpha_{k_1}^1\right) \cup \ldots \cup \left(\beta_{k_n}^n - \alpha_{k_n}^n\right)$ is in

$$
\begin{aligned}
W_{k-1} &= W_{k_1-1}^1 \cup \ldots \cup W_{k_n-1}^n; \\
S(\alpha_k; W_{k-1}) &= S(\beta_k; W_{k-1}) \\
&= p_k \\
&= S\left(\alpha_{k_1}^1, W_{k_1-1}^1\right) \cup \ldots \cup S\left(\alpha_{k_n}^n, W_{k_n-1}^n\right) \\
&= S\left(\beta_{k_1}^1, W_{k_1-1}^1\right) \cup \ldots \cup S\left(\beta_{k_n}^n, W_{k_n-1}^n\right) \\
&= p_{k_1}^1 \cup \ldots \cup p_{k_n}^n
\end{aligned}
$$

and since

$$p_k \alpha_k = 0 \cup \ldots \cup 0.$$

$$p_{k_1}^1 \alpha_{k_1}^1 \cup \ldots \cup p_{k_n}^n \alpha_{k_n}^n = 0 \cup \ldots \cup 0.$$

we have

$$W_{k-1} \cap Z(\alpha_k; T) = \{0\} \cup \ldots \cup \{0\}.$$

That is

$$W_{k_1-1}^1 \cap Z\left(\alpha_{k_1}^1; T_1\right) \cup \ldots \cup W_{k_1-1}^n \cap Z\left(\alpha_{k_n}^n; T_n\right)$$

$$= \{0\} \cup \ldots \cup \{0\}.$$

Because each $\alpha_k = \alpha_k^1 \cup \ldots \cup \alpha_k^n$ satisfies the above two equations just mentioned it follows that $W_k = W_0 \oplus Z(\alpha_1; T) \oplus \ldots \oplus Z(\alpha_k; T)$; i.e.,

$$W_{k_1}^1 \cup \ldots \cup W_{k_n}^n = \left\{W_0^1 \oplus Z\left(\alpha_1^1; T_1\right) \oplus \ldots \oplus Z\left(\alpha_{k_1}^1, T_1\right)\right\} \cup$$

$$\left\{W_0^2 \oplus Z\left(\alpha_1^2;T_2\right) \oplus \ldots \oplus Z\left(\alpha_{k_2}^2;T_2\right)\right\} \cup \ldots \cup$$
$$\left\{W_0^n \oplus Z\left(\alpha_1^n;T_n\right) \oplus \ldots \oplus Z\left(\alpha_{k_n}^n;T_n\right)\right\}$$

and that $p_k = p_{k_1}^1 \cup \ldots \cup p_{k_n}^n$ is the T-n-annihilator of $\alpha_k = \alpha_{k_1}^1 \cup \ldots \cup \alpha_{k_n}^n$. In other words, n-vectors $\left\{\alpha_1^1, \ldots, \alpha_{r_1}^1\right\}$, $\left\{\alpha_1^2, \ldots, \alpha_{r_2}^2\right\}$, …, $\left\{\alpha_1^n, \ldots, \alpha_{r_n}^n\right\}$ define the same n-sequence of n-subspaces. $W_1 = \left\{W_1^1 \cup \ldots \cup W_1^n\right\}$, $W_2 = \left\{W_2^1 \cup \ldots \cup W_2^n\right\}$ …. as do the n-vectors $\left\{\beta_1^1 \cup \ldots \cup \beta_{r_1}^1\right\}, \left\{\beta_1^2 \cup \ldots \cup \beta_{r_2}^2\right\}, \ldots, \left\{\beta_1^n \cup \ldots \cup \beta_{r_n}^n\right\}$ and the T-n-conductors $p_k = S(\alpha_k; W_{k-1})$ that is

$$\left\{p_{k_1}^1 \cup \ldots \cup p_{k_n}^n\right\} = S\left(\alpha_{k_1}^1; W_{k_1-1}^1\right) \cup \ldots \cup S\left(\alpha_{k_n}^n; W_{k_n-1}^n\right)$$

have the same maximality properties.The n-vectors $\left\{\alpha_1^1, \ldots, \alpha_{r_1}^1\right\}$ $\left\{\alpha_1^n, \ldots, \alpha_{r_n}^n\right\}$ have the additional property that the n-subspaces

$$W_0 = \left\{W_0^1 \cup \ldots \cup W_0^n\right\},$$
$$Z(\alpha_1; T) = Z\left(\alpha_1^1;T_1\right) \cup \ldots \cup Z\left(\alpha_1^n;T_n\right)$$
$$Z(\alpha_2; T) = Z\left(\alpha_2^1,T_2\right) \cup \ldots \cup Z\left(\alpha_2^n,T_n\right), \ldots$$

are n-independent. It is therefore easy to verify condition (ii) of the theorem. Since

$$p_i\alpha_i = \left(p_{i_1}^1\alpha_{i_1}^1\right) \cup \ldots \cup \left(p_{i_n}^n\alpha_{i_n}^n\right) = (0 \cup \ldots \cup 0)$$

we have the trivial relation

$$p_k\alpha_k = p_{k_1}^1\alpha_{k_1}^1 \cup \ldots \cup p_{k_n}^n\alpha_{k_n}^n$$
$$= \left(0 + p_1^1\alpha_1^1 + \ldots + p_{k_1-1}^1\alpha_{k_1-1}^1\right) \cup \left(0 + p_1^2\alpha_1^2 + \ldots + p_{k_2-1}^2\alpha_{k_2-1}^2\right) \cup \ldots$$
$$\cup \left(0 + p_1^n\alpha_1^n + \ldots + p_{k_n-1}^n\alpha_{k_n-1}^n\right).$$

Apply step 2 with $\left\{\beta_1^1, \ldots, \beta_{k_1}^1\right\}, \ldots, \left\{\beta_1^n, \ldots, \beta_{k_n}^n\right\}$ replaced by $\left\{\alpha_1^1, \ldots, \alpha_{k_1}^1\right\}, \ldots, \left\{\alpha_1^n, \ldots, \alpha_{k_n}^n\right\}$ and with $\beta = \beta_1 \cup \ldots \cup \beta_n$ $= \alpha^1_{k_1} \cup \ldots \cup \alpha^n_{k_n}$ ; $p_k$ n-divides each $p_i$, i < k that is $(i_1, \ldots, i_n) < (k_1, \ldots, k_n)$ i.e., $p^1_{k_1} \cup \ldots \cup p^n_{k_n}$ n-divides each $p^1_{i_1} \cup \ldots \cup p^n_{i_n}$ i.e., each $p^t_{k_t}$ divides $p^t_{i_t}$ for t = 1, 2, …, n.

STEP 4:
The number $r = (r_1, \ldots, r_n)$ and the n-polynomials $(p_1, \ldots, p_{r_1})$, …, $(p_1, \ldots, p_{r_n})$ are uniquely determined by the conditions of the theorem. Suppose that in addition to the n-vectors $\{\alpha^1_1, \ldots, \alpha^1_{r_1}\}, \{\alpha^2_1, \ldots, \alpha^2_{r_2}\}, \ldots, \{\alpha^n_1, \ldots, \alpha^n_{r_n}\}$ we have non zero n-vectors $\{\gamma^1_1, \ldots, \gamma^1_{s_1}\}, \ldots, \{\gamma^n_1, \ldots, \gamma^n_{s_n}\}$ with respective T-n-annihilators $\{g^1_1, \ldots, g^1_{s_1}\}, \ldots, \{g^n_1, \ldots, g^n_{s_n}\}$ such that

$$V = W_0 \oplus Z(\gamma_1; T) \oplus \ldots \oplus Z(\gamma_s; T)$$

i.e.,

$$V = V_1 \cup \ldots \cup V_n$$
$$\{W^1_0 \oplus Z(\gamma^1_1; T_1) \oplus \ldots \oplus Z(\gamma^1_{s_1}; T_1)\} \cup \ldots \cup$$
$$\{W^n_0 \oplus Z(\gamma^n_1; T_n) \oplus \ldots \oplus Z(\gamma^n_{s_n}; T_n)\}$$

$g^t_{k_t}$ divides $g^t_{k_t - 1}$ for t = 1, 2, …, n and $k_t = 2, \ldots, s_t$. We shall show that r = s i.e., $(r_1, \ldots, r_n) = (s_1, \ldots, s_n)$ and $p^t_i = g^t_i$; $1 \le t \le n$. i.e., $p^1_i \cup \ldots \cup p^n_i = g^1_i \cup \ldots \cup g^n_i$ for each i. We see that $p_1 = g_1 = p^1_1 \cup \ldots \cup p^n_1 = g^1_1 \cup \ldots \cup g^n_1$. The n-polynomial $g_1 = g^1_1 \cup \ldots \cup g^n_1$ is determined by the above equation as the T-n-conductor of V into $W_o$ i.e., $V = V_1 \cup \ldots \cup V_n$ into $W^1_0 \cup \ldots \cup W^n_0$. Let $S(V; W_o) = S(V_1; W^1_0) \cup \ldots \cup S(V_n; W^n_0)$ be the collection of all n-polynomials $f = f_1 \cup \ldots \cup f_n$ such that $f\beta = f_1\beta_1 \cup \ldots \cup f_n\beta_n$ is in $W_o$ for every in $\beta = \beta_1 \cup \ldots \cup \beta_n$ in $V = V_1 \cup \ldots \cup V_n$ i.e., n-polynomials f such that the n-range of f(T) = range of $f_1(T_1) \cup$ range of $f_2(T_2) \cup \ldots \cup$ range of $f_n(T_n)$ is contained in $W_0 = W^1_0 \cup \ldots \cup W^n_0$, i.e., range of $f_i(T_i)$ is in $W^i_0$ for i = 1, 2, …, n. Then $S(V_i; W_i)$ is a non zero ideal in the polynomial algebra so that we see S(V; $W_0$) = S($V_1$; $W_0^1$) ∪ … ∪ S($V_n$; $W_0^n$) is a non zero n-ideal in the n-polynomial algebra.

The polynomial $g_1^t$ is the monic generator of that ideal i.e., the n-monic polynomial $g_1$ = $g_1^1 \cup ... \cup g_1^n$ is the monic n-generator of that n-ideal; for this reason. Each β = $\beta_1$ ∪ … ∪ $\beta_n$ in V = $V_1$ ∪ … ∪ $V_n$ has the form β = $\left(\beta_0^1 + f_1^1\gamma_1^1 + ... + f_{s_1}^1\gamma_{s_1}^1\right)$ ∪ … ∪ $\left(\beta_0^n + f_1^n\gamma_1^n + ... + f_{s_n}^n\gamma_{s_n}^n\right)$ and so

$$g_1\beta = g_1\beta_0 + \sum_1^s g_1 f_i \gamma_i$$

i.e.,

$$g_1^1\beta_1 \cup ... \cup g_1^n\beta_n = \left(g_1^1\beta_0^1 + \sum_1^{s_1} g_1^1 f_{i_1}^1 \gamma_{i_1}^1\right) \cup \ldots \cup \left(g_1^n\beta_0^n + \sum_1^{s_n} g_1^n f_{i_n}^n \gamma_{i_n}^n\right).$$

Since each $g_i^t$ divides $g_1^t$ for t = 1, 2, …, n we have $g_1\gamma_i$ = 0 ∪ … ∪ 0 i.e., $g_1^1\gamma_{i_1}^1 \cup ... \cup g_1^n\gamma_{i_n}^n$ = 0 ∪ … ∪ 0 for all i = ($i_1$, …, $i_n$) and $g_0\beta = g_1\beta_0$ is in $W_0$ = $W_0^1$ ∪ … ∪ $W_0^n$. Thus $g_i^t$ is in S($V_t$; $W_0^t$) for t = 1, 2, …, n, so $g_1$ = $g_1^1$ ∪ …∪ $g_1^n$ is in

S(V; $W_0$) = S($V_1$; $W_0^1$) ∪ …∪ S($V_n$; $W_0^n$).

Since each $g_i^t$ is monic; $g_1$ is a monic n-polynomial of least n-degree which sends $\gamma_1^1$ into $W_0^t$ so that $\gamma_1$ = $\gamma_1^1 \cup ... \cup \gamma_1^n$ into $W_0 = W_0^1 \cup ... \cup W_0^n$; we see that $g_1$ = $g_1^1 \cup ... \cup g_1^n$ is the monic n-polynomial of least n-degree in the n-ideal S (V; $W_0$). By the same argument $p_i$ is the n-generator of that ideal so $p_1 = g_1$; i.e., $p_1^1 \cup ... \cup p_1^n$ = $g_1^1 \cup ... \cup g_1^n$.

If f = $f_1$ ∪ … ∪ $f_n$ is a n-polynomial and W = $W_1$ ∪ … ∪ $W_n$ is a n-subspace of V = $V_1$ ∪ … ∪ $V_n$ we shall employ the short hand fW for the set of all n-vectors $f\alpha = f_1\alpha_1 \cup \ldots \cup f_n\alpha_n$ with $\alpha = \alpha_1 \cup \ldots \cup \alpha_n$ in $W = W_1 \cup \ldots \cup W_n$.

The three facts can be proved by the reader:

1. $fZ(\alpha; T) = Z(f_\alpha; T)$ i.e., $f_1Z(\alpha_1; T_1) \cup \ldots \cup f_nZ(\alpha_n; T_n)$ $= Z(f_1\alpha_1; T_1) \cup \ldots \cup Z(f_n\alpha_n; T_n)$.
2. If $V = V_1 \oplus \ldots \oplus V_{k'}$ i.e., $V_1^1 \oplus \ldots \oplus V_{k_1}^1 \cup \ldots \cup V_1^n \oplus \ldots \oplus V_{k_n}^n$.

   where each $V_t$ is n-invariant under T that is each. $V_t^i$ is invariant under $T_i$ ; $1 \leq i < t$; $t = 1, 2, \ldots, n$ then $fV = fV_1 \oplus \ldots \oplus fV_n$; i.e., $f_1V_1 \cup \ldots \cup f_nV_n = f_1V_1 \oplus \ldots \oplus f_1 V_{k_1}^1 \cup \ldots \cup f_n V_1^n \oplus \ldots \oplus f_n V_{k_n}^n$.
3. If $\alpha = \alpha_1 \cup \ldots \cup \alpha_n$ and $\gamma = \gamma_1 \cup \ldots \cup \gamma_n$ have the same T-n-annihilator then $f\alpha$ and $f\gamma$ have the same T-n-annihilator and hence n-dim $Z(f\alpha; T)$ = n-dim $Z(f\gamma; T)$. i.e., $f\alpha = f_1\alpha_1 \cup \ldots \cup f_n\alpha_n$ and $f\gamma = f_1\gamma_1 \cup \ldots \cup f_n\gamma_n$ with $\dim Z(f_1\alpha_1; T_1) \cup \ldots \cup \dim Z(f_n\alpha_n; T_n) = \dim Z(f_1\gamma_1; T_1) \cup \ldots \cup \dim Z(f_n\gamma_n; T_n)$.

Now we proceed by induction to show that r = s and $p_i = g_i$ for i = 1, 2, …, r. The argument consists of counting n-dimensions in the proper way. We shall give the proof if $r = (r_1, r_2, \ldots, r_n) \geq (2, 2, \ldots, 2)$ then $p_2 = p_2^1 \cup \ldots \cup p_2^n = g_2 = g_2^1 \cup \ldots \cup g_2^n$ and from that the induction should be clear. Suppose that $r = (r_1 \cup \ldots \cup r_n) \geq (2, 2, \ldots, 2)$; then n-dim $W_0$ + n-dim $Z(\alpha_1; T)$ < n-dim V i.e., $(\dim W_0^1 \cup \ldots \cup \dim W_0^n) + \dim Z(\alpha_1^1; T_1) \cup \ldots \cup \dim Z(\alpha_1^n; T_n) < \dim V_1 \cup \ldots \cup \dim V_n$;

$$\dim W_0^1 + \dim Z(\alpha_1^1; T_1) \cup \ldots \cup \dim W_0^n + \dim Z(\alpha_1^n; T_n)$$
$$< \dim V_1 \cup \ldots \cup \dim V_n.$$

Since we know $p_1 = g_1$ we know that $Z(\alpha_1; T)$ and $Z(\gamma_1; T)$ have the same n-dimension. Therefore n-dim $W_0$ + n-dim $Z(\gamma_1; T)$ < n-dim V as before

$$\dim W_0^1 + \dim Z\left(\gamma_1^1;T_1\right) \cup \ldots \cup \dim W_0^n + \dim Z\left(\gamma_1^n;T_n\right) <$$
$$\dim V_1 \cup \ldots \cup \dim V_n,$$

which shows that s = ($s_1$, $s_2$, …, $s_n$) > (2, 2, …, 2). Now it makes sense to ask whether or not $p_2 = g_2$, $p_2^1 \cup ... \cup p_2^n = g_2^1 \cup ... \cup g_2^n$. From the two decompositions, of $V = V_1 \cup \ldots \cup V_n$, we obtain two decomposition of the n-subspace

$$p_2V = p_2^1V_1 \cup ... \cup p_2^nV_n$$
$$p_2V = p_2^1W_0 \oplus Z\left(p_2\alpha_1;T\right)$$

i.e.,

$$p_2^1V_1 \cup ... \cup p_2^nV_n =$$
$$p_2^1W_0^1 \oplus Z\left(p_2^1\alpha_1^1;T_1\right) \cup ... \cup p_2^nW_0^n \oplus Z\left(p_2^n\alpha_1^n;T_n\right)$$
$$p_2V = p_2W_0 \oplus Z\left(p_2\gamma_1;T\right) \oplus ... \oplus Z\left(p_2\gamma_s;T\right);$$

i.e.,

$$p_2^1V_1 \cup ... \cup p_2^nV_n = p_2^1W_0^1 \oplus Z\left(p_2^1\gamma_1^1;T_1\right) \oplus ... \oplus Z\left(p_2^1\gamma_{s_1}^1;T_1\right)$$
$$\cup ... \cup p_2^nW_0^n \oplus Z\left(p_2^n\gamma_1^n;T_n\right) \oplus ... \oplus Z\left(p_2^n\gamma_{s_n}^n;T_n\right).$$

We have made use of facts (1) and (2) above and we have used the fact $p_2\alpha_i = p_2^1\alpha_{i_1}^1 \cup ... \cup p_2^n\alpha_{i_n}^n = 0 \cup \ldots \cup 0$; i = ($i_1$, …, $i_n$) > (2, 2, …, 2). Since we know that $p_1 = g_1$ fact (3) above tell us that

$$Z\left(p_2\alpha_i;T\right) = Z\left(p_2^1\alpha_1^1;T_1\right) \cup ... \cup Z\left(p_2^n\alpha_1^n;T_n\right)$$

and $Z\left(p_2\gamma_1;T\right) = Z\left(p_2^1\gamma_1^1;T_1\right) \cup ... \cup Z\left(p_2^n\gamma_1^n;T_n\right)$, have the same n-dimension. Hence it is apparent from above equalities that n-dim $Z\left(p_2\gamma_i;T\right) = 0 \cup \ldots \cup 0$.

$$\dim Z\left(p_2^1\gamma_{i_1}^1;T_1\right) \cup ... \cup \dim Z\left(p_2^n\gamma_{i_n}^n;T_n\right) =$$
$$(0 \cup \ldots \cup 0);\ i = (i_1, \ldots, i_n) \geq (2, \ldots, 2).$$

We conclude $p_2\gamma_2 = \left(p_2^1\gamma_2^1\right) \cup \ldots \cup \left(p_2^n\gamma_2^n\right) = 0 \cup \ldots \cup 0$

and $g_2$ n-divides $p_2$ i.e., $g_2^t$ divides $p_2^t$ for each t; t = 1, 2, …, n. The argument can be reversed to show that $p_2$ n-divides $g_2$ i.e., $p_2^t$ divides $g_2^t$ for each t; t = 1, 2, …, n. Hence $p_2 = g_2$.

We leave the two corollaries for the reader to prove.

**COROLLARY 1.2.17:** *If $T = T_1 \cup \ldots \cup T_n$ is a n-linear operator on a finite $(n_1, \ldots, n_n)$ dimensional n-vector space $V = V_1 \cup \ldots \cup V_n$ then T-n-admissible n-subspace has a complementary n-subspace which is also invariant under T.*

**COROLLARY 1.2.18:** *Let $T = T_1 \cup \ldots \cup T_n$ be a n-linear operator on a finite $(n_1, \ldots, n_n)$ dimensional n-vector space $V = V_1 \cup \ldots \cup V_n$.*

a. *There exists n-vectors $\alpha = \alpha_1 \cup \ldots \cup \alpha_n$ in $V = V_1 \cup \ldots \cup V_n$ such that the T-n-annihilator of $\alpha$ is the n-minimal polynomial for T.*
b. *T has a n-cyclic n-vector if and only if the n-characteristic and n-minimal polynomial for T are identical.*

Now we proceed on to prove the Generalized Cayley-Hamilton theorem for n-vector spaces of finite n-dimension.

**THEOREM 1.2.58:** *(GENERALIZED CAYLEY HAMILTON THEOREM): Let $T = T_1 \cup \ldots \cup T_n$ be a n-linear operator on a finite $(n_1, n_2, \ldots, n_n)$ dimension n-vector space $V = V_1 \cup \ldots \cup V_n$. Let p and f be the n-minimal and n-characteristic polynomials for T, respectively*

i. *p n-divides f i.e., if $p = p_1 \cup \ldots \cup p_n$ and $f = f^1 \cup \ldots \cup f^n$ then $p_i$ divides $f^i$ = i = 1, 2, …, n.*
ii. *p and f have same prime factors expect for multiplicities.*

*iii. If p = $f_1^{r_1} \ldots f_k^{r_k}$ is the prime factorization of p then f = $f_1^{d_1} \cup \ldots \cup f_k^{d_k}$ where $d_i$ is the n-multiplicity of $f_i(T)^{r_i}$ n-divided by the n-degree of $f_i$.*

*That is if p = $p_1 \cup \ldots \cup p_n$*

$$= \left(f_1^1\right)^{r_1^1} \ldots \left(f_{k_1}^1\right)^{r_{k_1}^1} \cup \ldots \cup \left(f_1^n\right)^{r_1^n} \ldots \left(f_{k_n}^n\right)^{r_{k_n}^n};$$

*then f =* $\left(f_1^1\right)^{d_1^1} \ldots \left(f_{k_1}^1\right)^{d_{k_1}^1} \cup \ldots \cup \left(f_1^n\right)^{d_1^n} \ldots \left(f_{k_n}^n\right)^{d_{k_n}^n}$; $d_i^t = \left(d_1^t, \ldots, d_{k_t}^t\right)$ *is the nullity of* $f_i^t(T_t)^{r_i^t}$ *which is n-divided by the n-degree* $f_i^t$ *i.e.,* $d_i^t$; $1 \leq i \leq k_t$. *This is true for each t = 1, 2, …, n.*

*Proof:* The trivial case V = {0} ∪… ∪ {0} is obvious. To prove (i) and (ii) consider a n-cyclic decomposition

$$\begin{aligned} V &= Z(\alpha_1; T) \oplus \ldots \oplus Z(\alpha_r; T) \\ &= Z\left(\alpha_1^1; T_1\right) \oplus \ldots \oplus Z\left(\alpha_{r_1}^1; T_1\right) \cup \ldots \cup \\ &\quad Z\left(\alpha_1^n; T_n\right) \oplus \ldots \oplus Z\left(\alpha_{r_n}^n; T_n\right). \end{aligned}$$

By the second corollary $p_1$ = p. Let $U_1$ = $U_i^1 \cup \ldots \cup U_i^n$ be the n-restriction of T = $T_1 \cup \ldots \cup T_n$ i.e., each $U_i^s$ is the restriction of $T_s$ (for s = 1, 2, …, $r_s$) to $Z\left(\alpha_i^s; T_s\right)$. Then $U_i$ has a n-cyclic n-vector and so $p_i$ = $p_{i_1}^1 \cup \ldots \cup p_{i_n}^n$ is both n-minimal polynomial and the n-characteristic polynomial for $U_i$. Therefore the n-characteristic polynomial f = $f^1 \cup \ldots \cup f^n$ is the n-product f = $p_i^1 \ldots p_{r_1}^1 \cup \ldots \cup p_i^n \ldots p_{r_n}^n$. That is evident from earlier results that the n-matrix of T assumes in a suitable n-basis.

Clearly $p_1$ = p n-divides f; hence the claim (i). Obviously any prime n-divisor of p is a prime n-divisor of f. Conversely a prime n-divisor of f = $p_i^1 \ldots p_{r_1}^1 \cup \ldots \cup p_1^n \ldots p_{r_n}^n$ must n-divide one of the factor $p_i^t$ which in turn n-divides $p_1$.

Let p = $\left(f_1^1\right)^{r_1^1}\ldots\left(f_{k_1}^1\right)^{r_{k_1}^1}\cup\ldots\cup\left(f_1^n\right)^{r_1^n}\ldots\left(f_{k_n}^n\right)^{r_{k_n}^n}$ be the prime n-factorization of p. We employ the n-primary decomposition theorem, which tell us if $V_t^i = V_1^i \cup \ldots \cup V_n^i$ is the n-null space for $f_t^i\left(T_t\right)^{r_i^t}$ then

$$V = V_1 \oplus \ldots \oplus V_k = \left(V_1^1 \oplus \ldots \oplus V_1^{k_1}\right) \cup \ldots \cup \left(V_n^1 \oplus \ldots \oplus V_n^{k_n}\right)$$

and $\left(f_i^t\right)^{r_i^t}$ is the minimal polynomial of the operator $T_t^i$ restricting $T_t$ to the invariant subspace $V_t^i$. This is true for each t = 1, 2, …, n. Apply part (ii) of the present theorem to the n-operator $T_t^i$. Since its minimal polynomial is a power of the prime $f_i^t$ the characteristic polynomial for $T_t^i$ has the form $\left(f_i^t\right)^{r_i^t}$ where $d_i^t > r_i^t$, t = 1, 2, …, n.

We have $d_i^t = \dfrac{\dim V_t^i}{\deg f_i^t}$ for every t = 1, 2 ,…, n and dim $V_t^i$ = nullity $f_i^t\left(T_t\right)^{r_i^t}$ for every t = 1, 2, …, n. Since $T_t$ is the direct sum of the operators $T_t^1,\ldots,T_t^{k_t}$ the characteristic polynomial $f^t$ is the product, $f^t = \left(f_1^t\right)^{d_i^t}\cdots\left(f_{k_t}^t\right)^{d_{k_t}^t}$. Hence the claim.

The immediate corollary of this theorem is left as an exercise for the reader.

**COROLLARY 1.2.19:** *If T = $T_1 \cup \ldots \cup T_n$ is a n-nilpotent operator of the n-vector space of ($n_1$, …, $n_n$) dimension then the n-characteristic n-polynomial for T is $x^{n_1} \cup \ldots \cup x^{n_n}$.*

*Let us observe the n-matrix analogue of the n-cyclic decomposition theorem. If we have the n-operator T and the n-direct sum decomposition, let $B^i$ be the n-cyclic ordered basis.*

$$\left\{\alpha_{i_1}^1, T_1\alpha_{i_1}^1, \ldots, T_1^{k_{i_1}^1-1}\alpha_{i_1}^1\right\} \cup \ldots \cup \left\{\alpha_{i_n}^n, T_n\alpha_{i_n}^n, \ldots, T_n^{k_{i_1}^1-1}\alpha_{i_n}^n\right\}$$

*for* $Z(\alpha_i; T) = Z\left(\alpha^1_{i_1}; T_1\right) \cup \ldots \cup Z\left(\alpha^n_{i_n}; T_n\right)$. *Here* $\left(k^1_i, \ldots, k^1_{i_n}\right)$ *denotes n-dimension of* $Z(\alpha_i; T)$ *that is the n-degree of the n-annihilator* $p_i = p^1_{i_1} \cup \ldots \cup p^n_{i_n}$. *The n-matrix of the induced operator* $T_i$ *in the ordered n-basis* $B_i$ *is the n-companion n-matrix of the n-polynomial* $p_i$. *Thus if we let B to be the n-ordered basis for V which is the n-union of* $B^i$ *arranged in order* $\left\{B^1_1, \ldots, B^1_{r_1}\right\} \cup \ldots \cup \left\{B^n_1, \ldots, B^n_{r_n}\right\}$ *then the n-matrix of T in the ordered n-basis B will be* $A = A^1 \cup A^2 \cup \ldots \cup A^n$

$$= \begin{bmatrix} A^1_1 & 0 & \cdots & 0 \\ 0 & A^1_2 & \cdots & 0 \\ \vdots & \vdots & & \vdots \\ 0 & 0 & \cdots & A^1_{r_1} \end{bmatrix} \cup \ldots \cup \begin{bmatrix} A^n_1 & 0 & \cdots & 0 \\ 0 & A^n_2 & \cdots & 0 \\ \vdots & \vdots & & \vdots \\ 0 & 0 & \cdots & A^n_{r_n} \end{bmatrix}$$

where $A^t_i$ is the $k^t_i \times k^t_i$ companion matrix of $p^t_i$; for t = 1, 2, …, n. A $(n_1 \times n_1, \ldots, n_n \times n_n)$ n-matrix A which is the n-direct sum of the n-companion matrices of the non-scalar monic n-polynomial $\left\{p^1_i, \ldots, p^1_{r_1}\right\} \cup \ldots \cup \left\{p^n_1, \ldots, p^n_{r_n}\right\}$ such that $p^t_{i_t+1}$ divides $p^t_{i_t}$ for $i_t = 1, 2, \ldots, r_t - 1$ and t = 1, 2, …, n will be said to be the rational n-form or equivalently n-rational form.

**THEOREM 1.2.59:** *Let* $F = F_1 \cup \ldots \cup F_n$ *be a n-field. Let* $B = B_1 \cup \ldots \cup B_n$ *be a* $(n_1 \times n_1, \ldots, n_n \times n_n)$ *n-matrix over F. Then B is n-similar over the n-field F to one and only one matrix in the rational form.*

*Proof:* We know from the usual square matrix every square matrix over a fixed field is similar to one and only one matrix which is in the rational form.

So the n-matrix $B = B_1 \cup \ldots \cup B_n$ over the n-field $F = F_1 \cup \ldots \cup F_n$ is such that each $B_i$ is a $n_i \times n_i$ square matrix over $F_i$, is similar to one and only one matrix which is in the rational form, say $C_i$. This is true for every i, so $B = B_1 \cup \ldots \cup B_n$ is n-similar over the field to one and only one n-matrix C which is in the rational n-form.

The n-polynomials $\left\{\left(p_1^1,\ldots,p_{r_1}^1\right),\ldots,\left(p_1^n,\ldots,p_{r_n}^n\right)\right\}$ are called invariant n-factors or n-invariant factors for the n-matrix B = $B_1 \cup \ldots \cup B_n$.

We shall just introduce the notion of n-Jordan form or Jordan n-form. Suppose that N = $N_1 \cup \ldots \cup N_n$ be a nilpotent n-linear operator on a finite ($n_1, n_2,\ldots, n_n$) dimension space V = $V_1 \cup \ldots \cup V_n$. Let us look at the n-cyclic decomposition for N which we have depicted in the theorem. We have positive integers ($r_1, \ldots, r_n$) and non zero n-vectors $\left\{\alpha_1^n,\ldots,\alpha_{r_n}^n\right\}$ in V with n-annihilators $\left\{p_1^1,\ldots,p_{r_1}^1\right\}\cup\ldots\cup\left\{p_1^n,\ldots,p_{r_n}^n\right\}$ such that

$$\begin{aligned} V &= Z(\alpha_1;N)\oplus\cdots\oplus Z(\alpha_r;N) \\ &= Z(\alpha_1^1;N_1)\oplus\cdots\oplus Z(\alpha_{r_1}^1;N_1)\cup\cdots\cup \\ &\quad Z(\alpha_1^n;N_n)\oplus\cdots\oplus Z(\alpha_{r_n}^n;N_n) \end{aligned}$$

and $p_{i_t+1}^t$ divides $p_{i_t}^t$ for $i_t$ = 1, 2 ,…, $r_t$ – 1 and t = 1, 2, …, n. Since N is n-nilpotent the n-minimal polynomial is $x^{k_1}\cup\ldots\cup x^{k_n}$ with $k_t \le n_t$; t = 1, 2, …, n. Thus each $p_{i_t}^t$ is of the form $p_{i_t}^t = x^{k_i^t}$ and the n-divisibility condition simply says $k_1^t \ge k_2^t \ge \ldots \ge k_{r_t}^t$ ; t = 1, 2, …, n. Of course $k_1^t = k^t$ and $k_r^t \ge 1$. The n-companion n-matrix of $x^{k_{i_1}^1}\cup\ldots\cup x^{k_{i_n}^n}$ is the $k_{i_r} \times k_{i_r}$ n-matrix. A = $A_{i_1}^1 \cup\ldots\cup A_{i_n}^n$ with

$$A_{i_t}^t = \begin{bmatrix} 0 & 0 & \cdots & 0 & 0 \\ 1 & 0 & \ldots & 0 & 0 \\ 0 & 1 & & 0 & 0 \\ \vdots & \vdots & & \vdots & \vdots \\ 0 & 0 & \ldots & 1 & 0 \end{bmatrix}; t = 1, 2, \ldots, n.$$

Thus we from earlier results have an n-ordered n-basis for V = $V_1 \cup \ldots \cup V_n$ in which the n-matrix of N is the n-direct sum of the elementary nilpotent n-matrices the sizes $i_t$ of which decrease as $i_t$ increases. One sees from this that associated with a n-nilpotent ($n_1 \times n_1, \ldots, n_n \times n_n$), n-matrix is a positive n-tuple integers ($r_1, \ldots, r_n$) i.e.,

$$\left\{k_1^1, \ldots, k_{r_1}^1\right\} \cup \ldots \cup \left\{k_1^n, \ldots, k_{r_n}^n\right\}$$

such that

$$k_1^1 + \ldots + k_{r_1}^n = n_1$$

$$k_1^2 + \ldots + k_{r_2}^n = n_2$$

and

$$k_1^n + \ldots + k_{r_n}^n = n_n$$

and $k_{i_t}^t \geq k_{i_t+1}^t$; t = 1, 2, …, n and $1 \leq i$, $i + 1 \leq r_t$ and these n-sets of positive integers determine the n-rational form of the n-matrix, i.e., they determine the n-matrix up to similarity.

Here is one thing, we like to mention about the n-nilpotent, n-operator N above. The positive n-integer ($r_1, \ldots, r_n$) is precisely the n-nullity of N infact the n-null space has a n-basis with ($r_1, \ldots, r_n$) n-vectors $N_1^{k_{i_1}-1}\alpha_{i_1}^1 \cup \ldots \cup N_n^{k_{i_n}-1}\alpha_{i_n}^n$. For let $\alpha$ = $\alpha_1 \cup \ldots \cup \alpha_n$ be in n-nullspace of N we write

$$\alpha = \left(f_1^1\alpha_1^1 + \ldots + f_{r_1}^1\alpha_{r_1}^1\right) \cup \ldots \cup \left(f_1^n\alpha_1^n + \ldots + f_{r_1}^n\alpha_{r_n}^n\right)$$

where $\left(f_{i_1}^1, \ldots, f_{i_n}^n\right)$ is a n-polynomial the n- degree of which we may assume is less than $k_{i_1}, \ldots, k_{i_n}$. Since $N\alpha = 0 \cup \ldots \cup 0$; i.e., $N_1\alpha_1 \cup \ldots \cup N_n\alpha_n = 0 \cup \ldots \cup 0$ for each $i_r$ we have

$$\begin{aligned} 0 \cup 0 \cup \ldots \cup 0 &= N(f_i\, \alpha_i) \\ &= N_1\left(f_{i_1}\alpha_{i_1}\right) \cup \ldots \cup N_n\left(f_{i_n}\alpha_{i_n}\right) \\ &= N_1 f_{i_1}(N_1)\alpha_{i_1} \cup \cdots \cup N_n f_{i_n}\left(N_n\right)\alpha_{i_n} \\ &= \left(x f_{i_1}\right)\alpha_{i_1} \cup \ldots \cup (x f_{i_n})\alpha_{i_n}. \end{aligned}$$

Thus $xf_{i_1} \cup \ldots \cup xf_{i_n}$ is n-divisible by $x^{k_{i_1}} \cup \ldots \cup x^{k_{i_n}}$ and since n-deg $\left(f_{i_1}, \ldots, f_{i_n}\right) > \left(k_{i_1}, k_{i_2}, \ldots, k_{i_n}\right)$ this imply that

$$f_{i_1} \cup \ldots \cup f_{i_n} = c^1_{i_1} x^{k_{i_1-1}} \cup \ldots \cup c^n_{i_1} x^{k_{i_n-1}}$$

where $c^1_{i_1} \cup \ldots \cup c^n_{i_1}$ is some n-scalar. But then

$$\begin{aligned}
\alpha &= \alpha_1 \cup \ldots \cup \alpha_{n_s} \\
&= c^1_1\left(x^{k_{i_1-1}}\alpha^1_1\right) + \cdots + c^1_{r_1}\left(x^{k_{i_1-1}}\alpha^1_{r_1}\right) \cup \\
&\quad c^2_1\left(x^{k_{i_2-1}}\alpha^2_1\right) + \ldots + c^2_{r_2}\left(x^{k_{i_2-1}}\alpha^2_{r_2}\right) \cup \ldots \cup \\
&\quad c^n_1\left(x^{k_{i_n-1}}\alpha^n_1\right) + \ldots + c^n_{r_n}\left(x^{k_{i_n-1}}\alpha^n_{r_n}\right);
\end{aligned}$$

which shows that all the n-vectors form a n-basis for the n-null space of $N = N_1 \cup \ldots \cup N_n$. Suppose T is a n-linear operator on $V = V_1 \cup \ldots \cup V_n$ and that T factors over the n-field $F = F_1 \cup \ldots \cup F_n$ as

$$f = f_1 \cup \ldots \cup f_n$$

$$= \left(x - c^1_1\right)^{d^1_1} \ldots \left(x - c^1_{k_1}\right)^{d^1_{k_1}} \cup \ldots \cup \left(x - c^n_1\right)^{d^n_1} \cdots \left(x - c^n_{k_n}\right)^{d^1_{k_n}}$$

where $\left\{c^1_1, \ldots, c^1_{k_1}\right\} \cup \ldots \cup \left\{c^n_1, \ldots, c^n_{k_n}\right\}$ are n- distinct n-element of $F = F_1 \cup \ldots \cup F_n$ and $d^t_{i_t} \geq 1$; t = 1, 2, …, n.

Then the n-minimal polynomial for T will be

$$p = \left(x - c^1_1\right)^{r^1_1} \ldots \left(x - c^1_{k_1}\right)^{r^1_{k_1}} \cup \ldots \cup \left(x - c^n_1\right)^{r^n_1} \ldots \left(x - c^n_{k_n}\right)^{r^n_{k_n}}$$

where $1 < r^t_{i_t} \leq d^t_{i_t}$. t = 1, 2, …, n.

If $W^1_{i_1} \cup \ldots \cup W^n_{i_n}$ is the n-null space of

$$\left(T - c_i I\right)^{r_i} = \left(T_1 - C^1_{i_1} I_1\right)^{r^1_{i_1}} \cup \ldots \cup \left(T_n - c^n_{i_n} I_n\right)^{r^n_{i_n}}$$

then the n-primary decomposition theorem tells us that

$$V = V_1 \cup \ldots \cup V_n = \left(W^1_1 \oplus \ldots \oplus W^1_{k_1}\right) \cup \ldots \cup \left(W^n_1 \oplus \ldots \oplus W^n_{k_n}\right)$$

and that the operator $T^t_{i_t}$ induced on $W^t_{i_t}$ defined by $T^{i_t}_t$ has n-minimal polynomial $\left(x - c^t_{i_t}\right)^{r^t_i}$ for t = 1, 2, …, n; $1 \leq i_t \leq k_t$. Let $N^t_{i_t}$ be the n-linear operator on $W^t_{i_t}$ defined by $N^t_{i_t} = T^t_{i_t} - c^t_{i_t} I_t$. $1 \leq i_t \leq k_t$; then $N^t_{i_t}$ is n-nilpotent and has n-minimal polynomial $x^{r^t_{i_t}}_{i_t}$. On $W^t_{i_t}$, $T_t$ acts like $N^t_{i_t}$ plus the scalar $c^t_{i_t}$ times the identity operator. Suppose we choose a n-basis for the n-subspace $W^1_{i_1} \cup \ldots \cup W^n_{i_n}$ corresponding to the n-cyclic decomposition for the n-nilpotent $N^t_{i_t}$. Then the k-matrix $T^t_{i_t}$ in this ordered n-basis will be the n-direct sum of n-matrices.

$$\begin{bmatrix} c_1 & 0 & ... & 0 & 0 \\ 1 & c_1 & & 0 & 0 \\ \vdots & \vdots & & \vdots & \vdots \\ & & & c_1 & \\ 0 & 0 & \ldots & 1 & c_1 \end{bmatrix} \cup \begin{bmatrix} c_2 & 0 & ... & 0 & 0 \\ 1 & c_2 & & 0 & 0 \\ \vdots & \vdots & & \vdots & \\ & & & c_2 & \vdots \\ 0 & 0 & \ldots & 1 & c_2 \end{bmatrix} \cup \ldots \cup$$

$$\begin{bmatrix} c_n & 0 & ... & 0 & 0 \\ 1 & c_n & \ldots & 0 & 0 \\ \vdots & \vdots & & \vdots & \vdots \\ & & & c_n & \\ 0 & 0 & \ldots & 1 & c_n \end{bmatrix}$$

each with c = $c^t_{i_t}$ for t = 1, 2, …, n. Further more the sizes of these n-matrices will decrease as one reads from left to right. A n-matrix of the form described above is called an n-elementary Jordan n-matrix with n-characteristic values $c_1 \cup \ldots \cup c_n$.

Suppose we put all the n-basis for $W^1_{i_1} \cup \ldots \cup W^n_{i_n}$ together and we obtain an n-ordered n-basis for V. Let us describe the n-matrix A of T in the n-order basis.

The n-matrix A is the n-direct sum

$$A = \begin{bmatrix} A_1^1 & 0 & \dots & 0 \\ 0 & A_2^1 & \dots & 0 \\ \vdots & \vdots & & \vdots \\ 0 & 0 & \dots & A_{k_1}^1 \end{bmatrix} \cup \dots \cup \begin{bmatrix} A_1^n & 0 & \dots & 0 \\ 0 & A_2^n & \cdots & 0 \\ \vdots & \vdots & & \vdots \\ 0 & 0 & \cdots & A_{k_n}^n \end{bmatrix}$$

of the $k_i$ sets of n-matrices $\left\{A_1^1, ..., A_{k_1}^1\right\} \cup \dots \cup \left\{A_1^n, ..., A_{k_n}^n\right\}$. Each

$$A_{i_t}^t = \begin{bmatrix} J_{t1}^{(i_t)} & 0 & ... & 0 \\ 0 & J_{t2}^{(i_t)} & ... & 0 \\ & & & \vdots \\ 0 & 0 & ... & J_{tn}^{(i_t)} \end{bmatrix}$$

where each $J_{jt}^{(i_t)}$ is an elementary Jordan matrix with characteristic value $c_{i_t}^t$ ; $1 < i_t < k_t$; t = 1, 2, …, n. Also with in each $A_{i_t}^t$ the sizes of the matrices $J_{tj_t}^{i_t}$ decrease as $j_t$ increase $1 \leq j_t \leq n_t$ ; t = 1, 2, …, n. A $(n_1 \times n_1, \dots, n_n \times n_n)$ n-matrix A which satisfies all the conditions described so far for some n-sets of distinct $k_i$-scalars $\left\{c_1^1 \dots c_{k_1}^1\right\} \cup \dots \cup \left\{c_1^n \dots c_{k_n}^n\right\}$ will be said to be in Jordan n-form or n-Jordan form.

The interested reader can derive several interesting properties in this direction.

In the next chapter we move onto define the notion of n-inner product spaces for n-vector spaces of type II.

**Chapter Two**

# n-INNER PRODUCT SPACES OF TYPE II

Throughout this chapter we denote by $V = V_1 \cup \ldots \cup V_n$ only a n-vector space of type II over n-field $F = F_1 \cup \ldots \cup F_n$. The notion of n-inner product spaces for n-vector spaces of type II is introduced. We using these n-inner product define the new notion of n-best approximations n-orthogonal projection, quadratic n-form and discuss their properties.

**DEFINITION 2.1:** *Let $F = F_1 \cup \ldots \cup F_n$ be a n-field of real numbers and $V = V_1 \cup \ldots \cup V_n$ be a n-vector space over the n-field F. An n-inner product on V is a n-function which assigns to each n-ordered pair of n-vectors $\alpha = \alpha_1 \cup \ldots \cup \alpha_n$ and $\beta = \beta_1 \cup \ldots \cup \beta_n$ in V a n-scalar $(\alpha / \beta) = (\alpha_1 / \beta_1) \cup \ldots \cup (\alpha_n / \beta_n)$ in $F = F_1 \cup \ldots \cup F_n$ i.e., $(\alpha_i / \beta_i) \in F_i$, $i = 1, 2, \ldots, n$; in such a way that for all $\alpha = \alpha_1 \cup \ldots \cup \alpha_n$, $\beta = \beta_1 \cup \ldots \cup \beta_n$ and $\gamma = \gamma_1 \cup \ldots \cup \gamma_n$ in V and for all n-scalars $c = c_1 \cup \ldots \cup c_n$ in $F_1 \cup \ldots \cup F_n = F$.*

a. $(\alpha + \beta/\gamma) = (\alpha/\gamma) + (\beta/\gamma)$
i.e., $(\alpha_1+\beta_1/\gamma_1) \cup \ldots \cup (\alpha_n+\beta_n/\gamma_n) =$
$\left[(\alpha_1/\gamma_1)+(\beta_1/\gamma_1)\right]\cup\cdots\cup\left[(\alpha_n/\gamma_n)+(\beta_n/\gamma_n)\right]$

b. $(c\alpha/\beta) = (c\alpha/\beta)$ i.e.,
$(c_1\alpha_1/\beta_1) \cup \ldots \cup (c_n\alpha_n/\beta_n)$
$= c_1(\alpha_1/\beta_1) \cup \ldots \cup c_n(\alpha_n/\beta_n)$

c. $(\alpha/\beta) = (\beta/\alpha)$ i.e., $(\alpha_1/\beta_1) \cup \ldots \cup (\alpha_n/\beta_n)$
$= (\beta_1/\alpha_1) \cup \ldots \cup (\beta_n/\alpha_n)$

d. $(\alpha / \alpha) = (\alpha_1 / \alpha_1) \cup \ldots \cup (\alpha_n/\alpha_n) > (0 \cup \ldots \cup 0)$ if $\alpha_i \neq 0$ for $i = 1, 2, \ldots, n$.

*From the above conditions we have*

$(\alpha / c\beta + \gamma) = c\,(\alpha/\beta) + (\alpha / \gamma)$
$= (\alpha_1 / c_1\, \beta_1 + \gamma_1) \cup \ldots \cup (\alpha_n / c_n\, \beta_n + \gamma_n)$
$= [c_1\,(\alpha_1 / \beta_1) + (\alpha_1 / \gamma_1)] \cup \ldots \cup [c_n\,(\alpha_n / \beta_n) + (\alpha_n / \gamma_n)].$

*A n-vector space $V = V_1 \cup \ldots \cup V_n$ endowed with an n-inner product is called the n-inner product space. Let $F = F_1 \cup \ldots \cup F_n$ for $V = F_1^{n_1} \cup \ldots \cup F_n^{n_n}$ a n-vector space over F there is an n-inner product called the n-standard inner product. It is defined for $\alpha = \alpha_1 \cup \ldots \cup \alpha_n = \left(x_1^1 \ldots x_{n_1}^1\right) \cup \ldots \cup \left(x_1^n \ldots x_{n_n}^n\right)$ and $\beta = \beta_1 \cup \ldots \cup \beta_n = \left(y_1^1 \ldots y_{n_1}^1\right) \cup \ldots \cup \left(y_1^n \ldots y_{n_n}^n\right)$*

$$(\alpha/\beta) = \sum_{j_1} x_{j_1}^1 y_{j_1}^1 \cup \ldots \cup \sum_{j_n} x_{j_n}^n y_{j_n}^n .$$

*If $A = A_1 \cup \ldots \cup A_n$ is a n-matrix over the n-field $F = F_1 \cup \ldots \cup F_n$ where $A_i \in F_i^{n_i \times n_i}$ is a vector space over $F_i$ for $i = 1, 2, \ldots, n$ $V = F_1^{n_1 \times n_1} \cup \ldots \cup F_n^{n_n \times n_n}$ over the n-field $F = F_1 \cup \ldots \cup F_n$, then V is isomorphic to $F_1^{n_1^2} \cup \ldots \cup F_n^{n_n^2}$ in a natural way. It therefore follows*

$$(A / B) = \sum_{j_1 k_1} A_{j_1 k_1}^1 B_{j_1 k_1}^1 \cup \ldots \cup \sum_{j_n k_n} A_{j_n k_n}^n B_{j_n k_n}^n ,$$

*defines an n-inner product on V. A n-vector space over the n-field F is known as the n-inner product space.*

It is left as an exercise for the reader to verify.

**THEOREM 2.1:** *If $V = V_1 \cup \ldots \cup V_n$ is an n-inner product space then for any n-vectors $\alpha = \alpha_1 \cup \ldots \cup \alpha_n$, $\beta = \beta_1 \cup \ldots \cup \beta_n$ in V and any scalar $c = c_1 \cup \ldots \cup c_n$.*

(1) $|| c\alpha|| = | c | \; ||\alpha||$

i.e., $|| c\alpha || = || c_1\alpha_1|| \cup \ldots \cup ||c_n \alpha_n||$

$= |c_1| \; ||\alpha_1|| \cup \ldots \cup |c_n| \; ||\alpha_n||$

(2) $||\alpha|| > (0, \ldots, 0) = (0 \cup \ldots \cup 0)$ for $\alpha \neq 0$

i.e., $|| \alpha_1 || \cup \ldots \cup || \alpha_n|| > (0, 0, \ldots, 0) = (0 \cup 0 \cup \ldots \cup 0)$

(3) $|| (\alpha/\beta) || \leq || \alpha || \; || \beta ||$

$|| (\alpha_1/\beta_1) || \cup \ldots \cup || (\alpha_n / \beta_n) ||$

$\leq ||\alpha_1|| \; ||\beta_1 || \cup \ldots \cup || \alpha_n|| \; || \beta_n||.$

Proof as in case of usual inner produce space.

We wish to give some notation for n-inner product spaces

$$\gamma = \beta - \frac{(\beta / \alpha)}{|| \alpha ||^2}\alpha$$

$$\text{i.e., } \gamma_1 \cup \ldots \cup \gamma_n = \left(\beta_1 - \frac{(\beta_1 / \alpha_1)}{|| \alpha_1 ||^2}\alpha_1\right) \cup \ldots \cup \left(\beta_n - \frac{(\beta_n / \alpha_n)}{|| \alpha_n ||^2}\alpha_n\right).$$

It is left as an exercise for the reader to derive the Cauchy – Schwarz inequality in case of n-vector spaces. We now proceed on to define the notion of n-orthogonal n-set and n-orthonormal n-set.

**DEFINITION 2.2:** *Let $\alpha = \alpha_1 \cup \ldots \cup \alpha_n$, $\beta = \beta_1 \cup \ldots \cup \beta_n$ be n-vectors in an n-inner product space $V = V_1 \cup \ldots \cup V_n$. Then $\alpha$ is n-orthogonal to $\beta$ if $(\alpha/\beta) = 0 \cup \ldots \cup 0$.*

*i.e.,* $(\alpha_1 / \beta_1) \cup \ldots \cup (\alpha_n / \beta_n) = 0 \cup \ldots \cup 0$.

*Since this implies* $\beta = \beta_1 \cup \ldots \cup \beta_n$ *is n-orthogonal to* $\alpha = \alpha_1 \cup \ldots \cup \alpha_n$.

We often simply say that $\alpha$ and $\beta$ are n-orthogonal. Let $S = S_1 \cup \ldots \cup S_n$ be n-set of n-vectors in $V = V_1 \cup \ldots \cup V_n$. S is called an n-orthogonal n-set provided all n-pair of n-distinct n-vectors in S are n-orthogonal. An n-orthonormal n-set is an n-orthogonal n-set with addition property $\| \alpha \| = \| \alpha_1 \| \cup \ldots \cup \|\alpha_n\| = 1 \cup \ldots \cup 1$ for every $\alpha = \alpha_1 \cup \ldots \cup \alpha_n$ in $S = S_1 \cup \ldots \cup S_n$.

**THEOREM 2.2:** *A n-orthogonal n- set of non-zero n-vectors is n-linearly independent.*

The proof is left as an exercise for the reader.

**THEOREM 2.3:** *Let* $V = V_1 \cup \ldots \cup V_n$ *be an n-inner product space and let*

$$\left\{\beta_1^1, \ldots, \beta_{n_1}^1\right\} \cup \ldots \cup \left\{\beta_1^n, \ldots, \beta_{n_n}^n\right\}$$

*be any n-independent vectors in V. Then one way to construct n-orthogonal vectors*

$$\left\{\alpha_1^1, \ldots, \alpha_{n_1}^1\right\} \cup \ldots \cup \left\{\alpha_1^n, \ldots, \alpha_{n_n}^n\right\}$$

*in* $V = V_1 \cup \ldots \cup V_n$ *is such that for each k = 1, 2, …, n the n-set*

$$\left\{\alpha_1^1, \ldots, \alpha_{k_1}^1\right\} \cup \ldots \cup \left\{\alpha_1^n, \ldots, \alpha_{k_n}^n\right\}$$

*is a n-basis for the n-subspace spanned by*

$$\left\{\beta_1^1, \ldots, \beta_{k_1}^1\right\} \cup \ldots \cup \left\{\beta_1^n, \ldots, \beta_{k_n}^n\right\}.$$

*Proof:* The n-vectors $\left\{\alpha_1^1, \ldots, \alpha_{n_1}^1\right\} \cup \ldots \cup \left\{\alpha_1^n, \ldots, \alpha_{n_n}^n\right\}$ can be obtained by means of a construction analogous to Gram-Schmidt orthogonalization process called or defined as Gram-Schmidt n-orthogonalization process. First let $\alpha_1 = \alpha_1^1 \cup \ldots \cup \alpha_1^n$ and $\beta_1 = \beta_1^1 \cup \ldots \cup \beta_1^n$. The other n-vectors are then given inductively as follows:

Suppose $\left\{\alpha_1^1, \ldots, \alpha_{m_1}^1\right\} \cup \ldots \cup \left\{\alpha_1^n, \ldots, \alpha_{m_n}^n\right\}$ ($1 \le m_i \le n_i$; i = 1, 2, …, n) have been chosen so that for every $k_i$, $\left\{\alpha_1^1, \ldots, \alpha_{k_1}^1\right\} \cup \ldots \cup \left\{\alpha_1^n, \ldots, \alpha_{k_n}^n\right\}$; $1 \le k_i \le m_i$, i = 1, 2, …, n is an n-orthogonal basis for the n-subspace of $V = V_1 \cup \ldots \cup V_n$ that is spanned by $\left\{\beta_1^1, \ldots, \beta_{k_1}^1\right\} \cup \ldots \cup \left\{\beta_1^n, \ldots, \beta_{k_n}^n\right\}$. To construct next n-set of n-vectors, $\alpha_{m_1+1}^1, \ldots, \alpha_{m_n+1}^n$;

let

$$\alpha_{m_1+1}^1 = \beta_{m_1+1}^1 - \sum_{k_1=1}^{m_1} \frac{\left(\beta_{m_1+1}^1 / \alpha_{k_1}^1\right)}{\|\alpha_{k_1}^1\|^2}\ \alpha_k^1$$

and so on

$$\alpha_{m_n+1}^n = \beta_{m_n+1}^n - \sum_{k_n=1}^{m_n} \frac{\left(\beta_{m_n+1}^n / \alpha_{k_n}^n\right)}{\|\alpha_{k_n}^n\|^2}\ \alpha_{k_n}^n .$$

Thus $\left(\alpha_{m_1+1}^1 \cup \ldots \cup \alpha_{m_n+1}^n \mid \alpha_{j_1}^1 \cup \ldots \cup \alpha_{j_n}^n\right)$ = $(0 \cup \ldots \cup 0)$ for $1 \le j_t \le m_t$; t = 1, 2, …, n. For other wise $\left(\beta_{m_n+1}^1 \cup \ldots \cup \beta_{m_n+1}^n\right)$ will be a n-linear combination of $\left\{\alpha_1^1, \ldots, \alpha_{m_1}^1\right\} \cup \ldots \cup \left\{\alpha_1^n, \ldots, \alpha_{m_n}^n\right\}$. Further $1 \le j_t \le m_t$; t = 1, 2, …, n.

Then

$$\begin{aligned}
(\alpha_{m+1} / \alpha_j) &= \left(\alpha_{m_1+1}^1 / \alpha_{j_1}^1\right) \cup \ldots \cup \left(\alpha_{m_n+1}^n / \alpha_{j_n}^n\right) \\
&= \left[\left(\beta_{m_1+1}^1 / \alpha_{j_1}^1\right) - \sum_{k_1=1}^{m_1} \frac{\left(\beta_{m_1+1}^1 / \alpha_{k_1}^1\right)}{\|\alpha_{k_1}^1\|^2}\left(\alpha_{k_1}^1 / \alpha_{j_1}^1\right)\right] \cup \ldots \cup \\
&\quad \left[\left(\beta_{m_n+1}^n / \alpha_{j_n}^n\right) - \sum_{k_n=1}^{m_n} \frac{\left(\beta_{m_n+1}^n / \alpha_{k_n}^n\right)\left(\alpha_{k_n}^n / \alpha_{j_n}^n\right)}{\|\alpha_{k_n}^n\|^2}\right] \\
&= \left[\left(\beta_{m_1+1}^1 / \alpha_{j_1}^1\right) - \left(\beta_{m_1+1}^1 / \alpha_{j_1}^1\right)\right] \cup \ldots \cup \\
&= \left[\left(\beta_{m_n+1}^n / \alpha_{j_n}^n\right) - \left(\beta_{m_n+1}^n / \alpha_{j_n}^n\right)\right] \\
&= 0 \cup \ldots \cup 0.
\end{aligned}$$

Therefore $\left\{\alpha_1^1, \ldots, \alpha_{m_1+1}^1\right\} \cup \ldots \cup \left\{\alpha_1^n, \ldots, \alpha_{m_n+1}^n\right\}$ is an n-orthogonal n-set consisting of {$m_1$ + 1 ∪ … ∪ $m_n$ + 1} non zero n-vectors in the n-subspace spanned by $\left\{\beta_1^1, \ldots, \beta_{m_1+1}^1\right\} \cup \ldots \cup \left\{\beta_1^n, \ldots, \beta_{m_n+1}^n\right\}$. By theorem we have just proved it is a n-basis for this n-subspace. Thus the n-vectors $\left\{\alpha_1^1, \ldots, \alpha_{n_1}^1\right\} \cup \ldots \cup \left\{\alpha_1^n, \ldots, \alpha_{n_n}^n\right\}$ may be constructed one after the other given earlier. In particular when $\alpha_1 = \beta_1$.

i.e., $\left(\alpha_1^1 \cup \ldots \cup \alpha_{n_1}^1\right) = \left(\beta_1^1 \cup \ldots \cup \beta_{n_1}^1\right)$

$$\left(\alpha_1^2 \cup \ldots \cup \alpha_{n_2}^2\right) = \left[\beta_1^2 - \frac{\left(\beta_2^1/\alpha_1^1\right)}{\|\alpha_1^1\|^2}\alpha_1^1\right] \cup \ldots \cup \left[\beta_n^2 - \frac{\left(\beta_n^1/\alpha_n^1\right)}{\|\alpha_n^1\|^2}\alpha_n^1\right].$$

$$\left(\alpha_1^3 \cup \ldots \cup \alpha_n^3\right) = \left[\beta_1^3 - \frac{\left(\beta_1^3/\alpha_1^1\right)}{\|\alpha_1^1\|^2}\alpha_1^1 - \frac{\left(\beta_2^3/\alpha_2^2\right)}{\|\alpha_2^2\|^2}\alpha_2^2\right] \cup \ldots \cup \left[\beta_n^3 - \frac{\left(\beta_n^3/\alpha_n^1\right)}{\|\alpha_n^1\|^2}\alpha_n^1 - \frac{\left(\beta_n^3/\alpha_n^2\right)}{\|\alpha_n^2\|^2}\alpha_n^2\right].$$

The reader is expected to derive the following corollary.

**COROLLARY 2.1:** *Every finite ($n_1$ … $n_n$) dimensional n-inner product space has an n-orthonormal basis.*

Next we define the new notion of n-best approximations or best n-approximations.

**DEFINITION 2.3:** *Let V = $V_1$ ∪… ∪ $V_n$ be a n-inner product vector space over the n- field F = $F_1$ ∪… ∪ $F_n$ of type II. Let W = $W_1$ ∪… ∪ $W_n$ be a n-subspace of V. Let β = $\beta_1$ ∪… ∪ $\beta_n$ be* *a n-vector in V. To find the n-best approximation to* $\beta = \beta_1 \cup \ldots \cup \beta_n$ *by n-vectors in* $W = W_1 \cup \ldots \cup W_n$*. This means to find a vector* $\alpha = \alpha_1 \cup \ldots \cup \alpha_n$ *for which* $\|\beta - \alpha\| = \|\beta_1 - \alpha_1\| \cup \ldots \cup \|\beta_n - \alpha_n\|$ *is as small as possible subject to the restriction that* $\alpha = \alpha_1 \cup \ldots \cup \alpha_n$ *should belong to* $W = W_1 \cup \ldots \cup W_n$*. i.e., to be more precise.*

*A n-best approximation to* $\beta = \beta_1 \cup \ldots \cup \beta_n$ *in* $W = W_1 \cup \ldots \cup W_n$ *is a n-vector* $\alpha = \alpha_1 \cup \ldots \cup \alpha_n$ *in W such that*

$$\|\beta - \alpha\| \leq \|\beta - \gamma\|$$

*i.e.,*

$$\|\beta_1 - \alpha_1\| \cup \ldots \cup \|\beta_n - \alpha_n\| \leq \|\beta_1 - \gamma_1\| \cup \ldots \cup \|\beta_n - \gamma_n\|$$

*for every n-vector* $\gamma = \gamma_1 \cup \ldots \cup \gamma_n$ *in W.*

The following theorem is left as an exercise for the reader.

**THEOREM 2.4:** *Let* $W = W_1 \cup \ldots \cup W_n$ *be a n-subspace of an n-inner product space* $V = V_1 \cup \ldots \cup V_n$ *and* $\beta = \beta_1 \cup \ldots \cup \beta_n$ *be a n-vector in* $V = V_1 \cup \ldots \cup V_n$*.*

a. *The n-vector* $\alpha = \alpha_1 \cup \ldots \cup \alpha_n$ *in W is a n-best approximation to* $\beta = \beta_1 \cup \ldots \cup \beta_n$ *by n-vectors in* $W = W_1 \cup \ldots \cup W_n$ *if and only if* $\beta - \alpha = \beta_1 - \alpha_1 \cup \ldots \cup \beta_n - \alpha_n$ *is n-orthogonal to every vector in W. i.e., each* $\beta_i - \alpha_i$ *is orthogonal to every vector in* $W_i$*; true for i = 1, 2, …, n.*
b. *If a n-best approximation to* $\beta = \beta_1 \cup \ldots \cup \beta_n$ *by n-vectors in* $W = W_1 \cup \ldots \cup W_n$ *exists, it is unique.*
c. *If W is finite dimensional and* $\left(\alpha_1^1, \ldots, \alpha_{n_1}^1\right) \cup \ldots \cup \left(\alpha_1^n, \ldots, \alpha_{n_n}^n\right)$ *is any n-orthonormal n-basis for W then the n-vector*

$$\begin{aligned}\alpha &= \alpha_1 \cup \ldots \cup \alpha_n \\ &= \sum_{k_1} \frac{\left(\beta_1 / \alpha_{k_1}^1\right)}{\|\alpha_{k_1}^1\|^2} \alpha_{k_1}^1 \cup \ldots \cup \sum_{k_n} \frac{\left(\beta_n / \alpha_{k_n}^n\right)}{\|\alpha_{k_n}^n\|^2} \alpha_{k_n}^n\end{aligned}$$

*is the unique n-best approximation to* $\beta = \beta_1 \cup \ldots \cup \beta_n$ *by n-vectors in W.*

**DEFINITION 2.4:** *Let $V = V_1 \cup \ldots \cup V_n$ be n-inner product space and $S = S_1 \cup \ldots \cup S_n$ any set of n-vectors in V. The n-orthogonal complement of S denoted by $S^{\perp} = S_1^{\perp} \cup \ldots \cup S_n^{\perp}$ is the set of all n-vectors in V which are n-orthogonal to every n-vector in S.*

Whenever the n-vector $\alpha = \alpha_1 \cup \ldots \cup \alpha_n$ in the theorem 2.4 exists its is called the n-orthogonal projection to $\beta = \beta_1 \cup \ldots \cup \beta_n$ on $W = W_1 \cup \ldots \cup W_n$. If every n-vector has an n-orthogonal projection on $W = W_1 \cup \ldots \cup W_n$ the n-mapping that assigns to each n-vector in V its n-orthogonal projection on $W = W_1 \cup \ldots \cup W_n$ is called the n-orthogonal projection of V on W.

We leave the proof of this simple corollary as an exercise for the reader.

**COROLLARY 2.2:** *Let $V = V_1 \cup \ldots \cup V_n$ be an n-inner product space, $W = W_1 \cup \ldots \cup W_n$ a finite dimensional n-subspace and $E = E_1 \cup \ldots \cup E_n$ the n-orthogonal projection of V on W. Then the n-mapping $\beta \to \beta - E\beta$ i.e., $\beta_1 \cup \ldots \cup \beta_n \to (\beta_1 - E_1\beta_1) \cup \ldots \cup (\beta_n - E_n\beta_n)$ is the n-orthogonal projection of V on W.*

The following theorem is also direct and can be proved by the reader.

**THEOREM 2.5:** *Let W be a finite dimensional n-subspace of an n-inner product space $V = V_1 \cup \ldots \cup V_n$ and let $E = E_1 \cup \ldots \cup E_n$ be the n-orthogonal projection of V on W. Then $E = E_1 \cup \ldots \cup E_n$ is an idempotent n-linear transformation of V onto W, $W^{\perp}$ is the null space of E and $V = W \oplus W^{\perp}$ i.e. $V = V_1 \cup \ldots \cup V_n = W_1 \oplus W_1^{\perp} \cup \ldots \cup W_n \oplus W_n^{\perp}$.*

The following corollaries are direct and expect the reader to prove them.

**COROLLARY 2.3:** *Under the conditions of the above theorem $1 - E$ is the n-orthogonal n-projection of V on $W^{\perp}$. It is an n-*

*idempotent n-linear transformation of V onto $W^{\perp}$ with n-null space W.*

**COROLLARY 2.4:** *Let $\{\{\alpha_1^1, \ldots, \alpha_{n_1}^1\} \cup \{\alpha_1^2, \ldots, \alpha_{n_2}^2\} \cup \ldots \cup \{\alpha_1^n, \ldots, \alpha_{n_n}^n\}\}$ be an n-orthogonal set of nonzero n-vectors in an n-inner product space $V = V_1 \cup \ldots \cup V_n$. If $\beta = \beta_1 \cup \ldots \cup \beta_n$ is any n-vector in V then*

$$\sum_{K_1} \frac{|(\beta_1 \mid \alpha_{K_1}^1)|^2}{\|\alpha_{K_1}^1\|^2} \cup \sum_{K_2} \frac{|(\beta_2 \mid \alpha_{K_2}^2)|^2}{\|\alpha_{K_2}^2\|^2} \cup \ldots \cup \sum_{K_n} \frac{|(\beta_n \mid \alpha_{K_n}^n)|^2}{\|\alpha_{K_n}^n\|^2} \leq \|\beta_1\|^2 \cup \ldots \cup \|\beta_n\|^2$$

*and equality holds if and only if*

$$\beta = \sum_{K_1} \frac{(\beta_1 \mid \alpha_{K_1}^1)}{\|\alpha_{K_1}\|^2} \alpha_{K_1}^1 \cup \ldots \cup \sum \frac{(\beta_n \mid \alpha_{K_n}^n)}{\|\alpha_{Kn}^n\|^2} \alpha_{K_n}^n = \beta_1 \cup \ldots \cup \beta_n .$$

Now we proceed onto define the new notion of linear functionals and adjoints in case of n-vector spaces over n-fields of type II.

Let $V = V_1 \cup \ldots \cup V_n$ be a n-vector space over a n-field $F = F_1 \cup \ldots \cup F_n$ of type II.

Let $f = f_1 \cup \ldots \cup f_n$ be any n-linear functional on V. Define n-inner product on V, for every $\alpha = \alpha_1 \cup \ldots \cup \alpha_n \in V$; $f(\alpha) = (\alpha \mid \beta)$ for a fixed $\beta = \beta_1 \cup \ldots \cup \beta_n$ in V.

$f_1(\alpha_1) \cup \ldots \cup f_n(\alpha_n) = (\alpha_1 \mid \beta_1) \cup \ldots \cup (\alpha_n \mid \beta_n)$.

We use this result to prove the existence of the n-adjoint of a n-linear operator $T = T_1 \cup \ldots \cup T_n$ on V; this being a n-linear operator $T^* = T_1^* \cup \ldots \cup T_n^*$ such that $(T\alpha \mid \beta) = (\alpha \mid T^*\beta)$ i.e., $(T_1\alpha_1 \mid \beta_1) \cup \ldots \cup (T_n\alpha_n \mid \beta_n) = (\alpha_1 \mid T_1^* \beta_1) \cup \ldots \cup (\alpha^n \mid T_n^* \beta_n)$ for all α, β in V.

By the use of an n-orthonormal n-basis the n-adjoint operation on n-linear operators (passing from T to $T^*$) is identified with the operation of forming the n-conjugate transpose of a n-matrix.

Suppose $V = V_1 \cup \ldots \cup V_n$ is an n-inner product space over the n-field $F = F_1 \cup \ldots \cup F_n$ and let $\beta = \beta_1 \cup \ldots \cup \beta_n$ be some fixed n-vector in V. We define a n-function $f_\beta$ from V into the scalar n-field by $f_\beta(\alpha) = (\alpha|\beta)$;

i.e. $f_{1\beta_1}(\alpha_1) \cup \ldots \cup f_{n\beta_n}(\alpha_n) = (\alpha_1 \mid \beta_1) \cup \ldots \cup (\alpha_n \mid \beta_n)$

for $\alpha = \alpha_1 \cup \ldots \cup \alpha_n \in V = V_1 \cup \ldots \cup V_n$.

The n-function $f_\beta$ is a n-linear functional on V because by its very definition $(\alpha|\beta)$ is a n-linear n-function on V, arises in this way from some $\beta = \beta_1 \cup \ldots \cup \beta_n \in V$.

**THEOREM 2.6**: *Let V be a finite $(n_1, \ldots, n_n)$ dimensional n-inner product space and $f = f_1 \cup \ldots \cup f_n$, a n-linear functional on $V = V_1 \cup \ldots \cup V_n$. Then there exists a unique n-vector $\beta = \beta_1 \cup \ldots \cup \beta_n$ in V such that $f(\alpha) = (\alpha \mid \beta)$ i.e. $f_1(\alpha_1) \cup \ldots \cup f_n(\alpha_n) = (\alpha_1 \mid \beta_1) \cup \ldots \cup (\alpha_n|\beta_n)$ for all α in V.*

*Proof:* Let $\{\alpha_1^1 \ldots \alpha_{n_1}^1\} \cup \ldots \cup \{\alpha_1^n \ldots \alpha_{n_n}^n\}$ be a n-orthonormal n-basis for V. Put

$$\beta = \beta_1 \cup \ldots \cup \beta_n = \sum_{j_1=1}^{n_1} f_1(\alpha_{j_1}^1)\alpha_{j_1}^1 \cup \ldots \cup \sum_{j_n=1}^{n_n} f_n(\alpha_{j_n}^n)\alpha_{j_n}^n$$

and let $f_\beta = f_{1\beta_1} \cup \ldots \cup f_{n\beta_n}$ be the n-linear functional defined by $f_\beta(\alpha) = (\alpha|\beta)$; $f_{1\beta_1}(\alpha_1) \cup \ldots \cup f_{n\beta_n}(\alpha_n) = (\alpha_1|\beta_1) \cup \ldots \cup (\alpha_n|\beta_n)$; then

$f_\beta(\alpha_K) = (\alpha_K \mid \sum_j f(\alpha_j)\alpha_j) = f(\alpha_K)$ i.e., if $\alpha_K = (\alpha_K^1 \cup \ldots \cup \alpha_K^n)$

then

$$f_\beta(\alpha_K) = f_{1\beta_1}(\alpha_K^1) \cup \ldots \cup f_{n\beta_n}(\alpha_K^n)$$

$$= \left(\alpha_K^1 \mid \sum_{j_1} f_1(\alpha_{j_1}^1)\alpha_{j_1}^1\right) \cup \ldots \cup \left(\alpha_K^n \mid \sum_{j_n} f_n(\alpha_{j_n}^n)\alpha_{j_n}^n\right)$$
$$= f_1(\alpha_K^1) \cup \ldots \cup f_n(\alpha_K^n).$$

Since this is true for each $\alpha_K = \alpha_K^1 \cup \ldots \cup \alpha_K^n$, it follows that f = $(f_\beta)$. Now suppose $\gamma = \gamma_1 \cup \ldots \cup \gamma_n$ is a n-vector in V such that $(\alpha|\beta) = (\alpha|\gamma)$ for all $\alpha \in$ V; i.e.,$(\alpha_1|\beta_1) \cup \ldots \cup (\alpha_n \mid \beta_n) = (\alpha_1 \mid \gamma_1) \cup \ldots \cup (\alpha_n \mid \gamma_n)$. Then $(\beta - \gamma \mid \beta - \gamma) = 0$ and $\beta = \gamma$. Thus there is exactly one n-vector $\beta = \beta_1 \cup \ldots \cup \beta_n$ determining the n-linear functional f in the stated manner.

**THEOREM 2.7**: *For any n-linear operator $T = T_1 \cup \ldots \cup T_n$ on a finite $(n_1, \ldots, n_n)$ dimensional n-inner product space $V = V_1 \cup \ldots \cup V_n$, there exists a unique n-linear operator $T^* = T_1^* \cup \ldots \cup T_n^*$ on V such that $(T\alpha|\beta) = (\alpha \mid T^*\beta)$ for all α, β in V.*

*Proof:* Let $\beta = \beta_1 \cup \ldots \cup \beta_n$ be any n-vector in V. Then $\alpha \to (T\alpha|\beta)$ is a n-linear functional on V. By the earlier results there is a unique n-vector $\beta' = \beta_1' \cup \ldots \cup \beta_n'$ in V such that $(T\alpha|\beta) = (\alpha|\beta')$ for every α in V. Let $T^*$ denote the mapping $\beta \to \beta'$; $\beta' = T^*\beta$. We have $(T_1 \alpha_1| \beta_1) \cup \ldots \cup (T_n\alpha_n \mid \beta_n) = (\alpha_1^1 |T_1^*\beta_1) \cup \ldots \cup (\alpha_n \mid T_n^*\beta_n)$, so must verify that $T^* = T_1^* \cup \ldots \cup T_n^*$ is an n-linear operator. Let β, γ be in V and $c = c_1 \cup \ldots \cup c_n$ be a n-scalar. Then for any α,

$$\begin{aligned}
(\alpha \mid T^*(c\beta + \gamma)) &= (T\alpha \mid c\beta + \gamma) \\
&= (T\alpha|c\beta) + (T\alpha|\gamma) \\
&= c(T\alpha|\beta) + (T\alpha \mid \gamma) \\
&= c(\alpha|T^*\beta) + (\alpha \mid T^*\gamma) \\
&= (\alpha \mid cT^*\beta) + (\alpha \mid T^*\gamma) \\
&= (\alpha \mid cT^*\beta + T^*\gamma)
\end{aligned}$$

i.e. $(\alpha_1 \mid T_1^* (c_1\beta_1 + \gamma_1)) \cup \ldots \cup (\alpha_n \mid T_n^*(c_n\beta_n + \gamma_n)) =$
$(\alpha_1 \mid c_1 T_1^* \beta_1 + T_1^*\gamma_1) \cup \ldots \cup (\alpha_n \mid c_n T_n^* \beta_n + T_n^*\gamma_n)$.

The uniqueness of $T^*$ is clear. For any β in V the n-vector $T^*\beta$ is uniquely determined as the vector β′ such that (Tα|β) = (α|β′) for every α.

**THEOREM 2.8**: *Let $V = V_1 \cup ... \cup V_n$ be a finite ($n_1$, ..., $n_n$) dimensional n-inner product space and let $B = \{\alpha_1^1 \ldots \alpha_{n_1}^1\} \cup ... \cup \{\alpha_1^n \ldots \alpha_{n_n}^n\}$ be an n-ordered n-orthogonal basis for V. Let $T = T_1 \cup ... \cup T_n$ be a n-linear operator on V and let $A = A_1 \cup ... \cup A_n$ be the n-matrix of T in the ordered n-basis B. Then $A_{Kj} = (T\alpha_j| \alpha_K)$ i.e.*

$$A_{K_1 j_1}^1 \cup \ldots \cup A_{K_n j_n}^n = (T_1\alpha_{j_1}^1 \mid \alpha_{K_1}^1) \cup \ldots \cup (T_n\alpha_{j_n}^n \mid \alpha_{K_n}^n).$$

*Proof:* Since $B = \{\alpha_1^1 \ldots \alpha_{n_1}^1\} \cup \ldots \cup \{\alpha_1^n \ldots \alpha_{n_n}^n\}$ is an ordered n-basis we have

$$\alpha = \sum_{K_1=1}^{n_1} (\alpha_1 \mid \alpha_{K_1}^1)\, \alpha_K^1 \cup \ldots \cup \sum_{K_n=1}^{n_n} (\alpha_n \mid \alpha_{K_n}^n)\, \alpha_K^n$$

$$= \alpha_1 \cup \ldots \cup \alpha_n.$$

The n-matrix A is defined by

$$T\alpha_j = \sum_{K=1}^{n} A_{Kj}\,\alpha_K \text{ , i.e.,}$$

$$T_1\,\alpha_{j_1}^1 \cup \ldots \cup T_n \alpha_{j_n}^n = \sum_{K_1=1}^{n_1} A_{K_1 j_1}^1 \alpha_{K_1}^1 \cup \ldots \cup \sum_{K_n=1}^{n_n} A_{K_n j_n}^n \alpha_{K_n}^n$$

since

$$T\alpha_j = \sum_{K=1}^{n} (T\alpha_j \mid \alpha_K)\alpha_K \;;$$

i.e.,

$$T_1\alpha_{j_1}^1 \cup \ldots \cup T_n\alpha_{j_n}^n = \sum_{K_1=1}^{n_1} (T_1\alpha_{j_1}^1 \mid \alpha_{K_1}^1)\alpha_{K_1}^1$$

$$\cup \ldots \cup \sum_{K_n=1}^{n_n} (T_n\alpha_{j_n}^n \mid \alpha_{K_n}^n)\alpha_{K_n}^n$$

we have $A_{Kj} = (T\alpha_j \mid \alpha_K)$;

$$A_{K_1 J_1}^1 \cup \ldots \cup A_{K_n J_n}^n = (T_1\alpha_{j1}^1 \mid \alpha_{K_1}^1) \cup \ldots \cup (T_n\alpha_{jn}^n \mid \alpha_{K_n}^1).$$

The following corollary is immediate and left for the reader to prove.

**COROLLARY 2.5:** *Let V be a finite ($n_1$, ..., $n_n$) dimensional n-inner product space over the n-field F and set $T = T_1 \cup ... \cup T_n$ a n-linear operator on V. In any n-orthogonal basis for V, the n-matrix of $T^*$ is the n-conjugate transpose of the n-matrix of T.*

**DEFINITION 2.5:** *Let $T = T_1 \cup ... \cup T_n$ be a n-linear operator on an inner product space $V = V_1 \cup ... \cup V_n$, then we say that $T = T_1 \cup ... \cup T_n$ has an n-adjoint on V if there exists a n-linear operator $T^* = T_1^* \cup ... \cup T_n^*$ on V such that $(T\alpha|\beta) = (\alpha|T^*\beta)$ i.e. $(T_1\alpha_1|\beta_1) \cup ... \cup (T_n\alpha_n \mid \beta_n) = (\alpha_1 \mid T_1^* \beta_1) \cup ... \cup (\alpha_n \mid T_n^* \beta_n)$ for all $\alpha = \alpha_1 \cup ... \cup \alpha_n$ and $\beta = \beta_1 \cup ... \cup \beta_n$ in $V = V_1 \cup ... \cup V_n$.*

It is left for the reader to prove the following theorem.

**THEOREM 2.9**: *Let $V = V_1 \cup ... \cup V_n$ be a finite dimensional n-inner product space over the n-field $F = F_1 \cup ... \cup F_n$. If T and U are n-linear operators on V and c is n-scalar.*

i. $(T + U)^* = T^* + U^*$
   i.e., $(T_1 + U_1)^* \cup ... \cup (T_n + U_n)^*$
   $= T_1^* + U_1^* \cup ... \cup (T_n^* + U_n^*)$.
ii. $(cT)^* = cT^*$.
iii. $(TU)^* = U^*T^*$.
iv. $(T^*)^* = T$.

**DEFINITION 2.6:** *Let $V = V_1 \cup ... \cup V_n$ and $W = W_1 \cup ... \cup W_n$ be n-inner product spaces over the same n-field $F = F_1 \cup ... \cup F_n$ and let $T = T_1^* \cup ... \cup T_n^*$ be a n-linear transformation from V into W. We say that T-n-preserves inner products if $(T\alpha|T\beta) = (\alpha|\beta)$ i.e., $(T_1\alpha_1 \mid T_1\beta_1) \cup ... \cup (T_n\alpha_n \mid T_n\beta_n) = (\alpha_1 \mid \beta_1) \cup ... \cup (\alpha_n|\beta_n)$ for all α, β in V. An n-isomorphism of V onto W is a n-vector space n-isomorphism T of V onto W which also preserves inner products.*

**THEOREM 2.10:** *Let V and W be finite dimensional n-inner product spaces over the same n-field $F = F_1 \cup \ldots \cup F_n$. Both V and W are of same n-dimension equal to $(n_1,\ldots,n_n)$. If $T = T_1 \cup \ldots \cup T_n$ is a n-linear transformation from V into W, the following are equivalent.*

i. *$T = T_1 \cup \ldots \cup T_n$ preserves n-inner products.*
ii. *$T = T_1 \cup \ldots \cup T_n$ is an n-isomorphism.*
iii. *$T = T_1 \cup \ldots \cup T_n$ carries every n-orthonormal n-basis for V into an n-orthonormal n-basis for W.*
iv. *T carries some n-orthonormal n-basis for V onto an n-orthonormal n-basis for W.*

*Proof:* Clearly (i) → (ii) i.e., if $T = T_1 \cup \ldots \cup T_n$ preserves n-inner products, then $|T\alpha|| = ||\alpha||$ for all $\alpha = \alpha_1 \cup \ldots \cup \alpha_n$ in V, i.e., $|| T_1\alpha_1|| \cup \ldots \cup || T_n \alpha_n || = || \alpha_1 || \cup \ldots \cup || \alpha_n ||$. Thus T is n-non singular and since n-dim V = n-dim W = $(n_1, \ldots, n_n)$ we know that $T = T_1 \cup \ldots \cup T_n$ is a n-vector space n-isomorphism.

(ii) → (iii) Suppose $T = T_1 \cup \ldots \cup T_n$ is an n-isomorphism. Let $\{\alpha_1^1 \ldots \alpha_{n_1}^1\} \cup \ldots \cup \{\alpha_1^n \ldots \alpha_{n_n}^n\}$ be an n-orthonormal basis for V. Since $T = T_1 \cup \ldots \cup T_n$ is a n-vector space isomorphism and n-dim V = n-dim W, it follows that $\{T_1\alpha_1^1, \ldots, T_n\alpha_{n_1}^1\} \cup \{T_2\alpha_2^1, \ldots, T_2\alpha_{n_2}^2\} \cup \ldots \cup \{T_n\alpha_1^n, \ldots, T_n\alpha_{n_n}^n\}$ is a n-basis for W. Since T also n-preserves inner products $(T\alpha_j \mid T\alpha_K) = (\alpha_j \mid \alpha_K) = \delta_{jK}$, i.e.,

$$(T_1\ \alpha_{j_1}^1 \mid T_1\alpha_{K_1}^1) \cup \ldots \cup (T_n\alpha_{j_n}^1 \mid T_n\alpha_{K_n}^n)$$
$$(\alpha_{j_1}^1 \mid \alpha_{K_1}^1) \cup \ldots \cup (\alpha_{j_n}^n \mid \alpha_{K_n}^n) = \delta_{j_1K_1} \cup \ldots \cup \delta_{j_nK_n}.$$

(iii) → (iv) is obvious.

(iv) → (i).

Let $\{\alpha_1^1 \ldots \alpha_{n_1}^1\} \cup \ldots \cup \{\alpha_1^n \ldots \alpha_{n_n}^n\}$ be an n-orthonormal basis for V such that $\{T_1\alpha_1^1, \ldots, T_1\alpha_{n_1}^1\} \cup \ldots \cup \{T_n\alpha_n^n, \ldots, T_n\alpha_{n_n}^n\}$ is an n-orthonormal basis for W. Then $(T\alpha_j \mid T\alpha_K) = (\alpha_j \mid \alpha_K) = \delta_{jK}$;

i.e., $\{T_1\alpha^1_{j_1} \mid T_1\alpha^1_{K_1}\} \cup \ldots \cup \{T_n\alpha^n_{j_n} \mid T_n\alpha^n_{K_n}\}$

$= (\alpha^1_{j_1} \mid \alpha^1_{K_1}) \cup \ldots \cup (\alpha^n_{j_n} \mid \alpha^n_{K_n})$

$= \delta_{j_1K_1} \cup \ldots \cup \delta_{j_nK_n}$.

For any

$\alpha = \alpha_1 \cup \ldots \cup \alpha_n$

$= x^1_1\alpha^1_1 + \ldots + x^1_{n_1}\alpha^1_{n_1} \cup \ldots \cup x^n_1\alpha^n_1 + \ldots + x^n_{n_n}\alpha^n_{n_n}$

and

$\beta = y^1_1\alpha^1_1 + \ldots + y^1_{n_1}\alpha^1_{n_1} \cup \ldots \cup y^n_1\alpha^n_1 + \ldots + y^n_{n_n}\alpha^n_{n_n}$

in V we have

$$(\alpha|\beta) = \sum_{j=1}^{n} x_j y_j$$

that is

$$[(\alpha_1 \mid \beta_1) \cup \ldots \cup (\alpha_n \mid \beta_n)] = \sum_{j_1=1}^{n_1} x^1_{j_1} y^1_{j_1} \cup \ldots \cup \sum_{j_n=1}^{n_n} x^n_{j_n} y^n_{j_n} .$$

$(T\alpha \mid T\beta) = (T_1\alpha_1 \mid T_1\beta_1) \cup \ldots \cup (T_n\alpha_n \mid T_n\beta_n)$

$$= \left(\sum_j x_j T\alpha_j \,\middle|\, \sum_K y_K T\alpha_K\right)$$

$$= \left(\sum_{j_1} x^1_{j_1} T_1\alpha^1_{j_1} \,\middle|\, \sum_{K_1} y^1_{K_1} T_1\alpha^1_{K_1}\right) \cup \ldots \cup \left(\sum_{j_n} x^n_{j_n} T_n\alpha^n_{j_n} \,\middle|\, \sum_{K_n} y^n_{K_n} T_n\alpha^n_{K_n}\right)$$

$$= \sum_j \sum_K x_j y_K (T\alpha_j \mid T\alpha_K)$$

$$= \sum_{j_1} \sum_{K_1} x_{j_1} y_{K_1} (T_1\alpha^1_{j_1} \mid T_1\alpha^1_{K_1}) \cup \ldots \cup \sum_{j_n} \sum_{K_n} x^n_{j_n} y^n_{K_n} (T_n\alpha^n_{j_n} \mid T_n\alpha^n_{K_n})$$

$$= \sum_{j=1}^{n} x_j y_j$$

$$= \sum_{j_1=1}^{n} x^1_{j_1} y^1_{j_1} \cup \ldots \cup \sum_{j_n=1}^{n_n} x^n_{j_n} y^n_{j_n}$$

and so T-n-preserves n-inner products.

**COROLLARY 2.6:** *Let $V = V_1 \cup \ldots \cup V_n$ and $W = W_1 \cup \ldots \cup W_n$ be finite dimensional n-inner product spaces over the same n-field $F = F_1 \cup \ldots \cup F_n$. Then V and W are n-isomorphic if and only if they have the same n-dimension.*

*Proof:* If $\{\alpha^1_1 \ldots \alpha^1_{n_1}\} \cup \ldots \cup \{\alpha^n_1 \ldots \alpha^n_{n_n}\}$ is a n-orthonormal n-basis for $V = V_1 \cup \ldots \cup V_n$ and $\{\beta^1_1 \ldots \beta^1_{n_1}\} \cup \ldots \cup \{\beta^n_1 \ldots \beta^n_{n_n}\}$ is an n-orthonormal n-basis for W; let T be the n-linear transformation from V into W defined by $T\alpha_j = \beta_j$ i.e. $T_1\alpha^1_{j_1} \cup \ldots \cup T_n\alpha^n_{j_n} = \beta^1_{j_1} \cup \ldots \cup \beta^n_{j_n}$. Then T is an n-isomorphism of V onto W.

**THEOREM 2.11:** *Let $V = V_1 \cup \ldots \cup V_n$ and $W = W_1 \cup \ldots \cup W_n$ be two n-inner product spaces over the same n-field and let $T = T_1 \cup \ldots \cup T_n$ be a n-linear transformation from V into W. Then T-n-preserves n-inner products if and only if $||T\alpha|| = ||\alpha||$ for every $\alpha$ in V.*

*Proof:* If $T = T_1 \cup \ldots \cup T_n$, n-preserves inner products then T-n-preserves norms. Suppose $||T\alpha|| = ||\ \alpha\ ||$ for every $\alpha = \alpha_1 \cup \ldots \cup \alpha_n$ in V; i.e., $||T_1\alpha_1\ || \cup \ldots \cup ||T_n\alpha_n|| = ||\ \alpha_1\ || \cup \ldots \cup ||\ \alpha_n\ ||$. Now using the appropriate polarization identity for real space and the fact T is n-linear we see $(\alpha \mid \beta) = (T\alpha \mid T\beta)$ i.e. $(\alpha_1|\beta_1) \cup \ldots \cup (\alpha_n|\beta_n) = (T_1\alpha_1 \mid T_1\beta_1) \cup \ldots \cup (T_n\alpha_n \mid T_n\beta_n)$ for all $\alpha = \alpha_1 \cup \ldots \cup \alpha_n$ and $\beta = \beta_1 \cup \ldots \cup \beta_n$ in V.

Recall a unitary operator on an inner product space is an isomorphism of the space onto itself.

A n-unitary operator on an n-inner product space is an n-isomorphism of the n-space V onto itself.

It is left for the reader to verify that the product of two n-unitary operator is n-unitary.

**THEOREM 2.12:** *Let $U = U_1 \cup \ldots \cup U_n$ be a n-linear operator on an n-inner product space $V = V_1 \cup \ldots \cup V_n$. Then V is n-unitary if and only if the n-adjoint $U^*$ of U exists and $UU^* = U^*U = I$.*

*Proof:* Suppose $U = U_1 \cup \ldots \cup U_n$ is n-unitary. Then U is n-invertible and $(U\alpha \mid \beta) = (U\alpha|UU^{-1}\beta) = (\alpha \mid U^{-1} \beta)$ for all $\alpha = \alpha_1 \cup \ldots \cup \alpha_n$ and $\beta = \beta_1 \cup \ldots \cup \beta_n$ in V. Hence $U^{-1}$ is the n-adjoint of U for

$$(U_1\alpha_1 \mid \beta_1) \cup \ldots \cup (U_n \beta_n \mid \beta_n)$$
$$= (U_1\alpha_1 \mid U_1 U_1^{-1} \beta_1) \cup \ldots \cup (U_n\alpha_n \mid U_n U_n^{-1} \beta_n)$$
$$= (\alpha_1 \mid U_1^{-1} \beta_1) \cup \ldots \cup (\alpha_n \mid U_n^{-1} \beta_n).$$

Conversely suppose $U^*$ exists and $UU^* = U^*U = I$ i.e.,

$$U_1 U_1^* \cup \ldots \cup U_n U_n^*$$
$$= U_1^* U_1 \cup \ldots \cup U_n^* U_n$$
$$= I_1 \cup \ldots \cup I_n.$$

Then $U = U_1 \cup \ldots \cup U_n$ is n-invertible with $U^{-1} = U^*$; i.e., $U_1^{-1} \cup \ldots \cup U_n^{-1} = U_1^* \cup \ldots \cup U_n^*$. So we need only show that $U = U_1 \cup \ldots \cup U_n$ preserves n-inner products. We have $(U\alpha|U\beta) = (\alpha \mid U^*U\beta) = (\alpha|I\beta) = (\alpha \mid \beta)$;
i.e.,

$$(U_1\alpha_1 \mid U_1\beta_1) \cup \ldots \cup (U_n\alpha_n \mid U_n\beta_n)$$
$$= (\alpha_1| U_1^* U_1\beta_1) \cup \ldots \cup (\alpha_n| U_n^* U_n\beta_n)$$
$$= (\alpha_1|I_1\beta_1) \cup \ldots \cup (\alpha_n|I_n\beta_n) = (\alpha_1|\beta_1) \cup \ldots \cup (\alpha_n|\beta_n)$$

for all $\alpha, \beta \in V$.

We call a real n-mixed square matrix $A = A_1 \cup \ldots \cup A_n$ over the n-field $F = F_1 \cup \ldots \cup F_n$ to be n-orthogonal if $A^tA = I$ i.e., $A_1^t A_1 \cup \ldots \cup A_n^t A_n = I_1 \cup \ldots \cup I_n$. We say $A = A_1 \cup \ldots \cup$ $A_n$ to be n-anti orthogonal if $A^tA = -I$, i.e., $A_1^t A_1 \cup \ldots \cup A_n^t A_n = -I_1 \cup -I_2 \cup \ldots \cup -I_n$.

Let $V = V_1 \cup \ldots \cup V_n$ be a finite dimensional n-inner product space and $T = T_1 \cup \ldots \cup T_n$ be a n-linear operator on V. We say T is n-normal if it n-commutes with its n-adjoint i.e. $TT^* = T^*T$, i.e., $T_1 T_1^* \cup \ldots \cup T_n T_n^* = T_1^* T \cup \ldots \cup T_n^* T_n$.

The following theorem speaks about the properties enjoyed by n-self adjoint n-operators on a n-vector space over the n-field of type II.

**THEOREM 2.13:** *Let $V = V_1 \cup \ldots \cup V_n$ be an n-inner product space and $T = T_1 \cup \ldots \cup T_n$ be a n-self adjoint operator on V. Then each n-characteristic value of T is real and n-characteristic vectors of T associated with distinct n-characteristic values are n-orthogonal.*

*Proof:* Suppose $c = c_1 \cup \ldots \cup c_n$ is a n-characteristic value of T i.e., $T\alpha = c\alpha$ for some nonzero n-vector $\alpha = \alpha_1 \cup \ldots \cup \alpha_n$ i.e., $T_1\alpha_1 \cup \ldots \cup T_n\alpha_m = c_1\alpha_1 \cup \ldots \cup c_n\alpha_n$.

$$\begin{aligned} c(\alpha|\alpha) &= (c\alpha|\alpha) = (T\alpha \mid \alpha) \\ &= (\alpha \mid T\alpha) = (\alpha \mid c\alpha) \\ &= \bar{c}\,(\alpha|\alpha) \end{aligned}$$

i.e.,

$$\begin{aligned} c(\alpha|\alpha) &= c_1 (\alpha_1 |\alpha_1) \cup \ldots \cup c_n(\alpha|\alpha_n) \\ &= (c_1\alpha_1 \mid \alpha_1) \cup \ldots \cup (c_n\alpha_n|\beta_n) \\ &= (T_1\alpha_1 \mid \alpha_1) \cup \ldots \cup (T_n\alpha_n|\alpha_n) \\ &= (\alpha_1 \mid T_1\alpha_1) \cup \ldots \cup (\alpha_n|T_n\alpha_n) \\ &= (\alpha_1 \mid c_1\alpha_1) \cup \ldots \cup (\alpha_n|c_n\alpha_n) \\ &= \bar{c}_1 (\alpha_1|\alpha_1) \cup \ldots \cup \bar{c}_n (\alpha_n \mid \alpha_n). \end{aligned}$$

Since $(\alpha \mid \alpha) \neq 0 \cup \ldots \cup 0$; i.e., $(\alpha_1 \mid \alpha_1) \cup \ldots \cup (\alpha_n \mid \alpha_n) \neq (0 \cup \ldots \cup 0)$. We have $c = \bar{c}$, i.e., $c_i = \bar{c}_i$ for i = 1, 2, …, n. Suppose we also have $T\beta = d\beta$ with $\beta \neq 0$ i.e., $\beta = \beta_1 \cup \ldots \cup \beta_n \neq 0 \cup \ldots \cup 0$.

$$(c\alpha|\beta) = (T\alpha|\beta) = (\alpha \mid T\beta)$$

$= (\alpha \mid d\beta) \quad = \quad d(\alpha \mid \beta)$
$= d(\alpha \mid \beta)$.

If $c \neq d$ then $(\alpha \mid \beta) = 0 \cup \ldots \cup 0$; i.e., if $d_1 \cup \ldots \cup d_n \neq c_1 \cup \ldots \cup c_n$ i.e., $c_i \neq d_i$ for i = 1, 2, …, n then $(\alpha_1 \mid \beta_1) \cup \ldots \cup (\alpha_n \mid \beta_n) = 0 \cup \ldots \cup 0$.

The reader is advised to derive more properties in this direction. We derive the spectral theorem for n-inner product n-vector space over the n-real field $F = F_1 \cup \ldots \cup F_n$.

**THEOREM 2.14:** *Let $T = T_1 \cup \ldots \cup T_n$ be a n-self adjoint operator on a finite ($n_1$, ..., $n_n$) dimensional n-inner product vector space $V = V_1 \cup \ldots \cup V_n$. Let*

$$\{c_1^1 \ldots c_{K_1}^1\} \cup \ldots \cup \{c_1^n \ldots c_{K_n}^n\}$$

*be the n-distinct n-characteristic values of T. Let $W_j = W_{j_1}^1 \cup \ldots \cup W_{j_n}^n$ be a n-characteristic space associated with n-scalar $c_j = c_{j_1}^1 \cup \ldots \cup c_{j_n}^n$ and $E_j = E_{j_1}^1 \cup \ldots \cup E_{j_n}^n$ be the n-orthogonal projection of V on $W_j$ i.e., $V_t$ on $W_j^t$ for t = 1, 2, ..., n. Then $W_i$ is n-orthogonal to $W_j$ if $i \neq j$ i.e., if $W_i = W_{i_1}^1 \cup \ldots \cup W_{i_n}^n$ and $W_j = W_{j_1}^1 \cup \ldots \cup W_{j_n}^n$ then $W_{i_t}^t$ is orthogonal to $W_{j_t}^t$ if $i_t \neq j_t$ for t = 1, 2, ..., n. V is the n-direct sum of $\{W_1^1, \ldots, W_{K_1}^1\} \cup \ldots \cup \{W_1^n, \ldots, W_{K_n}^n\}$ and*

$$T = T_1 \cup \ldots \cup T_n$$
$$= [c_1^1 E_1^1 + \ldots + c_1^{K_1} E_1^{K_1}] \cup \ldots \cup [c_n^1 E_n^1 + \ldots + c_n^{K_n} E_n^{K_n}] .$$

*Proof:* Let $\alpha = \alpha_1 \cup \ldots \cup \alpha_n$ be a n-vector in

$$W_j = W_{j_1}^1 \cup \ldots \cup W_{j_n}^n .$$

$\beta = \beta_1 \cup \ldots \cup \beta_n$ be a n-vector in

$$W_i = W_{i_1}^1 \cup \ldots \cup W_{i_n}^n$$

and suppose $i \neq j$; $c_j (\alpha|\beta) = (T\alpha|\beta) = (\alpha|T^*\beta) = (\alpha|c_i\beta)$ i.e.,

$$c_{j_1}^1 (\alpha_1 \mid \beta_1) \cup \ldots \cup c_{j_n}^n (\alpha_n \mid \beta_n)$$
$$= (\alpha_1 \mid c_{i_1}^1 \beta_1) \cup \ldots \cup (\alpha_n \mid c_{i_n}^n \beta_n) .$$

Hence $(c_j - c_i)(\alpha|\beta) = 0 \cup \ldots \cup 0$. i.e.,

$$(c^1_{j_1} - c^1_{i_1})(\alpha_1 \mid \beta_1) \cup \ldots \cup (c^n_{j_n} - c^n_{i_n})(\alpha_n \mid \beta_n) = 0 \cup \ldots \cup 0;$$

since $c_j - c_i \neq 0 \cup \ldots \cup 0$. It follows $(\alpha|\beta) = (\alpha_1 \mid \beta_1) \cup \ldots \cup (\alpha_n|\beta_n) = 0 \cup \ldots \cup 0$. Thus $W_j = W^1_{j_1} \cup \ldots \cup W^n_{j_n}$ is n-orthogonal to $W_i = W^1_{i_1} \cup \ldots \cup W^n_{i_n}$ when $i \neq j$. From the fact that V has an n-orthonormal basis consisting of n-characteristic, it follows that $V = W_1 + \ldots + W_K$ i.e.,

$$V = V_1 \cup \ldots \cup V_n$$

$$= (W^1_1 + \ldots + W^1_{K_1}) \cup \ldots \cup (W^n_1 + \ldots + W^n_{K_n})$$

when $i_t \neq j_t$ ; $1 \leq i_t, j_t \leq K_t$ or $t = 1, 2, \ldots, n$. From the fact that $V = V_1 \cup \ldots \cup V_n$ has an n-orthonormal n-basis consisting of n-characteristic n-vectors it follows that $V = W_1 + \ldots + W_K$, If $\alpha^t_{j_t} \in W^t_{j_t}$ $(1 \leq j_t \leq K_t)$ and $(\alpha^1_1 + \ldots + \alpha^1_{K_1}) \cup \ldots \cup (\alpha^n_1 + \ldots + \alpha^n_{K_n}) = 0 \cup \ldots \cup 0$ then $0 \cup \ldots \cup 0$

$$= \left( (\alpha_i \mid \sum_j \alpha_j) \right)$$

$$= \sum_j (\alpha_i \mid \alpha_j)$$

$$\text{i.e. } \left( (\alpha^1_{i_1} \mid \sum_{j_1} \alpha^1_{j_1}) \right) \cup \ldots \cup \left( (\alpha^n_{i_n} \mid \sum \alpha^n_{j_n}) \right)$$

$$= \sum_{j_1} (\alpha^1_{i_1} \mid \alpha^1_{j_1}) \cup \ldots \cup \sum_{j_n} (\alpha^n_{i_n} \mid \alpha^n_{j_n})$$

$$= \|\alpha^1_{i_1}\|^2 \cup \ldots \cup \|\alpha^n_{i_n}\|^2$$

for every $1 \leq i_t \leq K_t$; $t = 1, 2, \ldots, n$ so that V is the n-direct sum of $(W^1_1, \ldots, W^1_{K_1}) \cup \ldots \cup (W^n_1, \ldots, W^n_{K_n})$. Therefore

$$E^1_1 + \ldots + E^1_{K_1} \cup \ldots \cup E^n_1 + \ldots + E^n_{K_n} = I_1 \cup \ldots \cup I_n$$

and

$$T = (T_1 E^1_1 + \ldots + T_1 E^1_{K_1}) \cup \ldots \cup (T_n E^n_1 + \cdots + T_n E^n_{K_n}) =$$
$$(c^1_1 E^1_1 + \ldots + c^1_{K_1} E^1_{K_1}) \cup \ldots \cup (c^n_1 E^n_1 + \ldots + c^n_{K_n} E^n_{K_n}).$$

This decomposition is called the n-spectral resolution of T. The following corollary is immediate however we sketch the proof of it.

**COROLLARY 2.7:** *If*

$$e_j = e^1_{j_1} \cup \ldots \cup e^n_{j_n}$$

$$= \prod_{i_1 \neq j_1} \left( \frac{x - c^1_{i_1}}{c^1_{j_1} - c^1_{i_1}} \right) \cup \ldots \cup \prod_{i_n \neq j_n} \left( \frac{x - c^n_{i_n}}{c^n_{j_n} - c^n_{i_n}} \right)$$

*then* $E_j = E^1_{j_1} \cup \ldots \cup E^n_{j_n}$ *for* $1 \leq j_t \leq K_t$; $t = 1, 2, \ldots, n$.

*Proof:* Since $E^t_{i_t} E^t_{j_t} = 0$ for every $1 \leq i_t, j_t \leq K_t$; $t = 1, 2, \ldots, n$ ($i_t \neq j_t$) it follows that $T^2 = (T_1 \cup \ldots \cup T_n)^2 = T_1^2 \cup \ldots \cup T_n^2$.

$$(c^1_1)^2 E^1_1 + \ldots + (c^1_{K_1})^2 E^1_{K_1} \cup \ldots \cup (c^n_1)^2 E^n_1 + \ldots + (c^n_{K_n})^2 E^n_{K_n}$$

and by an easy induction argument that

$$T^n = \left[ (c^1_1)^n E^1_1 + \ldots + (c^1_{K_1})^n E^1_{K_1} \right] \cup \ldots \cup \left[ (c^n_1)^n E^n_1 + \ldots + (c^n_{K_n})^n E^n_{K_n} \right]$$

for every n-integer $(n_1, \ldots, n_n) \geq (0, 0, \ldots, 0)$. For an arbitrary n-polynomial $f = \sum_{n=0}^{r} a_n x^n$;

$$f = f_1 \cup \ldots \cup f_n = \sum_{n_1=0}^{r_1} a^1_{n_1} x^{n_1} \cup \ldots \cup \sum_{n_n=0}^{r_n} a^n_{n_n} x^{n_n},$$

we have

$$\begin{aligned} f(T) &= f_1(T_1) \cup \ldots \cup f_n(T_n) \\ &= \sum_{n_1=0}^{r_1} a^1_{n_1} T_1^{n_1} \cup \ldots \cup \sum_{n_n=0}^{r_n} a^n_{n_n} T_n^{n_n} \\ &= \sum_{n_1=0}^{r_1} a^1_{n_1} \sum_{j_1=1}^{K_1} c^{n_1}_{j_1} E^1_{j_1} \cup \ldots \cup \\ &\quad \sum_{n_n=0}^{r_n} a^n_{n_n} \sum_{j_n=1}^{K_n} c^{n_n}_{j_n} E^n_{j_n} \end{aligned}$$

$$= \sum_{j_1=1}^{K_1}\left(\sum_{n_1=0}^{r_1} a^1_{n_1} c^{n_1}_{j_1}\right)E^1_{j_1} \cup \ldots \cup \sum_{j_n=1}^{K_n}\left(\sum_{n_n=0}^{r_n} a^n_{n_n} c^{n_n}_{j_n}\right)E^n_{j_n}$$

$$= \sum_{j_1=1}^{K_1} f_1(c^1_{j_1})E^1_{j_1} \cup \ldots \cup \sum_{j_n=1}^{K_n} f_n(c^n_{j_n})E^n_{j_n}\,.$$

Since $e^1_{j_1}(c^t_{m_1}) = \delta_{j_1 m_1}$ for t = 1, 2, …, n it follows $e^t_{j_t}(T_t) = E^t_{j_t}$ true for t = 1, 2, …, n.

Since $\{E^1_1 \ldots E^1_{K_1}\} \cup \ldots \cup \{E^n_1 \ldots E^n_{K_n}\}$ are canonically associated with T and $I_1 \cup \ldots \cup I_n = (E^1_1 + \ldots + E^1_{K_1}) \cup \ldots \cup (E^n_1 + \ldots + E^n_{K_n})$ the family of n-projections $\{E^1_1 \ldots E^1_{K_1}\} \cup \ldots \cup \{E^n_1 \ldots E^n_{K_n}\}$ is called the n-resolution of the n-identity defined by $T = T_1 \cup \ldots \cup T_n$.

Next we proceed onto define the notion of n-diagonalized and the notion of n-diagonalizable normal n-operators.

**DEFINITION 2.7**: *Let $T = T_1 \cup \ldots \cup T_n$ be a n-diagonalizable normal n-operator on a finite dimensional n-inner product space and*

$$T = T_1 \cup \ldots \cup T_n = \sum_{j_1=1}^{K_1} c^1_{j_1} E^1_{j_1} \cup \ldots \cup \sum_{j_n=1}^{K_n} c^n_{j_n} E^n_{j_n}$$

*is its spectral n-resolution.*

*Suppose $f = f_1 \cup \ldots \cup f_n$ is a n-function whose n-domain includes the n-spectrum of T that has n-values in the n-field of scalars. Then the n-linear operator $f(T) = f_1(T_1) \cup \ldots \cup f_n(T_n)$ is defined by the n-equation*

$$f(T) = f_1(T_1) \cup \ldots \cup f_n(T_n) = \sum_{j_1=1}^{K_1} f_1(c^1_{j_1})E^1_{j_1} \cup \ldots \cup \sum_{j_n=1}^{K_n} f_n(c^n_{j_n})E^n_{j_n}\,.$$

Now we prove the following interesting theorem for n-diagonalizable normal n-operator.

**THEOREM 2.15:** *Let $T = T_1 \cup \ldots \cup T_n$ be a n-diagonalizable normal n-operator with n-spectrum $S = S_1 \cup \ldots \cup S_n$ on a finite $(n_1, \ldots, n_n)$ dimensional n-inner product space $V = V_1 \cup \ldots \cup V_n$. Suppose $f = f_1 \cup \ldots \cup f_n$ is a n-function whose domain contains S that has n-values in the field of n-scalars. Then f(T) is a n-diagonalizable normal n-operator with n-spectrum $f(S) = f_1(S_1) \cup \ldots \cup f_n(S_n)$. If $U = U_1 \cup \ldots \cup U_n$ is a n-unitary map of V onto V′ and $T' = UTU^{-1} = U_1T_1U_1^{-1} \cup \ldots \cup U_nT_nU_n^{-1}$ then $S = S_1 \cup \ldots \cup S_n$ is the n-spectrum of T′ and $f(T') = f_1(T'_1) \cup \ldots \cup f_n(T'_n) = Uf(T)U^{-1} = U_1f_1(T_1)U_1^{-1} \cup \ldots \cup U_nf_n(T_n)U_n^{-1}$.*

*Proof:* The n-normality of $f(T) = f_i(T_i) \cup \ldots \cup f_n(T_n)$ follows by a simple computation from the earlier results and the fact

$$\begin{aligned} f(T)^* &= f_1(T_1)^* \cup \ldots \cup f_n(T_n)^* \\ &= \sum_{j_1} f(c^1_{j_1})E^1_{j_1} \cup \ldots \cup \sum_{j_n} f(c^n_{j_n})E^n_{j_n}. \end{aligned}$$

Moreover it is clear that for every $\alpha = \alpha_1 \cup \ldots \cup \alpha_n$ in $E_j(V) = E^1_{j_1}(V_1) \cup \ldots \cup E^n_{j_n}(V_n)$; $f(T)\alpha = f(c_j)\alpha$, i.e., $f_1(T_1)\alpha_1 \cup \ldots \cup f_n(T_n)\alpha_n = f_1(c^1_{j_1})\,\alpha_1 \cup \ldots \cup f_n(c^n_{j_n})\alpha_n$. Thus the set f(S) of all f(c) with c in S is contained in the n-spectrum of f(T).

Conversely, suppose $\alpha = \alpha_1 \cup \ldots \cup \alpha_n = 0 \cup \ldots \cup 0$ and that $f(T)\alpha = b\alpha$ i.e., $f_1(T_1)\alpha_1 \cup \ldots \cup f_n(T_n)\alpha_n = b_1\alpha_1 \cup \ldots \cup b_n\alpha_n$.

Then $\alpha = \sum_j E_j\alpha$ i.e.,

$$\alpha_1 \cup \ldots \cup \alpha_n = \sum_j E^1_{j_1}\,\alpha_1 \cup \ldots \cup \sum_{j_n} E^n_{j_n}\alpha_n$$

and

$$\begin{aligned} f(T)\alpha &= f_1(T_1)\alpha_1 \cup \ldots \cup f_n(T_n)\alpha_n \\ &= \sum_{j_1} f_1(T_1)\,E^1_{j_1}\alpha_1 \cup \ldots \cup \sum_{j_n} f_n(T_n)\,E^n_{j_n}\alpha_n \\ &= \sum_{j_1} f_1\,(c^1_{j_1})E^1_{j_1}\alpha_1 \cup \ldots \cup \sum_{j_n} f_n\,(c^n_{j_n})E^n_{j_n}\alpha_n \\ &= \sum_{j_1} b_1\,E^1_{j_1}\alpha_1 \cup \ldots \cup \sum_{j_n} b_n\,E^n_{j_n}\alpha_n. \end{aligned}$$

Hence

$$\left\|\sum_j (f(c_j)-b)E_j\alpha\right\|^2$$

$$= \left\|\sum_{j_1}(f_1(c^1_{j_1})-b_1)E^1_{j_1}\alpha_1\right\|^2 \cup \ldots \cup \left\|\sum_{j_n}(f_n(c^n_{j_n})-b_n)E^n_{j_n}\alpha_n\right\|^2$$

$$= \sum_{j_1}\left|f_1(c^1_{j_1}-b_1)\right|^2 \left\|E^1_{j_1}\alpha_1\right\|^2 \cup \ldots \cup \sum_{j_n}\left|f_n(c^n_{j_n}-b_n)\right|^2 \left\|E^n_{j_n}\alpha_n\right\|^2.$$

Therefore $f(c_j) = f_1(c^1_{j_1}) \cup \ldots \cup f_n(c^n_{j_n}) = f_1 \cup \ldots \cup f_n$ or $E^1_{j_1}\alpha_1 \cup \ldots \cup E^n_{j_n}\alpha_n = 0 \cup \ldots \cup 0$; by assumption $\alpha = \alpha_1 \cup \ldots \cup \alpha_n = 0 \cup \ldots \cup 0$ so there exists an n-index tuple $i = (i_1, \ldots, i_n)$ such that $E_i\alpha = E_{i_1}\alpha_1 \cup \ldots \cup E_{i_n}\alpha_n = 0 \cup \ldots \cup 0$. By assumption $\alpha = \alpha_1 \cup \ldots \cup \alpha_n \neq 0 \cup \ldots \cup 0$. It follows that $f(c_i) = b$ i.e., $f_1(c^i_1) \cup \ldots \cup f_n(c^i_n) = b_1 \cup \ldots \cup b_n$ and hence that f(S) is the n-spectrum of $f(T) = f_1(T_1) \cup \ldots \cup f_n(T_n)$. Infact that $f(S) = \{b^1_1, \ldots, b^1_{r_1}\} \cup \ldots \cup \{b^n_1, \ldots, b^n_{r_n}\} = f_1(S_1) \cup \ldots \cup f_n(S_n)$ where $b^t_{m_t} \neq b^t_{n_t}$; t = 1, 2, …, n; where $m_t \neq n_t$. Let $X_m = X_{m_1} \cup \ldots \cup X_{m_n}$ be the set of indices $i = (i_1, \ldots, i_n)$ such that

$f(c_i) = f_1(c^1_{i_1}) \cup \ldots \cup f_n(c^n_{i_n}) = b_{m_1} \cup \ldots \cup b_{m_n}$.

Let

$$\begin{aligned} P_m &= P_{m_1} \cup \ldots \cup P_{m_n} \\ &= \sum_i E_i \\ &= \sum_{i_1} E^1_{i_1} \cup \ldots \cup \sum_{i_n} E^n_{i_n} \end{aligned}$$

the sum being extended over the n-indices $(i_1, \ldots, i_n)$ in $X_m = X^1_{m_1} \cup \ldots \cup X^n_{m_n}$. Then $P_m = P^1_{m_1} \cup \ldots \cup P^n_{m_n}$ is n-orthogonal projection of $V = V_1 \cup \ldots \cup V_n$ on the n-subspace of n-characteristic n-vectors belonging to the n-characteristic value $b_m$ of f(T) i.e., $b^1_{m_1} \cup \ldots \cup b^n_{m_n}$ of $f_1(T_1) \cup \ldots \cup f_n(T_n)$ and

$$f(T) = f_1(T_1) \cup \ldots \cup f_n(T_n) = \sum_{m_1=1}^{r_1} b^1_{m_1} P^1_{m_1} \cup \ldots \cup \sum_{m_n=1}^{r_n} b^n_{m_n} P^n_{m_n}$$

is the n-spectral resolution of $f(T) = f_1(T_1) \cup \ldots \cup f_n(T_n)$. Now suppose $U = U_1 \cup \ldots \cup U_n$ is a n-unitary transformation of V onto V' and that $T' = U T U^{-1}$, i.e., $T'_1 \cup \ldots \cup T'_n = U_1 T_1 U_1^{-1} \cup \ldots \cup U_n T_n U_n^{-1}$. Then the equation $T\alpha = c\alpha$ i.e., $T_1\alpha_1 \cup \ldots \cup T_n\alpha_n = c_1\alpha_1 \cup \ldots \cup c_n\alpha_n$ holds good if and only if $T'U\alpha = cU\alpha$ i.e.,

$T'_1 U_1\alpha_1 \cup \ldots \cup T'_n U_n\alpha_n = c_1U_1\alpha_1 \cup \ldots \cup c_nU_n\alpha_n$. Thus $S = S_1 \cup \ldots \cup S_n$ is the n-spectrum of $T' = T'_1 \cup \ldots \cup T'_n$ and $U = U_1 \cup \ldots \cup U_n$ maps each n-characteristic n-subspace for $T = T_1 \cup \ldots \cup T_n$ onto the corresponding n-subspace for T′. In fact using earlier results we see that

$$T' = \sum_j c_j E'_j$$

$$T'_1 \cup \ldots \cup T'_n = \sum_{j_1} c^1_{j_1} (E^1_{j_1})' \cup \ldots \cup \sum_{j_n} c^n_{j_n} (E^n_{j_n})'$$

Here

$$E' = (E^1_{j_1})' \cup \ldots \cup (E^n_{j_n})'.$$

$U_1 E^1_{j_1} U_1^{-1} \cup \ldots \cup U_n E^n_{j_n} U_n^{-1}$ is the n-spectral resolution of $T' = T'_1 \cup \ldots \cup T'_n$. Hence

$$\begin{aligned} f(T') &= f_1(T'_1) \cup \ldots \cup f_n(T'_n) \\ &= \sum_{j_1} f_1(c^1_{j_1})(E^1_{j_1})' \cup \ldots \cup \sum_{j_n} f_n(c^n_{j_n})(E^n_{j_n})' \\ &= \sum_{j_1} f_1(c^1_{j_1}) U_1 E^1_{j_1} U_1^{-1} \cup \ldots \cup \sum_{j_n} f_n(c^n_{j_n}) U_n E^n_{j_n} U_n^{-1} \end{aligned}$$

is the n-spectral resolution of $T' = T'_1 \cup \ldots \cup T'_n$.

Hence

$$f(T') = f_1(T'_1) \cup \ldots \cup f_n(T'_n)$$

$$= \sum_{j_1} f_1(c^1_{j_1}) U_1 E^1_{j_1} U_1^{-1} \cup \ldots \cup$$

$$\sum_{j_n} f_n(c^n_{j_n}) U_n E^n_{j_n} U_n^{-1}$$

$$= U_1 \sum_{j_1} f_1(c^1_{j_1})(E^1_{j_1}) U_1^{-1} \cup \ldots \cup U_n \sum_{j_n} f_n(c^n_{j_n})(E^n_{j_n}) U_n^{-1}$$

$$= U_1 f_1(T_1) U_1^{-1} \cup \ldots \cup U_n f_n(T_n) U_n^{-1}$$

$$= Uf(T) U^{-1}.$$

In view of the above theorem we have the following corollary.

**COROLLARY 2.8:** *With the assumption of the above theorem suppose $T = T_1 \cup \ldots \cup T_n$ is represented in an n-ordered basis $B = B_1 \cup \ldots \cup B_n = \{\alpha^1_1 \ldots \alpha^1_{n_1}\} \cup \ldots \cup \{\alpha^n_1 \ldots \alpha^n_{n_n}\}$ by the n-diagonal matrix $D = D_1 \cup \ldots \cup D_n$ with entries $\{d^1_1 \ldots d^1_{n_1}\} \cup \ldots \cup \{d^n_1 \ldots d^n_{n_n}\}$. Then in the n-basis B, $f(T) = f_1(T_1) \cup \ldots \cup f_n(T_n)$ is represented by the n-diagonal matrix $f(D) = f_1(D_1) \cup \ldots \cup f_n(D_n)$ with entries $\{f_1(d^1_1), \ldots, f_{n_1}(d^1_{n_1})\} \cup \ldots \cup \{f_n(d^n_1), \ldots, f_{n_n}(d^n_{n_n})\}$. If $B' = \{(\alpha^1_1)', \ldots, (\alpha^1_n)'\} \cup \ldots \cup \{(\alpha^n_1)', \ldots, (\alpha^n_n)'\}$ is any other n-ordered n-basis and P the n-matrix such that*

$$\beta^t_{j_t} = \sum_{i_t} P^t_{i_t j_t} \alpha^t_{i_t}$$

*for t = 1, 2, …, n then $P^{-1}(f(D))P$ is the n-basis B′.*

*Proof:* For each n-index $i = (i_1, \ldots, i_n)$ there is a unique n-tuple, $j = (j_1, \ldots, j_n)$ such that $1 \leq j_t, n_p \leq K_t$; t = 1, 2, …, n. $\alpha_i = \alpha^1_{i_1} \cup \ldots \cup \alpha^n_{i_n}$ belongs to $E_{j_t}(V_t)$ and $d^t_{i_t} = c^t_{j_t}$ for every t = 1, 2, …, n. Hence $f(T)\alpha_i = f(d_i)\alpha_i$ for every $i = (i_1, \ldots, i_n)$ i.e. $f_1(T_1)\alpha^1_{i_1} \cup \ldots \cup f_n(T_n)\alpha^n_{i_n} = f_1(d^1_{i_1})\alpha^1_{i_1} \cup \ldots \cup f_n(d^n_{i_n})\alpha^n_{i_n}$ and $f(T)\alpha'_j = f_1(T_1)(\alpha^1_{j_1})' \cup \ldots \cup f_n(T_n)(\alpha^n_{j_n})'$.

$$\sum_i P_{ij} f(T)\alpha_i = \sum_{i_1} P^1_{i_1 j_1} f_1(T_1)\alpha^1_{i_1} \cup \ldots \cup \sum_{i_n} P^n_{i_n j_n} f_n(T_n)\alpha^n_{i_n}$$

$$
\begin{aligned}
\sum_{i} d_i P_{ij}\, \alpha_i &= \sum_{i_1} d^1_{i_1} P^1_{i_1 j_1}\, \alpha^1_{i_1} \cup \ldots \cup \sum_{i_n} d^n_{i_n} P^n_{i_n j_n}\, \alpha^n_{i_n} \\
&= \sum_{i} (DP)_{ij}\, \alpha_i \\
&= \sum_{i_1} (D^1P^1)_{i_1 j_1}\, \alpha_{i_1} \cup \ldots \cup \sum (D^nP^n)_{i_n j_n}\, \alpha_{i_n} \\
&= \sum_{i} (DP)_{ij} \sum_{K} P^{-1}_{Kt} \alpha'_K \\
&= \sum_{i_1} (DP)_{i_1 j_1} \sum_{K_1} P^{-1}_{K_1 t_1} (\alpha^1_{K_1})' \cup \ldots \cup \\
&\quad \sum_{i_n} (DP)_{i_n j_n} \sum_{K_n} P^{-1}_{K_n t_n} (\alpha^n_{K_n})' \\
&= \sum_{K_1} (P_1^{-1} D_1 P_1)_{K_1 j_1} (\alpha^1_{K_1})' \cup \ldots \cup \\
&\sum_{K_n} (P_n^{-1} D_n P_n^{-1})_{K_n j_n} (\alpha^n_{K_n})' \; .
\end{aligned}
$$

Under the above conditions we have $f(A) = f_1(A_1) \cup \ldots \cup f_n(A_n) = P^{-1} f(D)\, P = P_1^{-1} f_1(D_1) P_1 \cup \ldots \cup P_n^{-1} f_n(D_n) P_n$.

The reader is expected to derive other interesting analogous results to usual vector spaces for the n-vector spaces of type II.

Now we proceed onto define the notion of bilinear n-forms, for n-vector spaces of type II.

**DEFINITION 2.8:** *Let $V = V_1 \cup \ldots \cup V_n$ be a n-vector space over the n-field $F = F_1 \cup \ldots \cup F_n$. A bilinear n-form on V is a n-function $f = f_1 \cup \ldots \cup f_n$ which assigns to each set of ordered pairs of n-vectors $\alpha$, $\beta$ in V a n-scalar $f(\alpha, \beta) = f_1(\alpha_1, \beta_1) \cup \ldots \cup f_n(\alpha_n, \beta_n)$ in $F = F_1 \cup \ldots \cup F_n$ and which satisfies:*

$$f(c\alpha^1 + \alpha^2, \beta) = cf(\alpha^1, \beta) + f(\alpha^2, \beta)$$

$$f_1(c_1 \alpha_1^1 + \alpha_1^2, \beta_1) \cup \ldots \cup f_n(c_n \alpha_1^n + \alpha_n^2, \beta_n)$$

$$= c_1 f_1(\alpha_1^1, \beta_1) + f_2(\alpha_1^2, \beta_1) \cup \ldots \cup c_n f_n(\alpha_n^1, \beta_n) + f_n(\alpha_n^2, \beta_n).$$

If we let $V \times V = (V_1 \times V_1) \cup \ldots \cup (V_n \times V_n)$ denote the set of all n-ordered pairs of n-vectors in $V_1$ this definition can be rephrased as follows. A bilinear n-form on V is a n-function f from $V \times V$ into $F = F_1 \cup \ldots \cup F_n$ i.e., $(V_1 \times V_1) \cup \ldots \cup (V_n \times V_n)$ into $F_1 \cup \ldots \cup F_n$ i.e., $f: V \times V \to F$ is $f_1 \cup \ldots \cup f_n : (V_1 \times V_1) \cup \ldots \cup (V_n \times V_n) \to F_1 \cup \ldots \cup F_n$ with $f_t: V_t \times V_t \to F_t$ for $t = 1, 2, \ldots, n$ which is a n-linear function of either of its arguments when the other is fixed. The n-zero function or zero n-function from $V \times V$ into F is clearly a bilinear n-form. It is also true that any bilinear combination of bilinear n-forms on V is again a bilinear n-form.

Thus all bilinear n-forms on V is a n-subspace of all n-functions from $V \times V$ into F. We shall denote the space of bilinear n-forms on V by $L^n (V, V, F) = L(V_1, V_1, F_1) \cup \ldots \cup L(V_n, V_n, F_n)$.

**DEFINITION 2.9:** *Let $V = V_1 \cup \ldots \cup V_n$ be a finite $(n_1, \ldots, n_n)$ dimensional n-vector space and let $B = \{\alpha_1^1, \ldots, \alpha_{n_1}^1\} \cup \ldots \cup \{\alpha_1^n, \ldots, \alpha_{n_n}^n\} = (B_1 \cup \ldots \cup B_n)$ be a n-ordered n-basis of V. If f $= f_1 \cup \ldots \cup f_n$ is a bilinear n-form on V then the n-matrix of f in the ordered n-basis B, is $(n_1 \times n_1, \ldots, n_n \times n_n)$ n-matrix $A = A_1 \cup \ldots \cup A_n$ with entries $A^t_{i_t j_t} = f_t(\alpha^t_{i_t}, \alpha^t_{j_t})$. At times we shall denote the n-matrix by $[f]_B = [f_1]_{B_1} \cup \ldots \cup [f_n]_{B_n}$.*

**THEOREM 2.16:** *Let $V = V_1 \cup \ldots \cup V_n$ be a $(n_1, \ldots, n_n)$ dimensional n-vector space over the n-field $F = F_1 \cup \ldots \cup F_n$. For each ordered n-basis $B = B_1 \cup \ldots \cup B_n$ of V, the n-function which associates with each bilinear n-form on V its n-matrix in the ordered n-basis $B = B_1 \cup \ldots \cup B_n$ is an n-isomorphism of the n-space $L^n(V, V, F) = L(V_1, V_1, F_1) \cup \ldots \cup L(V_n, V_n, F)$ onto the n-space of $(n_1 \times n_1, \ldots, n_n \times n_n)$ n-matrices over the n-field $F = F_1 \cup \ldots \cup F_n$.*

*Proof:* We see $f \to [f]_B$ i.e., $f_1 \cup \ldots \cup f_n \to [f_1]_{B_1} \cup \ldots \cup [f_n]_{B_n}$ where $f_t \to [f_t]_{B_t}$ for every $t = 1, 2, \ldots, n$ which is a one to one n-correspondence between the set of bilinear n-forms on V and the set of all $(n_1 \times n_1, \ldots, n_n \times n_n)$ matrices over the n-field F. That this is a n-linear transformation is easy to see because

$$(cf + g)(\alpha_i, \alpha_j) = cf(\alpha_i, \alpha_j) + cg(\alpha_i, \alpha_j)$$

i.e., $(c^1f_1 + g_1)\ (\alpha^1_{i_1}, \alpha^1_{j_1}) \cup \ldots \cup (c^nf_n + g_n)\ (\alpha^n_{i_n}, \alpha^n_{j_n})$

$$= [c^1 f_1 (\alpha^1_{i_1}, \alpha^1_{j_1}) + g_1 (\alpha^1_{i_1}, \alpha^1_{j_1})] \cup \ldots \cup [c^n f_n (\alpha^n_{i_n}, \alpha^n_{j_n}) + g_n (\alpha^n_{i_n}, \alpha^n_{j_n})].$$

This simply says $[cf + g]_B = c[f]_B + [g]_B$ i.e.,

$$[c^1 f_1 + g_1]_{B_1} \cup \ldots \cup [c^n f_n + g_n]_{B_n}$$
$$= (c^1[f_1]_{B_1} + [g_1]_{B_1}) \cup \ldots \cup (c^n[f_n]_{B_n} + [g_n]_{B_n}).$$

We leave the following corollary for the reader to prove.

**COROLLARY 2.9:** *If* $B = \{\alpha^1_1, \ldots, \alpha^1_{n_1}\} \cup \ldots \cup \{\alpha^n_1, \ldots, \alpha^n_{n_n}\}$ *is an n-ordered n-basis for* $V = V_1 \cup \ldots \cup V_n$ *and*

$$B^* = \{L^1_1, \ldots, L^1_{n_1}\} \cup \ldots \cup \{L^n_1, \ldots, L^n_{n_n}\}$$

*is the dual n-basis for* $V^* = V^*_1 \cup \ldots \cup V^*_n$ *, then* $(n^2_1, \ldots, n^2_n)$ *bilinear n-forms* $f_{ij}(\alpha, \beta) = L_i(\alpha)\, L_j(\beta)$*; i.e.,*

$$f^1_{i_1 j_1}(\alpha_1, \beta_1) \cup \ldots \cup f^n_{i_n j_n}(\alpha_n, \beta_n) =$$
$$L^1_{i_1}(\alpha_i) L^1_{j_1}(\beta_1) \cup \ldots \cup L^n_{i_n}(\alpha_n) L^n_{j_n}(\beta_n);$$

$1 \le i_t, j_t \le n_t$*;* $t = 1, 2, \ldots, n$*, forms a n-basis for the n-space* $L^n(V, V, F) = L(V_1, V_1, F_1) \cup \ldots \cup L(V_n, V_n, F_n)$ *is* $(n^2_1, n^2_2, \ldots, n^n_2)$.

**THEOREM 2.17:** *Let* $f = f_1 \cup \ldots \cup f_n$ *be a bilinear n-form on the finite dimensional n-vector space V of n-dimension* $(n_1, n_2, \ldots, n_n)$*. Let* $R_f$ *and* $L_f$ *be a n-linear transformation from V into* $V^*$ *defined by* $(L_f\alpha)\, \beta = f(\alpha, \beta) = (R_f\beta)\alpha$ *i.e.,*

$$(L^1_{f_1}\alpha_1)\beta_1 \cup \ldots \cup (L^n_{f_n}\alpha_n)\beta_n = f_1(\alpha_1, \beta_1) \cup \ldots \cup f_n(\alpha_n, \beta_n)$$
$$= (R^1_{f_1}\beta_1)\alpha_1 \cup \ldots \cup (R^n_{f_n}\beta_n)\alpha_n.$$

*Then n-rank* $L_f$ *= n-rank* $R_f$*.*

The reader is expected to give the proof of the above theorem.

**DEFINITION 2.10**: *If $f = f_1 \cup \ldots \cup f_n$ is a bilinear n-form on the finite dimensional n-vectors space V and n-rank of f is the n-tuple of integers $(r_1, \ldots, r_n)$; $(r_1, \ldots, r_n)$ = n-rank $L_f$ = (rank $L^1_{f_1}$, ..., rank $L^n_{f_n}$) and n-rank $R_f$ = (rank $R^1_{f_1}$, ..., rank $R^n_{f_n}$).*

We state the following corollaries and the reader is expected to give the proof.

**COROLLARY 2.10:** *The n-rank of a bilinear n-form is equal to the n-rank of the n-matrix of the n-form in any ordered n-basis.*

**COROLLARY 2.11:** *If $f = f_1 \cup \ldots \cup f_n$ is a bilinear n-form on the $(n_1, \ldots, n_n)$ dimensional n-vector space $V = V_1 \cup \ldots \cup V_n$ the following are equivalent*

a. *n-rank $f = (n_1, \ldots, n_n)$ = (rank $f_1$, ..., rank $f_n$).*
b. *For each non zero $\alpha = \alpha_1 \cup \ldots \cup \alpha_n$ in V there is a $\beta = \beta_1 \cup \ldots \cup \beta_n$ in V such that $f(\alpha, \beta) \neq 0 \cup \ldots \cup 0$; $f_1(\alpha_1, \beta_1) \cup \ldots \cup f_n(\alpha_n, \beta_n) = 0 \cup \ldots \cup 0$.*
c. *For each non zero $\beta = \beta_1 \cup \ldots \cup \beta_n$ in V there is an $\alpha = \alpha_1 \cup \ldots \cup \alpha_n$ in V such that $f(\alpha,\beta) = 0 \cup \ldots \cup 0$ i.e., $f_1(\alpha_1, \beta_1) \cup \ldots \cup f_n(\alpha_n, \beta_n) = 0 \cup \ldots \cup 0$.*

Now we proceed onto define the non degenerate of a bilinear n-form.

**DEFINITION 2.11:** *A bilinear n-form $f = f_1 \cup \ldots \cup f_n$ on a n-vector space $V = V_1 \cup \ldots \cup V_n$ is called non-degenerate if it satisfies the conditions (b) and (c) of the above corollary.*

We now proceed onto define the new notion of symmetric bilinear n-forms.

**DEFINITION 2.12:** *Let $f = f_1 \cup \ldots \cup f_n$ be a bilinear n-form on a n-vector space $V = V_1 \cup \ldots \cup V_n$. We say that f is n-symmetric if $f(\alpha,\beta) = f(\beta,\alpha)$ for all $\alpha, \beta \in V$ i.e., $f_1(\alpha_1, \beta_1) \cup \ldots \cup f_n(\alpha_n, \beta_n) = f_1 (\beta_1, \alpha_1) \cup \ldots \cup f_n (\beta_n, \alpha_n)$.*

*If $f = f_1 \cup ... \cup f_n$ is a n-symmetric bilinear n-form the quadratic n-form associated with f is the function $q = q_1 \cup ... \cup q_n$ from $V = V_1 \cup ... \cup V_n$ onto $F = F_1 \cup ... \cup F_n$ defined by $q(\alpha) = q_1(\alpha_1) \cup ... \cup q_n(\alpha_n) = f(\alpha, \alpha) = f_1(\alpha_1, \alpha_1) \cup ... \cup f_n(\alpha_n, \alpha_n)$. If V is a real n-vector space an n-inner product on V is a n-symmetric bilinear form f on V which satisfies $f(\alpha, \alpha) > 0 \cup ... \cup 0$ i.e., $f_1(\alpha_1, \alpha_1) \cup ... \cup f_n(\alpha_n, \alpha_n) > (0 \cup ... \cup 0)$ where each $\alpha_i \neq 0$ for i = 1, 2, ..., t. A bilinear n-form in which $f_i(\alpha_i, \alpha_i) > 0$ for each i = 1, 2, ..., n is called n-positive definite.*

*So two n-vectors α, β in V are n-orthogonal with respect to a n-inner product $f = f_1 \cup f_2 \cup ... \cup f_n$ if $f(\alpha, \beta) = 0 \cup ... \cup 0$ i.e., $f(\alpha, \beta) = f_1(\alpha_1, \beta_1) \cup ... \cup f_n(\alpha_n, \beta_n) = 0 \cup ... \cup 0$.*

*The quadratic n-form $q(\alpha) = f(\alpha, \alpha)$ takes only non negative values.*

The following theorem is significant on its own.

**THEOREM 2.18:** *Let $V = V_1 \cup ... \cup V_n$ be a n-vector space over the n-field $F = F_1 \cup ... \cup F_n$ each $F_i$ of characteristic zero, i = 1, 2, ..., n and let $f = f_1 \cup ... \cup f_n$ be a n-symmetric bilinear n-form on V. Then there is an ordered n-basis for V in which f is represented by a diagonal n-matrix.*

*Proof:* What we need to find is an ordered n-basis $B = B_1 \cup \ldots \cup B_n = \{\alpha^1_1 \ldots \alpha^1_{n_1}\} \cup \ldots \cup \{\alpha^n_1 \ldots \alpha^n_{n_n}\}$ such that

$$\begin{aligned} f(\alpha_i, \alpha_j) &= f_1(\alpha^1_{i_1}, \alpha^1_{j_1}) \cup \ldots \cup f_n(\alpha^n_{i_n}, \alpha^n_{j_n}) \\ &= 0 \cup \ldots \cup 0 \end{aligned}$$

for $i_t \neq j_t$; t = 1, 2, …, n. If $f = f_1 \cup \ldots \cup f_n = 0 \cup \ldots \cup 0$ or $n_1 = n_2 = \ldots = n_n = 0$ the theorem is obviously true thus we suppose $f = f_1 \cup \ldots \cup f_n \neq 0 \cup \ldots \cup 0$ and $(n_1, \ldots, n_n) > (1, \ldots, 1)$. If $f(\alpha, \alpha) = 0 \cup \ldots \cup 0$ i.e., $f_1(\alpha_1, \alpha_1) \cup \ldots \cup f_n(\alpha_n, \alpha_n) = 0 \cup \ldots \cup 0$ for every $\alpha_i \in V_i$; i = 1, 2, …, n the associated n-quadratic form q is identically $0 \cup \ldots \cup 0$ and by the polarization n-identity

$$\begin{aligned} f(\alpha, \beta) &= f_1(\alpha_1, \beta_1) \cup \ldots \cup f_n(\alpha_n, \beta_n) \\ &= \left[\frac{1}{4} q_1(\alpha_1 + \beta_1) - \frac{1}{4} q_1(\alpha_1 - \beta_1)\right] \cup \ldots \cup \end{aligned}$$

$$\frac{1}{4}\left[ q_n(\alpha_n + \beta_n) \; - \; \frac{1}{4} q_n(\alpha_n - \beta_n) \right].$$

$$\begin{aligned} f &= f_1 \cup \ldots \cup f_n \\ &= 0 \cup \ldots \cup 0. \end{aligned}$$

Thus there is n-vector $\alpha = \alpha_1 \cup \ldots \cup \alpha_n$ in V such that $f(\alpha, \alpha) = q(\alpha) \neq 0 \cup \ldots \cup 0$ i.e.,

$$\begin{aligned} f_1(\alpha_1,\alpha_1) \cup \ldots \cup f_n(\alpha_n, \alpha_n) &= q_1(\alpha_1) \cup \ldots \cup q_n(\alpha_n) \\ &\neq 0 \cup \ldots \cup 0. \end{aligned}$$

Let W be the one-dimensional n-subspace of V which is spanned by $\alpha = \alpha_1 \cup \ldots \cup \alpha_n$ and let $W^{\perp}$ be the set of all n-vectors $\beta = \beta_1 \cup \ldots \cup \beta_n$ in V such that

$$\begin{aligned} f(\alpha, \beta) &= f_1(\alpha_1, \beta_1) \cup \ldots \cup f_n(\alpha_n, \beta_n) \\ &= 0 \cup \ldots \cup 0. \end{aligned}$$

Now we claim $W \oplus W^{\perp} = V$ i.e., $W_1 \oplus W_1^{\perp} \cup \ldots \cup W_n \oplus W_n^{\perp} = V = V_1 \cup \ldots \cup V_n$.

Certainly the n-subspaces W and $W^{\perp}$ are independent. A typical n-vector in V is $c\alpha$ where c is a n-scalar. If $c\alpha$ is in $W^{\perp}$ i.e.,

$$c\alpha = (c_1\alpha_1 \cup \ldots \cup c_n\alpha_n) \in W^{\perp} = W_1^{\perp} \cup \ldots \cup W_n^{\perp}$$

then

$$\begin{aligned} & f_1(c_1\alpha_1, c_1\alpha_1) \cup \ldots \cup f_n(c_n\alpha_n, c_n\alpha_n) \\ &= c_1^2 f_1(\alpha_1, \alpha_1) \cup \ldots \cup c_n^2 f_n(\alpha_n, \alpha_n) \\ &= 0 \cup \ldots \cup 0. \end{aligned}$$

But $f_i(\alpha_i, \alpha_i) \neq 0$ for every i, i = 1, 2, …, n thus each $c_i = 0$. Also each n-vector in V is the sum of a n-vector in W and a n-vector in $W^{\perp}$. For let $\gamma = \gamma_1 \cup \ldots \cup \gamma_n$ be any n-vector in V and put $\beta = \beta_1 \cup \ldots \cup \beta_n$

$$= \gamma - \frac{f(\gamma,\alpha)}{f(\alpha,\alpha)}\alpha$$

$$= \gamma_1 - \frac{f_1(\gamma_1,\alpha_1)}{f_1(\alpha_1,\alpha_1)}\alpha_1 \cup \ldots \cup \gamma_n - \frac{f_n(\gamma_n,\alpha_n)}{f_n(\alpha_n,\alpha_n)}\alpha_n.$$

Then

$$f(\alpha, \beta) = f(\alpha, \gamma) - \frac{f(\gamma,\alpha)}{f(\alpha,\alpha)} f(\alpha,\alpha);$$

and since f is n-symmetric, $f(\alpha, \beta) = 0 \cup \ldots \cup 0$. Thus $\beta$ is in the n-subspace $W^{\perp}$. The expression

$$\gamma = \frac{f(\gamma,\alpha)}{f(\alpha,\alpha)} \alpha + \beta$$

shows that $V = W + W^{\perp}$. Then n-restriction of f to $W^{\perp}$ is a n-symmetric bilinear n-form on $W^{\perp} = W_1^{\perp} \cup \ldots \cup W_n^{\perp}$. Since $W^{\perp}$ has $(n_1 - 1, \ldots, n_n - 1)$ dimension we may assume by induction that $W^{\perp}$ has a n-basis $\{\alpha_2^1, \ldots, \alpha_{n_1}^1\} \cup \ldots \cup \{\alpha_2^n, \ldots, \alpha_{n_n}^n\}$ such that $f_t(\alpha_{i_t}^t, \alpha_{j_t}^t) = 0$; $i_t \neq j_t$, $i_t \geq 2$, $j_t \geq 2$ and t = 1, 2, …, n. Putting $\alpha_1 = \alpha$ we obtain a n-basis $\{\alpha_1^1, \ldots, \alpha_{n_1}^1\} \cup \ldots \cup \{\alpha_1^n, \ldots, \alpha_{n_n}^n\}$ for V such that $f(\alpha_i, \alpha_j) = f_1(\alpha_{i_1}^1, \alpha_{j_1}^1) \cup \ldots \cup f_n(\alpha_{i_n}^n, \alpha_{j_n}^n) = 0 \cup \ldots \cup 0$ for $i_t \neq j_t$; t = 1, 2, …, n.

**THEOREM 2.19:** *Let $V = V_1 \cup \ldots \cup V_n$ be a $(n_1, \ldots, n_n)$ dimensional n-vector space over the n-field of real numbers and let $f = f_1 \cup \ldots \cup f_n$ be the n-symmetric bilinear n-form on V which has n-rank $(r_1, \ldots, r_n)$. Then there is an n-ordered n-basis $\{\beta_1^1 \ldots \beta_{n_1}^1\}$, $\{\beta_1^2 \ldots \beta_{n_2}^2\} \cup \ldots \cup \{\beta_1^n \ldots \beta_{n_n}^n\}$ for $V = V_1 \cup \ldots \cup V_n$ in which the n-matrix of $f = f_1 \cup \ldots \cup f_n$ is n-diagonal and such that*

$$\begin{aligned} f(\beta_j \beta_j) &= f_1(\beta_{j_1}, \beta_{j_1}) \cup \ldots \cup f_n(\beta_{j_n}, \beta_{j_n}); \\ &= \pm 1 \cup \ldots \cup \pm 1; \end{aligned}$$

*$1 \leq j_t \leq r_t$ and t = 1, 2, …, n. Further more the number of n-basis vectors $\beta_{j_t}$ for which $f_t(\beta_{j_t}, \beta_{j_t}) = 1$ is independent of the choice of the n-basis for t = 1, 2, …, n.*

The proof of the theorem is lengthy and left for the reader to prove.

Now for complex vector spaces we in general find it difficult to define n-vector spaces of type II as it happens to be the algebraically closed field. Further when n happens to be arbitrarily very large the problem of defining n-vector spaces of type II is very difficult.

First we shall call the field $F_1' \cup \ldots \cup F_n' = F'$ to be a special algebraically closed n-field of the n-field $F = F_1 \cup \ldots \cup F_n$ if and only if each $F_i'$ happens to be an algebraically closed field of $F_i$ for a specific and fixed characteristic polynomial $p_i$ of $V_i$ over $F_i$ relative to a fixed transformation $T_i$ of T, this is true of every i.

## Chapter Three

# SUGGESTED PROBLEMS

In this chapter we suggest nearly 120 problems about n-vector spaces of type II which will be useful for the reader to understand this concept.

1. Find a 4-basis of the 4-vector space V = $V_1 \cup V_2 \cup V_3 \cup V_4$ of (3, 4, 2, 5) dimension over the 4-field $Q(\sqrt{2}) \cup Q(\sqrt{3}) \cup Q(\sqrt{5}) \cup Q(\sqrt{7})$. Find a 4-subspace of V.

2. Given a 5-vector space V = $V_1 \cup V_2 \cup V_3 \cup V_4 \cup V_5$ over $Q(\sqrt{2}) \cup Q(\sqrt{3}, \sqrt{5}) \cup Q(\sqrt{7}) \cup Q(\sqrt{11}) \cup Q(\sqrt{13})$ of (5, 3, 2, 7, 4) dimension.

    a. Find a linearly independent 5-subset of V which is not a 5-basis.
    b. Find a linearly dependent 5-subset of V.
    c. Find a 5-basis of V.
    d. Does there exists a 5-subspace of (4, 2, 1, 6, 3) dimension in V?
    e. Find a 5-subspace of (4, 3, 1, 5, 2) dimension of V.

3. Let $V = V_1 \cup V_2 \cup V_3 \cup V_4$ be a 4-vector space over the 4-field $(Q(\sqrt{2}) \cup Q(\sqrt{7}) \cup Z_{11} \cup Z_2)$ of dimension (3, 4, 5, 6).

   a. Find a 4-subset of V which is a dependent 4-subset of V.
   b. Give an illustration of a 4-subset of V which is a independent 4-subset of V but is not a 4-basis of V.
   c. Give a 4-subset of V which is a 4-basis of V.
   d. Give a 4-subset of V which is semi n-dependent in V.
   e. Find a 4-subspace of V of dimension (2, 3, 4, 5) over the 4-field F where $F = (Q(\sqrt{2}) \cup Q(\sqrt{7}) \cup Z_{11} \cup Z_2)$ is a field.

4. Given $V = V_1 \cup V_2 \cup V_3$ is a 3-vector space of type II over the 3-field $F = Q \cup Z_2 \cup Z_5$ of dimension (3, 4, 5). Suppose

$$A = A_1 \cup A_2 \cup A_3$$

$$= \begin{bmatrix} 0 & -1 & 6 \\ 3 & 1 & 0 \\ 1 & 0 & -1 \end{bmatrix} \cup \begin{bmatrix} 2 & 0 & 3 & 0 \\ 1 & 2 & 0 & 1 \\ 0 & 0 & 1 & 3 \\ -1 & 1 & 0 & 0 \end{bmatrix} \cup \begin{bmatrix} 3 & 1 & 0 & 0 & 1 \\ 0 & 1 & 0 & -1 & 2 \\ -1 & 0 & 2 & 1 & 0 \\ 0 & -1 & 0 & 3 & 1 \\ 1 & 0 & 1 & 0 & 0 \end{bmatrix}$$

   is the 3-matrix find the 3-linear operator related with A.

   a. Is $A_3$ diagonalizable?
   b. Does A give rise to a 3-invertible 3-transformation?
   c. Find 3-nullspace associated with A.

5. Let $V = V_1 \cup V_2 \cup V_3 \cup V_4$ be a 4-vector space over the 4-field $F = (Q(\sqrt{3}) \cup Q(\sqrt{2}) \cup Z_5 \cup Z_2)$ of dimension (3, 2, 5, 4). Let $W = W_1 \cup W_2 \cup W_3 \cup W_4$ be a 4-vector space over the same F of dimension (4, 3, 2, 6).

Find a 4-linear transformation T from V into W, which will give way to a nontrivial 4-null subspace of V.

6. Let V = $V_1 \cup V_2 \cup V_3 \cup V_4 \cup V_5$ be a 5-vector space over the 5-field F = $Q(\sqrt{3}) \cup Q(\sqrt{2}) \cup Q(\sqrt{5}) \cup Z_5 \cup Z_2$ of type II of dimension (3, 4, 2, 5, 6).

   a. Find a 5-linear operator T on V such that T is 5-invertible.
   b. Find a 5-linear operator T on V such that T is non 5-invertible.
   c. Prove 5-rank T + 5 nullity T = (3, 4, 2, 5, 6) for any T.
   d. Find a T which is onto on V and find the 5-range of T.

7. Let V = $V_1 \cup V_2 \cup V_3 \cup V_4$ be a (5, 4, 3, 2)-4 vector space over the 4-field $Z_2 \cup Z_7 \cup Q(\sqrt{3}) \cup Q(\sqrt{2})$. Find V*. Obtain for any nontrivial 4-basis B its B* explicitly.

8. Given V* is a (3, 4, 5, 7) dimensional dual space over the 4-field F = $Z_2 \cup Z_7 \cup Z_3 \cup Z_5$. Find V. What is V** ?

9. Let V = $V_1 \cup V_2 \cup V_3 \cup V_4$ be a 4-vector space of (3, 4, 5, 6) dimension over the 4-field F = $Z_5 \cup Z_2 \cup Z_3 \cup Z_7$. Suppose W is a 4-subspace of (2, 3, 4, 5) dimension over the 4-field F. Prove dim W + dim $W^o$ = dim V. Find explicitly $W^o$. What is $W^{oo}$? Is $W^o$ a 4-subspace of V or V*? Justify your claim.

   Give a 4-basis for W and a 4-basis for $W^o$. Are these two sets of 4-basis related to each other in any other way?

10. Let V = $V_1 \cup V_2 \cup V_3$ be a 3-vector space over the 3-field F = $Z_7 \cup Z_5 \cup Z_2$ of (5, 4, 6) dimension over F. Let S = {(2 0 2 1 2), (1 0 1 1 1)} ∪ {(0 1 2 3), (4 2 0 4)} ∪ {(1 1 1 0 0 0), (1 1 1 1 1 0), (0 1 1 1 0 0)} = $S_1 \cup S_2 \cup S_3$ be a 3-set of V. Find $S^o$. Is $S^o$ a 3-subspace? Find the basis for $S^o$.

11. Given $V = V_1 \cup V_2 \cup V_3 \cup V_4 \cup V_5$ is a 5-space of (5, 4, 2, 3, 7) dimension over the 5-field $F = Q(\sqrt{3}) \cup Z_2 \cup Z_5 \cup Q(\sqrt{2}) \cup Z_7$. Find a 5-hypersubspace of V. Find V*? Prove W is a 5-hypersubspace of V. What is $W^o$? Find a 5-basis for V and its dual 5-basis. Find a 5-basis for W. What is its dual 5-basis?

12. Let $V = V_1 \cup V_2 \cup V_3 \cup V_4$ be a 4-space over the 4-field $F = F_1 \cup F_2 \cup F_3 \cup F_4$ of dimension (7, 6, 5, 4) over F. Let $W_1$ and $W_2$ be any two 4 subspaces of (4, 2, 3, 1) and (3, 5, 4, 2) dimensions respectively of V. Find $W_1^o$ and $W_2^o$. Is $W_1^o = W_2^o$? Find a 4-basis of $W_1$ and $W_2$ and their dual 4-basis. Is $W_1$ a 4-hypersubspace of V? Justify your claim! Find a 4-hyper subspace of V.

13. Let $V = V_1 \cup V_2 \cup V_3 \cup V_4 \cup V_5$ be a (2, 3, 4, 5, 6) dimensional 5-vector space over the 5-field $F = Z_2 \cup Z_3 \cup Z_5 \cup Z_7 \cup Q$. Find V*. Find a 5-basis and its dual 5-basis. Define a 5-isomorphism from V into V**. If $W = W_1 \cup W_2 \cup \ldots \cup W_5$ is a (1, 2, 3, 4, 5), 5-subspace of V find the n-annihilator space of V. Is W a 5-hyper space of V? Can V have any other 5-hyper space other than W?

14. Let $V = V_1 \cup V_2 \cup V_3 \cup V_4$ be a 4-space of (4, 3, 7, 2) – dimension over the 4-field $F = Z_3 \cup Z_5 \cup Z_7 \cup Q$. Find a 4-transformation T on V such that rank $(T^t)$ = rank T. (Assume T is a 4-linear transformation which is not a 4-isimorphism on V). Find 4-null space of T.

15. Give an example of 5-linear algebra over a 5-field which is not commutative 5-linear algebra over the 5-field.

16. Give an example of a 6-linear algebra over a 6-field which has no 6-identity.

17. Given A = $A_1 \cup A_2 \cup A_3$ where $A_1$ is a set of all 3 × 3 matrices with entries from $Q(\sqrt{3})$. $A_1$ is the linear algebra over $Q(\sqrt{3})$. $A_2$ = {All polynomials in the variable x with coefficient from $Z_7$}; $A_2$ is a linear algebra over $Z_7$ and $A_3$ = {set of all 5 × 5 matrices with entries from $Z_2$}. $A_3$ a linear algebra over $Z_2$. Is A a 3-linear algebra over the 3-field F = $Q(\sqrt{3}) \times Z_7 \times Z_2$? Is A a 3-commutative, 3-linear algebra over F? Does A contain the 3-identity?

18. Define a n-sublinear algebra of a n-linear algebra A over a n-field F.

19. Give an example of a 4-sublinear algebra of the 4-linear algebra over the 4-field.

20. Give an example of an n-vector space of type II which is not an n-linear algebra of type II.

21. Give an example of an n-commutative n-linear algebra over an n-field F (take n = 7).

22. Give an example of a 5-linear algebra which is not 5-commutative over the 5-field F.

23. Let A = $A_1 \cup A_2 \cup A_3$ where $A_1$ is a set of all 3 × 3 matrices over Q, $A_2$ = all polynomials in the variable x with coefficients from $Z_2$ and

$$A_3 = \left\{ \sum_{i=0}^{n} \alpha_{2i} x^{2i} \Big/ \alpha \in Z_3 \right\}.$$

Prove A is a 3-linear algebra over the 3-field F = $Q \cup Z_2 \cup Z_3$. Does A contain 3-identity find a 3-sublinear algebra of A over F. Is A a commutative 3-linear algebra?

24. Obtain Vandermode 4-matrix with (7 + 1, 6 + 1, 5 + 1, 3 + 1) from the 4-field F = $F_1 \cup F_2 \cup F_3 \cup F_4 = Q \cup Z_3 \cup Z_2 \cup Z_5$.

25. Using the 5-field $F = Z_2 \cup Z_7 \cup Z_3 \cup Z_5 \cup Z_{11}$, construct the 5-vector space = $V = Z_2 [x] \cup Z_7 [x] \cup Z_3 [x] \cup Z_5 [x] \cup Z_{11} [x]$. Verify

(i) $(cf + g)\alpha = cf (x) + g (\alpha)$

(ii) $(fg)\alpha = f(\alpha) g(\alpha)$

For $C = C^1 \cup C^2 \cup C^3 \cup C^4 \cup C^5$
$= 1 \cup 5 \cup 2 \cup 3 \cup 10$;
$f = (x^2 + x + 1) \cup (3x^3 + 5x + 1) \cup (2x^4 + x +1) \cup (4x + 2) \cup (10x^2 + 5)$ where $x^2 + x + 1 \in Z_2[x]$, $3x^3 + 5x + 1 \in Z_7[x]$, $2x^4 + x + 1 \in Z_3[x]$, $4x + 2 \in Z_5[x]$ and $10x^2 + 5 \in Z_{11}[x]$ and $g(x) = x^3 + 1 \cup x^2 + 1 \cup 3x^3 + x^2 + x + 1 \cup 4x^5 + x + 1 \cup 7x^2 + 4x + 5$ where $x^3 + 1 \in Z_2[x]$, $x^2 + 1 \in Z_7[x]$, $3x^3 + x^2 + x + 1 \in Z_3[x]$, $4x^5 + x + 1 \in Z_5[x]$ and $7x^2 + 4x + 5 \in Z_{11}[x]$.

26. Let $V = Z_2[x] \cup Z_7[x] \cup Q[x]$ be a 3-linear algebra over the 3-field $F = Z_2 \cup Z_7 \cup Q$. Let $M = M_1 \cup M_2 \cup M_2$ be the 3-ideal generated by $\{\langle x^2 + 1, x^5 + 2x + 1\rangle\} \cup \{\langle 3x^2 + x + 1, x + 5, 7x^3 + 1\rangle\} \cup \{\langle x^2 + 1, 7x^3 + 5x^2 + x + 3\rangle\}$. Is M a 3-principal ideal of V.

27. Prove in the 5-linear algebra of 5-polynomials $A = Z_3[x] \cup Z_2[x] \cup Q[x] \cup Z_7[x] \cup Z_{17}[x]$ over the 5-field $Z_3 \cup Z_2 \cup Q \cup Z_7 \cup Z_{17}$ every 5-polynomial $p = p^1 \cup p^2 \cup \ldots \cup p^5$ can be made monic. Find a nontrivial 5-ideal of A.

28. Let $A = A_1 \cup A_2 \cup \ldots \cup A_6 = Z_3[x] \cup Z_7[x] \cup Z_2[x] \cup Q[x] \cup Z_{11}[x] \cup Z_{13}[x]$ be a 6-linear algebra over the 6-field $F = Z_3 \cup Z_7 \cup Z_2 \cup Q \cup Z_{11} \cup Z_{13}$. Find a 6-minimal ideal of A. Give an example of a 6-maximal ideal of A.
*Hint:* We say in any n-polynomial n-linear algebra $A = A_1 \cup A_2 \cup \ldots \cup A_n$ over the n-field $F = F_1 \cup \ldots \cup F_n$ where $A_i = F_i[x]$; i = 1, 2, …, n. An n-ideal $M = M_1 \cup \ldots \cup M_n$ is said to be a n-maximal ideal of A if and only if each ideal $M_i$ of $A_i$ is maximal in $A_i$ for i = 1, 2, …, n. We say the n-ideal $N = N_1 \cup \ldots \cup N_n$ of A is n-minimal ideal of A if and only if each ideal $N_i$ of $A_i$ is minimal in $A_i$, for i

= 1, 2, …, n. An ideal $M = M_1 \cup \ldots \cup M_n$ is said to be n-semi maximal if and only if there exists atleast m number of ideals in M, $m \leq n$ which are maximal and none of them are minimal. Similarly we say $M = M_1 \cup \ldots \cup M_n$ is n-semi minimal if and only if M contains atleast some p-ideals which are minimal, $p \leq n$ and none of the ideals in M are maximal.

29. Give an example of an n-maximal ideal.

30. Give an example of an n-semi maximal ideal which is not n-maximal.

31. Give an example of an n-semi minimal ideal which is not an n-minimal ideal.

32. Let $A = Z_3[x] \cup Z_2[x] \cup Z_7[x] \cup Z_5[x]$ be a 4-linear algebra over the 4-field $F = Z_3 \cup Z_2 \cup Z_7 \cup Z_5$. Give an example of a 4-maximal ideal of A. Can A have a 4-minimal ideal? Justify your claim. Does A have a 4-semi maximal ideal? Can A have a 4-semi minimal ideal?

33. Give an example of a n-linear algebra which has both n-semi maximal and n-semi minimal ideals.

34. Give an example of a n-linear algebra which has n-ideal which is not a n-principal ideal? Is this possible if $A = A_1 \cup \ldots \cup A_n$ where each $A_i$ is $F_i[x]$ where A is defined over the n-field $F = F_1 \cup \ldots \cup F_n$; for i = 1, 2, …, n?

35. Let $A = Z_3[x] \cup Q[x] \cup Z_2[x]$ be a 3-linear algebra over the 3-field $Z_3 \cup Q \cup Z_2$. Find the 3-gcd of $\{x^2 + 1, x^2 + 2x + x^3 + 1\}$ $\{x^2 + 2, x + 2, x^2 + 8x + 16\}$ $\{x + 1, x^3 + x^2 + 1\}$ $= P = P_1 \cup P_2 \cup P_3$. Find the 3-ideal generated by P. Is P a 3-principle ideal of A? Justify your claim. Can P generate an n-maximal or n-minimal ideal? Substantiate your answer.

36. Find the 4-characteristic polynomial and 4-minimal polynomial for any 4-linear operator $T = T_1 \cup T_2 \cup T_3 \cup T_4$ defined on the 4-vector space $V = V_1 \cup V_2 \cup V_3 \cup V_4$ where V is a (3, 4, 2, 5) dimensional over the 4-field $F = Z_5 \cup Q \cup Z_7 \cup Z_2$.

a. Find a T on V so that the 4-characteristic polynomial is the same as 4-minmal polynomial.
b. Define a T on the V so that T is not a 4-diagonializable operator on V.
c. For a T in which the 4-chracteristic polynomial is different from 4-minimal polynomial find the 4-ideal of polynomials over the 4-field F which 4-annihilate T.
d. For every T can T be 4-diagonalizable. Justify your claim.

37. Let

$$A = \begin{bmatrix} 0 & 1 & 0 \\ 1 & 0 & 1 \\ 0 & 0 & 2 \end{bmatrix} \cup \begin{bmatrix} 4 & 0 & 1 & 1 \\ 0 & 1 & 1 & 3 \\ 1 & 0 & 0 & 0 \\ 2 & 1 & 0 & 0 \end{bmatrix} \cup \begin{bmatrix} 1 & 0 & 0 & 0 & 1 \\ 0 & 1 & 0 & 1 & 0 \\ 1 & 0 & 1 & 0 & 1 \\ 0 & 0 & 1 & 1 & 0 \\ 1 & 1 & 0 & 0 & 1 \end{bmatrix} \cup$$

$$\begin{bmatrix} 1 & 0 \\ 2 & 6 \end{bmatrix} \cup \begin{bmatrix} 0 & 1 & 0 & 2 & 3 & 5 \\ 7 & 0 & 0 & -1 & 0 & 0 \\ -1 & 0 & 2 & 0 & 0 & 1 \\ 0 & 0 & 1 & 1 & 0 & 0 \\ 0 & -1 & 0 & 0 & -1 & 0 \\ 1 & 0 & 0 & 1 & 1 & 1 \end{bmatrix}$$

be a 5-matrix over the five field; $F = Z_3 \cup Z_5 \cup Z_2 \cup Z_7 \cup Q$. Find the 5-characteristic value of A. Find the 5-characteristic vectors of A. Find the 5-characteristic polynomial and 5-minimal polynomial of A. Find the 5-characteristic space of all $\alpha$ such that $T\alpha = c\alpha$, T related to A.

38. Let $A = A_1 \cup \ldots \cup A_5$ be a $(5 \times 5, 3 \times 3, 2 \times 2, 4 \times 4, 6 \times 6)$, 5-matrix over the 5-field $F = Z_2 \cup Z_3 \cup Z_5 \cup Z_7 \cup Q$. If $A^2 = A$ i.e., $A_i^2 = A_i$ for i = 1, 2, …, 5. Prove A is 5-similar to a 5-diagonal matrix. What can you say about the 5-characteristic vectors and 5-characteristic values associated with A.

39. Let $V = V_1 \cup \ldots \cup V_6$ be a 6-vector space over a 6-field, $Z_2 \cup Q \cup Z_7 \cup Z_5 \cup Z_3 \cup Z_{11}$. Let T be a 6-diagonalizable operator on V. Let $W = W_1 \cup \ldots \cup W_6$ be a 6-invariant subspace of V under T. Prove the 6-restriction operator $T_W$ is 6-diagonalizable.

40. Let $V = V_1 \cup V_2 \cup V_3$ be a 3-vector space over the 3-field, $F = Q \cup Z_2 \cup Z_5$ of (5, 3, 4) dimension over F. Let $W = W_1 \cup W_2 \cup W_3$ be a 3-invariant 3-subspace of T (*Hint*: Find T and find its 3-subspace W which is invariant under T). Prove that the 3-minimal 3-polynomial for the 3-restriction 3-operator $T_W$ divides the 3-minimal polynomial for T. Do this without referring to 3-matrices.

41. Let $V = V_1 \cup \ldots \cup V_n$ be a $(3 \times 3, 2 \times 2, 4 \times 4, 5 \times 5)$, 4-space of matrices over the 4-field $F = Z_2 \cup Z_3 \cup Z_5 \cup Q$. Let T and U be 4-linear operators on V defined by

$$T (B) = AB$$

i.e., $T(B_1 \cup B_2 \cup B_3 \cup B_4) = A_1B_1 \cup \ldots \cup A_4B_4$

where $A = (A_1 \cup A_2 \cup A_3 \cup A_4)$ is a fixed chosen $(3 \times 3, 2 \times 2, 4 \times 4, 5 \times 5)$ matrix over

$$F = Z_2 \cup Z_3 \cup Z_5 \cup Q.$$

$$\begin{aligned} U(B) &= AB - BA \\ &= (A_1B_1 \cup A_2B_2 \cup A_3B_3 \cup A_4B_4) \\ &\quad - (B_1A_1 \cup B_2A_2 \cup B_3A_3 \cup B_4A_4) \\ &= (A_1B_1 - B_1A_1) \cup (A_2B_2 - B_2A_2) \cup \\ &\quad (A_3B_3 - B_3A_3) \cup (A_4B_4 - B_4A_4). \end{aligned}$$

a. If $A = (A_1 \cup A_2 \cup A_3 \cup A_4)$, choose fixed 4-matrix 4-diagonalizable; then T is 4-diagonalizable. Is this statement true or false?

b. If A is 4-diagonalizable then U is 4-diagonalizable. Is this statement true or false?

42. Suppose $A = A_1 \cup A_2 \cup A_3 \cup A_4 \cup A_5$ is a 5-triangular n-matrix similar to a 5-diagonal matrix then is A a 5-diagonal matrix? Justify your claim.

43. Let $V = V_1 \cup V_2 \cup V_3 \cup V_4$ be a (3, 4, 5, 6) dimensional vector space over the four field $F = F_1 \cup F_2 \cup F_3 \cup F_4$. Let $T = T_1 \cup T_2 \cup T_3 \cup T_4$ be a 4-linear operator on V. Suppose there exists positive integers $(k_1, k_2, k_3, k_4)$ such that $T_1^{k_1} = 0$, $T_2^{k_2} = 0$, $T_3^{k_3} = 0$ and $T_4^{k_4} = 0$. Will $T^{(3, 4, 5, 6)} = T_1^3 \cup T_2^4 \cup T_3^5 \cup T_4^6 = 0 \cup 0 \cup 0 \cup 0$ ?

44. Let $V = V_1 \cup \ldots \cup V_n$ be a n-vector space over the n-field $F = F_1 \cup F_2 \cup \ldots \cup F_n$. What is the n-minimal polynomial for the n-identity operator on V? What is the n-minimal polynomial for the n-zero operator?

45. Find a 3-matrix $A = A_1 \cup A_2 \cup A_3$ of order $(3 \times 3, 2 \times 2, 4 \times 4)$ such that the 3-minimal polynomial is $x^2 \cup x^2 \cup x^3$ over any suitable 3-field F.

46. Let $A = A_1 \cup \ldots \cup A_n$ be a $(n_1 \times n_1, n_2 \times n_2, \ldots, n_n \times n_n)$ matrix over the n-field $F = F_1 \cup \ldots \cup F_n$ with the n-characteristic n-polynomial $f = f_1 \cup \ldots \cup f_n =$

$$\left(x - c_1^1\right)^{d_1} \ldots \left(x - c_{k_1}^1\right)^{d_{k_1}^1} \cup \ldots \cup \left(x - c_1^n\right) \ldots \left(x - c_{k_n}^n\right)^{d_{k_n}^n} .$$

Show that trace A = trace $A_1 \cup \ldots \cup$ trace $A_n$

$$= c_1^1 d_1^1 + \ldots + c_{k_1}^1 d_{k_1}^1 \cup \ldots \cup c_1^n d_1^n + \ldots + c_{k_n}^n d_{k_n}^n .$$

47. Let $A = A_1 \cup \ldots \cup A_n$ be a n-matrix over the n-field $F = F_1 \cup \ldots \cup F_n$ of n-order $(n_1 \times n_1, \ldots, n_n \times n_n)$ with the n-characteristic polynomial

$$\left(x-c_1^1\right)^{d_1^1}\ldots\left(x-c_{k_1}^1\right)^{d_{k_1}^1}\cup\ldots\cup\left(x-c_1^n\right)^{d_1^n}\ldots\left(x-c_{k_n}^n\right)^{d_{k_n}^n}$$

where $\left(c_1^1,\ldots,c_{k_1}^1\right)\cup\ldots\cup\left(c_1^n,\ldots,c_{k_n}^n\right)$ are n-distinct n-characteristic values. Let V = $V_1 \cup \ldots \cup V_n$ be a n-vector space of ($n_1 \times n_1$, …, $n_n \times n_n$) matrices B = $B_1 \cup \ldots \cup B_n$ be such that AB = BA. Prove

$$\text{n dimV} = \left(\left(d_1^1\right)^2+\ldots+\left(d_{k_1}^1\right)^2,\ \ldots,\ \left(d_1^n\right)^2+\ldots+\left(d_{k_n}^n\right)^2\right).$$

48. Let V = $V_1 \cup \ldots \cup V_n$ be the n-space of ($n_1 \times n_1$, …, $n_n \times n_n$), n-matrices over the n-field F = $F_1 \cup \ldots \cup F_n$. Let A = $A_1 \cup \ldots \cup A_n$ be a fixed n-matrix of n-order ($n_1 \times n_1$, …, $n_n \times n_n$). Let T be a n-linear operator 'n-left multiplication by A' on V. Is it true that A and T have the same n-characteristic values?

49. Let A = $A_1 \cup \ldots \cup A_n$ and B = $B_1 \cup \ldots \cup B_n$ be two n-matrices of same n-order over the n-field F = $F_1 \cup \ldots \cup F_n$. Let n-order of A and B be ($n_1 \times n_1$, …, $n_n \times n_n$). Prove if (I – AB) is n-invertible then I-BA is n-invertible and $(I - BA)^{-1} = I + B(I - AB)^{-1}A$. Using this result prove both AB and BA have the same n-characteristic values in F = $F_1 \cup \ldots \cup F_n$ .

50. Let T is a n-linear operator of a ($n_1 \times n_1$, …, $n_n \times n_n$) dimensional n-vector space over the n-filed F = $F_1 \cup \ldots \cup F_n$ and suppose T has ($n_1$, …, $n_n$) distinct n-characteristic values. Prove T is n-diagonalizable.

51. Let A = $A_1 \cup \ldots \cup A_n$ be a ($n_1 \times n_1$, …, $n_n \times n_n$) triangular n-matrix over the n-field F. Prove that the n-characteristic values of A are the diagonal entries $\left(A_{ii}^1, A_{ii}^2, \ldots, A_{ii}^n\right)$.

52. Let $A = A_1 \cup \ldots \cup A_n$ be a n-matrix which is n-diagonal over the n-field $F = F_1 \cup \ldots \cup F_n$ of $(n_1 \times n_1, \ldots, n_n \times n_n)$ order i.e., if $A_k = \left(A^k_{ij}\right)$; $A^k_{ij} = 0$ if $i \neq j$ for $k = 1, 2, \ldots, n$.

Let $f = f_1 \cup \ldots \cup f_n$ be the n-polynomial over the n-field F defined by F =

$$\left(x - A^1_{11}\right)\ldots\left(x - A^1_{n_1 n_1}\right) \cup \left(x - A^2_{11}\right)\ldots\left(x - A^2_{n_2 n_2}\right) \cup \ldots \cup$$
$$\left(x - A^n_{11}\right)\ldots\left(x - A^n_{n_n n_n}\right).$$

What is the n-matrix of $f(A) = f_1(A_1) \cup \ldots \cup f_n(A_n)$?

53. Let $F = F_1 \cup \ldots \cup F_n$ be a n-field $F[x] = F_1[x] \cup \ldots \cup F_n[x]$ be the n-polynomial in the variable x. Show that the intersection of any number of n-ideals in F[x] is a minimal n-ideal.

54. Let $A = A_1 \cup \ldots \cup A_n$ be a n-matrix of $(n_1 \times n_1, n_2 \times n_2, \ldots, n_n \times n_n)$ order over a n-field $F = F_1 \cup \ldots \cup F_n$. Show that the set of n-polynomials $f = f_1 \cup \ldots \cup f_n$ in F[x] is such that $f(A) = f_1(A_1) \cup \ldots \cup f_n(A_n) = 0 \cup 0 \cup \ldots \cup 0$.

55. Let $F = F_1 \cup \ldots \cup F_n$ be a n-field. Show that the n-ideal generated by a finite number of n-polynomial $f^1, \ldots, f^n$. where $f^i = f^i_1 \cup \ldots \cup f^i_n$; $i = 1, 2, \ldots, n$ in $F[x] = F_1[x] \cup \ldots \cup F_n[x]$ is the intersection of all n-ideals in F[x] is an n-ideal.

56. Let $(n_1, \ldots, n_n)$ be a n-set of positive integers and $F = F_1 \cup \ldots \cup F_n$ be a n-field, let W be the set of all n-vectors

$$\left(x^1_1 \ldots x^1_{n_1}\right) \cup \left(x^2_1 \ldots x^2_{n_2}\right) \cup \ldots \cup \left(x^n_1 \ldots x^n_{n_n}\right)$$

in $F_1^{n_1} \cup F_2^{n_2} \cup \ldots \cup F_n^{n_n}$ such that

$$\left(x^1_1 + \ldots + x^1_{n_1}\right) = 0, \left(x^2_1 + \ldots + x^2_{n_2}\right) = 0, \ldots,$$
$$\left(x^n_1 + \ldots + x^n_{n_n}\right) = 0.$$

a. Prove $W^o = W^o_1 \cup W^o_2 \cup \ldots \cup W^o_n$ consists of all n-linear functionals $f = f_1 \cup \ldots \cup f_n$ of the form

$$f_1\left(x_1^1 \ldots x_{n_1}^1\right) \cup \ldots \cup f_n\left(x_1^n \ldots x_{n_n}^n\right) = c_i \sum_{j=1}^{n_1} x_j^1 \cup \ldots \cup c_n \sum_{j=1}^{n_n} x_j^n.$$

b. Show that the n-dual space W* of W can be naturally identified with n-linear functionals.

$$f_1\left(x_1^1 \ldots x_{n_1}^1\right) \cup \ldots \cup f_n\left(x_1^n \ldots x_{n_n}^n\right)$$
$$= c_1^1 x_1^1 + \ldots + c_{n_1}^1 x_{n_1}^1 \cup \ldots \cup c_1^n x_1^n + \ldots + c_{n_n}^n x_{n_n}^n$$

on $F_1^{n_1} \cup F_2^{n_2} \cup \ldots \cup F_n^{n_n}$ which satisfy $c_1^i + \ldots + c_{n_i}^i = 0$ for i = 1, 2, …, n.

57. Let W = $W_1 \cup \ldots \cup W_n$ be a n-subspace of a finite ($n_1$, …, $n_n$) dimensional n-vector space over V = $V_1 \cup \ldots \cup V_n$ and if

$$\left\{g_1^1 \ldots g_{r_1}^1\right\} \cup \left\{g_1^2 \ldots g_{r_2}^2\right\} \cup \ldots \cup \left\{g_1^n \ldots g_{r_n}^n\right\}$$

is a basis for $W^o = W_1^o \cup \ldots \cup W_n^o$ then

$$W = \bigcap_i N_{g_i} = \bigcap_{i_1=1}^{r_1} N_{g_{i_1}}^1 \cup \ldots \cup \bigcap_{i_n=1}^{r_n} N_{g_{i_n}}^n$$

where $\left\{N_1^1 \ldots N_{r_1}^1\right\} \cup \ldots \cup \left\{N_1^n \ldots N_{r_n}^n\right\}$ is the n-set of n-null space of the n-linear functionals

$$f = f^1 \cup \ldots \cup f^n = \left\{f_1^1 \ldots f_{r_1}^1\right\} \cup \ldots \cup \left\{f_1^n \ldots f_{r_n}^n\right\}$$

and

$$\left\{g_1^1 \ldots g_{r_1}^1\right\} \cup \ldots \cup \left\{g_1^n \ldots g_{r_n}^n\right\}$$

is the n-linear combination of the n-linear functions f = $f^1 \cup \ldots \cup f^n$.

58. Let S = $S_1 \cup \ldots \cup S_n$ be a n-set, F = $F_1 \cup \ldots \cup F_n$ a n-field. Let V(S; F) = $V_1(S_1; F_1) \cup V_2(S_2; F_2) \cup \ldots \cup V_n(S_n; F_n)$ the n-space of all n-functions from S into F;
i.e., $S_1 \cup \ldots \cup S_n$ into = $F_1 \cup \ldots \cup F_n$. (f + g) x = f (x) + g(x) and (cf)(x) = cf(x).
If f = $f_1 \cup \ldots \cup f_n$ and g = $g_1 \cup \ldots \cup \cup g_n$; then (f + g)(x)
= $f_1(x) + g_1(x) \cup \ldots \cup f_n(x) + g_n(x)$,

$cf(x) = c_1f_1(x) \cup \ldots \cup c_nf_n(x)$ where $c = c_1 \cup \ldots \cup c_n \in F = F_1 \cup \ldots \cup F_n$ . Let $W = W_1 \cup \ldots \cup W_n$ be any $(n_1, \ldots, n_n)$ dimensional space of $V(S; F) = V_1(S_1; F_1) \cup \ldots \cup V_n(S_n, F_n)$.
Prove that there exists $(n_1, \ldots, n_n)$ points

$$\left(x_1^1 \ldots x_{n_1}^1\right) \cup \left(x_1^2 \ldots x_{n_2}^2\right) \cup \ldots \cup \left(x_1^n \ldots x_{n_n}^n\right)$$

in $S = S_1 \cup \ldots \cup S_n$ and $(n_1, \ldots, n_n)$ functions

$$\left(f_1^1 \ldots f_{n_1}^1\right) \cup \left(f_1^2 \ldots f_{n_2}^2\right) \cup \ldots \cup \left(f_1^n \ldots f_{n_n}^n\right)$$

in W such that $f_i^t\left(x_j^t\right) = \delta_{ij}^t$; i = 1, 2, …, n.

59. Let $V = V_1 \cup \ldots \cup V_n$ be a n-vector space over the n-field $F = F_1 \cup \ldots \cup F_n$ and let $T = T_1 \cup \ldots \cup T_n$ be a n-linear operator on V. Let $C = C_1 \cup \ldots \cup C_n$ be a n-scalar suppose there is a non zero n-vector $\alpha = \alpha_1 \cup \ldots \cup \alpha_n$ in V such that $T\alpha = C\alpha$ i.e., $T_1\alpha_1 = C_1\alpha_1, \ldots, T_n\alpha_n = C_n\alpha_n$.
Prove that there is a non-zero n-linear functional $f = f_1 \cup \ldots \cup f_n$ on V such that $T^tf = Cf$ i.e., $T_1^tf_1 \cup \ldots \cup T_n^tf_n = C_1f_1 \cup \ldots \cup c_nf_n$.

60. Let $A = A_1 \cup \ldots \cup A_n$ be a $(m_1 \times m_1, \ldots, m_n \times m_n)$ n-mixed rectangular matrix over the n-field $F = F_1 \cup \ldots \cup F_n$ with real entries. Prove that $A = 0 \cup \ldots \cup 0$ (i.e., each $A_i = (0)$ for i = 1, 2, …, n) if and only if the trace $A^tA = (0)$ i.e., $A_1^tA_1 \cup \ldots \cup A_n^tA_n = 0 \cup \ldots \cup 0$.

61. Let $(n_1, \ldots, n_n)$ be n-tuple of positive integers and let V be the n-space of all n-polynomials functions over the n-field of reals which have n-degrees atmost $(n_1, n_2, \ldots, n_n)$ i.e., n-functions of the form

$$f(x) = f_1(x) \cup \ldots \cup f_n(x) =$$
$$c_0^1 + c_1^1x + \ldots + c_{n_1}^1x^{n_1} \cup \ldots \cup c_0^n + c_1^nx + \ldots + c_{n_n}^nx^{n_n}.$$

Let $D = D_1 \cup \ldots \cup D_n$ be the n-differential operator on $V = V_1 \cup \ldots \cup V_n$ over the n-field $F = F_1 \cup \ldots \cup F_n$. Find an n-basis for the n-null space of the n-transpose operator $D^t$.

62. Let $V = V_1 \cup \ldots \cup V_n$ be a n-vector space over the n-field $F = F_1 \cup \ldots \cup F_n$. Show that for a n-linear transformation $T = T_1 \cup \ldots \cup T_n$ on $T \to T^t$ is an n-isomorphism of $L^n(V, V)$ onto $L^n(V^*, V^*)$;
i.e., $T = T_1 \cup \ldots \cup T_n \to T^t = T_1^t \cup \ldots \cup T_n^t$ is an n-isomorphism of $L^n(V, V) = L(V_1, V_1) \cup \ldots \cup L(V_n, V_n)$ onto $L^n(V^*, V^*) = L\left(V_1^*, V_1^*\right) \cup \ldots \cup L\left(V_n^*, V_n^*\right)$.

63. Let $V = V_1 \cup \ldots \cup V_n$ be a n-vector space over the n-field $F = F_1 \cup \ldots \cup F_n$, where V is a n-space of $(n_1 \times n_1, n_2 \times n_2, \ldots, n_n \times n_n)$; n-matrices with entries from the n-field $F = F_1 \cup \ldots \cup F_n$.

a. If $B = B_1 \cup \ldots \cup B_n$ is a fixed $(n_1 \times n_1, n_2 \times n_2, \ldots, n_n \times n_n)$ n-matrix define a n-function $f_B = f_1B_1 \cup \ldots \cup f_nB_n$ on V by $f_B(A) = \text{trace}(B^tA)$;

$$f_B(A) = f_{B_1}(A_1) \cup \ldots \cup f_{B_n}(A_n)$$

$$= \text{trace}\left(B_1^t A_1\right) \cup \ldots \cup \text{trace}\left(B_n^t A_n^t\right)$$

where $A = A_1 \cup \ldots \cup A_n \in V = V_1 \cup \ldots \cup V_n$. Show $f_B$ is a n-linear functional on V.

b. Show that every n-linear functional on V of the above form i.e., is $f_B$ for some B.

c. Show that $B \to f_{1B_1} \cup \ldots \cup f_{nB_n}$ i.e., $B_1 \cup \ldots \cup B_n \to f_{1B_1} \cup \ldots \cup f_{nB_n}$ i.e., $B_i \to f_{iB_i}$; i = 1, 2, …, n is an n-isomorphism of V onto V*.

64. Let $F = F_1 \cup \ldots \cup F_n$ be a n-field. We have considered certain special n-linear functionals $F[x] = F_1[x] \cup \ldots \cup F_n[x]$ obtained via "evaluation at t" $t = t_1 \cup \ldots \cup t_n$ given by $L(f) = f(t)$.

$$L^1(f_1) \cup \ldots \cup L^n(f_n) = f_1(t_1) \cup \ldots \cup f_n(t_n).$$

Such n-functional are not only n-linear but also have the property that

$$L(fg) = L^1(f_1g_1) \cup \ldots \cup L^n(f_ng_n) =$$

$L^1(f_1)\, L^1(g_1) \cup \ldots \cup L^n(f_n)\, L^n(g_n)$.

Prove that if $L = L^1 \cup \ldots \cup L^n$ is any n-linear functional on $F[x] = F_1[x] \cup \ldots \cup F_n[x]$ such that $L(fg) = L(f)L(g)$ for all f, g then $L = 0 \cup \ldots \cup (0)$ or there is a $t = t_1 \cup \ldots \cup t_n$ such that $L(f) = f(t)$ for all f. i.e., $L^1(f_1) \cup \ldots \cup L^n(f_n) = f_1(t) \cup \ldots \cup f_n(t)$.

65. Let $K = K_1 \cup \ldots \cup K_n$ be the n-subfield of a n-field $F = F_1 \cup \ldots \cup F_n$ . Suppose $f = f_1 \cup \ldots \cup f_n$ and $g = g_1 \cup \ldots \cup g_n$ be n-polynomials in $K[x] = K_1[x] \cup \ldots \cup K_n[x]$. Let $M_K = M_{K_1} \cup \ldots \cup M_{K_n}$ be the n-ideal generated by f and g in K[x]. Let $M_F = M_{F_1} \cup \ldots \cup M_{F_n}$ be the n-ideal, n-generated in $F[x] = F_1[x] \cup \ldots \cup F_n[x]$. Show that $M_k$ and $M_F$ have the same n-monic generator. Suppose f, g are n-polynomial in $F[x] = F_1[x] \cup \ldots \cup F_n[x]$, if $M_F$ is the n-ideal of F[x] i.e., $M_F = M_{F_1} \cup \ldots \cup M_{F_n}$ . Find conditions under which the n-ideal $M_k$ of $K[x] = K_1[x] \cup \ldots \cup K_n[x]$ can be formed.

66. Let $F = F_1 \cup \ldots \cup F_n$ be a n-field. $F[x] = F_1[x] \cup \ldots \cup F_n[x]$. Show that the intersection of any number of n-ideals in F[x] is again a n-ideal.

67. Prove the following generalization of the Taylors formula for n-polynomials. Let f, g, h be n-polynomials over the n-subfield of complex numbers with n-deg f = (deg $f_1, \ldots,$ deg $f_n) \leq (n_1, \ldots, n_n)$, where $f = f_1 \cup \ldots \cup f_n$.

Then $f(g) = f_1(g_1) \cup \ldots \cup f_n(g_n)$

$$= \sum_{k_1=0}^{n_1} \frac{1}{\lfloor k_1} f_1^{(k_1)}(h_1)(g_1 - h_1)^{k_1} \cup \ldots \cup \sum_{k_n=0}^{n_n} \frac{1}{\lfloor k_n} f_n^{(k_n)}(h_n)(g_n - h_n)^{k_n}$$

68. Let $T = T_1 \cup \ldots \cup T_n$ be a n-linear operator on $(n_1, \ldots, n_n)$ dimensional space and suppose that T has $(n_1, \ldots, n_n)$ distinct n-characteristic values. Prove that any n-linear operator which commutes with T is a n-polynomial in T.

69. Let $V = V_1 \cup \ldots \cup V_n$ be a $(n_1, n_2, \ldots, n_n)$ dimensional n-vector space. Let $W_1 = \left(W_1^1 \cup \ldots \cup W_n^1\right)$ be any n-subspace of V. Prove that there exists a n-subspace $W_2 = \left(W_2^1 \cup \ldots \cup W_2^n\right)$ of V such that $V = W_1 \oplus W_2$ i.e.,

$$V = W_1 \oplus W_2$$
$$= \left(W_1^1 \oplus W_2^1\right) \cup \left(W_1^2 \oplus W_2^2\right) \cup \ldots \cup \left(W_1^n \oplus W_2^n\right).$$

70. Let $V = V_1 \cup \ldots \cup V_n$ be a $(n_1, \ldots, n_n)$ finite dimensional n-vector space and let $W_1, \ldots, W_n$ be n-subspaces of V such that $V = W_1 + \ldots + W_n$ i.e., $V = V_1 \cup \ldots \cup V_n = \left(W_1^1 + \ldots + W_{k_1}^1\right) \cup \left(W_1^2 + \ldots + W_{k_2}^2\right) \cup \ldots \cup \left(W_1^n + \ldots + W_{k_n}^n\right)$ and n dim V = (n dim $W_1$ + … + n dim $W_n$) = (dim $W_1^1$ + … + dim $W_{k_1}^1$, dim $W_1^2$ + … + dim $W_{k_2}^2$, …, dim $W_1^n$ + … + dim $W_{k_n}^n$). Prove that $V = W_1 + \ldots + W_n = \left(W_1^1 \oplus \ldots \oplus W_{k_1}^1\right) \cup \left(W_1^2 \oplus \ldots \oplus W_{k_2}^2\right) \cup \ldots \cup \left(W_1^n \oplus \ldots \oplus W_{k_n}^n\right)$.

71. If $E_1$ and $E_2$ are n-projections $E_1 = \left(E_1^1 \cup \ldots \cup E_1^n\right)$ and $E_2 = \left(E_2^1 \cup \ldots \cup E_2^n\right)$ on independent n-subspaces then is $E_1 + E_2$ an n-projection? Justify your claim.

72. If $E_1, \ldots, E_n$ are n-projectors of a n-vector space V such that $E_1 + \ldots + E_n = I$.
i.e., $\left(E_1^1 + \ldots + E_1^{k_1}\right) \cup \ldots \cup \left(E_n^1 + \ldots + E_n^{k_n}\right) = I_1 \cup \ldots \cup I_n$.
Then prove $E_t^i . E_t^j = 0$ if $i \neq j$ for t = 1, 2, …, n then $E_t^{p_t^2} = E_t^{p_t}$ for every t = 1, 2, …, n where $p_t = 1, 2, \ldots, k_t$.

73. Let $V = V_1 \cup \ldots \cup V_n$ be a real n-vector space and $E = E_1 \cup \ldots \cup E_n$ be an n-idempotent n-linear operator on V i.e., n-projection. Prove that $1 + E = 1 + E_1 \cup \ldots \cup 1 + E_n$ is n-invertible find $(1 + E)^{-1} = (1 + E_1)^{-1} \cup \ldots \cup (1 + E_n)^{-1}$.

74. Obtain some interesting properties about n-vector spaces of type II which are not true in case of n-vector spaces of type I.

75. Let T be a n-linear operator on the n-vector space V of type II which n-commutes with every n-projection operator on V. What can you say about T?

76. If N = $N_1 \cup \ldots \cup N_n$ is a n-nilpotent linear operator on a $(n_1, \ldots, n_n)$ dimensional n-vector space V = $V_1 \cup \ldots \cup V_n$ then the n-characteristic polynomial for N is ( $x^{n_1} \cup \ldots \cup x^{n_1}$ ).

77. Let V = $V_1 \cup \ldots \cup V_n$ be a $(n_1, n_2, \ldots, n_n)$ dimensional n-vector space over the n-field F = $F_1 \cup \ldots \cup F_n$ and T = $T_1 \cup \ldots \cup T_n$ be a n-linear operator on V such that n-rank T = (rank $T_1$, …, rank $T_n$) = (1, 1, …, 1). Prove that either T is n-diagonalizable or T is n-nilpotent, not both simultaneously.

78. Let V = $V_1 \cup \ldots \cup V_n$ be a $(n_1, n_2, \ldots, n_n)$ dimensional n-vector space over the n-field F = $F_1 \cup \ldots \cup F_n$ . Let T = $T_1 \cup \ldots \cup T_n$ be a n-operator on V. Suppose T commutes with every n-diagonalizable operator on V, i.e., each $T_i$ commutes with every $n_i$-diagonalizable operator on $V_i$ for i = 1, 2, …, n then prove T is a n-scalar multiple of the n-identity operator on V.

79. Let T be a n-linear operator on the $(n_1, n_2, \ldots, n_n)$ dimensional n-vector space V = $V_1 \cup \ldots \cup V_n$ over the n-field F = $F_1 \cup \ldots \cup F_{n.}$ Let p = $p_1 \cup \ldots \cup p_n$ be the n minimal polynomial for T i.e.,

$$p_1 = p_{11}^{r_1} \ldots p_{1k_1}^{r_{k_1}} \cup p_{21}^{r_2} \ldots p_{2k_2}^{r_{k_2}} \cup \ldots \cup p_{n1}^{r_n} \ldots p_{nk_n}^{r_{k_n}}$$

be the n-minimal polynomial for T.

Let V = $\left(W_1^1 \oplus \ldots \oplus W_{k_1}^1\right) \cup \ldots \cup \left(W_1^n \oplus \ldots \oplus W_{k_n}^n\right)$ = $W_1 \cup \ldots \cup W_n$ be the n-primary decomposition for T, i.e.,

$W_j^t$ is the null space of $p_{tj}(T_t)^{r_j}$, j = 1, 2, …, $k_t$. and t = 1, 2, …, n. Let $W^r$ = $W_1^r \cup \ldots \cup W_n^r$ be any n-subspace of V which is n-invariant under T. Prove that

$$W^r = \left(W_1^r \cap W_1\right) \oplus \left(W_2^r \cap W_2\right) \oplus \ldots \oplus \left(W_n^r \cap W_n\right).$$

80. Define some new properties on the n-vector spaces of type II relating to the n-linear operators.

81. Compare the n-linear operators on n-vector spaces V of type II and n-vector space of type I. Is every n-linear transformation of type I always be a n-linear operator of a type I n-vector space?

82. State and prove Bessel's inequality in case of n-vector spaces.

83. Derive Gram-Schmidt orthogonalization process for n-vector spaces of type II.

84. Define for a 3-vector space V = ($V_1 \cup V_2 \cup V_3$) of (7, 2, 5) dimension over the 3-field F = Q $\cup Z_3 \cup Z_2$ two distinct 3-inner products.

85. Let V = $V_1 \cup V_2 \cup V_3$ be a 3-spaces of (3, 4, 5) dimension over the 3-field F = $Q(\sqrt{2}) \cup Q(\sqrt{3}) \cup Q(\sqrt{7})$. Find $L^3$(V, V, F) = L($V_1$, $V_1$, $Q(\sqrt{2})$ ) $\cup$ L($V_2$, $V_2$, $Q(\sqrt{3})$) $\cup$ L ($V_3$, $V_3$, $Q(\sqrt{7})$).

86. Does their exists a skew-symmetric bilinear 5-forms on $R^{n_1} \cup R^{n_2} \cup \ldots \cup R^{n_5}$ ; $n_i \neq n_j$, $1 \leq i, j \leq 5$? Justify your claim.

87. Prove that 5-equation (Pf)($\alpha$, $\beta$) = ½ [f($\alpha$, $\beta$) – f($\beta$,$\alpha$)] defines a 5-linear operator P on $L^5$(V,V,F) where V = $V_1 \cup V_2 \cup \ldots \cup V_5$ is a 5-dimension vector space over a special 5-subfield $F = F_1 \cup F_2 \cup F_3 \cup F_4 \cup F_5$ of the complex field (such that $F_i \neq F_j$ if $i \neq j$, $1 \leq i, j \leq 5$.

a. V is of (7, 3, 4, 5, 6) dimension over $F = F_1 \cup \ldots \cup F_5$.
b. Prove $P^2 = P$ is a 5-projection.
c. Find 5-rank P and 5-nullily P. Is 5 rank P = (21, 3, 6, 10, 15) and 5-nullity P = (28, 6, 10, 15, 21).

88. Let $V = V_1 \cup \ldots \cup V_4$ be a finite (4, 5, 3, 6) dimension vector space over the 4-field $Q(\sqrt{2}) \cup Q(\sqrt{3}) \cup Q(\sqrt{5}) \cup Q(\sqrt{7})$ and f a symmetric bilinear 4-form on V. For each 4-subspace W of V. Let $W^{\perp}$ be the set of all 4-vector $\alpha = \alpha_1 \cup \alpha_2 \cup \alpha_3 \cup \alpha_4$ in V such that $f(\alpha, \beta) = 0 \cup 0 \cup 0 \cup 0$ for $\beta$ in W. Show that

a. $W^{\perp}$ is a 4-subspace.
b. $V = \{0\}^{\perp} \cup \{0\}^{\perp} \cup \{0\}^{\perp} \cup \{0\}^{\perp}$.
c. $V^{\perp} = \{0\} \cup \{0\} \cup \{0\} \cup \{0\}$ if and only if f is a non degenerate! Can this occur?
d. 4-rank f = 4-dim V – 4-dim $V^{\perp}$ .
e. If 4-dim V = $(n_1, n_2, n_3, n_4)$ and 4-dim W = $(m_1, m_2, m_3, m_4)$ then 4-dim $W^{\perp} \geq (n_1 - m_1, n_2 - m_2, n_3 - m_3, n_4 - m_4)$.
f. Can the 4-restriction of f to W be a non-degenerate if $W \cap W^{\perp} = \{0\} \cup \{0\} \cup \{0\} \cup \{0\}$?

89. Prove if U and T any two normal n operators which commute on a n-vector space over a n-field of type II prove U + T and UT are also normal n-operators ($n \geq 2$).

90. Define positive n-operator for a n-vector space of type II. Prove if S and T are positive n-operators every n-characteristic value of ST is positive.

91. Let $V = V_1 \cup \ldots \cup V_n$ be a $(n_1, n_2, \ldots, n_n)$ inner n-product space over a n-field $F = F_1 \cup \ldots \cup F_n$ . If T and U are positive linear n-operators on V prove that (T + U) is positive. Show by an example TU need not be positive.

92. Prove that every positive n-matrix is the square of a positive n-matrix.

93. Prove that a normal and nilpotent n-operator is the zero n-operator.

94. If $T = T_1 \cup \ldots \cup T_n$ is a normal n-operator prove that the n-characteristic n-vectors for T which are associated with distinct n-characteristic values are n-orthogonal.

95. Let $V = V_1 \cup \ldots \cup V_n$ be a $(n_1, n_2, \ldots, n_n)$ dimensional n-inner product space over a n-field $F = F_1 \cup \ldots \cup F_n$ for each n-vector $\alpha$, $\beta$ in V let $T_{\alpha,\beta} = T^1_{\alpha_1,\beta_1} \cup \ldots \cup T^n_{\alpha_n,\beta_n}$ (where $T = T^1 \cup \ldots \cup T^n$, $\alpha = \alpha_1 \cup \ldots \cup \alpha_n$ and $\beta = \beta_1 \cup \ldots \cup \beta_n$) be a linear n-operator on V defined by $T_{\alpha,\beta}(\gamma) = (\gamma / \beta)\alpha$ i.e., $T_{\alpha,\beta}(\gamma) = T^1_{\alpha_1,\beta_1}(\gamma_1) \cup \ldots \cup T^n_{\alpha_n,\beta_n}(\gamma_n) = (\gamma_1 / \beta_1)\alpha_1 \cup \ldots \cup (\gamma_n / \beta_n)\alpha_n$. Show that

a. $T^*_{\alpha,\beta} = T_{\beta,\alpha}$.

b. Trace $(T_{\alpha,\beta}) = (\alpha/\beta)$.

c. $(T_{\alpha,\beta})(T_{\gamma,\delta}) = T_{\alpha, (\beta/\gamma)\delta}$.

d. Under what conditions is $T_{\alpha,\beta}$ n-self adjoint?

96. Let $V = V_1 \cup \ldots \cup V_n$ be a $(n_1, n_2, \ldots, n_n)$ dimensional n-inner product space over the n-field $F = F_1 \cup \ldots \cup F_n$ and let $L^n(V, V) = L(V_1, V_1) \cup \ldots \cup L(V_n, V_n)$ be the n-space of linear n-operators on V. Show that there is a unique n-inner product on $L^n(V, V)$ with the property that $||T_{\alpha,\beta}||^2 = ||\alpha||^2\, ||\beta||^2$ for all $\alpha, \beta \in V$ i.e.,

$$\left\|T^1_{\alpha_1\beta_1}\right\|^2 \cup \ldots \cup \left\|T^n_{\alpha_n\beta_n}\right\|^2 = \left\|\alpha_1\right\|^2 \left\|\beta_1\right\|^2 \cup \ldots \cup \left\|\alpha_n\right\|^2 \left\|\beta_n\right\|^2.$$

$T_{\alpha,\beta}$ is an n-linear operator defined in the above problem. Find an n-isomorphism between $L^n(V,V)$ with this n-inner space of $(n_1 \times n_1, \ldots, n_n \times n_n)$, n-matrix over the n-field $F = F_1 \cup \ldots \cup F_n$ with the n-inner product $(A/B) = tr(AB^*)$ i.e., $(A_1/ B_1) \cup \ldots \cup (A_n / B_n) = tr(A_1B_1^*) \cup \ldots \cup$ $tr(A_nB_n^*)$ where $A = A_1 \cup ... \cup A_n$ and $B = B_1 \cup ... \cup B_n$ are $(n_1 \times n_1, \ldots, n_n \times n_n)$ n-matrix.

97. Let $V = V_1 \cup \ldots \cup V_n$ be a n-inner product space and let $E = E_1 \cup \ldots \cup E_n$ be an idempotent linear n-operator on V. i.e., $E^2 = E$ prove E is n-self adjoint if and only if $EE^* = E^*E$ i.e., $E_1 E_1^* \cup \ldots \cup E_n E_n^* = E_1^*E_1 \cup \ldots \cup E_n^*E_n$.

98. Show that the product of two self n-adjoint operators is self n-adjoint if and only if the two operators commute.

99. Let $V = V_1 \cup \ldots \cup V_n$ be a finite $(n_1, n_2, \ldots, n_n)$ dimensional n-inner product vector space of a n-field F. Let $T = T_1 \cup \ldots \cup T_n$ be a n-linear operator on V. Show that the n-range of T* is the n-orthogonal complement of the n-nullspace of T.

100. Let V be a finite $(n_1, n_2, \ldots, n_n)$ dimensional inner product space over the n-field F and T a n-linear operator on V. If T is n-invertible show that T* is n-invertible and $(T^*)^{-1} = (T^{-1})^*$.

101. Let $V = V_1 \cup \ldots \cup V_n$ be a n-inner product space over the n-field $F = F_1 \cup \ldots \cup F_n$. $\beta$ and $\gamma$ be fixed n-vectors in V. Show that $T\alpha = (\alpha/\beta)\gamma$ defines a n-linear operator on V. Show that T has an n-adjoint and describe T* explicitly.

102. Let $V = V_1 \cup \ldots \cup V_n$ be an n-inner product space of n-polynomials of degree less than or equal to $(n_1, \ldots, n_n)$ over the n-field $F = F_1 \cup \ldots \cup F_n$ i.e., each $V_i = F_i[x]$ for i = 1, 2, …, n. Let D be the differentiation on V. Find D*.

103. Let $V = V_1 \cup \ldots \cup V_n$ be a n-real vector space over the real n-field $F = F_1 \cup \ldots \cup F_n$. Show that the quadratic n-form determined by the n-inner product satisfies the n-parallelogram law.

$$||\alpha + \beta||^2 + ||\alpha - \beta||^2 = 2||\alpha||^2 + 2||\beta||^2$$

For $\alpha = \alpha_1 \cup \ldots \cup \alpha_n$ and $\beta = \beta_1 \cup \ldots \cup \beta_n$ in V.
i.e., $||\alpha_1 + \beta_1||^2 + ||\alpha_1 - \beta_1||^2 \cup \ldots \cup ||\alpha_n+\beta_n||^2 + ||\alpha_n - \beta_n||^2$
$= 2||\alpha_1||^2 + 2||\beta_1||^2 \cup \ldots \cup 2||\alpha_n||^2 + 2||\beta_n||^2$.

104. Let $V = V_1 \cup \ldots \cup V_n$ be a n-vector space over the n-field $F = F_1 \cup \ldots \cup F_n$. Show that the sum of two n-inner product on V is an n-inner product on V. Is the difference of two n-inner products an n-inner product? Show that a positive multiple of an n-inner product is an n-inner product.

105. Derive n-polarization identity for a n-vector space $V = V_1 \cup \ldots \cup V_n$ over the n-field $F = F_1 \cup \ldots \cup F_n$ for the standard n-inner product on V.

106. Let $A = A_1 \cup \ldots \cup A_n$ be a $(n_1 \times n_1, \ldots, n_n \times n_n)$ n-matrix with entries from n-field $F = F_1 \cup \ldots \cup F_n$ . Let $\left\{f_1^1 \ldots f_{n_1}^1\right\} \cup \ldots \cup \left\{f_1^n \ldots f_{n_n}^n\right\}$ be the n-diagonal entries of the n-normal form of $xI - A = xI_1 - A_1 \cup \ldots \cup xI_n - A_n$. For which n-matrix A is $\left(f_1^1, \ldots, f_1^n\right) \neq (1, 1, \ldots 1)$?

107. Let $T = T_1 \cup \ldots \cup T_n$ be a linear n-operator on a finite $(n_1, \ldots, n_n)$ dimensional vector space over the n-field $F = F_1 \cup \ldots \cup F_n$ and $A = A_1 \cup \ldots \cup A_n$ be a n-matrix associated with T in some ordered n-basis. Then T has a n-cyclic vector if and only if the n-determinants of $(n_1 - 1) \times (n_1 - 1), \ldots, (n_n - 1) \times (n_n - 1)$ n-submatrices of $xI - A$ are relatively prime.

108. Derive some interesting properties about n-Jordan forms or Jordan n-form (Just we call it as Jordan n-form or n-Jordan forms and both mean one and the same notion).

109.
    a. Let $T = T_1 \cup \ldots \cup T_n$ be a n-linear operator on the n space V of n-dim $(n_1, \ldots, n_n)$ . Let $R = R_1 \cup \ldots \cup R_n$ be the n-range of T. Prove that R has a n-

complementary T-n-invariant n-subspace if and only if R is n-independent of the n-null space $N = N_1 \cup \ldots \cup N_n$ of T.

b. Prove if R and N are n-independent N is the unique T-n-variant n-subspace complementary to R.

110. Let $T = T_1 \cup \ldots \cup T_n$ be a n-linear operator on the n-space $V = V_1 \cup \ldots \cup V_n$. If $f = f_1 \cup \ldots \cup f_n$ is a n-polynomial over the n-field $F = F_1 \cup \ldots \cup F_n$ and $\alpha_1 \cup \ldots \cup \alpha_n = \alpha$ and let $f(\alpha) = f(T)\alpha$

$$\text{i.e., } f_1(\alpha_1) \cup \ldots \cup f_n(\alpha_n) = f_1(T_1)\,\alpha_1 \cup \ldots \cup f_n(T_n)\,\alpha_n.$$

If $\left\{V_1^1 \ldots V_{k_1}^1\right\}, \ldots, \left\{V_1^n \ldots V_{k_n}^n\right\}$ are T-n-invariant n- sub space of $V = \left(V_1^1 \oplus \ldots \oplus V_{k_1}^1\right) \cup \ldots \cup \left(V_1^n \oplus \ldots \oplus V_{k_n}^n\right)$ show that $fV = \left(f_1 V_1^1 \oplus \ldots \oplus f_1 V_{k_1}^1\right) \cup \ldots \cup \left(f_n V_1^n \oplus \ldots \oplus f_n V_{k_n}^n\right)$.

111. Let T, V and F are as in the above problem (110). Suppose $\alpha_1 \cup \ldots \cup \alpha_n$ and $\beta = \beta_1 \cup \ldots \cup \beta_n$ are n-vectors in V which have the same T n-annihilator. Prove that for any n-polynomial $f = f_1 \cup \ldots \cup f_n$ the n-vectors $f\alpha = f_1\alpha_1 \cup \ldots \cup f_n\alpha_n$ and $f\beta = f_1\beta_1 \cup \ldots \cup f_n\beta_n$ have the same n-annihilator.

112. If $T = T_1 \cup \ldots \cup T_n$ is a n- diagonalizable operator on a n-vector space then every T n-invariant n-subspace has a n-complementary T-n-invariant subspace.

113. Let $T = T_1 \cup \ldots \cup T_n$ be a n-operator on a finite $(n_1, n_2, \ldots, n_n)$ dimensional n-vector space $V = V_1 \cup \ldots \cup V_n$ over the n-field $F = F_1 \cup \ldots \cup F_n$. Prove that T has a n cyclic vector if and only if every n-linear operator $U = U_1 \cup \ldots \cup U_n$ which commutes with T is a n-polynomial in T.

114. Let $V = V_1 \cup \ldots \cup V_n$ be a finite $(n_1, n_2, \ldots, n_n)$ n-vector space over the n-field $F = F_1 \cup \ldots \cup F_n$ and let $T = T_1 \cup$

… ∪ $T_n$ be a n-linear operator on V. When is every non zero n-vector in V a n-cyclic vector for T? Prove that this is the case if and only if the n-characteristic n-polynomial for T is n-irreducible over F.

115. Let T = $T_1 \cup \ldots \cup T_n$ be a n-linear operator on the ($n_1$, $n_2$, …, $n_n$) dimensional n-vector space over the n-field of type II. Prove that there exists a n-vector $\alpha = (\alpha_1 \cup \ldots \cup \alpha_n)$ in V with this property. If f is a n-polynomial i.e., $f = f_1 \cup \ldots \cup f_n$ and $f(T)\alpha = 0 \cup \ldots \cup 0$ i.e., $f_1(T_1)\alpha_1 \cup \ldots \cup f_n(T_n)\alpha_n = 0 \cup \ldots \cup 0$ then $f(T) = f_1(T_1) \cup \ldots \cup f_n(T_n) = 0 \cup \ldots \cup 0$. (such a n-vector is called a separating n-vector for the algebra of n-polynomials in T). When T has a n-cyclic vector give a direct proof that any n-cyclic n-vector is a separating n-vector for the algebra of n-polynomials in T.

116. Let T = $T_1 \cup \ldots \cup T_n$ be a n-linear operator on the n-vector space V = $V_1 \cup \ldots \cup V_n$ of ($n_1$, $n_2$, …, $n_n$) dimension over the n-field F = $F_1 \cup \ldots \cup F_n$ suppose that
   a. The n-minimal polynomial for T is a power of an irreducible n-polynomial.
   b. The minimal n-polynomial is equal to the characteristic n-polynomial. Then show that no nontrivial T-n-invariant n-subspace has an n-complementary T-n-invariant n-subspace.

117. Let A = $A_1 \cup \ldots \cup A_n$ be a ($n_1 \times n_1$, …, $n_n \times n_n$) n-matrix with real entries such that $A^2 + I = A_1^2 + I_1 \cup \cdots \cup A_n^2 + I_n$ = $0 \cup \ldots \cup 0$. Prove that ($n_1$, …, $n_2$) are even and if ($n_1$, …, $n_n$) = ($2k_1$, $2k_2$, …, $2k_n$) then A is n-similar over the n-field of real numbers to a n-matrix of the n-block form B = $B_1 \cup \ldots \cup B_n$ where B = $\begin{bmatrix} 0 & -I_{k_1} \\ I_k & 0 \end{bmatrix} \cup \cdots \cup \begin{bmatrix} 0 & -I_{k_n} \\ I_{k_n} & 0 \end{bmatrix}$ where $I_{k_t \times k_t}$ is a $k_t \times k_t$ identity matrix for t = 1, 2, …, n.

118. Let $A = A_1 \cup \ldots \cup A_n$ be a $(m_1 \times n_1, \ldots, m_n \times n_n)$ n-matrix over the n-field $F = F_1 \cup \ldots \cup F_n$ and consider the n-system of n-equation AX = Y i.e., $A_1X_1 \cup \ldots \cup A_nX_n = Y_1 \cup \ldots \cup Y_n$. Prove that this n-system of equations has a n-solution if and only if the row n-rank of A is equal to the row n-rank of the augmented n-matrix of the n-system.

119. If $P^n = P_1^n \cup \ldots \cup P_n^{n_n}$ is a stochastic n-matrix, is P a stochastic n-matrix? Show if $A = A_1 \cup \ldots \cup A_n$ is a stochastic n-matrix then (1, …, 1) is an n-eigen value of A.

120. Derive Chapman Kolmogorov equation for

$$p_{ij}^{n} = p_{i_1j_1}^{(n_1)} \cup \ldots \cup p_{i_nj_n}^{(n_n)}$$

$$= \sum_{k_1 \in s_1} p_{i_1j_1}^{(n_1-1)} p_{k_1j_1} \cup \cdots \cup \sum_{k_n \in s_n} p_{i_nj_n}^{(n_n-1)} p_{k_nj_n}$$

$S = S_1 \cup \ldots \cup S_n$ is an associated n-set.

# FURTHER READING

# INDEX

**M**

**N**

**P**

# ABOUT THE AUTHORS

**Dr.W.B.Vasantha Kandasamy** is an Associate Professor in the Department of Mathematics, Indian Institute of Technology Madras, Chennai. In the past decade she has guided 12 Ph.D. scholars in the different fields of non-associative algebras, algebraic coding theory, transportation theory, fuzzy groups, and applications of fuzzy theory of the problems faced in chemical industries and cement industries.

She has to her credit 646 research papers. She has guided over 68 M.Sc. and M.Tech. projects. She has worked in collaboration projects with the Indian Space Research Organization and with the Tamil Nadu State AIDS Control Society. This is her 39th book.

On India's 60th Independence Day, Dr.Vasantha was conferred the Kalpana Chawla Award for Courage and Daring Enterprise by the State Government of Tamil Nadu in recognition of her sustained fight for social justice in the Indian Institute of Technology (IIT) Madras and for her contribution to mathematics. (The award, instituted in the memory of Indian-American astronaut Kalpana Chawla who died aboard Space Shuttle Columbia). The award carried a cash prize of five lakh rupees (the highest prize-money for any Indian award) and a gold medal.

She can be contacted at vasanthakandasamy@gmail.com
You can visit her on the web at: http://mat.iitm.ac.in/~wbv

---

**Dr. Florentin Smarandache** is a Professor of Mathematics and Chair of Math & Sciences Department at the University of New Mexico in USA. He published over 75 books and 150 articles and notes in mathematics, physics, philosophy, psychology, rebus, literature.

In mathematics his research is in number theory, non-Euclidean geometry, synthetic geometry, algebraic structures, statistics, neutrosophic logic and set (generalizations of fuzzy logic and set respectively), neutrosophic probability (generalization of classical and imprecise probability). Also, small contributions to nuclear and particle physics, information fusion, neutrosophy (a generalization of dialectics), law of sensations and stimuli, etc. He can be contacted at smarand@unm.edu